S.-D. Poisson

**Researches into the Probabilities of Judgements
in Criminal and Civil Cases**

Recherches sur la probabilité des jugements
en matière criminelle et en matière civile

Paris, 1837

Translated by Oscar Sheynin

Berlin, 2013

# Contents

**From the Translator**

I have translated the main work of Poisson on probability theory but deleted many intermediate transformations not directly belonging to that discipline and full of formulas difficult to copy. Neither have I checked the derivations, sometimes insufficiently explained.

Poisson (Preamble, § 13) indicated (and obviously understated, see below) the aim of his first four chapters dealing with probability proper: they

*Include the general formulas of the calculus of probability which disperses with searching for them elsewhere and enables to treat other problems alien to the special aim of these researches but proper for the calculus to explain.*

Many authors have described the work of Poisson in the theory of probability, e. g., Bru (1981). I myself published a pertinent paper (1978), too long-winded but useful and was co-author of Gnedenko & Sheynin (1978) which contained a section on Poisson.

Here are some points contained in those four chapters.

**1. Random variable.** In § 53 Poisson defined a discrete random variable first appearing a few years earlier (Poisson 1829, § 2.5). Such variables were known and applied by Simpson, Laplace and Gauss and even before them, but Poisson introduced them formally although calling them by a provisional term. True, he still invariably considered extractions of balls of two or more colours from urns.

**2. Chances and probabilities.** Throughout his book, Poisson introduces both chances and probabilities. I am not at all sure that he was wholly consistent and anyway later authors have been applying much more appropriate terms, objective and subjective probabilities respectively. In many instances Poisson dealt with subjective probabilities and in a few cases his conclusions were therefore barely useful; one of his examples in § 11 was meaningless, see also Sheynin (2002). The probability of the studied event was 1/2, which signified *complete perplexity* (Poisson, § 4 of Chapter 1) and which information theory is telling us. The same is true about the celebrated Bertrand's problem about the length of a random chord.

**3. The law of large numbers and the central limit theorem.** For Poisson, that law was rather a loose principle, and for many decades statisticians had understood it just as Poisson. De Moivre had introduced the simplest form of the central limit theorem later developed by Laplace. Then several authors, Laplace certainly included, provided non-rigorous proofs of a few of its other versions, whereas Poisson's proof, on which he based his law of large numbers in a proper restricted sense, was methodically defective.

The law of large numbers (if understood in a loose sense) has been applied in the inverse case in which the theoretical probability was not known (or even did not exist) and had to be estimated (or effectively defined) by the observed frequencies. This case is less precise, but only Bayes rather than Jakob Bernoulli or De Moivre noted this circumstance. For his part, Poisson at the very least did not study it quantitatively.

**4. Causes of events.** Poisson consistently introduced causes of events, probabilities of their action and chances which they provide for the appearance of the events. Later authors did not follow him; Quetelet, for example, studied causes of events quite independently.

**5. Mathematical statistics.** When estimating the significance of empirical discrepancies, a problem which now would have been attributed to it, Poisson effectively applied consistent estimators. In the same context he introduced null-hypotheses both as main assumptions and tests of no difference.

Poisson's former student Gavarret (1840), who understood null-hypotheses in the second sense, declared that they were indispensable. Both Poisson and he restricted the application of statistics (a term that Poisson never mentioned!) to the availability of numerous observations. See Sheynin (2012).

**6. Theory of errors.** The main point here is that, because of the ill-fated priority strife between Legendre and Gauss, French mathematicians including Poisson had been ignoring the latter's pertinent fundamental work. Add to this that Laplace's theory of errors was barely useful since it demanded a large number of observations and his measure of error involved computations only possible in the case of the normal distribution. It was the existence of his own version of the error theory rather than the mentioned strife that prevented Laplace from following Gauss.

In several places of his book Poisson discussed geodetic measurements, but his deliberations were useless the more so since he was obviously ignorant of practical requirements and circumstances of field work.

**7. Three more points.** Poisson (§ 103) introduced functions now called after Dirac. They first appeared in Poisson (1811/1833, p. 637). For the sake of completeness I also mention that Poisson (1824, p. 278) introduced the distribution later called after Cauchy and effectively discovered its stability. And he (1829, § 1; 1837, pp. 63 and 80) also introduced distribution functions whose essential application did not begin until the 20$^{th}$ century.

As a most general conclusion, I ought to add, first, that for Poisson the theory of probability remained, just like for Laplace, an applied mathematical discipline; and second, that his book is not a *Lehrbuch*, although this is how its translator into German had called it.

**A few words about the fifth chapter.** This is Poisson's indirect explanation (Preamble, § 8) of its essence:

*The precise aim of the theory is to calculate for a jury panel composed of a certain number of people and judging a very large number of cases by an also known majority verdict, the future rates of acquittal and conviction likely to take place and the chance of mistaken judgement for those who were or will be randomly selected as jurymen.*

Poisson introduced both the probability of mistakes made by jurymen and judges and, unlike Laplace, the prior probability of the defendant's guilt. He studied the consequences of changes in criminal legislation on the conviction rate and on mistakes in convictions and acquittals but resolutely replaced the concepts of *guilt* and *innocence* by notions of *subject to be convicted, to be acquitted*.

Poisson extracted statistical data, including the particulars of verdicts (unanimous decisions, majority verdicts with various definite numbers of votes pro and con), by studying several sources. He never mentioned either the dependence between the decisions reached by jurors or judges or the inevitable changes with time in their attitudes to cases of the same kind. However, many of his conclusions were possibly more or less proper.

In spite of Laplace's opinion, the application of probability beyond natural sciences had been strongly criticized, for example by Poinsot (Poisson 1836, p. 380 of the Discussion), and, after the appearance of Poisson's book, by Mill (1843/1886, p. 353) and Poincaré (1896/1912, p. 20) who stated that people act like the *moutons de Panurge*. And Leibniz, in his letters to Jakob Bernoulli, insisted that a proper consideration of the pertinent circumstances was more important than calculations.

So how to estimate the influence of Poisson's study? Heyde & Seneta (1977, p. 31) concluded that *there was a surge of activity stimulated by Poisson*. They mentioned the preceding work of Poisson in the same field and the ensuing discussions and inserted a sketch of the pertinent work of earlier authors, but did not justify their conclusion.

I can only name Cournot (1843). He included a chapter on criminal statistics and, like earlier (1838), attempted to investigate the dependence between jurors or judges. In several places Cournot criticized Poisson, sometimes indirectly, but severely (§§ 61, 93, 149 Note, 225 and 237).

**Unpleasant circumstances**

**1.** The book was badly printed which is especially true with respect to formulas.

**2.** There are many misprints/mistakes unnoticed by Poisson. I have indicated at least some of them.

**3.** In the first two chapters, formulas are not numbered at all; in the next two, the enumeration is very incomplete. Poisson's references to previous formulas are therefore somewhat awkward, and in a few cases barely understandable. I have additionally numbered many formulas and adopted a convenient form of their enumeration.

**4.** The Preamble was not separated into sections and I myself accomplished this task.

**5.** Poisson obviously intended his work for a wider circle of readers; some explanations are too detailed, and conclusions are illustrated by possible bets on one or another circumstance. However, in many places the explanations are patently insufficient. Furthermore, at least the discussion of the Petersburg paradox (§ 25) is unsatisfactory and the *Bayes principle* is superficially described in a few words (§1 of Preamble).

**6.** For a modern reader Poisson's mention of *magnitude A* or *fraction p* is superfluous, but I did not dare to omit these words. Note that by *fraction* he always meant a proper fraction. Understandably, his terminology is dated but left without change. Some modern notation, for example, $n!$, $C_m^n$, $\bar{p}$, is lacking

I conclude that in an ideal world Poisson's book should have been rewritten anew rather than reprinted (Paris, 2003).

In 1841 Poisson's book was translated into German (see Bibliography) with four apparently forgotten supplements about 150 pages long by the translator (annuities and life insurance; moral expectation; probabilities of mean values, this being a translation of Poisson (1824 – 1829); application of probability to natural sciences). The print quality is good which allowed me to avoid at least some mistakes in copying formulas, and I also checked my translation of some phrases against their German rendition.

I am citing my website www.sheynin.de which is being diligently copied by Google: Oscar Sheynin, Home.

# Contents

# Chapter 5. Application of the General Rules of Probability to the Decisions of Jury Panels and the Judgements of the Tribunals.

Determining the probability that an accused person will be convicted or acquitted by a certain majority of the jurymen, each having a given probability of faultlessness, when considering, in addition, the given probability of the defendant's guilt existing before the judgement. Also determined, according to the rule about the probability of causes or hypotheses, are the probabilities of the accused, having been convicted or acquitted, being guilty or innocent, **§§ 114 – 117**

Formulas pertaining to the case of a certain number of jurymen, all of them having the same chance of faultlessness, whose judgement took or will take place, either by a given majority verdict, or by a verdict with only its minimal majority given. It is seen that the probability of a convictive verdict is always lower than the prior probability of guilt. The probability of a correct judgement only depends, other things being equal, on the majority of the verdict, but not on the total number of the jurymen[8] if their chance of faultlessness is given beforehand. This statement does not persist if that chance has to be determined *a posteriori* by issuing from that majority, **§§ 118 – 120**

Applying these formulas to the case in which the jury panel is very large which greatly lowers the probability of a conviction returned by a small majority, **§ 121**

A theorem about the jury panel composed of some number of jurymen each of them having many differing and unequally probable chances of faultlessness. An example of calculating the mean chance when the number of possible chances becomes infinite and their law of probabilities is given. This mean chance is the same for all the jurymen if they should be randomly selected from the same general list. Formulas which determine in that case the probabilities of a convictive verdict; of the convicted accused to be guilty; and of the chance of a mistake made by the jurymen to be contained within given limits, **§§ 122 – 127**

Application pf these formulas to a panel composed of a very large number of jurymen, **§§ 128 − 131**

In all cases, the use of those formulas require a formulated hypothesis about the law of probabilities of the jurymen's chances of error. Examining the Laplace's hypothesis. Its resulting consequences render it inadmissible. Since no appropriately justified hypothesis based on [some previously adopted] law of probabilities can be formulated, it is equally impossible to determine the probability of the correctness of an isolated judgement given the number of the jurymen and the majority of the verdict returned. The need to apply the results of a very large number of judgements for deriving the two special elements [parameters] included in the preceding formulas, the chance $u$ of faultlessness common to all the jurymen randomly selected from the same general list and the probability $k$ of the guilt of the accused resulting from the procedures preceding the pleadings at the assize courts, **§§ 132 – 133**

Probability of the difference between the rate of conviction derived from a series of trials and the special value, which this rate attains if the numbers of the convicted and accused become infinite, being contained between given limits. Probability of the difference between the former rate and the rate resulting from another series of trials to become contained within limits given as well, **§ 134**

Observations which serve for determining the numerical values of $u$ and $k$ extracted from the *Comptes généraux de l'administration de la justice criminelle*. These values are the various ratios taken into account before applying the preceding probability formulas. Influence of the successive changes in the legislation concerning the jury panels in France on the magnitude of those ratios. Separating crimes into two distinctive categories. At present, we are obliged to suppose that the values of $u$ and $k$, very different in those categories of crime, are respectively almost the same for all France, **§§ 135 – 138**

Calculating these values either for all France or for the département de la Seine in particular. Probability that, when issuing from those values, judgements of conviction or acquittal are returned unanimously, **§§ 139 – 141**

The meaning to be attached to the words *convictable* and *not guilty*. This is explained in more detail in the Preamble, **§ 142**

Formulas which indicate the measure of danger to the accused of being convicted although not being convictable, and to the society of acquitting an accused who should have been convicted, **§ 143**

Calculating these measures and the probabilities of innocence and guilt of the convicted accused in periods of differing legislation, **§§ 144 – 145**

Indicating a similar calculation, impossible owing to the lack of necessary observations, concerning the judgements of the police courts[9] and of the military justice, **§ 146**

Formulas concerning the correctness, more or less probable, of judgements in civil cases returned in courts of first instance and courts of appeal, **§§ 147 – 149**

Lacking the observational data necessary for determining the two distinctive elements [parameters] included in those formulas, we are obliged to suppose that the chances of error are the same for all the judges in both these courts. Calculating that [common] chance by issuing from the ratio, given by observations, of the number of judgements confirmed by the royal courts to that of the judgements returned by the courts of first instance and yearly submitted to them. During three consecutive years, that ratio varied but little which is a very remarkable proof of the general law of large numbers. Issuing from those observational data, the probabilities of correct judgements of both types of courts either coinciding or not are derived, **§§ 150 – 151**

## Notes

**1.** The numbering of the sections was mistaken, so that § 13 follows after § 11. Poisson

**2.** *Thing A* (chose A) was the author's term for *random variable*. Poisson (1829, p. 3) even before 1837 formally introduced the random variable and, although naming it by a provisional term, applied it in probability theory at large. He hardly connected the letter A with *aléa* (chance, risk) since he (1811/1833, pp. 141, 146) denoted a constant by the same letter. See also my note *From the Translator*.

**3.** Not the *rule*, but the definition of mathematical expectation.

**4.** In § 2 of Chapter 1 Poisson mentioned that irrational values of probability were also possible.

**5.** When applying this general proposition, for example, to therapeutics, it also conforms to simple common sense. If a medicine was successfully taken by a very large number of similar patients, so that the number of cases when it failed to help was very small as compared with the total number of these experiments, it will likely be successful in a new trial. Medicine will not become either a science or an art if not based on numerous observations, on the tact and proper experience of the physicists who judge the similarity of cases and take into account exceptional circumstances. Poisson

**6.** The author suggests a possible cause of those glowing bodies.

**7.** The author's term was *équations de condition*, but it contradicts the terminology adopted in the classical theory of errors.

**8.** Ostrogradsky made the same conclusion a few years earlier, see Note 1 to Chapter 5.

**9.** Police courts dealt with minor offences by summary jurisdiction.

**Preamble**

**[1]** A problem about games of chance proposed by a man about town to an austere Jansenist[1] became the origin of the calculus of probability. Its aim was to determine the proportion according to which the *stakes* should be shared by the gamblers if they agree to stop playing and to base their decision on the unequal number of points gained by them[2]. Pascal was the first to solve that problem although only in the case of two gamblers. Fermat solved it somewhat later in the general case of any number of gamblers.

Nevertheless, those geometers of the 17th century who had been occupying themselves with the calculus of probability merely determined chances in various contemporary games[3]. Only in the next century did that calculus wholly extend and become one of the main branches of mathematics owing both to the number and utility of its applications and to the type of analysis which it beget.

One of the most important among those applications concerns the probability of judgements, or, generally, of decisions reached by majority vote. Condorcet (1785) made the first attempt to determine it. He wrote his book during the lifetime, and at the request of Minister Turgot who conceived all the advantages which the moral sciences and the public administration could have elicited from the calculus of probability, whose indications are always precise even when lacking sufficient observational data and not therefore leading to complete solutions of problems.

Condorcet's work includes a preliminary discourse in which the author describes his results without the aid of analytical formulas and thoroughly develops suitable considerations for proving the utility of research of that kind. In his *Traité des probabilités* [*Théorie analytique des probabilités*] Laplace also occupied himself with calculating the chances of the error to be feared in judgements of tribunals or jury panels of a given strength by known majority vote returned against an accused. His solution, one of the most delicate in the theory of probability[4], is based on the principle that serves to determine the probabilities of various causes which can be attributed to the observed facts. At first, Bayes presented it in a slightly different form, and Laplace, in his [earlier] memoirs and treatise, discovered its most successful application for calculating the probabilities of future events by issuing from observation of the past events.

However, with regard to the problem of the probability of judgements, it is fair to say, that the ingenious idea of subordinating the solution to the Bayes principle by successively considering the guilt and innocence of the accused as an unknown cause of the received judgement which is the observed fact, is due to Condorcet. By issuing from this fact it is required to deduce the probability of that cause. The exactitude of this principle is demonstrated with full rigor, and the applications to the question with which we are occupied can not at all be doubtful[5]. Nevertheless, for applying it, Laplace introduced a hypothesis which is not at all incontestable. He supposed that the probability of a juryman's faultlessness can take all equally possible degrees from certitude to indifference which corresponds to 1/2 in the calculus and means an equality of chances of error and verity. The illustrious geometer[6] based his hypothesis on the opinion that a juryman [a judge] doubtlessly tends rather to verity than to error, and this should really be admitted in general. However, there exist infinitely many different laws of probabilities of error satisfying the stated condition without requiring the supposition either that the chance of a juryman to be faultless can never descend lower than 1/2, or that above that limit all its values are equally possible. Therefore, Laplace's particular presumption can not be justified beforehand. Either because of that hypothesis or owing to its consequences, which seem to me inadmissible, the solutions of the problem about the probability of judgements contained in Laplace's treatise (1812/1886, pp. 469 – 470) and its *Premier supplément* (1816/1886, p. 528) differing one from another invariably left me very doubtful.

**[2]** I would have applied to the illustrious author had I occupied myself with this problem during his lifetime. The authority of his name would have obliged me to do so and his friendliness with which he always glorified me would have facilitated its realization. It is easy to conceive that only after a long reflection I have decided to consider the question from another viewpoint and that before going ahead I may allow myself to describe the main reasons which impelled me to abandon his last solution whose numerical results he inserted in his *Essai* (1814/1996, Chapter *On the probability of Judicial Decisions*).

Laplace's formula expressing the probability of an erroneous judgement only depends on the majority vote and the total number of judges but does not include anything about their more or less extensive knowledge of the case under their consideration. It follows, however, that the probability of a mistaken decision reached, for example, by a majority of 7 votes against 5, will be the same as delivered by 12 jurymen selected form any group of people. For me, this consequence alone seems sufficient for definitely rejecting the formula from which it was derived.

That same formula supposes that prior to the decision of the jury panel there is no presumption of the accused being guilty; the more or less high probability of his guilt should only be deduced from the convicting decision. This, however, is once more inadmissible. The accused, upon arriving at the assize courts[7], had already been under detention awaiting trial, then committed for trial. This establishes against him a probability higher than 1/2 that he is guilty. And in a fair game certainly no one will hesitate to stake more for his guilt than for his innocence.

Now, the rules for establishing the probability of an observed event given the probability of its cause, which are the basis of the theory under consideration, require taking into account all the presumptions prior to the observation, if only they are thought to exist, or if proven that they are not absent. On the contrary [contrary to Laplace], such a presumption is evident owing to the criminal proceedings and I am obliged to consider it when solving that problem. And it is actually seen that otherwise it will be impossible to reconcile the consequences of calculations with the invariable results of observations.

This presumption is similar to that which takes place in civil cases when one of the litigants appeals against the first judgement to a superior court: there exists the assumption contrary to his cause, and it will be a grave mistake to disregard this circumstance when calculating the probability of error to be feared in the definitive decision. Finally, Laplace restricted his attention to considering the probability of judgement received by a given majority vote. However, the danger that the accused, when brought before the jury panel, will suddenly be mistakenly convicted by that majority, depends not only on that probability; in addition, it depends on the chance that such a conviction will be pronounced. And, admitting for a moment that the probability of an erroneous judgement returned by a majority of 7 votes against 5 is expressed by a fraction almost equal to 2/7, as it results from Laplace's formula, it should also be noted that, according to experience[8], each year the jury panels in France convict only 7/100 of the accused by that majority. The danger for the accused to be falsely judged by 7 votes against 5 is therefore measured by the product of 2/7 and 7/100, or by 1/50. Indeed, for all eventual things the danger of a loss or the hope of a gain is expressed by the product of the value of that thing which is dreaded or expected, multiplied by the probability that it will take place.

   This consideration alone reduces the proportion of the innocent accused being yearly convicted by the least possible majority vote of the jury panel to 1/50. Still, this is doubtless too much if all those accused were really innocent, but it is convenient to expound right here the real sense to be attributed in this (? - O.S.)[i] theory to the words *guilty* and *innocent*, the sense that Laplace and Condorcet [Condorcet and Laplace] had actually attributed to them. We will never be able to find out mathematically whether an accused is guilty. Even his confession cannot be regarded as a probability very nearly equal to certainty. The most enlightened and humane juryman will only pronounce conviction when having a strong probability of guilt, [although] often lower than that resulting from a confession of guilt.

   There exists an essential difference between him and a judge in a civil case. When the judge, after deeply examining a case, is only able, owing to its difficulty, to discern a feeble probability in favour of one of the two parties, it will suffice for convicting the other party [for deciding against it]. The juryman, however, should only vote for conviction if, to his eyes, the probability that the accused is guilty, reaches a certain limit and much exceeds the probability of his innocence.

   **[3]** However hard you attempt, it is impossible to avoid every chance of error in criminal judgements, so how much it could be reduced for ensuring the innocent man the greatest possible guarantee? It is difficult to answer this question in a general way. Condorcet thought that the chance to be unjustly convicted can be equivalent to the chance of such a danger which we believe to be sufficiently slim for attempting to avoid it in habitual life. Because, as he states, for its security the society can rightfully expose its member to a danger whose chance is for him, so to say, indifferent. For such a grave problem that consideration is, however, much too subtle.

Laplace offered a definition much more suitable for clearing up this question of the chance of error which we are compelled to admit in judgements in criminal cases. He maintained that that probability should be such that more danger for public security results from acquitting a guilty man than occurs from the fear of convicting an innocent. And, as he expressively stated, it was this question rather than the guilt itself of the accused that each juryman was called for to decide, in his own way, according to his enlightenment and opinion. An error in his vote, whether he convicted or acquitted, can originate from two different causes. Either he wrongly appreciates the arguments contrary to, or favourable for the accused, or his established limit of probability necessary for conviction is too high or too low. Not only is that limit different for all those selected for judging, it also changes with the nature of the accusations and even depends on the circumstances under which the proceedings are taking place. In the army, face to face with the enemy, or when trying cases of espionage, that limit is doubtlessly much lower than in ordinary cases. It lowers, and the number of convictions increases in cases of such crimes that become more frequent and more dangerous for the society.

The decisions of jury panels therefore concern the appropriateness of conviction or acquittal. The language will become more exact if the really true word *convictable* is substituted instead of *guilty* which needs to be explained and which we continue to use so as to conform to custom. And so, if we find out that in a very large number of judgements there was a certain proportion of mistaken convictions, it should not be understood that all the thus convicted were innocent. That proportion concerns those convicted for whom the probability was too low not for establishing their guilt rather than their innocence, but for their conviction to become necessary for public security. It is not the aim of our calculations to determine the number of those who, among the convicted, were really not guilty although there is ground to believe that that number is fortunately very small, at least beyond political cases.

For the ordinary cases this conclusion can be confirmed by the very small number of convictions pronounced by jury panels but protested by public opinion; by the small number of freely pardoned; by the number, once more very small, of cases in which an assize court, having decided that the oral debate at the court of first instance had destroyed the accusation and that the convicted was not guilty, had enjoyed the right granted them by law to annul a conviction returned by a jury panel and to return the accused to be judged by other jurymen.

The results concerning the chances of error in criminal judgements, to which Laplace had arrived, seem exorbitant, in disagreement with general ideas which contradict his own phrase that *the theory of probability is basically only common sense reduced to a calculus* (Laplace 1814/1995, p. 124). They [the results] were wrongly interpreted and it is too hasty to conclude that mathematical analysis is not applicable at all to such problems, or generally to things called moral.

**[4]** This is a prejudice which, as I see with regret, is shared by those well disposed. For destroying it, I believe it useful to recall here some general considerations which, in addition, are proper for comparing it with other problems. No one contests that for them the calculus is legitimate and necessary, and such comparison is also needed for thoroughly understanding the aim of the problem which I specially proposed in this work.

Things of every kind obey a universal law that we may call the law of large numbers. Its essence is that if we observe a very large number of events of the same nature, depending on constant causes and on causes varying irregularly, sometimes in one manner, sometimes in another, i.e., not progressively in any determined sense, then almost constant proportions will be found among these numbers. For things of every kind these ratios will have a special value from which they deviate ever less as the series of events increases, and which they reach if that series can be prolonged to infinity.

As the amplitudes of the variation of the irregular causes become more or less large, the number of events should also be more or less large for their ratio almost to attain permanence. Observations themselves will show, in each problem, whether the series of trials is sufficiently long. And, according to the number of established facts and the magnitude of the deviation still remaining between their ratios, the calculus will provide definite rules for determining the probability that the special value to which those ratio converge is contained between arbitrarily confined limits. If new trials are made and it is established that those ratios considerably deviate from their final values determined by the previous observations, it can be concluded that the causes on which depend the observed facts had experienced a progressive variation or even some abrupt change in the interval between the two series of observations.

Without the aid of the calculus of probability you run a great risk of being mistaken about the necessity of that conclusion. However, the calculus leaves nothing vague here and in addition provides necessary rules for determining the chance of the change of the causes indicated by comparing the observed facts at different times.

The law of large numbers is noted in events which are attributed to pure chance since we do not know their causes or because they are too complicated[9]. Thus, games, in which the circumstances determining the occurrence of a certain card or certain number of points on a die infinitely vary, can not be subjected to any calculus. If the series of trials is continued for a long time, the different outcomes nevertheless appear in constant ratios. Then, if calculations according to the rules of a game are possible, the respective probabilities of eventual outcomes conform to the known Jakob Bernoulli theorem. However, in most problems of contingency a prior determination of chances of the various events is impossible and, on the contrary, they are calculated from the observed result.

**[5]** For example, it is impossible to calculate beforehand the probability of a loss of a ship during a long voyage, and it is determined by comparing the number of shipwrecks and voyages and when the latter is very large, the ratio of those numbers is almost

constant, at least for each sea and each nation in particular. Its value can be assumed as the probability of future shipwrecks and it is this natural consequence of the law of large numbers that serves as a basis for marine insurance. If the insurer only deals with an insignificant number of cases, it is a simple bet with no values for computation, but if he operates a very large number of them, it will be a speculation whose success is almost certain.

Just the same, the law of large numbers governs phenomena produced by known forces acting together with accidental causes whose effect lacks any regularity. Successive elevations and abatements of the sea at harbours and sea-coasts offer an example of remarkable precision. In spite of the inequalities produced by the winds, which destroy the laws of the mentioned phenomenon in isolated or not numerous observations, the mean of a large number of observed tides at the same place reveals that they almost conform to the laws of *ebb* and *flow* of the sea resulting from the attraction of the Moon and the Sun as though the accidental winds had no influence. What can be the effect of the winds blowing in the same direction for some part of the year on the tides during that time is not yet determined at all.

The small difference between these means calculated for observations made at the beginning and end of the last [the $18^{th}$] century, and therefore separated by [about] a hundred years from each other, can be attributed to some changes having occurred in the localities. To provide [another] example of the law I am considering, I also refer to the mean length of human life. Among a considerable number of infants born each year in very near places and about the same time some die early, others live longer and still others reach the limits of longevity. And, in spite of the vicissitudes of human life, which produce such great differences in the ages of the dying, when dividing the sum of those ages by their number supposed to be very large, the quotient, or what is called the *mean life*, will not depend on that number.

Its duration can be different for the two sexes, in different countries and at different times because it depends on the climate and doubtless on the well-being of the people[10]. It increases if a disease disappears like smallpox did due to the blessing of the vaccine. And in all cases the calculus of probability shows whether those variations revealed in that duration are sufficiently large and manifested in a sufficiently large number of observations for necessarily attributing them to some changes having taken place in general causes.

The yearly sex ratio at birth in a large country has a constant value just as well. It does not seem to depend on the climate, but, owing to some singularity, to which it will possibly be easy to assign a likely cause, it appears to be different for infants born in and out of wedlock.

The constitution of bodies formed by disjointed molecules separated by spaces devoid of ponderable matter also offers a special application of the law of large numbers, Draw a straight line in a certain direction from a point in the interior of a body. The distance to the point where the line encounters the first molecule, although very short in any sense, nevertheless greatly varies with the direction. It can be 10, 20, 100 times longer in one than in another. Around each point the distribution of the molecules can be very irregular and very different from one point to another. Due to internal oscillations of the molecules it even incessantly changes since a body in repose is just an assemblage of molecules undergoing continual vibrations whose amplitudes are imperceptible but comparable to the distances between molecules. And when dividing each imperceptibly small portion of the volume by the number of the molecules contained there, by an extremely large number because of their excessive minuteness, and extracting the cubic root of the quotient, the result will be the *mean interval* between the molecules. When disregarding the unequal compression of its parts produced by its own weight, it is independent from the irregularity of their distribution and constant throughout a homogeneous body having everywhere the same temperature. Similar considerations served as a foundation for calculating molecular forces and calorific radiation in the interior of bodies in my other works[11].

All these diverse examples of the law of large numbers were taken from the category of physical things, and, if necessary, we can still multiply them. And it is not more difficult to cite other examples pertaining to things of the moral category. Among those, we may indicate the income from indirect taxes, constant if not yearly then at least during a few consecutive years. Such, along with others, is the judicial duty which adds almost the same sum to the annual revenue and which depends on the number and the importance of the cases, i. e., on the opposing and variable interests of the people and their greater or lesser aptitude of pleading.

Such are also the incomes from the Lottery of France before it was luckily suppressed and from the games of Paris[12] whose suppression is no less desirable. These games present constant magnitudes of two different kinds: each year, or during each few years, the sum of the stakes is almost the same; and, on the other hand, the banker's gain is palpably [sufficiently] proportional to this sum. This proportionality is a natural effect of randomness which provides the banker favourable outcomes in a constant proportion calculable beforehand according to the rules of the game. The constancy of the sum of the stakes is a fact belonging to the moral category because they depend on the gamblers' number and volition. It is good that those two elements scarcely vary; otherwise, the contractor of the games would have been unable to evaluate in advance the pledged yearly payment to the government by issuing from the profit he was able to receive during his previous lease.

**[6]** I will soon describe below the results of the experience on which I base myself when considering the problem of the probability of judgements and present additional decisive examples of the law of large numbers observed in things of the moral category. It will be seen that under the authority of the same legislation the conviction rate in all France scarcely varied from one year to another. It is therefore sufficient to consider about seven thousand cases, which is the number of judgements returned yearly by the jury panels, for that rate to become sufficiently permanent. In other problems, for example

concerning the mean life (see above), such a number, however, would have been far from sufficient for leading to a constant result. It is also strikingly seen that the influence of general causes varies each time that the legislation changes.

It is therefore impossible to doubt that the law of large numbers suits the moral things depending on human volition, interests, enlightenment and passions on a par with physical things. Actually, not the nature of causes that has to do with this (? - O.S.) problem, but rather the variation of their isolated effects and the number of cases necessary for the irregularity of the observed facts to be balanced out in the mean results.

The magnitude of those numbers can not be assigned in advance, it differs in different problems, and, as stated above, they are the larger, in general, the greater is the amplitude of those irregularities. And on this point it should not be thought that the effects of spontaneous volition, infatuation with passion, lack of enlightenment vary greater than human life with some babies dying at birth and other becoming centenarians. It is more difficult to foresee them than the circumstances leading to the loss of a ship during her long voyage; they are more capricious than the chance that determines a certain card or a certain number of points on a die. No, we do not attach these ideas to those effects and their causes. Only the calculus and observations can establish probable limits of their variation in a very large number of trials.

It follows from those examples of all kinds that for us the universal law of large numbers is already a general and incontestable fact resulting [derived] from experience that will never be contradicted. And in addition it is the foundation of all the applications of the calculus of probability, and we understand now their independence from the nature of the pertinent problems, and their perfect similitude, whether when concerning physical or moral things, if only the special data required by the calculus is provided by observations. However, because of the importance of the law of large numbers it is necessary to demonstrate it directly, and this is what I have attempted to do, and I believe that I have finally succeeded, as seen farther in this contribution.

The Jakob Bernoulli theorem cited above coincides with that law in the particular case in which the chances of the events remain constant during the series of trials. The author, who is known to have pondered for twenty years, essentially presumed this requirement in his demonstration. His theorem does not therefore suffice for considering problems concerning the repetition of moral things or physical phenomena whose chances are in general continuously varying, and most often without any regularity. For supplementing that theorem it was necessary to study the problem more generally and more completely than the state of mathematical analysis then permitted.

**[7]** The constancy of ratios between the number of times that an event had occurred and the very large number of trials, which establishes itself and is manifested in spite of the variations of the chance of that event during these trials, tempts us to attribute this remarkable regularity to some ceaselessly acting occult cause.

However, the theory of probability determines that the constancy of those ratios is a natural state of things belonging to physical and moral categories and maintains all by itself without any aid by some alien cause. On the contrary, it can only be hindered or disturbed by an intervention of a similar [alien] cause.

The government has published the *Comptes généraux de l'administration de la justice criminelle* for the nine years 1825 – 1833. It is from this authentic compilation that I have drawn all the documents made use of[13]. The number of yearly cases judged by the assize courts of the kingdom was very near to 5000 with about 7000 accused. From 1825 to 1830 inclusive the criminal legislation did not change and convictions by jury panels had been returned by majority verdicts of at least 7 votes against 5 but not when the cases in which that majority verdict was the least possible and the court had to intervene. In 1831, such interventions were suppressed and the required least majority verdict of 8 votes against 4 was introduced.

Acquittals became more frequent; neglecting the third significant digit, the acquittal rate during the six first years had been equal to 0.39. In one year only it lowered to 0.38 and another year it heightened to 0.40. This means that during those years it only deviated from the mean value by 0.01[14]. For the period of legislation valid until 1831 the value of that rate can be therefore taken as 0.39 with 0.61 being the conviction rate. At the same time the ratio of the number of convictions pronounced by the minimal majority verdict of 7 votes against 5 to the total number of the accused equalled 0.07 and it also scarcely varied from year to year. When subtracting this value from 0.61, the rate of convictions returned by a majority verdict of more than 7 votes against less than 5 becomes equal to 0.54. The acquittal rate was therefore 0.46.

And so it actually occurred because during 1831 the difference between this rate derived for the previous years and observed in 1831 was barely 0.005. In 1832 the minimal majority verdict of 1831 was preserved and the law prescribed that *mitigating circumstances*[15] be considered so that in affirmative cases the penalties were reduced. The effect of this novelty should have made conviction by jury panels easier, but by how much? Only experience could have allowed us to evaluate it. It could not have been calculated in advance, as was the increase in the number of acquittals which took place due to changing the least possible majority verdict.

Experience showed that in 1832 the acquittal rate had lowered to 0.41. During 1833, under the same legislation, it only changed by 0.001. The conviction rate in 1831 and later was 0.61, 0.54, 0.59, so that after lowering by 0.7 [0.07] because the minimal majority verdict was increased by one vote, it increased only by 0.5 [0.05] under the influence of the *mitigating circumstances* on the mind of the jurymen, see the *Comptes* for 1834[16].

During 1832 and 1833, the number of political cases submitted to the assize courts was considerable. When estimating the acquittal rate at 0.41, those cases were subtracted from the total number of criminal cases. If these are taken into account, that rate becomes equal to 0.43 which already proved the influence of the kind of cases on the number of acquittals pronounced by jury panels.

This influence is distinctly evident in the *Comptes généraux*. Criminal cases are there separated into two main categories: crimes against property (theft, robbery, or other encroachments), and outrages against the person numbering about 1/3 of the former, or 1/4 of the total number of cases. From 1825 to 1830 the rate of acquittals only amounted to 0.34 in the first category and increased to 0.52 in the latter. In other words, it exceeded the conviction rate by 0.04 [0.52 − 0.48 = 0.04]. The yearly values of each of these rates only varied not more than by 0.02 in either direction of those fractions, 0.34 and 0.52. It should also be noted that the number of convictions pronounced by the minimal majority verdict of 7 votes against 5 only amounted to 0.05 of the number of accused of crimes against property but rose to 0.11 in the other category. And, not only were the convictions proportionally more numerous in the first case, they were in general returned by a stronger majority.

Those differences can partly depend on a less severity of jurymen in cases of outrages against the person than encroachments on property which are doubtless considered more dangerous for the society because of being more frequent. However, this cause is not sufficient for producing the large inequality in the rates of acquittals as provided by experience. And calculation proves that such inequalities also originate from the greater presumption of guilt of the accused of theft or robbery than of the other accused which resulted from the information gained prior to the judgement.

**[8]** The *Comptes généraux* provide other ratios which, owing to large numbers, are almost invariable, but which I do not make use of at all. Thus, from 1826, when the sex of the accused began to be shown there, and until 1833 there were among the accused almost exactly 18/100 women yearly. Once only their share rose to ca. 0.20, and only once it lowered to 0.16. And it was always higher when concerned with theft than in crimes against the person. The acquittal rate was considerably higher for women than for men and reached almost 0.43 whereas it only was 0.39 for both sexes together.

However, the constancy of these diverse proportions observed each year for the entire France did not take place for separate assize courts. In the same départements and under the same legislation the acquittal rate remarkably varies from one year to another. This proves that for such courts the yearly number of criminal cases is not at all large enough for balancing out the irregularities of the judges' votes or for the acquittal rate to become permanent. That rate varies still more from one département to another and the number of cases in each assize court is not considerable enough for deciding with sufficient probability in which part of France the judges (jurys) tend more or less to severity.

Only in the département de la Seine are criminal cases sufficiently numerous for the observed acquittal rate not to be too variable and therefore allow to compare it with that rate for France as a whole. About 800 accused, or about 1/9 of the total for the entire kingdom, are yearly brought to the Paris assize court. From 1825 to 1830 the acquittal rate there varied between 0.33 and 0.40 and its mean value was only 0.35 whereas for France it was 0.39, or 0.04 larger. The rate of convictions returned by the least majority verdict of 7 votes against 5 is also a little less for Paris; it only reached 0.065 instead of 0.07 for the whole of France when calculated without distinguishing the kind of crimes.

Such are the data concerning the assize courts provided until now by experience. The precise aim of the theory is to calculate for a jury panel composed of a certain number of people and judging a very large number of cases by an also known majority verdict the future rates of acquittals and convictions likely to take place and the chance of mistaken judgement for those who were or will be randomly selected as jurymen.

I believe that determining that chance for convictions and acquittals at a given and isolated case is impossible if not basing the calculus on precarious presumptions, leading according to the adopted hypotheses almost to the desired results. However, for the society and the accused it is not the chance for a particular judgement that is most important, but the chance which concerns all the cases submitted to the assize courts during a year or many years, and which is determined by observation and calculation.

The probability of an error of some convicting judgement multiplied by the chance that it occurs is the veritable measure of the danger to which the society exposes an innocent accused. The product of the chance of an error in acquitting and the probability that such a judgement is pronounced is the measure of the danger to which the society itself is subjected and which should be known just as well. Indeed, it is the magnitude alone of this danger that can substantiate a possible unjust conviction.

In this important problem of humaneness and public order nothing can replace analytical formulas expressing these various probabilities. Without their aid, if a change of the number of jurymen is considered, or two countries where it differs are compared, how can we know whether a jury panel of 12 judging by a majority verdict of 8 votes against 4 offers more or less guarantee to the accused and the society than another panel of, for example, 9 jurymen selected from the same list as the former panel was, and judging by some majority verdict? How was it decided whether the situation existing in France before 1831 of a least majority verdict by 7 votes against 5 supplemented with an intervention of judges in cases of that minimal majority had been more advantageous or, on the contrary, less favourable than what we have today, the same majority and the influence of *mitigating circumstances*? We can not know it since the observational data pertaining to our time is lacking.

**[9]** The formulas which define this aim and which are included in this work are derived without introducing any hypotheses by issuing from the general and known laws of the calculus of probability. They include two special magnitudes which depend on the moral state of the nation, the methods of the adopted criminal proceedings and the ability of those charged with directing them. One of those formulas expresses the probability that a juryman randomly selected from a list covering

the jurisdiction of an assize court does not err in his vote. The other one is the probability that the accused is guilty existing before the beginning of the pleading.

These are the two essential elements in the problem of criminal judgements. Their numerical values should be determined by the observational data just as the constants included in astronomical formulas are deduced from observations. The entire solution of the problem proposed in these researches requires an interaction of theory and experience.

The observational figures which I made use of, two in number, just equal to the number of elements to be determined, are the number of convictions by least majority verdicts of 7 votes against 5, and of the accused convicted by that least majority, each figure divided by the total number of the accused. The obtained ratios are very different for crimes against the person and encroachments on property, and I consider them separately. They are not the same in all the départements, but the need to deduce them from very large numbers compelled me to unite, again separately, the judgements of all the assize courts of the kingdom. The values that I derived, as though these elements did not vary much from one département to another, are therefore only approximate.

However, the new law re-established the least majority convictive verdicts of 7 votes against 5 and stipulated that the jury panels should make known the convictions returned with that least majority. For each département we know therefore the considerable enough number of convictions decided by that, and, separately, by any other majority which is necessary for deducing our two elements. And we will find out whether the chance of the juryman's error remarkably varies from place to place. For the département de la Seine taken all by itself calculation had already established that that chance is a bit less than for the rest of France. Following are the main numerical results contained in this work; it seemed useful to present here their summary.

Before 1831, for the entire France the probability of a juryman's faultless vote was a bit higher than 2/3 for crimes against the person, and almost equal to 13/17 for crimes against property. Without distinguishing between these categories of crime, that chance was very little lower than 3/4, once more for the whole kingdom, and a bit higher for the département de la Seine. At the same time, the other element of criminal judgement, the preliminary probability of the defendant's guilt for the entire France and crimes against the person did not much exceed 1/2, was constant and contained between 0.53 and 0.54. It somewhat exceeded 2/3 for crimes against property; without distinguishing between these categories of crime it was almost equal to 0.64, but amounted to 0.68 for the jurisdiction of the Paris assize court. When subtracting these diverse fractions from unity we will obtain the probabilities corresponding to the juryman's error and to an erroneous conviction.

**[10]** It can be remarked that the preliminary probability of the defendant's guilt was always higher than the conviction rate. For example, in the case in which that probability was the lowest, and only exceeded 1/2 by 0.03 or 0.04, that rate, as stated above, was about 0.02

lower than 1/2. This is a general result, and the formulas of probability show that it always takes place whichever is the ability of the jurymen, the chance of their error, and the required majority convicting verdict. It should also be noted that that prior estimate only expresses the probability that the accused will be convicted by a jury panel according to their method of judging; that is, according to the unknown level of probability which they require for conviction and which is doubtless lower than the probability that an accused is really guilty as results from preliminary information.

Actually, no one will hesitate to stake more than even money that an accused is guilty when brought to an assize court on charges of crime against the person, and this in spite of the prior probability only a bit higher than 1/2 as found out for those crimes.

In 1831, only the majority required for conviction was changed and the two elements which we consider should have remained the same. In the following years the problem of *mitigating circumstances* doubtless influenced the values of those elements. However, for 1832 and 1833 we only know the conviction rate which is insufficient for determining the two elements, and we do not know whether the chance of a juryman's faultlessness had increased or decreased. We are unable to decide this without imposing on the other element a hypothesis and risking to deviate considerably from the truth. Just the same, the present legislation had imposed a secret vote on the jurymen[17], and we do not know whether this chance of faultlessness changed once more.

When it, as well as the chance of guilt, resulting from the information gathered before the judgement, will be determined by a sufficient number of observations, we will also know, by repeating our calculations for ever more distant times, whether these two elements have been varying in France progressively in one or another sense, which will provide an important document about the moral state of our country.

Although the judges are doubtless more experienced in criminal cases, their chance of faultlessly voting seem to differ little from that of the jurymen. Actually, from 1826 to 1830 in the entire France convicting verdicts by majority verdicts of 7 votes against 5 were returned 1911 times. The assize courts, then consisting of 5 judges, had appropriately been called for intervening, and in 314 cases accepted the minority vote. Assuming that the probability of faultlessness was the same for judges and jurymen, calculation will prove that that should have happened ca. 282 times. These two numbers, 314 and 282, are not sufficiently considerable for deciding with a very high probability just how did this assumption deviate from verity.

The small difference between them is a cause for believing that there should also be a very small difference between the chances of error by judges and jurymen. For the jurymen, that chance does not therefore originate, as it is possible to think, from their want of skill. All other things being equal, it is evident that the conviction rate diminishes when a greater majority convicting verdict is demanded.

If, as in England, a unanimous vote of the 12 jurymen is required both for convictions and acquittals, and if adopting for England the values of the two elements of criminal justice pertaining to the entire France without distinguishing between the categories of crime, the probability of such convicting verdicts will little differ from 1/50 and be about a half lower for acquittals. This makes decisions very difficult, at least if some arrangement between the jurymen is not reached with a part of them sacrificing their opinion. It is also seen that otherwise unanimous acquittals will be twice more difficult and twice rarer that convictions. It is only possible to bet even money on a unanimous conviction or acquittal happening once in 22 randomly selected cases[18].

**[11]** After the judgement is pronounced, the probability of the defendant's guilt becomes much higher or much lower when he is convicted or acquitted respectively. The formulas in this work, if the two included there elements are determined by observation, provide the value [of change] according to the majority by which the judgement was returned. When the majority of a convicting verdict is at least 8 votes against 4, the probability of guilt of a convicted accused is a bit higher than 0.98 in cases of crimes against the person and a bit higher than 0.998 in cases of crimes against property. It reduces the chances of an erroneous judgement to a little less than 0.02 and 0.002 respectively.

When taking into account the probability of not being acquitted, the chance of an erroneous conviction for crimes of the first category is about 1/150 and only 4/10,000 for the other crimes[19]. At the same time the probabilities of innocence of an acquitted accused are approximately equal to 0.72 and 0.82 respectively. When taking into account the probability of not being convicted, we can also establish that the chances of a guilty accused to be acquitted in these two cases are about 0.18 and 0.07 respectively. Therefore, among a very large number of those acquitted there will be more than 1/6 and about 1/14 who should have been convicted.

During the seven years, 1825 – 1831, in the entire France the accused, convicted by that majority verdict of at least 8 votes against 4, numbered about 6000 for crimes against the person and about 22,000 for crimes against property. According to the chances of an erroneous conviction cited above, it can be thought that about 40 and 9 were innocent[20] which amounts to 7 people yearly [for both categories of crime taken together]. At the same time, the number of acquitted but guilty accused should have been 50 times larger[21], or about 360 each year. However, we should not lose sight of the sense attached, as explained above, to the word *guilty*. According to that sense, the number 18 (? - O.S.) is only the superior limit for the really innocent but convicted, whereas 360 is, on the contrary, the inferior limit of the number of those acquitted who were not at all innocent.

These results of calculations, far from injuring our proper respect for a judged thing or from diminishing our confidence in the decisions of jury panels, prevent, on the contrary, any exaggeration of the error to be feared in convictions. In truth, it is not in essence possible to verify them by experience; however, this circumstance is common to many other applications of mathematics whose certainty only rests, just as it does here, upon the rigour of demonstration and precision of the observational data.

For the years preceding 1831 and for the entire France the probability of an erroneous conviction pronounced by the least majority verdict of 7 votes against 5 was about 0.16 and 0.04 for crimes against the person and against property respectively. Without distinguishing these categories, it was 0.06. According to the Laplace formula, that chance should have been the same for all cases and about five times greater than 0.06. However, it should also be noted that the intervention of a court was then necessary in cases of least majority verdicts. If decisions of jury panels were confirmed, it resulted in reducing that chance of error, 0.06, a little less than by 0.01. This means that [yearly] during 1826 – 1830 about 15 or 16 out of the 1597 convictions by the majority verdicts mentioned can be supposed erroneous; the accused should not have been convicted which does not imply their innocence.

**[12]** The distinctive character of this new theory of probability of criminal judgements therefore consisted of ascertaining first of all, by issuing from a very large number of cases of the same nature, the chance of error in the jurymen's voting and of the guilt of the accused existing before the beginning of pleadings.

It should be applicable to all the numerous kinds of judgements, to those made in police courts, by military justice and to the decisions in civil courts, if only there always exist sufficient data for determining those two elements. It should also be appropriate to judgements pronounced in a very large number by extraordinary tribunals during the ill-fated years of the revolution. On this point, however, it is necessary to enter into some explanation so that there will remain no doubt about the generality and exactitude of the theory. The difficulty that that exceptional case presents will not escape those who wish to hear attentively and are interested in the results of my work.

An accused can be convicted either because he is guilty and the judges did not err, or since he is innocent and the judges erred. The conviction rate does not vary if the prior probability of the defendant's guilt and the probability of each judge's faulty vote are changed into their complements to unity. It remains the same if, for example, those two probabilities are 2/3 and 3/4 or 1/3 and 1/4. It also retains its value when both probabilities are almost certain or very near to unity, and when they are almost zero. In these extreme cases the number of convictions will very little deviate from the number of accusations [of the accused].

The equation to be solved for determining those probabilities have two real roots making up unity, but each of those solutions has a distinguishing feature: when adopting one of them, the probability of guilt of a convicted accused will be higher than that of his innocence and lower otherwise. In ordinary cases we should therefore choose the first solution since it is unreasonable to suppose that the tribunals are in general unjust and most often judge contrary to common sense. This, however, is not so when the judgements are returned under the influence of passions; it is not the reasonable root of those equations, it is the other solution which should be chosen and which attaches such a high probability to unjust convictions. A great proportion of convictions pronounced by the revolutionary tribunals had not been sufficiently justified by proving the legal guilt of the accused. Issuing from the laws which the tribunals had to apply, we can not at all find out how many convicted accused were guilty or innocent.

It should always be borne in mind that, according to that [new probability] theory, the injustice of judges and passions of prosecutors, just as great pity or extreme indulgence are considered as chances of error and that calculations are based on the results of votes whichever motives dictated them.

In the police courts, the mean conviction rate for nine consecutive years and for the entire France was contained between 0.86 and 0.85. This indication is not, however, sufficient for determining the probability of the accused guilt existing before the judgement and the probability of a faultless vote for a judge of those tribunals. Presuming that the judgements were returned by 3 judges, which seems to have generally taken place, we should also know the proportion of convictions decided unanimously or by simple majority verdicts of 2 votes against 1. This proportion is not given by observations, and can only be supplied by some unjustified hypothesis.

In the case of military tribunals we also lack the two necessary data for determining the values of the two special elements included in the probability formulas. A court martial consisted of 7 judges and convictions could have only be decided by a majority verdict of at least 5 votes against 2. The overall conviction rate is estimated as 2/3 but we do not know the proportion of convictions returned unanimously or by a simple majority. Lacking this indication, we can not precisely compare military justice and the assize courts with respect to the chance of erroneous convictions and acquittals which, however, would have been very interesting.

**[13]** When describing civil matters, the probability formulas contain only one special magnitude rather than two of them, − that which expresses the chance of a judge's faultless vote. According to the information I had been given, judgements in tribunals of the first instance are generally returned by 3 judges. We do not know, however, the ratio of the number of cases in which decisions had been made unanimously and by a simple majority ruling of 2 votes against 1, and it is therefore impossible to determine directly the chance of an erroneous voting.

We are able to calculate that chance for the judgements appealed to the royal courts by comparing the number of confirmed and not confirmed cases and assuming an equal chance of error for the judges of those two tribunals. Although this hypothesis deviates perhaps considerably from the truth, I admit it for providing an example of calculating the error to be feared in civil judgements. The truth or justice results from a decision, necessarily unanimous, of judges having no chance of error. In each case, that *absolute justice* is an unknown thing; nevertheless, we reveal those cases which are decided contrary to it by the erroneous votes and judgements. The problem consists in determining their probabilities and, consequently, the proportion according to which they likely and almost exactly occur in a sufficiently large number of cases.

The *Compte général de l'administration de la justice civile* published by the government cites the number of judgements in courts of first instance confirmed or disaffirmed by the royal courts during the three last months of 1831 and the two following years, 1832 and 1833. For the entire France, the ratio of the first of those two numbers to their sum is almost equal to 0.68 and only varied from year to year by 1/70 of its value.

In spite of the diversity of cases that should have existed and doubtless of the kingdom's magistrates unequal enlightenment, about 8000 yearly decisions were sufficient for that ratio to attain an almost constant value. This presents one more very remarkable example of the universal law of large numbers. In the jurisdiction of the Paris royal court that ratio was considerably larger amounting to about 0.76.

When applying a value for France as a whole, and assuming that in each royal court of appeal there were 7 counsellors pronouncing decisions in civil matters, we find that 0.68 is the probability that one of them or of the judges in the courts of first instance randomly selected from those in all the kingdom will not err when opining about a case, again randomly selected among those yearly submitted to two degrees of jurisdiction. That probability is possibly different for the cases judged in courts of first instance without being appealed by either of the two sides.

Issuing from that fraction, 0.68, we find, when neglecting the thousandths, 0.76 for the probability of virtue of a judgement in courts of first instance. In courts of appeal it is 0.95 when its judgement confirms the decision of a court of first instance and 0.64 otherwise. Finally, that probability is 0.75 for a decision of a royal court, whether confirming or disaffirming the judgement in a court of first instance, confirmed by a second royal court issuing from the same information as the first royal court did.

The approximate probabilities that the court of first instance and first court of appeal judged properly; that the former judged properly and the latter inappropriately; that it was otherwise for both courts; and that both courts judged improperly, are 0.649, 0.203, 0.113, and 0.035 whose sum is unity.

The problems concerning the probability of judgements whose basis I describe and whose results I provide are in the fifth and last chapter of this work. The four preceding chapters include the rules and general formulas of the calculus of probability which disperses with searching for then elsewhere and enables the treatment of other problems alien to the special aim of these researches but proper for that calculus to explain. There also is the solution of a problem which proves how the majority of an elected assembly can be completely changed after a new election or to a larger extent than the change in the vote of those distributed among electoral colleges and voting by simple majority in each of them.

### Notes

**1.** The *man about town* was De Méré (1610 – 1685), and *an austere Jansenist*, Pascal. Jansen (1585 – 1638) was a theologian.

**2.** More precisely: on the unequal number of points still lacking them.

   **3.** In the 17th century, political arithmetic and population statistics in particular as well as insurance of life were born. In addition, Huygens had studied mortality, but his pertinent manuscripts were only published by the end of the 19th century.

   **4.** Laplace considered that problem in Chapter 11 of, and the First Supplement to his *Théorie*. Unlike Laplace, Poisson had as a rule applied the term *calculus of probability*, especially in the main text of his work, and Bertrand, Poincaré and Markov followed suit. However, in his mimeographed lectures the last-mentioned had until 1900 written *theory of probability*.

   **5.** This is wrong. The Bayes approach became the subject of debates and is still possibly argued about.

   **6.** From 1774, Laplace (e. g., 1812/1886, p. 365) separated himself from the *geometers* and actually sided with applied mathematicians.

   **7.** Assizes: it would have been important to indicate who voted there, judges, jurymen, or, as I found out about later times, both. However, Poisson did not clearly supply such information, and it appears that in different places he had in mind differing methods of voting. Then, in Chapter 5 he mentioned several kinds of courts dealing with civil cases without sufficiently explaining how those cases had passed from one of them to another.

   **8.** Poisson invariably wrote *experience* and never mentioned statistics.

   **9.** Poincaré is known to have repeated that interpretation of randomness. For him, however, the main pattern of the action of chance was small causes leading to large consequences.

   **10.** Separate mortality tables for men and women began to appear even before 1832 (Quetelet & Smit 1832, p. 33) whereas Corbaux (1833) followed this new practice and even somehow distinguished between several strata of population.

   **11.** These *other works* were hardly studied in the stochastic sense.

   **12.** At the end of 1837 they were suppressed. At that time, there were seven gambling houses in Paris (*La Grande Enc.*, t. 21, p. 152), but I did not find any mention of the *jeux de Paris*.

   **13.** Poisson actually made use of other sources as well, see his Chapter 5.

   **14.** Poisson many times applied the same method of estimation.

   **15.** Concerning these *circumstances*, see Chapter 5. O. S.

We see there that during that year, when legislation remained as it was in the two preceding years, the conviction rate amounted to 0.60 and thus only exceeded the rate for those years by 0.01. After the example of our country, the Belgian government publishes its own *Compte général de l'administration de la justice criminelle*. Jury panels were re-established there about mid-1831 and the required majority convicting verdict was at least 7 votes against 5. For the years 1832, 1833 and 1834, the rate of acquittals was 0.41, 0.40 and 0.39. The mean rate, 0.40, remarkably only differed by 0.01 from its value in France for the same majority verdicts.

Before jury panels were re-established, criminal tribunals in Belgium consisted of 5 judges and convictions were returned by simple majority of 3 votes against 2. From year to year, the rate of acquittals also varied very little, but it only reached about 0.17, less than a half of the value taking place for judgements by jury panels. That difference of more than twice did not only result from the difference, either of the numbers 5 and 12 (judges and jurymen) or of the minimal votes in majority verdicts (3 against 2 or 7 against 5). It also suggests, as seen here below, that for convictions judges require a considerably lower probability [of guilt] than the jurymen whichever the chances of their error are. Poisson

   **16.** Such references are not definite enough.

   **17.** Although, after participating in a secret ballot, jurymen can not go back on their decision, there is a special case which can sometimes occur and is appropriate to be indicated. Two men, call them Pierre and Paul, are accused of theft. When asked whether Pierre was guilty of that theft, 4 jurymen answer *yes*, 3 others also say *yes*, and the 5 others, *no*. The accused is found guilty by 7 votes against 5. When asked about Paul, the first 4 jurymen answer *yes*, the same 3 of them say *no*, and the 5 last ones, *yes*. Paul is declared guilty by 9 votes against 3.

Then comes the next question, whether the theft was done *by many*; when answered in the affirmative, it leads to a more severe punishment. [Necessarily]

keeping to their preceding votes, the first 4 jurymen answer *yes*, and the 8 others, who decided that either Pierre or Paul were innocent, say *no*. Without contradicting those votes, the decision of the panel will be that both the accused are guilty of the theft, but that it was not perpetrated *by many*. Poisson

**18.** According to documents published in England, and seem to earn confidence, the number of people yearly tried by jury panels had recently increased from year to year and the conviction rate had also increased progressively (Porter). Here are the results extracted from those documents; they can be compared with the situation in our country. The numbers only pertain to England and Wales.

Three periods are tabulated, each of them lasting 7 years and ending in 1818, 1825 and 1832. During the first period which ended in 1817 (? - O.S.) the [yearly] number of the accused did not amount to 35,000 and the conviction rate was a bit lower than 0.60 (? - O.S.). Only once, in 1832, the last year of those periods, did the number of the accused reach 20,829 of whom 14,947, or about 3/4 were convicted. I do not know whether this number had increased or decreased in the following years.

In England and France, the proportions of the feeblest punishments little differ. The appended Table shows that almost 2/3 of the total number of convictions consisted in imprisonment. In France, the fist number exceeded a half of the second one. The Table also proved that during the last of the three periods the mean yearly number of the executed, the least of them, amounted to 60. In France, it is now twice less; the yearly number is not larger than 30.

**Table**
showing the appropriate numbers in periods 1 – 3

1. Accused 64,538; 93,718; 127,910
2. Convicted 41,054; 63,418; 90,240
3. Conviction rate 0.636; 0.677; 0.705
4. Sentences to death 5802; 7770; 9729
5. Executed 636; 579; 414
6. Sentenced to imprisonment for two years or less 27,168; 42,713; 58,757

**19.** Someone inserted the word *dix* before *millièmes*.
**20.** Someone crossed out 88 and wrote 9.
**21.** Someone crossed out *vingt* and wrote *plus de cinquante*.

## Chapter 1. General Rules of Probabilities

**Misprints/Mistakes Unnoticed by the Author**

1. In § 8, p. 41 of the original text, the signs of the terms in both series for $(1 – p)^n$ are wrong.

**2.** At the very end of § 9, p. 43 of the original text, Poisson indicated in words that some magnitude was a bit less than $1/60·10^6$; instead of $1/60·10^7$; a similar mistake was made a few lines below.

**3.** In § 17, p. 58 of the original text, line 7 from bottom. Urn U should be urn A.

**4.** In § 21, p. 66. The right side of the second equality on last line should be *b/m* rather than *a/m*.

**1.** The *probability* of an event is our reason to believe that it will occur or occurred. Suppose that it happened in one case, and that in another one it is only probable. For us, however, all other things being equal, its probability in those cases, so different by themselves, is the same. A ball is extracted from an urn containing some white and black balls whose numbers are known. Or [in another case] the colour of the extracted ball remains unknown to me. Evidently I have the same reason to believe that in both cases that ball is white.

Probability depends on our knowledge about an event; for the same event it can differ for different persons. Thus, if a person only knows that an urn contains white and black balls, whereas another person also knows that there are more white balls than black ones, the latter has more grounds to believe in the extraction of a white ball. In other words, for him, that event has a higher probability than for the former. It is for this reason that two persons, A and B, having different knowledge about the same event, sometimes judge contrary to each other about it. If A knows everything known to B and something else, his judgement is more competent. It is therefore reasonable to adopt his opinion when having to choose between the contrary judgements of A and B although it is possibly based on a lower probability than that which motivates B's opinion. This means: even if A is less justified to believe in his own opinion than B is with respect to his.

In ordinary life, the words *chance* and *probability* are almost synonymous and most often used indifferently. However, if necessary to distinguish their meaning, we attach here the word *chance* to events taken independently from our knowledge, and retain its previous definition [!] for the word *probability*. Thus, by its nature an event has a greater or lesser chance, known or unknown, whereas its probability is relative to our knowledge about it[1].

For example, in the game of *heads* or *tails*, the arrival of *heads* results from the constitution of the tossed coin. It can be regarded as physically impossible that the chances of both outcomes are the same; however, if that constitution is unknown to us, and we did not yet try out the coin, the probability of *heads* is for us absolutely the same as that of *tails*. Actually, we have no reason to believe in one of these events rather than in the other one. This will not be the same after many tosses of the coin: the chance of each side does not change during the trials, but for someone who knows their results, the probability of the future occurrence of *heads* and *tails* varies in accord with the number of times they happened.

   **[2]** The measure of the probability of an event is the ratio of the number of cases favourable for it to the total number of favourable and contrary cases, all of them equally possible or having the same chance². That proposition signifies that when this ratio is the same for two events, we have the same reason to believe in the occurrence of either of them. Otherwise, we have more reason to believe in the arrival of that event for which it is larger.

   Suppose for example that an urn contains 4 white balls and 6 black ones, and that another urn contains 10 and 15 respectively. The numbers of favourable cases for the occurrence of a white ball and of all the possible cases are 4 and 10 for the first urn, and 10 and 25 for the second one. For each urn the ratio of the first number to the second is 2/5; it should be proved first of all, that the probability of extracting a white ball is the same for both urns which means that if we are somewhat interested in the arrival of a white ball, we will have absolutely no reason to choose rather one urn than the other one. [A long proof follows.]

   Now I will suppose that urn A contains 4 white balls and 3 black ones, and urn B, 3 white and 2 black balls. The same ratios are 4/7 and 3/5. The second fraction exceeds the first one by 1/35 so that there is a greater reason to believe that a white ball will be extracted from B rather than from A. When reducing both these fractions to a common denominator, they will become 20/35 and 21/35, and, according to the proven statement, the probability of extracting a white ball is the same as when extracting it from an urn C containing 35 balls, 21 of them white, and 14 black. Just the same can be stated about urns B and D the latter containing 35 balls, 21 of them white and 15 black. The urns C and D contain the same number of balls but D has more white balls than C and there evidently is a greater reason to believe that the white ball was extracted from B rather than C, and, just the same, from D rather than A. And this concludes the demonstration of the proposition stated above.

   It seems that the adopted measure of probability always leads to commensurable fractions³. However, if the numbers of all the possible cases and of those favourable to an event are infinite, the ratio of the second number to the first one can be an incommensurable magnitude. Suppose, for example, that $s$ and $\sigma$ are the extents of a plane surface and of one of its parts. When tossing a round coin whose centre will equally fall on any point of $s$, the probability of its landing on some point of $\sigma$ will evidently be equal to the ratio of $\sigma$ to $s$ whose values can be incommensurable.

   **[3]** In the preceding demonstration we selected some definite numbers of balls, but it is easy to see that the reasoning was general and independent from those particular numbers. It was also supposed that the event, whose probability had been considered, was an extraction of a white ball from an urn containing white and black balls so that their numbers represented favourable and contrary cases for that event. For simplifying the reasoning, it is always possible to

substitute [to adopt] a similar hypothesis in each problem concerning eventuality and things of quite another nature.

Thus, let E be an event of some kind; represent the number of cases favourable to its occurrence by $a$, of contrary cases by $b$, and by $p$, the probability of E. According to the demonstrated, the measure or numerical value of that last-mentioned magnitude will be

$p = a/(a + b)$.

At the same time, if F is an event contrary to E, then only one of these two events will necessarily arrive, just as in the preceding examples either a white, or a black ball will be extracted. And, denoting the probability of F by $q$, we will also have

$q = b/(a + b)$

since the $b$ cases contrary to E are favourable for F, and $p + q = 1$. The sum of the probabilities of two contrary events, as defined above, is always equal to unity.

If we have no greater reason to believe in the occurrence of E rather than F, their probabilities will be equal, $p = q = 1/2$. This[4] takes place when tossing a coin, whose physical constitution is unknown to us, for the first time. E and F are here the occurrences of one or another of its two sides. Instead of being an event which can arrive or not, E can be a thing whose veracity or falsity we wish to know. Then, $a$ and $b$ will be the numbers of cases in which we believe it to be true or false, and $p$ and $q$ will express the probabilities of these alternatives.

In each example of either eventuality or doubt and criticism, when evaluating the numbers of favourable and contrary cases for E and F, and being certain that these numbers are $a$ and $b$, the fractions $p$ and $q$ will be the chances of E and F. If, however, that evaluation only results from our knowledge about these two things, $p$ and $q$ will only be their probabilities which can differ, as we have explained, from their unknown chances. It is always necessary for the favourable and contrary cases to be equally possible either by themselves or according to our knowledge.

**[4]** In the theory of chances [!] *certitude* is considered as a particular case of probability, the case of an event having no contrary chances. In calculations, it is represented by unity whereas probabilities are fractions less than unity. Complete *perplexity* of our mind when selecting one of two contrary things is represented by 1/2, and *impossibility*, by 0. For us, that notion of certainty is sufficient; we do not need to define it by itself and in an absolute way which is also impossible.

Absolute certainty belongs to things about which we are only able to provide examples. Among those which are called *certain* there are very few rigorously such, like our own existence. Then, some axioms are not only certain, but evident; and propositions like geometrical theorems whose veracity we demonstrate, or whose contrariety we disprove. Things not contrary to the general laws of nature and attested by numerous testimonies, and those, confirmed by everyday experience, are nevertheless only very likely, but so probable that there is no need, either in usual life or even in physical or historical sciences, to distinguish their probability from certainty.

The aim of the calculus of probability is to determine, in each problem of eventuality or doubt, the ratio of the number of cases favourable for the occurrence of an event or for the veracity of a thing, to the total number of all possible cases[5]. We will then be able to know precisely, according to that fraction more or less approaching unity, our reason to believe that that thing is real, or that that event occurred or will occur, and it will also be possible for us without any illusions to compare that reason in any of the two completely different problems of nature.

It [the calculus of probability] is based on a small number of rules which we will provide and which are demonstrated in full rigour as was shown in an example concerning a proposition from § 2. Those principles (? - O.S.) should be regarded as a necessary supplement of logic[6] since there are so many problems in which the art of reasoning can not lead us to entire certainty. No other branch of mathematics is susceptible of more immediately useful applications. As shown in Chapter 2, it includes abstract and controversial problems of general philosophy and provides their clear and incontestable solutions.

**[5]** If $p$ and $p'$ are the probabilities of two events, E and E′, independent from each other[7], the probability of their concurrence or of an event consisting of those two, is $pp'$. Actually, let the event E be the extraction of a white ball from urn A containing $c$ balls, $a$ of them white, and $b$, black. And let E′ be a similar event concerning another urn, A′, containing $c'$ balls, $a'$ of then white, and $b'$, black. Then, according to the preceding,

$$p = a/c, p' = a'/c'$$

are the probabilities of E and E′. The compound event will then be the arrival of two white balls, one of them extracted from A, the other one, from A′. When randomly extracting a ball from each of those two urns, each ball from A can occur with each from A′ so that there will be $cc'$ equally possible cases. From this total number, those favourable for the compound event result from the combinations of each white ball in A with each of them in A′, and the number of these favourable cases is $aa'$. Therefore, the probability of the compound event will be (§ 2) $aa'/cc'$, or, which is the same, $pp'$.

Just the same, if $p, p', p'', \ldots$ are the probabilities of some number of events E, E′, E″, … independent from each other, the probability of their concurrence will be $pp'p'' \ldots$ This general proposition can also be derived from the particular case of an event consisting of two others. Indeed, the product of $p$ and $p'$, or $pp'$, is the probability of the concurrence of E and E′; the probability of the concurrence of that compound event and E″ will be […] etc.

All those fractions, $p, p', p'' \ldots$ are less than unity, at least when no event E, or E′, or E″, … is certain. It follows that the probability of a compound event is also lower than that of each of those events. It lowers ever more as their number increases; generally, it tends to zero and is exactly equal to it or becomes infinitely low if that number becomes infinite. The only exception is the case in which the infinite series of the probabilities $p, p', p'', \ldots$ is composed of terms infinitely approaching unity, or certitude. Their product can have a finite value less than unity. Thus, denote a positive magnitude $\alpha \leq 1$ and let

$p = \alpha, p' = 1 - \alpha^2, p'' = 1 - \alpha^2/4, p''' = 1 - \alpha^2/9, \ldots$

According to a known formula, their product, or the probability of the compound event, will be equal to $(1/\pi)\sin\alpha\pi$; as usual, $\pi$ is the ratio of a circumference to its diameter.

**[6]** Here is a problem about a probability of a compound event providing an example of the preceding rule. I suppose that one of the two randomly chosen numbers is subtracted from the other one. It is required to determine the probability that that operation is possible without having to increase any digit of the minuend.

Each of the corresponding digits of the minuend and subtrahend can take 10 different values from 0 to 9, and it follows that in each partial subtraction there are 100 distinct and equally possible cases. For being able to accomplish our task, the minuend's digits should exceed the subtrahend's corresponding digits or be equal to them.

This takes place in 55 of those 100 cases: only once if the minuend's digit is 0, twice if that digit is unity, …, and 10 times if that is 9. These numbers form a geometric progression having 10 terms, and its sum is $(1/2)10(1 + 10) = 55$. In each subtraction there is a probability of 55/100 so that the whole operation done without increasing the minuend's digits has probability $0.55^i$ where $i$ is the number of the minuend's (or subtrahend's) digits. When subtracting, for example, the decimal parts [the mantissas] of the logarithms in the Callet's table[8], we have

$i = 7, 0.55^i = 0.0152243 \ldots$[9],

that is, a probability contained between 1/66 and 1/65.

We also derive the same probability when adding two numbers containing $i$ digits each without having to memorize a unity.

**[7]** If one and the same event E occurs $m$ times successively, the product $pp'p'' \ldots$ becomes $p^m$ and expresses the probability that E will arrive $m$ times in as many trials during which its probability $p$ remains constant. Similarly, if the contrary events E and F have probabilities $p$ and $q$, so that $p + q = 1$ (§ 3), and if their chances remain constant during the $m + n$ trials, the product $p^m q^n$ will be the probability of E arriving $m$ times and F, $n$ times at the other $n$ trials, in a determined order. This follows from the rule of § 5 when supposing that the number of events E, E′, E″, … is equal to $m + n$ and that E is substituted for $m$ of them, and F, for the other $n$. The order in which E and F should occur does not influence that probability $p^m q^n$ of the compound event. It remains the same whether E arrives at the first $m$ trials, and F, at the last $n$ of them or vice versa, or even if these events should happen in a determined irregular way.

However, if the order which E and F should follow is not stated, and we only desire that in $m + n$ trials E occurs $m$ times, and F, $n$ times in whichever order, it is evident that the probability of that other compound event will exceed the one corresponding to some determined order. Actually, it will be a multiple of $p^m q^n$ and I will later provide its general expression.

If the chances of E and F are equal, $p = q = 1/2$; and if $m + n = \mu$, the probability of the arrival of these events in a determined order $m$ and $n$ times will become $(1/2)^\mu$. It is thus independent not only from the order of their occurrences, but from the proportion of these as well, and only depends on the total number $\mu$ of the trials. Such is the case of an urn containing an equal number of white and black balls from which $\mu$ successive extractions with replacement are made. The probability of extracting $\mu$ white balls is the same as of having $m$ white and $n$ black balls in a determined order. Both probabilities should be very low, but none lower than the other one. Before the trials we have the same reason to believe in the occurrence of a sequence of balls of the same colour, or of the same number of balls, some of them white and the others black in any arbitrary order.

However, if we see a sequence of, say, 30 balls of the same colour and are quite certain that white and black balls are invariably equal in number; or see a very different event presenting some symmetry such as the appearance of 30 balls, alternately white and black; or 15 black balls following 15 white balls, we are led to believe that these regular events are not due to randomness and that the person who extracted the 30 balls knew the colour of each and in selecting them had pursued a particular aim. In such cases, as will be seen below, an intervention of a non-random cause actually has probability very near to certainty[10].

**[8]** The power $q^n$ is the probability that event F will arrive $n$ times in succession, without interruption. Subtracting it from unity, we will therefore obtain the probability of the contrary event E occurring at least once in $n$ consecutive trials. Therefore, when denoting the probability of that latter compound event by $r$, and substituting $(1 − p)$ instead of $q$, we get

$$r = 1 − (1 − p)^n. \qquad (8.1)$$

When equating this value to 1/2, we will determine the number of trials necessary for having the same reason to believe that E either happens or not; that is, to believe that we can bet equal money on the arrival of E at least once. And so,

$$(1 − p)^2 = 1/2, \ n = − \lg 2/\lg(1 − p).$$

If E is the occurrence of a *six* or another determined face in a throw of die with six faces, then $p = 1/6$, and $n = 3.8018 \dots$ and it is advantageous to bet that a *six* will arrive at least once in 4 throws but disadvantageous to bet on 3 throws. If two dice are thrown, and E is the occurrence of a *double-six*, $p = 1/36$ and $n = 24.614 \dots$ so that it is advantageous[11] to bet that […].

The general expression (8.1) shows that however feeble is the chance *p* of event E, if only not exactly zero, it is always possible to choose a sufficiently large number *n* of trials for the probability of E arriving at least once to approach certitude as closely as desired. Indeed, however little (1 – *p*) differs from unity, it is always possible to choose a sufficiently large exponent *n* for (1 – *p*)$^n$ to become smaller than a given magnitude. This leads to an essential difference between an absolutely impossible thing and an event E whose chance is extremely slim. An impossible thing never happens, whereas an event having an arbitrarily low probability will always likely occur in a sufficiently long series of trials.

According to the binomial formula

$$(1-p)^n = 1 - np + C_n^2 p^2 - C_n^3 p^3 + ...,$$

and, if *n* is a very large number, *n* − 1, *n* − 2, … can be replaced by *n*, so that

$$(1-p)^n = 1 - np + C_n^2 p^2 - C_n^3 p^3 + ...$$

which is the series for $e^{-np}$ where *e* is the base of the Naperian logarithms. Therefore, the approximate value of *r* is

$$r = 1 - e^{-np}.$$

If *p* = 1/*n*, it will be equal to (*e* − 1)/*e* so that, if the chance of some event E is 1/*n* where *n* is a very large number, that number will be sufficient for E to arrive at least once with probability (*e* − 1)/*e* ≈ 2/3.

**9.** If two events, E и E$_1$, are not at all independent, which means that the occurrence of one of them influences the chance of the other event, the probability of the compound event consisting of E и E$_1$ will be *pp*$_1$, where *p* is the probability of E, which should occur before E$_1$ does and *p*$_1$ expresses the probability that E$_1$ will follow.

And so, if *a* and *b* denote the numbers of white and black balls in urn A, and *c* is their sum; if E is the occurrence of a white ball in the first trial, and E$_1$ is the same event happening at the second drawing, then *p* = *a*/*c*, and *p*$_1$ = (*a* − 1)/(*c* − 1) because before the second drawing *c* and *a* became *c* − 1 and *a* − 1. Therefore,

$$pp_1 = \frac{a(a-1)}{c(c-1)}$$

According to the same rule,

$$pp_1 = \frac{ab}{c(c-1)}.$$

for the probability of the extraction of a white and a black ball in a determined order and without replacing the first of them.

Generally, the probability $w$ of obtaining $m$ white and $n$ black balls in $(m + n)$ successive drawings without replacement and in any determined order is

$$w = \frac{a(a-1)...(a-m+1)b(b-1)...(b-n+1)}{c(c-1)...(c-m-n+1)}.$$

Actually, if after $m_1 + n_1$ first drawings, $m_1$ white balls and $n_1$ black ones are extracted, there will remain $(a - m_1)$ white and $(b - n_1)$ black balls. In a new drawing the probabilities of extracting balls of those colours will therefore be

$$\frac{a-m_1}{c-m_1-n_1}, \quad \frac{b-n_1}{c-m_1-n_1}.$$

Assume for $m_1$ all the [natural] numbers from 0 to $m - 1$, and from 0 to $n - 1$ for $n_1$; then the product of the $(m + n)$ thus obtained magnitudes will evidently make up the value of $w$ coinciding with the formula above.

When replacing the extracted balls of either colour, their chances will remain constant and equal to $a/c$ and $b/c$. The probability of getting $m$ white and $n$ black balls in a determined order will be $a^m b^n/c^{m+n}$. This is how the expression for $w$ is actually transformed when $a$ and $b$ are extremely large and can be considered infinite with respect to $m$ and $n$ so that during the trials the chances of the balls of each colour remain invariable. If $n = 0$, we get the probability of extracting $m$ white balls in succession

$$w = \frac{a(a-1)...(a-m+1)}{c(c-1)...(c-m+1)}.$$

Instead of an urn [problem], consider a game with 16 red cards and the same number of black ones. If required to determine the probability of extracting all the red cards in succession in 16 drawings, we will have $a = 16$, $c = 32$, $m = 16$ so that

$$w = \frac{16!}{17 \cdot 18 \cdot ... \cdot 31 \cdot 32} = \frac{1}{601,080,390}$$

which is a little less than $1/60 \cdot 10^7$. It will therefore be necessary to make somewhat more than $60 \cdot 10^7$ trials for obtaining probability 2/3 or for betting about 2 against 1 on the extraction of the 16 red cards at least once without interruption.

**10.** Suppose that an event E can occur in many distinct ways independent one from another with probabilities $p_1, p_2, ...$ Its composite probability $p$ will be the sum of all these partial probabilities:

$p = p_1 + p_2 + ...$

For the sake of definiteness let there be a given number $i$ of urns A containing white and black balls; the total number of balls and the number of white balls among them being $c_1$ and $a_1$ in the first urn, $c_2$ and $a_2$ in the second urn etc. Suppose also that E is the extraction of a white ball from a randomly chosen urn which can therefore happen in $i$ different ways, and that the probability of choosing one of these urns is the same for all of them and thus equal to $1/i$. The chances of extracting a white ball are $a_1/c_1$, $a_2/c_2$ etc. According to the rule of § 5, the probabilities $p_1, p_2, \ldots$ of the different ways in which E can occur will be

$$p_1 = \frac{a_1}{ic_1}, \; p_2 = \frac{a_2}{ic_2}, \; \ldots$$

and it is required to prove that the composite probability $p$ of extracting a white ball from any of the urns A will be

$$p = \frac{1}{i}\left[\frac{a_1}{c_1} + \frac{a_2}{c_2} + \ldots\right].$$

The demonstration is based on a lemma which will be equally useful on other occasions.

Consider some number $i$ of urns C each containing $\mu$ white and black balls in diverse proportions. The probability of extracting a white ball from any of them will not change if we combine all the $i\mu$ balls in one single urn B. Actually, they will form somehow arranged groups containing one and the same number $\mu$ of balls taken from the same urn. This will suffice for the chance of randomly choosing a group to be the same for all of them and equal to $1/i$ as though each group were contained in one of the urns C. The chance of extracting a white ball from a randomly chosen group will not change and the probability of extracting a white ball will be the same for urn B and for the system of urns C.

That conclusion will not persist if the numbers of balls contained in urns C are unequal. For any number of them the chance of randomly choosing an urn remains the same and equals $1/i$, but when all the balls are contained in B, the groups there will consist of unequal numbers of balls and the chances of randomly selecting a group will differ; it will evidently be larger for larger groups. And so, we will reduce all the fractions $a_1/c_1, a_2/c_2, \ldots$ to the same denominator $\mu$. Let $\alpha_1, \alpha_2, \ldots$ be their numerators, so that

$$\frac{a_1}{c_1} = \frac{\alpha_1}{\mu_1}, \; \frac{a_2}{c_2} = \frac{\alpha_2}{\mu_2}, \ldots$$

The chance of extracting a white ball from each of the urns A and, therefore, from all of them taken together, will not change if each number $c_1, c_2, \ldots$ of the balls of either colour is replaced by one and the same number $\mu$, and the numbers of white balls, $a_1, a_2, \ldots$, by $\alpha_1, \alpha_2, \ldots$ The probability of extracting a white ball will not change either if we then unite all the balls in urn C. And so, C contains the total number $i\mu$ of balls among which $\alpha_1 + \alpha_2 + \ldots$ are white. That probability will be

$$\frac{1}{i}[\frac{a_1}{\mu} + \frac{a_2}{\mu} + \ldots].$$

By virtue of the preceding equations, this magnitude coincides with the value of $p$, Q. E. D.

**11.** For applying this rule to examples let us suppose first of all that a person knows that a ball was extracted from urn A containing 5 white balls and 1 black ball, or from urn B with 3 and 4 of those balls, and that he has no reason to believe that that ball was extracted from one of those urns rather than from another. For him, the probability that the extracted ball is white is

$w = (1/2)(5/6) + (1/2)(3/7) = 53/84$

since that event could have happened in two different ways whose probabilities are the two terms written above.

For another person, who knows that the ball was extracted from B, the probability that it is black, is $p = 4/7 = 48/84$. Both fractions exceed 1/2, so the first person should think that the extracted ball is white, and the second, that it is black. Although 53/84 is larger than 48/84, we should choose the latter opinion since the second person is better informed. This is a very simple example and it is easy to provide many others concerning the previous statement in § 1 about contrary opinions on the same question formed by differently informed persons.

Suppose also that we know that an urn A contains a given number of white and black balls whose proportion is absolutely unknown to us. We can formulate $(n + 1)$ different and equally possible hypotheses about that proportion, and, at the same time, about the different ways in which a white ball can be extracted. These hypotheses are: $n$ white balls; $(n - 1)$ white balls and 1 black ball; $(n - 2)$ white and 2 black balls; ... $n$ black balls. All these assumptions are equally possible and the probability of each is $1/(n + 1)$. The partial probabilities of extracting a white ball are, accordingly,

$$p_1 = \frac{1}{n+1} \cdot \frac{n}{n}, \; p_2 = \frac{1}{n+1} \cdot \frac{n-1}{n}, \; p_3 = \frac{1}{n+1} \cdot \frac{n-2}{n}, \ldots$$

and the composite probability is

$$w = \frac{1}{n+1}(\frac{n}{n} + \frac{n-1}{n} + \frac{n-2}{n} + \ldots + \frac{n-n}{n}) = \frac{1}{2}$$

as it should be since we have no reason to believe in the arrival of a white ball rather than a black ball.

If, however, we certainly know that there are more white balls in the urn than black ones, then $w > 1/2$. For determining it, we should distinguish the cases of an odd and an even *n*. Let *i* be a natural number. If $n = 2i + 1$, we can only formulate $i + 1$ different and equally possible hypotheses: $2i + 1$ white balls; $2i$ white balls and 1 black ball; … $(i + 1)$ white and *i* black balls. In this first case, the composite value of *w* is

$$w = \frac{1}{(i+1)(2i+1)}[(2i+1) + 2i + (2i-1) + ... + (i+1)] = \frac{1}{2}\frac{3i+2}{2i+1}.$$

It is unity, as it should be, if $i = 0$, and it indefinitely decreases and approaches 3/4 when *i* increases ever more.

If $n = 2i + 2$, we can also formulate $(i + 1)$ equally possible hypotheses: A contains $(2i + 2)$ white balls; $(2i + 1)$ white balls and 1 black ball; … $(i + 2)$ white and *i* black balls. The composite value of *w* is

$$w = \frac{1}{(i+1)(2i+2)}[(2i+2) + (2i+1) + 2i + ... + (i+2)] = \frac{1}{2}\frac{3i+4}{2i+2}.$$

Just like in the former case, $w = 1$ and 3/4 if *i* takes the extreme values 0 and ∞. For any other natural number *i* the probability exceeds the preceding by $i/[4(i + 1)(2i + 1)]$ whose maximal value is 1/24 at $i = 1$.

An urn A contains *c* balls, *a* of them white. They are collected in groups the first of which has $c_1$ balls, $a_1$ of them white, the second has $c_2$ balls, $a_2$ of them white, etc. so that

$c_1 + c_2 + + ... = c$, $a_1 + a_2 + + ... = a$.

Let *p* be the probability of a white ball appearing from that urn. It should be equal to *a/c*, which only verifies the rule of § 10. A white ball can come from the first group, the chance of which is $(c_1/c) \cdot (a_1/a)$ and the composite value of *p*, equal to *a/c* according to the second of the two preceding equations, is

$$p = \frac{c_1}{c} \cdot \frac{a_1}{c_1} + \frac{c_2}{c} \cdot \frac{a_2}{c_2} + ...$$

However, if we place all these groups in different urns $A_1$, $A_2$, …, and if all the numbers $c_1$, $c_2$, … are unequal, the chance of extracting a white ball will not be *a/c* anymore. Generally, that chance depends on the way that the white and black balls are distributed among $A_1$, $A_2$, … Without knowing it, calculation is impossible. Nevertheless, for someone who does not know it, the reason to believe in the appearance of a white ball when a ball is randomly extracted from those urns is evidently the same as believing in an appearance of the ball of that colour from A. Therefore, for that person, the probability of such an extraction, unlike its proper chance, will be equal to *a/c*. I assume, for example, that A contains 2 white balls and 1 black ball, 2 of them going to $A_1$, and the third one, to $A_2$. For that person there will be 3 equally possible distributions: both white balls from $A_1$ and the black ball from $A_2$; 1 white, and 1 black ball from $A_1$ and the other white ball from $A_2$; that second white ball and the black ball, from $A_1$ and the first white ball from $A_2$.

For these three cases the probability of extracting a white ball from either urn will be

(1/2) (1 + 0); (1/2) [(1/2) + 1)]; (1/2) [(1/2) + 1)].

Taking their sum and dividing it by three we have 2/3 for the composite probability of that extraction, just as it is for an appearance of a white ball from A.

Consider finally a system of urns $D_1, D_2, \ldots$ containing $c_1$ balls, $a_1$ of them white, in $D_1$; $c_2$ balls, $a_2$ of them white, in $D_2 \ldots$ Suppose that for some reason the chances of choosing an urn for extracting a ball differ and are equal to $k_1, k_2, \ldots$ By the rule of § 5 the probabilities of extracting a white ball will be $k_1 a_1/c_1, k_2 a_2/c_2, \ldots$ These products express the partial probabilities $p_1, p_2, \ldots$ concerning the different ways of the extraction. The composite probability will therefore be

$$w = \frac{k_1 a_1}{c_1} + \frac{k_2 a_2}{c_2} + \ldots$$

Regarding a system $A_1, A_2, \ldots$ of urns for which all the probabilities $k_1, k_2, \ldots$ are the same, will suffice for demonstrating the rule of § 10 in all generality. And that rule being thus proved, its application to the other urns $D_1, D_2, \ldots$ for which the chances $k_1, k_2, \ldots$ take some values, leads, as is seen, to the expression of $w$ of the general case.

**13.** Let now E and F be contrary events excluding each other, one of which should always occur. Their probabilities are $p$ and $q$, $p + q = 1$. Suppose that each of these events can happen in various ways with probabilities $p_1, p_2 \ldots$ and $q_1, q_2 \ldots$ respectively. Successively applying the preceding rule, we will have

$$p = p_1 + p_2 + \ldots, q = q_1 + q_2 + \ldots, p_1 + p_2 + \ldots + q_1 + q_2 + \ldots = 1.$$

In problems of eventuality the terms on the left side of the last-written equation are the probabilities of the diverse favourable and unfavourable combinations for the appearance of E. That equation therefore shows that their sum is always equal to unity or certitude which should indeed take place if all the possible combinations are exhausted.

It follows from that same equation that $p$ can be written as

$$p = \frac{lp_1 + lp_2 + \ldots}{lp_1 + lp_2 + \ldots + lq_1 + lq_2 + \ldots}$$

where $l$ is an arbitrary magnitude. The terms of this fraction are proportional to the chances $p_1, p_2 \ldots, q_1, q_2 \ldots$ of the favourable and unfavourable cases for the arrival of E. And if we suppose that among the terms of the numerator $a'$ terms are equal one to another and equal $\alpha'$, $a''$ are equal one to another and equal $\alpha''$, etc; that among the terms of the denominator $c'$ terms are equal one to another and equal $\gamma'$, $c''$ terms are equal one to another and equal $\gamma''$, etc, then

$$p = \frac{\alpha'a' + \alpha''a'' + \ldots}{\gamma'c' + \gamma''c'' + \ldots}.$$

Thus, if those favourable and unfavourable cases do not have the same chance, the probability of E will be expressed by multiplying the numbers of equally probable cases by magnitudes proportional to their respective probabilities and then dividing by the sum of those products for all the possible cases. This rule is more general and often more handy for applying than that of § 2 since in no problem about the occurrence of an event whose probability we wish to know does it require an equality of chances of all the favourable and unfavourable cases.

**14.** The rules of §§ 5 and 10 suffice for deriving the formulas about the repetition of an event with known constant or variable chances during a series of trials. As always, we denote contrary events of some nature, one of which should occur at each trial, by E and F. Suppose first of all that their probabilities are constant and given and denote their chances by $p$ and $q$ respectively and also let the total number of trials be $\mu$ with E arriving during them $m$ times and F, $n$ times so that $p + q = 1$, $m + n = \mu$.

The probability that E and F will arrive $m$ and $n$ times in some determined order is independent from that particular requirement and equal to $p^m q^n$ (§ 7). Therefore, denoting the probability of such a result by $\Pi$ and the number of different ways in which it is possible by $K$, we will have in accordance to the rule of § 10

$$\Pi = Kp^m q^n.$$

For determining $K$, I suppose first of all that the $\mu$ events A, B, C, … which should take place are all different. Then $K$ will be the number of *permutations* of all the letters $\mu$ which is equal to $\mu!$.

Then, if some $m$ of those letters A, B, C, … represent one and the same event E, the number of different permutations will be $\mu!/m!$. And if the other $\mu - m = n$ letters represent one and the same event F, the number of the appropriate permutations will be $\mu!/n!$. Consequently, the number of various permutations with event E occurring $m$ times, and F, $n$ times, or the required value of $K$ will be

$$K = \mu!/m!n!.$$

Since $\mu = m + n$, this magnitude will be symmetric with respect to $m$ and $n$:

$$K = C_\mu^m = C_\mu^n.$$

This proves that Π is the $(m + 1)$-th term of the expansion of $(p + q)^\mu$ arranged in increasing powers of $p$, or its $(n + 1)$-th term arranged in increasing powers of $q$. It follows that in our case when the chances $p$ and $q$ of the contrary events E and F are constant, the chances of all the compound events that can occur in $\mu$ trials are expressed by different terms of $(p + q)^\mu$.

Those events are $(\mu + 1)$ in number and they are unequally probable either because of the multiplicity $K$ of the combinations which can lead to them, or since the chances $p$ and $q$ are unequal. If $p = q$, the most probable event corresponds to $m = n$ if $\mu$ is even, and to one of the two cases $m - n = \pm 1$ for odd values of $\mu$.

**15.** Let $P$ be the probability that E occurs at least $m$ times in $\mu$ trials. This compound event can take place in $(m + 1)$ different ways, when E arrives $\mu, \mu - 1, \ldots, \mu - n = m$ times. The corresponding probabilities can be derived from the preceding expression of Π by consecutively substituting $\mu$ and 0, $\mu - 1$ and 1, $\ldots$, $m$ and $n$ instead of these two last-mentioned numbers. By the rule of § 10 the composite probability of $P$ will be the sum of these $(n + 1)$ partial probabilities:

$$P = p^\mu + \mu p^{\mu-1} q + C_\mu^2 p^{\mu-2} q^2 + \ldots + C_\mu^n p^m q^n.$$

In other words, $P$ will be the sum of the $(n + 1)$ first terms of $(p + q)^\mu$ arranged in increasing powers of $q$. For $m = 0$ or $n = \mu$, $P = (p + q)^\mu = 1$, as it should be since this compound event includes all the possible combinations of E and F arriving in all the trials so that its probability should become certainty. For $m = 1$ the corresponding event is contrary to F occurring at all trials; actually, the value of $P$ will then be the entire expansion of $(p + q)^\mu$ less its last term, $q^\mu$, which accords with the value of $r$ from § 8.

If $\mu = 2i + 1$ is odd, and it is required to determine the probability that E occurs oftener than F, it can be derived from the general expression of $P$ when taking $m = i + 1$ and $n = i$. If, however, $\mu = 2i$ is even, the probability that E arrives at least as many times as F does, is obtained from the same expression by taking $m = n = i$.

**16.** The solution of the first probability problem, mentioned here at the very beginning of the Preamble and known as the *problem of points*, can also be provided from that formula. Gamblers A and B are playing some game in which one of them gains a point after each set with probabilities $p$ and $q$. To win the game, they lack $a$ and $b$ points, and it is required to determine their probabilities, α and β, of achieving this aim. One of these contrary events will necessarily occur, so that α + β = 1 and it is only sufficient to determine α.

Note first of all that the game ends at most after $(a + b - 1)$ sets since A will then win at least $a$ points, or B, at least $b$ points. Furthermore, without changing at all their respective chances of winning the game, the gamblers can agree to play $(a + b - 1)$ sets since, whatever happens after that, only one of them can win his lacking points: either A wins $a$ points before B wins $b$, or B wins $b$ points before A wins $a$.

For determining α and β, we may therefore suppose that the number of sets is always $(a + b - 1)$. And so, α is the probability that event E having chance *p* at each trial will appear at least *a* times in that number of trials so that its value can be determined from the preceding expression of *P* by taking

μ = *a* + *b* − 1, *m* = *a,* and *n* = *b* − 1.

Let for example *p* = 2/3, *q* = 1/3, *a* = 4 and *b* = 2. Then α = 112/243, β = 131/243 and β > α. It follows that a gambler A, twice more skilful than B, i. e., having a double chance of winning each set[12], still can not bet without disadvantage on winning 4 points before B wins 2.

As will be shown below, if the gamblers agree to stop playing without ending, A and B will receive the stakes multiplied by the chances α and β of winning them which means that they ought to share the stakes in the proportion α:β.

**17.** Suppose that instead of events E and F there is a larger number of them, three for example, and call them E, F, and G with only one of them arriving at each trial. Let their constant probabilities be *p, q* and *r*, and μ, the number of trials. It is easy to generalize the method of § 14 and arrive at the probability that E, F and G occur *m, n* and *o* times:

$$\frac{\mu! p^m q^n r^o}{m!n!o!}, \; p + q + r = 1, \; m + n + o = \mu,$$

so that the derived probability is the general term of the expansion of $(p + q + r)^\mu$.

Such is the case of an urn containing balls of three different colours in proportion *p:q:r* and the events E, F and G being extractions with replacement of these three kinds of balls. Taking the sum of the terms of the expansion of $(p + q + r)^\mu$ containing *p* to the power equal or larger than *m*, we will have the probability that E occurs at least *m* times in μ trials. Whichever is the number of the events E, F, G, …, only one of which arrives at each trial, we can immediately derive that probability by issuing from the preceding expression for *P*.

As always, denote the constant chances *p, q, r*, … of the occurrence of one of the events E, F, G, … in each trial. This pattern can be considered as a compound event[13], call it F′, and its probability, *q*′. Then

*q*′ = *q* + *r* + …, *p* + *q*′ = 1.

Events E and F′ are therefore contrary and only one of them will arrive at each trial. Therefore, the probability that E occurs at least *m* times in μ trials will be obtained by substituting *q*′ instead of *q* in the expression of *P*.

To provide an example of this rule based on the expansion of a power of a polynomial, I suppose that an urn A contains *m* balls numbered 1, 2, …, *m*. When extracting a ball with replacement μ times in succession, the chance of the arrival of a ball with a fixed number is equal to $1/m$. Let

$$n_1, n_2, \ldots, n_m \qquad (17.1)$$

be given numbers including zero or not and

$$n_1 + n_2 + \ldots + n_m = \mu.$$

Then denote by $U$ the probability that ball No. 1 will occur $n_1$ times in some particular order, ball No. 2, $n_2$ times, ... and let

$$(t_1 + t_2 + \ldots + t_m)^\mu = \theta. \qquad (17.2)$$

When expanding $\theta$ in powers and products of the undetermined magnitudes $t_1, t_2, \ldots, t_m$, $U$ will be the term of that expansion containing

$$t_1^{n_1} t_2^{n_2} \ldots t_m^{n_m}$$

with all these magnitudes made equal to $1/m$. Denote the numerical coefficient of that product by $N$, then

$$U = N/m^\mu.$$

$N$ is a natural number depending on $\mu$ and numbers (17.1); it is equal to

$$\frac{\mu!}{n_1! n_2! \ldots n_m!}$$

with $0! = 1$ in all such products.

Let now $s$ be the sum of the numbers extracted in $\mu$ drawings:

$$s = n_1 + 2n_2 + \ldots + mn_m.$$

If $s$ is given and all the natural numbers and zero satisfying that equation and adding up to $\mu$ are consecutively substituted for $n_1, n_2, \ldots, n_m$; if $N'$, $N''$, ... are the corresponding values of $N$; and if $V$ is the sum of the corresponding values of $U$, − then

$$V = (1/m^\mu)(N' + N'' + \ldots)$$

will be the probability of obtaining the given sum $s$.

It is easier to calculate $V$ by replacing the undetermined magnitudes $t_1, t_2, \ldots, t_m$ by the powers $t, t^2, \ldots, t^m$ of the same magnitude $t$. Denote the new value of $\theta$, see (17.2), by $T$, then[14]

$$T = (t + t^2 + \ldots + t^m)^\mu$$

and it will be easy to see that the sum $(N' + N'' + …)$ is just the numerical coefficient $M_s$ of $t^s$ in the expansion of $T$ so that

$V = M_s/m^\mu$.

That coefficient depends on the given numbers $\mu$, $m$ and $s$ and is easily calculated in each example [in each problem].

Instead of a single urn A we may suppose that there are $\mu$ urns $A_1$, $A_2$, …, $A_\mu$ each containing $m$ balls numbered 1, 2, …, $m$, and that a ball is extracted from each. Or, we may also replace these urns by the same number of ordinary dice with six faces numbered 1, 2, …, 5, 6. So $m = 6$ and $V$ expresses the probability that in a throw of $\mu$ dice the sum of the obtained points will be $s$. Let $\mu = 3$, then

$T = t^3(1 + t + t^2 + … + t^5)^3$, $V = M_s/6^3$.

The expansion of $T$ consists of 16 terms with their coefficients equally remote from the extremes, such as $M_3$ and $M_{18}$, $M_4$, and $M_{17}$, …, $M_{10}$ and $M_{11}$, being equal. The sum of all the coefficients will be equal to the value of $T$ at $t = 1$, or $6^3$. The sum of the first 8 coefficients, $M_3$, $M_4$, …, $M_{10}$, just as the sum of the last 8 of them, $M_{11}$, $M_{12}$, …, $M_{18}$, will be equal to $6^3/2$.

We conclude that in a throw of three dice the probability of gaining 10 points or less is 1/2, the same as when obtaining 11 or more. It will therefore be a fair bet that the sum of the arrived three numbers is smaller or larger than 10. It is on this result that the game *passe-dix* is based. Without being aided by calculation, it is easy to become assured in the equality of chances of the two gamblers when noting that the sum of the points on each pair of opposing faces of the same die equals 7 (6 and 1, 5 and 2, 4 and 3)[15]. It also follows that the sum of the three superior and inferior numbers is always 21, so that, when the three dice are falling on the gaming table, and the first sum is larger than 10, the second [the next] sum is less, and *vice versa*. […]

For finding the chances of the diverse values of $s$ from $s = 3$ to 18, it is necessary to turn to the expansion of $T$. We get

$M_3 = M_{18} = 1$, $M_4 = M_{17} = 3$, $M_5 = M_{16} = 6$, $M_6 = M_{15} = 10$,
$M_7 = M_{14} = 15$, $M_8 = M_{13} = 21$, $M_9 = M_{12} = 25$, $M_{10} = M_{11} = 27$

These are the numbers of the combinations of three throws for achieving the appropriate sums. By dividing those numbers by $6^3 = 216$, we get the chances of those sums.

**18.** When the chance of event E varies during the trials, the probability of its being repeated a given number of times, depends on the law of that variation. Suppose, as it was done in § 9, that E is the extraction of a white ball without replacement from an urn A containing $a$ white and $b$ black balls. Denote the number of drawings by $\mu$ and the probability of the arrival of $m$ white and $n$ black balls in some determined order by $w$. The value of $w$ is provided by the formula in that § 9 and is independent from the order of the appearance of those balls.

Denote the probability that they arrive in some fixed order by Π, then

Π = Kw

where K is the same as in § 14. Let, as always,

$m + n = \mu, a + b = c; a - m = a', b - n = b', c - \mu = a' + b' = c'.$

Then $a', b', c'$ will be the numbers of balls of both colours and their sum, which were initially $a, b$ and $c,$ and which are left after the drawings. When considering the expressions of $K$ and $w$, we will have

$$\Pi = \frac{\mu!a!b!c'!}{m!n!a'!b'!c!}.$$

This formula can easily be extended to the case in which A includes balls of three or more different colours. After suppressing common factors in the numerator and denominator that formula becomes simpler[16]:

$$\Pi = C_\mu^n \frac{a(a-1)...(a-m+1)b(b-1)...(b-n+1)}{c(c-1)...(\mu+1)}.$$

The probability of extracting not less than $m$ white balls in $\mu$ drawings will be the sum of $(n + 1)$ values of Π obtained by replacing $m$ and $n$ by $\mu$ и 0; by $\mu − 1$ and 1; …, $\mu − n$ and $n$. Denoting that probability by $P$, we will get

$Pc(c-1)...(\mu-1) + (\in a\ d)...(\mu-1) + $
$\mu ba(a-1)...(a-\mu+2) + C_\mu^2 b(b-1)a(a-1)...(\mu-3) + $
$C_\mu^3 b(b-1)(b-2)a(a-1)...(\mu-4) + $
$... + C_\mu^n b(b-1)...(b-n+1)a(a-1)...(a-m+1).$

If $m = 0$ and $n = \mu$, then $P = 1$ and we conclude that

$c(c-1)...(\mu-1) + (\in a\ d)...(\mu-1) + $
$\mu ba(a-1)...(a-\mu+2) + C_\mu^2 b(b-1)a(a-1)...(\mu-3) + $
$C_\mu^3 b(b-1)(b-2)a(a-1)...(\mu-4) + ... + $
$\mu b(b-1) ... (b - \mu + 2a) + b(b-1) ... (b - \mu + 1)$

which coincides with the generally known formula similar to that of the binomial. In this [in my] formula and in all others of the same kind, each magnitude such as $a(a − 1)(a − 2) … (a − m + 1)$, being a product of $m$ factors, should be supposed to equal unity at $m = 0$. Therefore, that magnitude is not suited for the case in which $\mu = 0$ and the same exception is just as applicable to a binomial to the power of zero.

**19.** When drawing μ times without replacement, the probability of obtaining *m* white and *n* black balls will be the same whether drawing consecutively or all at once. This can be verified in the following way[17].

I denote the number of groups, each consisting of μ balls, by $G_\mu = C_c^\mu$ which can be formed given the *c* balls contained in the urn. Actually, to form those groups from groups of (μ − 1) balls it is necessary to combine each of them with the (*c* − μ + 1) balls not contained there. The number of groups containing μ balls will be (*c* − μ + 1)$G_{\mu-1}$ but μ groups of (μ − 1) balls will produce one and the same group of μ balls, and we should therefore divide the product (*c* − μ + 1)$G_{\mu-1}$ by μ, to get the number of different groups of μ balls. And so,

$G_\mu = G_{\mu-1}(c − \mu + 1)/\mu$.

For μ = 1, evidently $G_1 = c$. Taking μ = 2, 3, 4, … we will obtain

$G_2 = [(c-1)/2]G_1 = C_c^2, \ G_3 = [(c-2)/3]G_2 = C_c^3,$
$G_4 = [(c-3)/4]G_3 = C_c^4, \ …$

and finally, $G_\mu$.

Denote now the value of $G_m$ when *c* and μ are replaced by *a,* and *m* by $G'_m$; and by $G''_n$ when those two magnitudes are replaced by *b* and *n*. Then

$G'_m = C_a^m, \ G''_n = C_b^n.$

The product $G'_m G''_n$ is the number of groups consisting of *m* + *n* = μ balls, each containing *m* white and *n* black balls, which can be formed from *a* + *b* = *c* balls. The probability of obtaining one of those groups when drawing μ balls at once is equal to their number divided by the number of all the groups of μ balls in the urn

Π = $G'_m G''_n / G_\mu$

and coincides with Π of § 18. The expression of *P* in the same section is also the probability of obtaining at least *m* white balls when drawing μ balls at once.

**20.** In the example of § 18 the chance of event E varied during the trials since at each new trial it depended on the number of the previous arrivals of that E and of the contrary event F. However, there are other problems in which the proper chances of these two events of some nature are independent from previous trials but vary from one trial to another.

Generally, in a series of μ past or future trials let $p_1$ and $q_1$ be the chances of E and F at the first trial, $p_2$ and $q_2$, at the second trial, …, and $p_\mu$ and $q_\mu$ at the last trial. Then $p_1 + q_1 = p_2 + q_2 = … = p_\mu + q_\mu = 1$. It is required to determine the probability that E will arrive or arrived *m* times, and F, *n* = μ − *m* times in some definite order. I denote the product of *m* of the fractions $p_1 p_2 \ldots p_\mu$ by $P_m$, and, by $Q_n$, the product of *n* fractions $q_1 q_2 \ldots q_\mu$. Here, if $P_m$ includes or not the fraction $p_i$, then $q_i$ is not, is included in $Q_n$.

I then multiply $P_m$ by $Q_n$ and form the sum of all the possible magnitudes $P_m Q_n$ whose number was denoted by *K* in § 14. And that sum expresses the required probability. It is also possible to formulate that rule in another way which will be useful below.

Let *u* and *v* be two undetermined magnitudes and

$R = (up_1 + vq_1)(up_2 + vq_2) \ldots (up_\mu + vq_\mu)$.

That *R* is the product of $\mu = m + n$ factors[18], and when actually multiplying them we will get a polynomial of ($\mu + 1$) terms arranged in powers of *u* and *v*. The coefficient of $u^m v^n$ in that polynomial will be the required probability of the arrival of E and F *m* and *n* times respectively in some definite order. Suppose for example that *m* = 3, then

$R = u^3 p_1 p_2 p_3 + u^2 v (p_1 p_2 q_3 + p_1 p_3 q_2 + p_2 p_3 q_1) +$
$\quad uv^2 (p_1 q_2 q_3 + p_2 q_1 q_3 + p_3 q_1 q_2) + v^3 q_1 q_2 q_3$.

The coefficient of $u^3$ is evidently the probability that E occurred three times; of $u^2 v$, the probability that E arrives at the first two trials, and F, at the third; or, F at the second trial, and E, at the two others; or, F at the first, and E, at the last two trials. The coefficient of $uv^2$ expresses the probability that F occurs twice, and E, once; and, finally, the coefficient of $v^3$ is evidently the probability of the appearance of F three times.

Now, if at each trial E can occur in many equally possible ways, we adopt the sum of the [partial] probabilities of those diverse ways divided by their number as the chance of E at that trial, and note that this procedure conforms to the rule of § 10. Subtracting that mean chance of E from unity, we get that chance of F. Those mean chances for each trial should be applied for calculating the probability that E and F will arrive *m* and *n* times in *m* + *n* trials, and the same holds for any other compound event consisting of E and F. However the partial chances of E and F vary in number and magnitude from one trial to another, if their mean chances remain constant, the probabilities of the compound events obey the same laws as when the chances of those events are invariable.

**21.** One of the most frequent applications of the calculus of probability aims at determining the advantages and damages attached to possible things as against the gains and losses which they produce and the chances of their occurrence. That aim is based on the following rule.

Suppose that one of the events E, F, G, … should appear and denote their probabilities by *p, q, r,* … so that $p + q + r + \ldots = 1$. Suppose also that for some person gain *g* is attached to the arrival of E, that for another person, it is attached to F, etc. If all those interested agree to share *g* before chance will decide the problem, or if they are obliged to share that gain by some cause, it should be shared between them proportionally to their probabilities of winning it, so that *gp* will be the share of the first person, *gq*, the share of the second person, etc.

Actually, if *m* is the number of all the equally possible cases among which *a, b, c*, … are favourable for E, F, G, …, then

$a + b + c + ... = m, p = a/m, q = b/m, r = c/m, …$

If there are *m* people and each of them is to gain by the appearance of one of the *m* possible causes, then evidently *g* should be equally shared among all of them so that each will get *g/m*. And the person who enjoys probability *p* of gaining or who has *a* possible chances for gaining, should then receive *ag/m = pg*. The shares of the other people will similarly be *qg, rg*, …

If a game has begun, that rule will indicate what will each gambler get according to his probability of winning the game, should those gamblers agree to separate without concluding it. And the share of each gambler before the game begins should be proportional to his chance of winning. Indeed, if the gamblers decide to separate instead of playing, each of them should get his share back. By the preceding rule, the sum returned to him should also be equal to the sum of the stakes multiplied by the probability of [his] winning the game. In games of chance that probability depends on the rules of the game and can be calculated in advance if only they are not too complicated. In games in which success depends on the skill of each participant, the probability of winning it is usually based on their reputations and can only be determined with some precision by long experience.

Let the probabilities of two contrary events E and F be *p* and *q*, $p + q = 1$. If A bets α on the arrival of E, and B bets β on F, those sums α and β should be in the same ratio as *p* and *q*, $pβ = qα$. However, it should not be forgotten that those probabilities *p* and *q* differ in general from the proper chances of E and F and depend on the knowledge which A and B could have obtained about those events. If the probabilities mentioned are based on the same knowledge of both gamblers, the bet will be equitable even if it favours one of them at the expense of the other. Otherwise, the ratio α:β will not be equal to the ratio of their probabilities and there will be no means equitably to regulate the bet.

**22.** The formulas of § 19 allow to calculate without difficulty the diverse chances of the Lottery of France, happily suppressed by a recent law. Comparing them with the multiples of the stakes paid by the Lottery for the winning tickets, we see that their multiples are much smaller than those for a fair game. For the Lottery, this resulted in an exorbitant advantage at the gamblers' expense which the law should have punished as being illegal.

Let in general *n* be the number of the numbers included in a lottery, *m* of them being extracted at each drawing, *l*, the number of winning numbers chosen by a gambler, and λ, the probability that those last-mentioned numbers will be drawn. According to the formulas of § 19, the quantities of possible groups consisting of *l* numbers formed either from the *n* numbers of the lottery, or from the *m* numbers extracted at each drawing, are $C_n^l$ and $C_m^l$ and the probability λ will be the ratio of the latter to the former:

$$\lambda = C_m^l \div C_n^l = \frac{m(m-1)...(m-l+1)}{n(n-1)...(n-l+1)}.$$

Suppose that a gambler's stake is unity, then the stake of the lottery should be (1 − λ)/λ, and if a gambler wins, the lottery should also give him back his stake. Let μ be the just multiple of the gambler's stake, then

μ = [(1 − λ)/λ] + 1 = 1/λ.

Let also *x* be the number of consecutive drawings for betting even money on the numbers chosen by the gambler to win at least once. By the rule of § 8

(1 − λ)$^x$ = 1/2,

and, if λ is a very small fraction, $x \approx 0.69315/\lambda$ where the coefficient is the Naperian logarithm of *x*.
In the Lottery of France *n* = 90 and *m* = 5. For gambling on three winning numbers *l* = 3 and

$$\lambda = \frac{5 \cdot 4 \cdot 3}{90 \cdot 89 \cdot 88}, \mu = 11,748, x = 8143.13...$$

In case of his winning, the Lottery should have paid the gambler 11,748 times more than his stake, but actually paid 5500 times more, i. e., less than a half of the indicated. The disproportion was even greater for betting on 4 or 5 numbers, and less when gambling on 2 or 1.
It was advantageous to bet even money on three chosen numbers to be extracted at least once in 8144 drawings and disadvantageous to bet on 8143. With respect to one number chosen in advance

(1 − 1/18)$^x$ = 1/2, *x* = lg2/(lg18 − lg17) = 12.137 …

so it was disadvantageous to bet even money on the extraction of that number at least once in 12 drawings and advantageous to bet on 13. It would have also been advantageous to bet even money on the extraction of the 90 numbers at least once in 85 or 84 drawings (Laplace 1812/1886, pp. 200 – 201)[19].
Some gamblers choose numbers which had not been drawn for a long time; others, on the contrary, choose those which had been most often drawn[20]. Both these preferences are equally ill-considered. For example, suppose that there exists a probability [1 − (17/18)]$^{100}$ ≈ 0.997 very nearly a certitude, that some fixed number will appear at least once in a 100 consecutive drawings, but that it was not extracted in the first 88 of them. The probability of its being drawn in the 12 last extractions will be almost 1/2, just as for any other determined number. As to the most often drawn numbers, that circumstance should not be thought to be due to chance compatible with the obvious equality of chances of all the numbers at each drawing. In no games of chance in which equal or unequal chances are known in a certain manner, past events have any influence on the probability of future events[21], and no combinations imagined by gamblers can augment the gain or diminish the loss resulting from chances as described by the rule of § 21.

In the Paris public games[22], the advantages of the banker at each set is inconsiderable. Thus, for the game *thirty-and-forty* it is a little less than 0.011 of the stake, see my paper (1825). However, because of the rapidity of that game and the large number of sets played during a few hours, the banker enjoys assured benefit almost constant from year to year for which he has to pay yearly five or six millions to the public administration which granted him the monopoly. Those games are more harmful than the Lottery since only in the capital the moneys in those games amount to many hundred millions, much more than in the Lottery over all France.

Any discussion of the reasons put forward to preserve those public games is here out of place. I had never been able to consider them good enough; suffice it to mention that these games have caused many misfortunes and perhaps crimes so that the administration should suppress them rather than share the benefits they provide with those whom they sell the privilege[23].

**23.** The product of a gain and the probability of obtaining it is what is called *mathematical expectation* of anyone interested in some speculation[24]. If the gain is 60,000 francs, say, and 1/3 is the chance of the event to which it is attached, the person who is to receive that sum under some circumstances may consider a third of that sum as his property and ought to include it in the list of his actual fortunes.

In general, if someone is to get $g$ at the arrival of an event E, $g'$ at the occurrence of E′, … and if the chances of those events are $p, p'$, …, his mathematical expectation will be $gp + g'p' + \ldots$ And if one or many magnitudes $g, g'$, … express losses to be feared by that person, those magnitudes should be reckoned with a minus. According to the total expectation being positive or negative, it will represent an increase or decrease of fortune. And if the man does not wish to wait for the events to occur, that expectation should be actually included in his credits or debts.

It is well known that when the gains or losses will only happen after a long time, they should be *discounted* and thus converted into actual values independent from eventualities. If $g$ should only be paid to the person who evaluates his fortune after $n$ years, $g'$, after $n'$ years, … those sums will be today worth $g/(1 + \theta)^n$, $g'/(1\theta)^{n'}$, … where $\theta$ is the yearly interest. Therefore, when denoting by ε the part of the fortune resulting from the mathematical expectation of that person, we will have

$$\varepsilon = gp/(1 + \theta)^n + g'p'/(1\theta)^{n'} + \ldots$$

When taking account of the gains and losses following from those events, ε is the sum that another person should pay or receive depending on ε being positive or negative. Calculations of annuities on one or more lives, of life insurance and pensions are based on that formula and on *mortality tables* as can be seen in special contributions treating those problems.

**24.** Since the advantage procured for someone by a gain depends on his fortune, we distinguish that relative advantage from mathematical expectation and call it *moral expectation*. If it is an infinitely small magnitude, its ratio to the actual fortune of the person is admitted as the measure of moral expectation which can be either positive or negative depending on whether that fortune eventually increases or decreases.

By issuing from that measure, integral calculus provides consequences, according with the rules that prudence indicates about the way in which everyone should direct his speculations. By that calculus we also discover the reason for abstaining from gambling even when the game is fair. True, that reason is perhaps not the best possible. An irrefutable argument against games of chance, if they cease to be a simple amusement, is that they do not create values and that winning gamblers can only attain advantage when causing misfortune to, and sometimes ruin of the losers.

Commerce is also a game in the sense that success of most prudent speculations is always only likely and leaves room for chances of loss which ability and prudence can only lessen. However, commerce increases the value of things by transporting them from one place to another, and it is in this increase that the merchant finds his benefit and at the same time advantages the consumer[25].

**25.** The rule of § 21, simple and natural as it is, is fraught with a difficulty which we sometimes have to deal with. A and B play *heads* and *tails* under the following conditions: **1**. A game ends when *heads* appear. **2**. B gives A 2 francs if *heads* appear at once, 4 if it only arrives at the second toss, …, and in general $2^n$ francs if *heads* [first] occurs at the $n$-th toss. **3**. The game ends in a draw if *heads* do not happen in the first $m$ tosses. Without this restriction the game could have lasted for an infinite time.

It is supposed that the coin does not tend to fall on one side rather than on the other so that at each toss the chances of both *heads* and *tails* are 1/2. It follows that the probability of *heads* first appearing at the $n$-th throw is $1/2^n$. Indeed, for that event to happen, *tails* should have occurred ($n$ – 1) times consecutively with probability $1/2^{n-1}$. After that, the appearance of *heads* has probability 1/2 and $(1/2^{n-1})(1/2) = 1/2^n$. In such a case, A gets $2^n$ francs and increases his mathematical expectation by 1 franc[26]. This happens, however, at each of the $m$ tosses of the game so that the entire value of his mathematical expectation is 1 franc multiplied by $m$.

For the game to be fair, A ought to give $m$ francs to B; that is, give a thousand, a million francs, even an infinite sum if the game can be indefinitely prolonged. However, no one will risk a thousand francs, say, on such a game. Here, the rule of mathematical expectation seems to be defective. For eliminating this difficulty, which we indicate, the rule of moral expectation and its measure were introduced. We ought to remark, however, that that difficulty was caused by the conditions of the game which disregarded the possibility of B to pay all the moneys that the chance can deliver to A. No matter how great is B's fortune, it is necessarily limited; let it be *b* francs, then A will never receive a greater sum from B which decreases his mathematical expectation to a very large extent. Actually, we always have

$$b = 2^\beta(1 + h) \qquad (25.1)$$

where β is a natural number and *h* is a positive number less than unity. If $\beta \geq m$, B can pay any sum due A; otherwise, however, B will be unable to pay if *heads* first appears after more than β tosses. For those first tosses, the mathematical expectation of A is therefore equal to β, but for the rest $m - \beta$ tosses it becomes constant and equal to *b*, see formula (25.1), multiplied by their respective probabilities from $1/2^{\beta+1}$ to $1/2^m$. Denoting the composite value of A's mathematical expectation by ε, which is what he should give B for the game to become fair, we have

$$\varepsilon = \beta + \frac{1}{2}(1+h)[1 + \frac{1}{2} + \frac{1}{4} + \ldots + \frac{1}{2^{m-\beta-1}}] = \beta + (1+h)[1 - \frac{1}{2^{m-\beta}}].$$

This magnitude does not increase with *m*; on the contrary (? - O.S.), it is almost independent from that number and, when being very large, it is appreciably reduced to

$$\varepsilon = \beta + 1 + h.$$

B's fortune can never be sufficiently great for β to become considerable, and A should only risk a rather small sum between $\beta + 1$ and $\beta + 2$. If B is a banker having a hundred million francs, we will find out that the largest power of 2 in formula (25.1) is 26 so that A will really be at a disadvantage to bet 28 francs or more against the proprietor of that great fortune.

The rule of moral expectation applied to that problem (Laplace 1812/1886, p. 451) leads to a different sum which A can risk and which depends on the fortune of A rather than of B, but it seems to me that it is the possibility of B to pay in full that should restrict the sum which A must give him in advance.

**26.** I conclude this chapter by some remarks on the influence of a chance favourable for an event without knowing which one[27]. It always increases, as will be seen, the probability of repeated events in a series of trials. Thus, when playing *heads* and *tails*, it is always possible to believe that the coin, according to its physical constitution, tends to fall on one side rather than on the other. We do not know in advance whether that circumstance favours the occurrence of *heads* or *tails* which does not prevent the heightening of the probability of the same side to arrive many times in succession.

For showing it, let us denote the chance of the favoured side by $(1 + δ)/2$; consequently, the chance of the other side will be $(1 − δ)/2$. Here, δ is a small positive fraction whose value is unknown, and we do not know which of these two chances corresponds to *heads* or *tails*. If the coin should be tossed only once, there will be no reason to believe that the side chosen by one of the gamblers is more favoured or less, and the probability of its occurrence is 1/2, as though δ = 0. If, however, there ought to be two tosses, it is advantageous to bet on the coincidence of the arrived sides. Four combinations are possible, two for coincidence (*heads, heads*; or, *tails, tails*) and two for dissimilitude. The chances of the first two outcomes are $[(1 + δ)/2]^2$ and $[(1 − δ)/2]^2$, so that the probability that one of them takes place is, by the rule of § 10, their sum, $(1 + δ^2)/2$. The chances of the two other combinations are the same, each of them being expressed by $[(1 + δ)/2][(1 − δ)/2]$, and their sum is $(1 − δ^2)/2$ which is less than the previous sum in the ratio of

$$\frac{1 - δ^2}{1 + δ^2} = 1 - \frac{2δ^2}{1 + δ^2}.$$

If, in a fair game, A bets a franc against B on coincidence, B should bet 1 franc less $2δ^2/(1 + δ^2)$; and if δ = 1/10, less by about 2 *centimes*. When the coin should be tossed three times in succession, 8 different combinations will be possible: three times *heads*, and three times *tails*, and 6 more, three of them consisting of 2 *heads* and 1 *tails*, and the three last ones, of 2 *tails* and 1 *heads*. Supposing that δ is exactly zero, the chances of those 8 combinations will be equal one to another. Therefore, A, betting always on coincidence, should stake three times less than B; however, doubtless δ ≠ 0 and that proportion of chances will even be more advantageous for A than in the case of two tosses. The probability of coincidence will be

$$[\frac{1}{2}(1+δ)]^3 + [\frac{1}{2}(1-δ)]^3 = \frac{1}{4}(1 + 3δ^2).$$

[…] This reasoning can easily be generalized on more than three trials, and, if desired, on other games in which more than two events are possible with their unknown chances possibly being unequal.

Suppose that in a game for two the gamblers' skills somewhat influence the result. Equality of skills is unlikely, and if it is unknown who plays better, we should bet on the same gambler winning the two first games. But even if we know who is more skilful, it is not always advantageous to bet on his winning both these games. Indeed, four combinations can take place, three of which are unfavourable to him and only one is favourable. And although that last-mentioned is the most probable, his chance probably will not offset the three others taken together.

Let *p* be the known probability of event E of some nature, and *q*, the probability of the contrary event F, so that $p + q = 1$. Suppose that some cause can increase the chance of one of those events, which remains unknown, and at the same time decrease the chance of the other by an unknown fraction α. Denote by *w* the probability that in *m* trials one of those events will always appear. If E is the favoured event, the probability of coincidences of an outcome during the *m* trials will be, according to the rule of § 10,

$$(p + \alpha)^m + (q - \alpha)^m$$

since the coincidence can happen in two different ways with either E or F arriving all the time. If, however, the favoured event is F, then the similar probability will be

$$(p - \alpha)^m + (q + \alpha)^m.$$

We do not know the chance of which event had increased or decreased and for us those two different values are equally possible with the probabilities of each being 1/2. By the rule of § 10 the composite probability of coincidence will be

$$w = (1/2)[(p + \alpha)^m + (q - \alpha)^m + (p - \alpha)^m + (q + \alpha)^m] = P + Q,$$
$$P = p^m + C_m^2 p^{m-2}\alpha^2 + C_m^4 p^{m-4}\alpha^4 + ...,$$
$$Q = q^m + C_m^2 q^{m-2}\alpha^2 + C_m^4 q^{m-4}\alpha^4 + ...$$

If, however, $\alpha = 0$, the probability of coincidence will simply be $(p^m + q^m)$ so that any cause that increases the chance of one of the two contrary events E and F without it being known which one of them, increases the probability of coincidences of the events in a series of trials since the value of *w* evidently becomes larger than $(p^m + q^m)$.

### Notes

**1.** Poisson invariably applied both terms, *chance* and *probability*, and in a letter of 1836 to Cournot (1843/1984, p. 6) indicated that he *beaucoup insisté* on distinguishing between them. Later authors have however adopted much more convenient terms, *theoretical* and *statistical* probabilities whereas subjective probabilities have much less significance than Poisson attached to them. Moreover, his attitude was fraught with difficulties, see Sheynin (2002).
**2.** Poisson repeated this restriction at the end of § 3. De Moivre (1712/1984, p. 237) introduced a more general definition adopted by Poisson (§ 13) as well.
**3.** Here and in some other instances Poisson's terminology is dated, but understandable.
**4.** In the theory of information, probabilities equal to 1/2 are also tantamount to complete ignorance.
**5.** According to Chebyshev (1845/1951, p. 29), the theory of probability aims at determining the probability of an event given its connection with other events whose probabilities are known. Laplace did not offer any similar definitions. Just below, Poisson indirectly stated that really important are cases in which a probability *more or less approaches unity*.
**6.** This statement ought to be noted.
**7.** Poisson indirectly defined independence of events in the beginning of § 9. The direct definition is due to De Moivre (1718/1738, p. 6).
**8.** François Callet published a table of logarithms with seven digits (Paris, 1795).
**9.** Calculations with superfluous digits had been universally practised up to the mid-20th century.
**10.** A similar statement was due to Laplace (1776/1891, p. 152; 1814/1995, p. 9).
**11.** De Moivre (1718/1756, pp. 37 – 38) had stated the same.
**12.** Poisson discussed a game of chance, but introduced skill!

**13.** More precisely, the compound event should be E + (F + G + …) so that F′ = F + G + … without the E.

**14.** The sequence $a_0 + a_1 t + a_2 t^2 + … + a_m t^m$ is called the generating function of $\{a\}$.

**15.** Poisson introduces a sufficient condition.

**16.** After μ drawings with $m$ white and $n$ black balls having arrived, the chance of extracting a white ball in a new drawing will depend on those numbers and equal $a'/c'$. However, for a person who only knows that μ balls were extracted and does not know the ratio of the arrived white and black balls, that chance is quite different. According to a note sent me by Mondesir, a graduate of the Ecole Polytechnique, the required probability is independent from $m$ and $n$ and equals $a/c$ as it was previous to the drawings.

For verifying this proposition by an example, let

$a = 4$, $b = 3$, $c = 7$, μ = 2 and $c' = 5$.

Concerning $m$ and $n$ there are three possible but unequally probable cases: $m = 2$ and $n = 0$; $m = 1$ and $n = 1$; and $m = 0$ and $n = 2$. Their probabilities derived from the expression of Π are 2/7, 4/7 and 1/7. And the chances of extracting a white ball in the next drawing will be 2/5, 3/5 and 4/5. By the rules of §§ 5 and 10 the composite probability of such an extraction will be

$$\frac{2}{7} \cdot \frac{2}{5} + \frac{4}{7} \cdot \frac{3}{5} + \frac{1}{7} \cdot \frac{4}{5} = \frac{4}{7} = \frac{a}{c}.$$

For the general demonstration I refer to the note of Mondesir [see Bibliography]. That [Mondesir's] proposition is evident if $a = b$ since in that case, for a person who does not know which balls had arrived previously, there is no reason to believe, either before or after that result, in the occurrence of a white rather than a black ball. Consequently, the probability of drawing a white ball always remains 1/2. We can also note that, if $a$ and $b$ are infinite, that proposition is in agreement with another one which will be proved below and according to which it is certain that $m:n = a:b$. And so, we are assured in that $a'$ and $b'$, the numbers of balls left in the urn, are also in the same ratio; the chance and the probability of the arrival of a new white ball do not anymore differ and $a'/c' = a/c$. Poisson

Poisson considered subjective probabilities, see Sheynin (2002). O. S.

**17.** Was it really necessary to prove this statement?

**18.** The bivariate generating function of the probabilities of E and F.

**19.** Laplace's calculations were extremely complicated. He solved the same problem in his early memoirs of 1774 and 1786. De Moivre (1712/1984, Problems 18 and 19; 1718/1756, Problem 39) investigated the appearance of any outcomes in a throw of an arbitrary number of dice with an arbitrary number of faces.

**20.** Montmort (1708, Préface) already mentioned both these useless *systems*.

**21.** Bertrand (1888, p. XXII) formulated this idea in the best possible way: *la roulette n'a ni conscience, ni mémoire*.

**22.** *Jeux publics de Paris*, as might be thought, were regularly played games in registered casinos.

**23.** This section was written before the latest issued law on finances suppressed games of chance beginning 1 January 1838. Poisson

**24.** Laplace (1812/1886, p. 189) replaced the term *expectation* by mathematical expectation to distinguish it from *moral expectation* which had come into vogue. His specification took root in French and Russian literature, but it became unnecessary long ago.

**25.** This statement is obviously lame.

**26.** A strange demonstration! The described *Petersburg game* (Montmort 1708/1713, p. 402) was due to Niklaus Bernoulli. It became generally known after Daniel Bernoulli (1738), who published his memoir in Petersburg, suggested to solve its paradox by means of the moral expectation.

Poisson considered it superficially. Condorcet and Lacroix suggested that even an infinite Petersburg game was only one single trial, and that a set of such games was necessary for studying the game. Freudenthal (1951) proposed the same and additionally recommended that at each game the gamblers should choose their roles by lot.

**27.** Laplace (1812, Chapter 7) and earlier, in his memoirs, considered that problem in lesser detail.

## Chapter 2. General Rules, Continued. Probabilities of Causes and Future Events Derived from Observing Past Events

### Misprints/Mistakes Unnoticed by the Author

**1.** In § 34, p. 95 of the original text, line 2 from bottom. Probability $p_1$ should be $p_n$.

**2.** In § 37, p. 101 of the original text, formulas for $y_1$ and $y_2$ in a single displayed line. In the second formula, in the second term of the denominator $y_2$ is printed instead of $y_1$. There also, when passing from $y_x$ to $y_{x-1}$, $(x-1)$, rather than $(y-1)$ as printed, was substituted instead of $x$.

**3.** In § 38, p. 104 of the original text. On line 10 of the second paragraph $n'$ is printed instead of $n$. In the same section, p. 106 of the original text, formulas for $w_n$ and $w_i$ in a single displayed line. In the second formula, in the first term of the denominator, $q_n p_n$ is printed instead of $q_i p_i$.

4. In § 39, p. 109 of the original text, line 6. The number of balls is stated as $a_n$ and $a_x$ instead of $a_n$ and $a_i$.

**5.** In § 40, p. 113 of the original text. The particular case of the displayed formula coincides with $p_n$ of § 38 rather than of the preceding section, as printed.

**6.** At the end of § 45, p. 124 of the original text, the reference to § 27 is printed instead of to § 26.

**7.** In § 46, p. 125 of the original text. The integrand $x^m(1-x)^m$ should be $x^m(1-x)^n$.

**8.** In § 53, p. 141 of the original text. In proving the *third main proposition*, magnitude A is stated as taking μ, instead of taking λ values.

**9.** In § 55, p. 145 of the original text. Extraction of balls from some urns. The number of drawings is printed as being equal to $m$ instead of μ.

**10.** In § 56, p. 149 of the original text. Two inequalities are printed twice instead of changing their signs in the second case.

**11.** In § 60, p. 156 of the original text. Just above the last displayed line $α + u$ is substituted instead of $z$ rather than instead of 2, as printed.

All those mistakes/misprints are corrected in the translation.

**27.** The rules provided in the preceding chapter assumed that the chances of certain events were known and aimed at deducing the probabilities of other events composed of those given[1]. Here, I put forward rules for calculating probabilities of causes given observed events and, after that, probabilities of future events. However, before going on, it is convenient to elucidate the exact sense which we attach to the word *cause* and which differs from that in ordinary language.

Usually, when stating that a thing is the *cause* of another thing, we attribute to the former a *possibility* of necessarily producing the latter, without, however, wishing either to imply that we know the nature of that power or how is it exercised. At the end of this chapter we return to this notion of *causality*, but at present suffice it to say that for the calculus of probability the word *cause* has a more general meaning.

We consider a *cause* C relative to some event E as a thing that attaches to its arrival a determined chance properly belonging to it[2]. In ordinary language, C is the cause of that chance, not of the event itself. And when E actually happens, it occurs because C coincided with other causes or circumstances which do not influence E's proper chance.

Let *p* be that known or unknown chance generally differing from probability; at the same time, C provides chance $1 - p$ to the contrary event F. When $p = 1$, the thing C necessarily produces the event E and is its cause in the restricted sense. When $p = 0$, C is a similar cause of F.

The set of causes which combine in the production of an event without influencing the magnitude of its chance, or the ratio of cases favourable for its appearance to all possible cases, is what should be understood as *hazard* (hasard)[3]. Thus, in a dice game the event occurring at each throw is the consequence of the number of the faces of the die, of the possible irregularities of form and density of the dice, and of the numerous shakes to which they are subjected before being thrown. These shakes are the causes which do not at all influence the chance of the occurrence of any fixed face; they aim at eliminating the influence of the dices' preliminary position in their box, so that that position will not be known anymore to the gamblers. And if that aim is attained, the relative chance of the arrival of each face will only depend on their number and the dice's defects which can lead to an inequality of the chances of different faces.

It is said that a thing happened randomly if it was executed without at all changing the respective chances of various possible events. An urn contains white and black balls; a ball is drawn at random if their arrangement in the urn was not taken into account. Supposing that all the balls had the same diameter[4], the chance of extracting a white ball can evidently only depend on the numbers of white and black balls and proved to be equal to the ratio of the former to the sum of both numbers.

A cause C can be a physical or a moral thing. At the game of *heads* and *tails*, it is the physical constitution of the coin that results in the chances of *heads* and *tails* generally differing a little from 1/2. When deciding a criminal case in court, the chance of a correct or mistaken vote of each juryman is determined by his moral principles, i. e., his own ability and conscience as applied to that case. Sometimes a cause C results from the coincidence of a moral and physical thing. Thus, in each kind of measurement or observation the chance of making an error of a given magnitude depends on the ability of the observer and the more or less perfect construction of his instrument[5]. However, in the calculus of probabilities the diverse causes of events are invariably considered independently from their particular nature, only regarding the magnitude of the chances they produce. And it is for this reason that that calculus is equally applicable to moral and physical things.

However, in most problems the chance determined by a given cause C is not known in advance and the cause itself or its chance is sometimes unknown either. If the chance is constant, it is determined, as is seen below, by a sufficiently long series of trials, but when jurymen vote the chance of error varies from one of them to another and doubtless differs for the same juryman in different cases. Repetition of trials concerning each juryman and each kind of cases is impossible, so that it is not the chance of error proper to each juryman that can be derived from observation, but, as it will be seen below, a certain probability concerning the set of all of them in the jurisdiction of an assize court, and this knowledge is sufficient for solving the problems which constitute the special goal of this work.

There often are many different causes which, when combined with randomness, can bring about a given event E or the contrary event F. Before one of those events takes place, each of these causes has a certain probability that changes according to which of the events E or F was observed. Supposing that the chance which each of those possible causes, if they are certain, provides to the arrival of E or F is known, we determine, first of all, the probabilities of all those causes existing after the observation, and then the probability of any other future event depending on the same causes as E and F.

**28.** Let E be the observed event. We suppose that its arrival can be attributed to $m$ distinct and only possible incompatible causes, equally probable before the observation. The occurrence of E renders those hypothetical causes unequally probable, and it is required to determine the probability of each resulting from the observation. This is achieved by means of the following theorem.

The probability of each of the possible causes of an observed event, if being certain, is equal to that which it provides to the event divided by the sum of the probabilities of that event resulting from all the causes which can be attributed to it.

And so, denote the $m$ possible causes of event E by $C_1, C_2, \ldots, C_m$. Let $p_1, p_2, \ldots, p_m$ be the known probabilities of its occurrence relative to those diverse causes so that $p_n$ is the probability that E will take place if the cause $C_n$ is unique; or, which is the same, if it is certain, so that all other causes are excluded. Also denote by $w_1, w_2, \ldots, w_m$ the unknown probabilities of those same causes so that $w_n$ is the probability of $C_n$; or, in other words, the probability that the arrival of E was due to that cause. It is required to prove that

$$w_n = p_n/(p_1 + p_2 + \ldots + p_n + \ldots + p_m).$$

No matter what is the nature of E, we may liken it to an extraction of a white ball from an urn containing white and black balls. Suppose then that there are $m$ such urns, $A_1, A_2, \ldots, A_m$, from which a white ball can be drawn with the ratio of white balls to the total number of balls in urn $A_n$ being $p_n$. A white ball can be extracted from each of those randomly chosen urns, and each of them represents one of the causes of its arrival; urn $A_n$ corresponds to cause $C_n$. The problem consists in determining the probability that the white ball was drawn from $A_n$.

Suppose that all fractions $p_1, p_2, \ldots, p_m$ are reduced to one and the same denominator:

$$p_1 = \alpha_1/\mu, \ p_2 = \alpha_2/\mu, \ \ldots$$

Here, μ and α₁, α₂, … are natural numbers. The chance of extracting a white ball from urn $A_n$ will not change at all if the balls contained there are replaced [represented] by numbers $α_n$ of white, and μ of all the balls, and a similar statement holds for all other urns. The total number of balls in all the urns will not change. It follows from the lemma of § 10 that if all the balls are united in one and the same urn A, and if those that had come from urn $A_1$ are provided with number 1; from urn $A_2$, with number 2 etc., the probability $w_n$ that a white ball drawn from the set $A_1, A_2, …, A_m$ originated in $A_n$ is the same as when it was extracted from A bearing the number $n$. This probability is equal to the ratio of $α_n$ to the sum of $α_1, α_2, …, α_m$ since that sum is the total number of white balls contained in A, and $α_n$ balls had number $n$. Therefore,

$$w_n = α_n/(α_1 + α_2 + … + α_m)$$

and coincides, according to the preceding equations, with $w_n$ which is what we should have derived.

**29.** When calculating the probability of many consecutive events, it is necessary to take into account the possible influence of the arrival of one of them on the chances of the next one (§ 9). And sometimes, when evaluating that chance, the probabilities of the diverse causes of the preceding event, or of the different ways in which it can take place, should also be allowed for. This will be seen for example in the following problem.

I suppose that there are $m$ urns A, B, C, … containing white and black balls and that the chances of drawing a white ball from them are $a, b, c, …$ A ball is extracted at random from one of those urns; then a second ball from one of the urns excepting that from which the first ball was drawn; then a third ball from an urn differing from those two, etc. This means that, after each drawing, the urn from which a ball was extracted is excluded. It is required to determine the probability of drawing $n$ white balls in $n$ extractions, $n ≤ m$.

Denote for brevity

$$a + b + c + d + … = s_1, ab + ac + ad + bc + bd + cd + … = s_2,$$
$$abc + abd + bcd + … = s_3, abcd + … = s_4, …$$

Here, $s_1$ is the sum of $a, b, c, …$; $s_2$ is the sum of their products with $a, b, c, …$ taken two at a time with their number being $C_m^2$; $s_3$, the sum of their products with $a, b, c, …$ taken three at a time with their number being $C_m^3$ etc. The probability of drawing a white ball at the first extraction is $s_1/m$. If the ball, be it white or black, was drawn from A, the probability of extracting a white ball at the second trial will be $(s_1 − a)/(m − 1)$, or $(s_1 − b)/(m − 1)$ had the first ball come from B etc. By the rules of §§ 9 и 10

$$\frac{1}{m-1}[α(s_1 - a) + β(s_1 - b) + γ(s_1 - c) + ...]$$

is the composite probability of the occurrence of a white ball at the second drawing. Here, α, β, γ, … are the probabilities that it was first extracted from A, B, C, … These probabilities are certainly different[6]. According to what was shown in § 28, we have $α = a/s_1$, $β = b/s_1$, $γ = c/s_1$, … and

$$a(s_1 - a) + b(s_1 - b) + c(s_1 - c) + \ldots = 2s_2,$$

so that the probability of extracting a white ball at the second drawing is $2s_2/(m-1)s_1$. Similarly, the probability for the same occurring at the third trial is $(s_1 - a - b)/(m-2)$ if the first two white or black balls came from A and B, $(s_1 - a - c)/(m-2)$ if they were extracted from A and C, etc.

Therefore, the composite probability of a white ball arriving at the third trial is

$$\frac{1}{m-2}[g(s_1 - a - b) + h(s_1 - a - c) + k(s_1 - b - c) + \ldots].$$

Here, $g, h, k, \ldots$ are the probabilities that the two first white balls came from A and B; from A and C; from B and C, etc. According to § 28, these probabilities are

$g = ab/s_2$, $h = ac/s_2$, $k = bc/s_2$, …

and also

$$ab(s_1 - a - b) + ac(s_1 - a - c) + bc(s_1 - b - c) + \ldots = 3s_3.$$

Therefore, the probability of the occurrence of a white ball at the third drawing is $3s_3/[(m-3)s_2]$. This reasoning can easily be continued as far as desired. And the result is that

$$\frac{s_1}{m}, \frac{2s_2}{(m-1)s_1}, \frac{3s_3}{(m-2)s_2}, \ldots, \frac{ns_n}{(m-n+1)s_{n-1}}$$

are the probabilities of the arrival of white balls in each of the $n$ first drawings. The required probability is therefore the product of these $n$ fractions, see § 5, $s_n/C_m^n$.

This answer can be verified by noting that each of the products is the probability of extracting $n$ white balls from $n$ fixed urns randomly chosen among A, B, C, … so that the sum of all those products divided by their number is the probability of extracting $n$ white balls from $n$ randomly chosen urns. That probability is evidently the same as the required. If $n = m$, then $C_m^m = 1$ and the probability is $s_m$ which immediately follows from the rule of § 5.

**30.** Let E′ be another event differing from E but depending on the same causes $C_1, C_2, \ldots$ and denote the chances of E′ resulting from those diverse causes by $p_1', p_2', \ldots, p_m'$, so that $p_n'$ is the known probability of the arrival of E′ if $C_n$ is certain. That cause becomes the only probable and its probability is $w_n$. The arrival of E′ due to cause $C_n$ will be a compound event whose chance is the product of those two probabilities (§ 5). Then, also, the composite probability $w'$ of E′ will be the sum of those chances relative to the *m* different ways in which it can occur (§ 10), i. e., the sum of $p_n'w_n$,

$$w' = p'_1 w_1 + p'_2 w_2 + \ldots + p'_n w_n + \ldots + p'_m w_m.$$

Substituting the values of $w_1, w_2, \ldots$, we get

$$w' = \frac{p_1 p'_1 + p_2 p'_2 + \ldots + p_n p'_n + \ldots + p_m p'_m}{p_1 + p_2 + \ldots + p_n + \ldots + p_m}.$$

This is the formula for calculating the probability of future events given observations of past events. The same expression can be derived without introducing common causes of E and E′ by considering them as compound events depending on the same simple event. And our reasoning is equally applicable to that other manner of contemplating our problem. However, if desired, it is possible to return immediately to the preceding.

Actually, suppose that E and E′ are two events composed from the same event G susceptible of chances $g_1, g_2, \ldots, g_n, \ldots, g_m$, all of them being equally probable before E is observed. Then they can be considered as so many distinct causes of E and E′. Therefore regard $g_n$ as the cause called $C_n$ in the reasoning above and the probability of $g_n$ will be the value of $w_n$ as already deduced. In other words, $w_n$ will be the probability that the chance of G is $g_n$, and the preceding expression of $w'$ will be the probability of the arrival of E′ resulting from the *m* possible values of the chance of G. In that formula [for $w'$], $p_n$ and $p_n'$ express the given probabilities of the arrival of E and E′ provided that $g_n$ is certainly the chance of G.

**31.** That determination of the probability of E′ given the observations of E should not be confused with some influence of the arrival of past events on future events, which is absurd to suppose. If I am sure, for example, that an urn A contains 3 white balls and 1 black ball, it is certain for me that the chance of drawing a white ball is 3/4. Therefore, if E′ is the extraction of two white balls from A with replacement of the first one, the chance of E′ will be $(3/4)^2 = 9/16$ independently from the possibly observed event E. And supposing that E is the appearance of a certain number of white and black balls successively drawn from A and returned there each time, I should always, without taking into account the ratio of those two numbers, bet 9 against 7 on the appearance of E′.

However, if the chance of the simple event G is unknown, and I only know that it is susceptible of certain values, the observation of E will provide the probabilities of each and I can then derive the probability of E′. That observation increases or decreases the reason I have to believe in the arrival of E′ without influencing at all that future event or the chance proper to it. For someone who observes another event $E_1$ depending on the same simple event G the reason to believe in the arrival of E′ can be much stronger or much weaker than for me, but this will not at all change the chance proper to E′.

Concerning that case of two people who observed E and $E_1$ respectively, both these events composed of the same event G, we should not forget that if $E_1$ includes E and something else as well, the opinion of the second person about the arrival of a new event E′, also depending on G, will be more informed than that of the first one and should be preferably adopted (§ 1).

Supposing that the observation of $E_1$ and E provide probabilities $k$ and $h$ to the future event E′, the second person will be more justified to bet $k$ against $(1 - k)$ than the first one to bet $h$ against $(1 - h)$ on the arrival of E′ whether the fractions $h$ and $k$ be larger or smaller than 1/2, and the difference $(h - k)$ positive or negative.

**32.** Before going ahead, it is appropriate to provide some simple examples of the use of the preceding expressions of $w_n$ and $w'$ which we will write down in an abbreviated form

$$w_n = p_n / \sum p_n, \quad w' = \sum p_n p_n' / \sum p_n$$

where $\sum$ indicates a sum from $n = 1$ to $m$.

It is known that an urn B contains $m$ white or black balls. A white ball is drawn, and it is required to determine the probability that the urn contains $n$ white balls. We can formulate $m$ different hypotheses about the number of those balls in the urn: $m$ white balls; $(m - 1)$ white balls and 1 black ball; $(m - 2)$ white and 2 black balls; …; 1 white ball and $(m - 1)$ black balls. All of them are equally possible and mutually incompatible; they can be considered as $m$ causes $C_1, C_2, \ldots$ of event E, the drawing of a white ball from B.

If $n$ white balls are among the $m$ balls in B, then the probability of such a drawing will be $p = n/m$ so that

$$\sum p_n = \frac{m+1}{2}, \quad w_n = \frac{2n}{m(m+1)}.$$

That probability can only be 1/2 if $m = n = 3$. In general, the probability that B contains only white balls, or that $n = m$, after one such ball had arrived, is $2/(m + 1)$. If E′ is the extraction of a new white ball, its probability $w'$ will differ depending on whether the first one was, or was not replaced. In the first case

$$p_n' = p_n = \frac{n}{m}, \quad \sum p_n p_n' = \frac{1}{m^2} \sum n^2,$$

but, as is known,

$$\sum \frac{n(n+1)}{1 \cdot 2} = \frac{m(m+1)(m+2)}{3!}, \quad \sum n = \frac{m(m+1)}{1 \cdot 2}$$

and therefore

$$\sum n^2 = 2\sum \frac{n(n+1)}{1\cdot 2} - \sum n = \frac{m(m+1)(2m+1)}{3!}, \quad w' = \frac{2m+1}{3m}.$$

In the second case, the number of the white balls and the total number of balls in B had decreased by unity, so that

$$p'_n = \frac{n-1}{m-1}, \quad \sum p_n p'_n = \frac{1}{m(m-1)} \sum n(n-1).$$

However, we always have $p_n = n/m$, $\sum p_n = (m+1)/2$ and $w' = 2/3$ since

$$\sum \frac{n(n-1)}{1\cdot 2} = \frac{(m-1)m(m+1)}{3!}.$$

The probability of drawing a white ball from an urn after a ball of that colour had already been extracted from it without replacement is therefore independent from the number $m$ of white or black balls in the urn and always equal to 2/3. The value of $w'$ in the first case is also reduced to 2/3 as it should be if $m$ is a very large number considered infinite.

If it is known that $(m-1)$ white balls have been drawn out of the $m$ white or black balls initially contained in B, there will be probability $m/(m+1)$ that the last ball is also white. We can only formulate two hypotheses, $C_1$ and $C_2$, that all the $m$ balls are white, or that one is black. According to $C_1$, the probability of the observed event is unity; according to $C_2$, it, or the chance of extracting $(m-1)$ white balls from B, is the same as when one black ball is left; and since that last ball could have initially been any of the $m$ balls in B, the probability that it is black equals $1/m$. Therefore, $p_1 = 1$, $p_2 = 1/m$, $w_1 = p_1/(p_1+p_2) = m/(m+1)$ for the probability of the first hypothesis; that is, for the probability that the last ball is white like all those extracted. The case of $m = 1$ is not included in that value of $w_1$ which is then equal to 1/2, as was evident from the very beginning. [The author's correction of the phrase in the text is taken account of in the translation.]

**33.** Here is one more immediate application of the preceding formulas in which we do not know the total number of white or black balls in an urn B. We only know, for example, that that number can not exceed three. The observed event is the appearance of $x$ white balls in a series of $n$ drawings with replacement. If $x \neq 0$ and $\neq n$, we can only formulate three hypotheses about the balls contained in B: $C_1$, there is 1 white and 1 black ball; $C_2$, 2 white balls and 1 black ball; $C_3$, 2 black balls and 1 white ball. The probabilities of E with respect to these three distinct causes are

$$p_1 = (1/2)^x (1/2)^{n-x} = 1/2^n, \quad p_2 = (2/3)^x (1/3)^{n-x} = 2^x/3^n,$$
$$p_3 = (1/3)^x (2/3)^{n-x} = 2^{n-x}/3^n.$$

Denote for the sake of brevity

$$3^n + 2^{n+x} + 2^{2n-x} = \mu,$$

then $w_1 = 3^n/\mu$, $w_2 = 2^{n+x}/\mu$, $w_3 = 2^{2n-x}/\mu$ will be the probabilities of $C_1$, $C_2$, $C_3$. Let the future event be $E'$, an extraction of a new white ball. Its probabilities relative to those three hypotheses will be $p'_1 = 1/2$, $p'_2 = 2/3$, $p'_3 = 1/3$, and the composite probability

$$w' = \frac{(1/2)3^n + (2/3)2^{n+x} + (1/3)2^{2n-x}}{3^n + 2^{n+x} + 2^{2n-x}}.$$

If $n = 2x$

$$w_1 = \frac{9^x}{9^x + 2 \cdot 8^x}, \quad w_2 = w_3 = \frac{8^x}{9^x + 2 \cdot 8^x}, \quad w' = \frac{(1/2)9^x + 8^x}{9^x + 2 \cdot 8^x} = \frac{1}{2}.$$

So $w' = 1/2$ as it should have been. Indeed, white and black balls have arrived the same number of times, and there is no reason to believe in the occurrence, in a new drawing, of a white rather than of a black ball. Nevertheless, the ratio $9^x/8^x$, should be larger than 2, or $x > 5$ for being able to bet more than one against one on the equality of the numbers of white and black balls in the urn or that it contains 1 white and 1 black ball. If $x$ is a very large number, the probability $w_1$ of that hypothesis will very little differ from certainty.

If $i$ is a natural number, $x = 2i$ and $n = 3i$, then

$$w' = \frac{(1/2)27^i + (2/3)32^i + (1/3)16^i}{27^i + 32^i + 16^i}.$$

For very large values of $i$ this expression very little differs from 2/3. At the same time, probability $w_2$ that B contains 2 white balls and 1 black ball also very little differs from certitude. Suppose now that $n = 3x$. Then

$$w' = \frac{(1/2)27^x + (2/3)16^x + (1/3)32^x}{27^x + 16^x + 32^x}.$$

If $x$ is very large, this expression will be almost equal to 1/3 and the probability $w_3$ that B includes 1 white ball and 2 black balls is almost unity. In the third case, when the number of drawings is supposed to be very large, the probability $w'$ of the arrival of a new white ball becomes very near to the ratio of the number of white balls extracted from B to the total number of trials, and at the same time with probability closely approaching certainty, that ratio is also equal to the ratio of the number of white balls to the total number of balls in B, i. e., to the chance proper to the drawing of a white ball from that urn.

It will be seen below, that when an event of some nature is observed a certain number of times in a very large number of trials, the ratio of the first number to the second is a very close and likely value of the known or unknown chance of that event. In our example, this chance can only be 1/2, 2/3 or 1/3, so that $x/n = 1/2$, 2/3 or 1/3 are the only values which can be thought likely when $x$ and $n$ are very large.

**34.** In the preceding, we have supposed that before the arrival of E all the causes $C_1, C_2, \ldots$ which can be attributed to that event were equally possible. However, if there exists some prior reason to believe in the existence of one of those causes rather than of another one, then, when evaluating the probabilities which those diverse causes will acquire after the arrival of E, it will be necessary to account for that inequality of the chances of $C_1, C_2, \ldots$ That necessity is an important point of the theory of probability[7], especially, as described in the Preamble, in problems pertaining to the judgements of tribunals. And the demonstration in § 28 is in addition easy to extend on the general case in which before the observation the causes of E had some known probabilities. Actually, like in that section, let us replace the event E by a drawing of a white ball from one of the urns $A_1, A_2, \ldots$ and first of all suppose that its extraction from each of them is equally possible. The probability that the ball has arrived from $A_n$ will be $p_n/\sum p_n$ where, as always, $p_n$ is the ratio of the number of white balls to the total number of balls in $A_n$ and the sum extends over all the urns. Similar expressions take place for all the other urns and, by the rule of § 10, the probability that the white ball is drawn from one of them is the sum of all such expressions. When the partial probabilities are equal one to another, that sum reduces to one of such fractions multiplied by their number.

Let now $p_1$ be the same ratio as above for $a_1$ urns $A_1$; $p_2$, the same ratio for $a_2$ urns $A_2, \ldots, p_i$, the same ratio for $a_i$ urns $A_i$ where $i$ expresses the number of groups of similar urns. If the number of all the urns is $s$, then

$$s = a_1 + a_2 + \ldots + a_i.$$

The sum $\sum p_n$ extending over all the urns can be replaced by $\sum a_n p_n$ extending over all the groups from $n = 1$ to $n = i$. And if some group consists of $a_n$ urns, the probability that the white ball arrived from it will be

$$w_n = a_n p_n / \sum a_n p_n.$$

However, before the observation the probability of drawing a ball from that group was evidently $q_n = a_n/s$, so that $a_n = sq_n$. Substituting this value of $a_n$ in $w_n$, and suppressing the common factors in the numerator and denominator, we will therefore have

$$w_n = q_n p_n / \sum q_n p_n.$$

Now, the different groups of urns which we consider represent all the $i$ possible and initially unequally probable causes $C_1, C_2, \ldots$ of event E. The fraction $q_n$ expresses the probability existing before the observation that the event will be due to cause $C_n$. After the observation, its arrival due to the same cause will have probability $w_n$. Since the causes $C_1, C_2, \ldots$ are mutually incompatible, $q_n$ and $w_n$ are the probabilities that that cause exists. The latter proves that the probability of each of the possible causes of the observed event is equal to the product of $q_n$ and the probability $p_n$ which it provides, if being certain, divided by the sum $\sum q_n p_n$ of such products relative to all the causes to which the event can be attributed.

The probability $w'$ of a future event E′ depending on the same causes as E, will be, as above, $\sum w_n p'_n$ where $w_n$ is its determined value, so that

$$w' = \frac{\sum q_n p_n p'_n}{\sum q_n p_n}.$$

Suppose that E′ was also observed after E, that E″ is the third event depending on those same causes, and that $p''_n$ is the chance that the cause $C_n$, if certain, provides to the future arrival of E″. According to the preceding rule, the probability of that cause, $w_n$, existing between observations of E and E′, becomes

$$\frac{w_n p'_n}{\sum w_n p'_n} = \frac{q_n p_n p'_n}{\sum q_n p_n p'_n}.$$

Multiplying it by $p''_n$, we get the probability of the arrival of E″ due to the cause $C_n$. The composite probability of that event is

$$w'' = \frac{\sum q_n p_n p'_n p''_n}{\sum q_n p_n p'_n}.$$

This expression can also be derived like $w'$ was by substituting $p''_n$ and $p_n p'_n$ instead of $p'_n$ and $p_n$. Actually, this product $p_n p'_n$ is the chance of the observed event, i. e., of the succession of E and E′, relative to the cause $C_n$.

**35.** For providing a very simple example of the preceding rule, which can also verify its exactitude and necessity, I suppose that two cards whose colour is unknown are lying on a table. Turning over one of them, I see that it is red. Only two hypotheses can be formulated about their colour: red, red; and red, black. If I do not at all know their provenance, prior to observation those hypotheses are equally probable. After the observation, the probability of the first one is 2/3, as seen in one of the examples in § 32. We can bet 2 against 1 on the second card to be also red. But it will not be the same anymore if it is known that, for example, these cards were chosen at random from [a pack for] the game of *piquet* consisting of 16 red and 16 black cards.

Before the observation we had (§ 18) the probabilities of the two hypotheses $q_1 = 16 \cdot 15/32 \cdot 31$, $q_2 = 2 \cdot 16 \cdot 16/32 \cdot 31$, and at the same time $p_1 = 1$, $p_2 = 1/2$ so that the probability of the first hypothesis after the observation is

$$w_1 = \frac{q_1 p_1}{q_1 p_1 + q_2 p_2} = \frac{15}{16}.$$

On the contrary, instead of betting 2 against 1, we can only bet 15 against 16 on the same event. This value of $w_1$ can immediately be verified since it is evident that the problem is tantamount to requiring the probability of drawing one more red card from the pack after one such had been extracted and only 15 red cards had been left.

In general, if *n* cards, *a* of them red and *b*, black, are randomly extracted from a pack of *m* cards, and (*n* − 1) of them turned over, *a′* of them red, and *b′* black, then, by the preceding rule, we have

$$w_1 = \frac{a - a'}{m - n + 1}, \quad w_2 = \frac{b - b'}{m - n + 1}$$

for the probabilities of the *n*-th card to be red and black.

That value of $w_1$ is also, as it should be, the probability of drawing a red card from the initial pack reduced to (*m* − *n* + 1) cards, (*a* − *a′*) of them red, by extracting *a′* red and *b′* black. Since *a* + *b* = *m*, and *a′* + *b′* = *n* − 1, we have $w_1 + w_2 = 1$ so that $w_2$ can also be verified.

**36.** The general consequence of the rule of § 34 is that when two events, E and E′, depend on the same cause, the probability of the future event E′ does not only result from the observed E; when evaluating it, we should account for the possible prior knowledge about that common cause. The probability of E′ can differ for two people who observed the same event E but had different previous knowledge about that problem.

It is just the same in problems of doubt or criticism to which the calculus of probability is also applicable (§ 3). When required to find out whether a fact attested by a witness is true or false, we should allow for the chance of his error and, in addition, for our knowledge prior to his testimony. Let us denote the probability that the witness is not deceiving us, involuntarily or otherwise, by *p*, and, by *q*, the probability of the truth of the attested fact existing before his testimony. The probability of the fact after the testimony depends on *p* and *q* and is determined in the following way.

Here, the observed event will be the attesting of a fact, not incontestable at all. Supposing that it is true, and the witness did not deceive us, its probability is *p*. And it is (1 − *p*) if the fact is false because the witness deceived us. Before the testimony, *q* was the probability of the first hypothesis, and (1 − *q*), of the second. Denote the probability of the former, or the truth of the fact existing after the testimony, by *r*. Then, by the rule of § 34,

$$r = \frac{pq}{pq + (1-p)(1-q)}, \quad r - q = \frac{q(1-q)(2p-1)}{pq + (1-p)(1-q)},$$

which proves that the difference (*r* − *q*) has the same sign as (*p* − 1/2) so that the witness heightens or lowers the prior probability of the fact according to the inequalities *p* > 1/2 or < 1/2. That difference disappears if *p* = 1/2, if the witness did not change the prior probability at all. We can bet even money on whether he tells the truth or not. If initially there was no reason to believe in the truth of the attested fact rather than in its falsity, the probability *q* = 1/2 and *r* = *p*. In this case, the probability of the fact being true only depends on the truthfulness and knowledge of the witness.

It is impossible to suppose that one of the two magnitudes, $p$ and $q$, is unity, and the other, zero, but if $p$ closely approaches certainty, and $q$, more closely approaches impossibility, so that $q/(1 - p)$ becomes a very small fraction, the probability $r$ will also be very low, almost equal to that ratio. Such is a case of a fact contrary to the general laws of nature but attested by a witness to which, without taking into account that opposition, a high degree of confidence is attached.

For us, those general laws are the result of a long series of experiences which provides them, if not absolute certainty, a very high probability, still more strengthened by the harmony which they present and which no testimony can offset. So, if the attested fact is contrary to those laws, the prior probability that it is exactly described, is almost zero. And, even supposing that the witness is honest, suffice it that he is not at all infallible so that his chance of error $(1 - p)$ is extremely high as compared with that prior probability $q$, and the probability $r$ existing after the testimony can still be considered insensible. In such cases, it is reasonable to reject our own testimony and believe that our senses had deceived us by presenting as true a thing contrary to the laws of nature.

**37.** Suppose that the fact, whose probability we are considering, is additionally attested by a second witness[8]. Let the probability that he is not deceiving us by $p'$, and by $r'$, the probability that the fact, as testified by both witnesses, is true. Noting that, independently from the second testimony, that probability was already equal to $r$, we conclude that $r'$ should be derived from $r$ by substituting $p'$ and $r$ instead of $p$ and $q$:

$$r' = \frac{p'r}{p'r + (1-p')(1-r)} = \frac{qpp'}{qpp' + (1-q)(1-p)(1-p')}.$$

Suppose that the second witness attests that the fact is false and thus contradicts his predecessor and note that, independently from that, the falsity of the fact already had probability $(1 - r)$. Denote by $r_1$ the probability that the fact is false as results from both contrary testimonies. It should be deduced from $r$ of § 36 by substituting $p'$ and $(1 - r)$ instead of $p$ and $q$:

$$r_1 = \frac{p'(1-r)}{p'(1-r) + r(1-p')} = \frac{p'(1-p)(1-q)}{p'(1-p)(1-q) + qp(1-p')}.$$

If $p = p'$, $r_1 = 1 - q$. Two contrary testimonies of the same weight nullify each other, and the probability of the fact's falsity remains as it was from the beginning.

The probability of the fact being true or false is similarly and easily determined when some witnesses contradict a number of others. But when the fact is unanimously attested by all the witnesses, the probability of its truthfulness takes the following form.

As before, let $q$ be the probability of the fact's truthfulness that existed prior to all the testimonies, and $y_x$ and $y_{x-1}$ be that probability after $x$ and $(x-1)$ witnesses had attested the same. Denote also by $p^{(x-1)}$ the probability that another witness, not included in those former, does not deceive us by stating the truthfulness of the fact as well. Then $y_x$ is derived from $r$ of § 36 by substituting $p^{(x-1)}$ and $y_{x-1}$ instead of $p$ and $q$:

$$y_x = \frac{p^{(x-1)} y_{x-1}}{p^{(x-1)} y_{x-1} + (1 - p^{(x-1)})(1 - y_{x-1})},$$

and $y_0 = q$. Taking $x = 1, 2, 3, \ldots$, we deduce the following formulas

$$y_1 = \frac{pq}{pq + (1-p)(1-q)}, \quad y_2 = \frac{p'q}{p'y_1 + (1-p')(q - y_1)}, \ldots$$

Thus, $y_2$ is calculated by excluding $y_1$; $y_3$, by excluding $y_2$; … Denote for the sake of brevity

$$\frac{1 - p^{(x-1)}}{p^{(x-1)}} = \rho_x,$$

and the preceding difference equation of the first order (and its complete solution) will become

$$y_x = \frac{y_{x-1}}{y_{x-1} + \rho_x(1 - y_{x-1})} = \frac{c}{c + (1-c)\rho_1 \rho_2 \ldots \rho_x}.$$

Here, $c$ is an arbitrary constant. Substituting $(x - 1)$ instead of $x$, we have

$$y_{x-1} = \frac{c}{c + (1-c)\rho_1 \rho_2 \ldots \rho_{x-1}}, \quad 1 - y_{x-1} = \frac{(1-c)\rho_1 \rho_2 \ldots \rho_{x-1}}{c + (1-c)\rho_1 \rho_2 \ldots \rho_{x-1}}.$$

Together with $y_x$, these values render the given equation identical and $c$ is determined by a particular value of $y_x$ by taking for example $x = 0$. Accordingly, assuming that $\rho_1 \rho_2 \ldots \rho_x = 1$, we get $y_0 = q = c$.

For $x$ witnesses

$$y_x = \frac{q}{q + (1-q)\rho_1 \rho_2 \ldots \rho_x}$$

and for witness $i$ the magnitude $\rho_i$ is equal to the ratio of the probabilities that he deceives or not deceives us so that $\rho_i > 1$ or $< 1$ when the former is higher or lower than the latter. And $\rho_i = 1$ when they are equal. When a very large number of witnesses are considered infinite and $\rho_i > 1$ for all of them, the probability $y_x$ of the truth of the attested fact becomes zero with an exceptional case. On the contrary, again when *x* is infinite and $\rho_i < 1$ for all of them, that probability is unity or certitude, also with an exception.

That exception takes place when $\rho_1, \rho_2, \ldots$ continuously decrease or increase and indefinitely approach unity. For example, let

$$\rho_i = 1 - \frac{4g^2}{(2i-1)^2 \pi^2},$$

where π is the ratio of a circumference to its diameter and *g* is a given constant not exceeding unity so that no magnitude $\rho_i$ becomes negative. According to a known formula, their product is equal to cos*g*,

$$\prod_i (1 - \frac{4g^2}{(2i-1)^2 \pi^2}) = \cos g \text{ and } y_\infty = \frac{q}{q + (1-q)\cos g}$$

and will differ much from unity if *g* differs from π/2. Let $g = h\sqrt{-1}$, then the new constant *h* can be smaller or larger than unity. Denoting the base of the Naperian logarithms by *e*, we will have

$$y_\infty = \frac{2q}{2q + (1-q)(e^h + e^{-h})},$$

and if $h \leq 1$, or at least is not very large, that probability $y_x$ will not be very low. Still, it is easy to be assured in that the first value of $y_x$ will always be higher than the probability *q* existing before the testimonies, and its second value, always lower.

These formulas assume that all the testimonies are direct, but we will now examine the case in which only one of them is direct, and all the others are traditional [traditionally prompted].

**38.** Suppose that a witness does not at all restrict himself by saying whether a thing is true or false, but attests to the arrival of an event in the case in which many of them were possible. The event which he is able to announce, if mistaken or wishing to deceive, is not unique and only ought to be one of the non-existent or those in which he does not believe. We will show that, independently from the prior probability, that circumstance influences the probability of the event existing after the testimony.

To fix the ideas, I suppose that an urn A contains μ balls, $a_1$ of them having number 1, $a_2$, number 2, …, $a_m$, number *m*, so that

μ = $a_1 + a_2 + \ldots + a_m$.

If a ball is extracted, we can formulate *m* different hypotheses, $C_1, C_2, \ldots, C_m$ about its number. Let their probabilities before any testimony be

$q_1 = a_1/μ, q_2 = a_2/μ, \ldots, q_m = a_m/μ$.

If a witness announces that the number of the drawn ball was *n*, the probabilities of those hypotheses become $w_1, w_2, \ldots, w_m$, and it is required to determine them according to the rule of § 34. Here, the observed event is the announcement of the arrival of number *n*. Each hypothesis provides a certain probability $p_n$, whose expression we ought to form first of all, to that event, so that $C_n, q_n, w_n, p_n$ correspond to the announced number *n*. I denote by *u* the probability that the witness was not mistaken, and by *v*, the probability that he did not wish to deceive. Then $(1 - u)$ will be the probability that he was mistaken, and $(1 - v)$, the probability that he wished to deceive. According to the *n*-th hypothesis, the witness will announce the drawing of that number provided he was not mistaken and did not wish to deceive.

That combination of circumstances has probability $uv$ (§ 5). If he was mistaken, he believed that the arrived ball had number $n' \ne n$; at the same time, when wishing to deceive, he will announce a number differing from $n'$ and taken from the $(m - 1)$ other numbers. The chance that the witness chooses the number *n* is $1/(m - 1)$, assuming however that he does not prefer one number rather than another. By the cited rule it follows that the probability that that number will be announced by a witness who was mistaken and wished to deceive, will be the product of three factors, $(1 - u)$, $(1 - v)$, and $1/(m - 1)$. If the witness was mistaken and did not wish to deceive, or was not mistaken and wished to deceive, he will not announce the extraction of number *n*. Indeed, in the first case he wishes to announce the number which he believes to be drawn, and which is not *n*; in the second case, he knows that that number was extracted and does not wish to announce it. All that discussion taken together with the rule of § 10 leads to the composite probability

$$p_n = uv + \frac{(1-u)(1-v)}{m-1}$$

which hypothesis $C_n$, if certain, attaches to the observed event. According to hypothesis $C_i$, corresponding to the arrival of number $i \ne n$, the witness, if not mistaken and not wishing to deceive, does not announce number *n*. If he is not mistaken and wishes to deceive, he knows that *i* has arrived but announces a number among the $(m - 1)$ others. The chance that that number will be *n*, is $1/(m - 1)$, so that $u(1 - v)/(m - 1)$ is the probability that he actually announces number *n*. If he is mistaken and does not wish to deceive that probability will be $v(1 - u)/(m - 1)$. Finally, if the witness is mistaken and wishes to deceive he will first of all believe that the arrived ball was one of the $(m - 1)$ numbers differing from those which he announces. Fraction $1/(m - 1)$ will be the probability that he announces *n*.

Therefore, the probability that he believes that the drawn ball had number $n'$ but announces *n* is $(m - 1)^2$. The resulting chance for that number *n* to be announced is therefore $1/(m - 1)^2$ multiplied by the number of such numbers like $n'$ which the witness can believe to have arrived. That number is only $(m - 2)$ since the mistaken witness wishing to deceive can not believe that the drawn number was *i* (which really happened) or the number *n* which he announces. On the other hand, the probability of this double error is the product $(1 - u)(1 - v)$. Therefore, the probability that *n* will be actually announced by that witness is $(1 - u)(1 - v)(m - 2)/(m - 1)^2$.

I combine the probabilities of that announcement with the three distinct possible cases:

$$p_i = \frac{u(1-v)}{m-1} + \frac{v(1-u)}{m-1} + \frac{(m-2)(1-u)(1-v)}{(m-1)^2}.$$

This is the composite probability of the observed event according to one of the $(m - 1)$ hypotheses contrary to the truth of that event. It is connected with $p_n$ by the equation

$$p_n + (m - 1) p_i = 1, \tag{38.1}$$

since the sum of the probabilities corresponding to the announcement of number *n* according to the *m* hypotheses $C_1, C_2, \ldots, C_m$ should equal unity. Then, by the rule of § 34, we have

$$w_n = \frac{q_n p_n}{q_n p_n + \sum q_i p_i}, \quad w_i = \frac{q_i p_i}{q_i p_i + \sum q_i p_i},$$

where the sums extend from $i = 1$ to *m* except *n*. And since $p_i$ is independent from *i*, and the sum of $q_i$ less its value at $i = n$ is $(\mu - a_n)/\mu$, after substituting the values of $p_n, q_n, p_i, q_i$ and multiplying the numerator and denominator by $\mu(m - 1)^2$, the expression of $w_n$ becomes

$$w_n = \frac{\alpha}{\alpha + \beta(\mu - a_n)}, \quad \alpha = (m - 1)[uv + (1 - u)(1 - v)]a_n,$$
$$\beta = (m - 1)(1 - v)u + (m - 1)(1 - u)v + (m - 2)(1 - u)(1 - v).$$

So this is the probability that number *n* announced by the witness had really come from A; the probability that it did not, is $(1 - w_n)$. In particular, the probability of the drawing of any other fixed number *i* is derived from $(1 - w_n)$ when multiplying it by $q_i p_i / \sum q_i p_i$ or by $a_i/(\mu - a_n)$, so that

$$w_i = \frac{(1 - w_n)a_i}{\mu - a_n}.$$

It should be noted that, for deriving this result, we assumed that when the witness is mistaken or wishes to deceive, the number he announces is only determined by chance rather than some particular reason. It will not be the same either if he wishes to deceive because he has some reason to believe in the arrival of one number rather than another, or when he is mistaken because, for example, he believed that the announced number arrived, but that it was similar to the number really extracted. It is difficult to evaluate such circumstances, and we leave them aside although they can seriously influence the probability of the announced number.

Instead of balls carrying different numbers the urn could have contained balls of the same number of different colours. If there are only white and black balls in the ratio $a/(\mu - a)$, and the witness announces the arrival of a white ball, then, in the expression representing $w_n$, $m = 2$, $a_n = a$. Denoting the result by $r$, we have

$$r = \frac{\eta}{\eta + [(1-v)u + (1-u)v](\mu - a)}, \eta = [uv + (1-u)(1-v)]a.$$

We may liken this particular case with the case of a true or false fact attested by a witness by supposing that the drawing of the white ball is that fact. The probability of its truthfulness will be $r$, and its expression should coincide with that of § 36.

We have first of all the probability that the witness did not deceive

$$p = uv + (1-u)(1-v).$$

That can happen either if he is not mistaken and does not wish to deceive, or if he is mistaken and wishes to deceive. Only in these two possible cases the drawing of a white and a black ball represent the truth and falsity of the attested fact; the witness believes the contrary to what happened, and testifies contrary to what he believes. At the same time, the probability that he deceives us is

$$1 - p = (1-v)u + (1-u)v,$$

which can be derived from the value of $p$ or obtained directly when noting that the witness can deceive us either when he is not mistaken but wishes to deceive, or when mistaken and not wishing to deceive. And in addition $q = a/\mu$ and $1 - q = (\mu - a)/\mu$ are the probabilities of the truth and falsity of the testified fact existing before the testimony. These diverse values actually identify the expression of $r$ in § 36 with the just obtained formula.

If the urn contains only one ball of each number from 1 to $m$, then $a_n = 1$ and $\mu = m$, and the general expression of $w_n$ will become much simpler:

$$w_n = uv + \frac{(1-u)(1-v)}{m-1}.$$

The probability that the number $n$ announced by the witness is indeed extracted, does not then differ from what was above denoted by $p_n$; that is, from the probability that the witness announcing that number supposes the same. It lowers as the number $m$ of balls in the urn increases and, if $m$ can become infinite, becomes equal to the probability that the witness was not mistaken and did not wish to deceive.

**39.** It remains to consider the general case of many witnesses, some of whom have direct knowledge of the fact they attest, and the others only know about it by tradition. However, for narrowing the extent of this digression from the probability of testimonies, we restrict our attention to resolving one particular problem of that kind.

Denote the $(x + 1)$ witnesses by $T, T_1, T_2, \ldots, T_{x-1}, T_x$. Just like in the preceding problem, a ball is extracted from urn A, and T directly knows its number whereas each of the other witnesses repeats after his predecessor that its number is $n$. The information is thus transmitted from T to $T_x$ and from him to us by a traditional and uninterrupted chain. So $T_x$ is our only witness, and he testifies, as borrowed from $T_{x-1}$, that the arrived ball had number $n$. It is required to determine the probability that that number was indeed extracted.

Let $y_x$ and $y'_x$ be the probabilities of the observed events, the drawing of number $n$ according to hypothesis $C_n$, and of number $i \neq n$ according to hypothesis $C_i$. Like previously, denote also the number of balls having numbers $n$ and $i$ by $a_n$ and $a_i$ and let $\mu$ be the total number of balls in the urn. Then the fractions $a_n/\mu$ and $a_i/\mu$ will be the prior chances of the arrival of numbers $n$ and $i$. By the rule of § 34 we have the probability of the hypothesis $C_n$

$$w_n = \frac{a_n y_x}{a_n y_x + \sum a_i y'_x}.$$

The sum extends over all values of $i$ from 1 to $m$ except $n$. We see at once that $y'_x$ is independent from $i$, and that the sum of the values of $a_i$ except $a_n$ equals $(\mu - a_n)$. And so,

$$w_n = \frac{a_n y_x}{a_n y_x + (\mu - a_n) y'_x}.$$

Probability $w_i$ of any other hypothesis $C_i$ is obtained when multiplying $(1 - w_n)$ by $a_i/(\mu - a_n)$. Our problem is now reduced to determining the unknowns $y_x$ and $y'_x$ as functions of $x$. I therefore represent by $k_x$ the probability that witness $T_x$ does not deceive; then $(1 - k_x)$ will be the probability of deception, involuntary or intended. That witness announces the extraction of number $n$ if he does not deceive and $T_{x-1}$ had stated the same. According to hypothesis $C_n$ that combination has probability $k_x y_{x-1}$. With regard to $T_{x-1}, y_{x-1}$ expresses the same as $y_x$ with regard to $T_x$. Then, he can announce the drawing of $n$ if he deceives and $T_{x-1}$ testified to another number. By the hypothesis $C_n$ the probability of that combination is $(1 - k_x)(1 - y_{x-1})$. But the chance that $T_x$ announces number $n$ taken from $(m - 1)$ numbers, which he does not believe to be announced by $T_{x-1}$, is $1/(m - 1)$, so the probability that $n$ will be announced should be reduced to $(1 - k_x)(1 - y_{x-1})/(m - 1)$. Finally, $T_x$ does not announce the arrival of that number, either because he deceives us and $T_{x-1}$ did announce it, or if he does not deceive us and $T_{x-1}$ announced the arrival of another number. So we have the composite probabilities of the observed event according to hypothesis $C_n$ and all the other hypotheses $C_i$

$$y_x = k_x y_{x-1} + \frac{(1-k_x)(1-y_{x-1})}{m-1}, \quad y'_x = k_x y'_{x-1} + \frac{(1-k_x)(1-y'_{x-1})}{m-1}.$$

The two unknowns, $y_x$ and $y'_x$, are included in the same difference equation of the first order and only differ one from another in the arbitrary constant.

When considering $y_x$ and denoting that constant by $c$, the complete solution of that equation is

$$y_x = \frac{1}{m} + \frac{c(mk_1 - 1)(mk_2 - 1)...(mk_{x-1} - 1)(mk_x - 1)}{(m-1)^x},$$

and, substituting $(x - 1)$ instead of $x$, we will derive

$$y_{x-1} = \frac{1}{m} + \frac{c(mk_1 - 1)(mk_2 - 1)...(mk_{x-1} - 1)}{(m-1)^{x-1}},$$

$$1 - y_{x-1} = \frac{m-1}{m} - \frac{c(mk_1 - 1)(mk_2 - 1)...(mk_{x-1} - 1)}{(m-1)^{x-1}}.$$

These values together with $y_x$ transform the given equation into an identity. For determining $c$, I take $x = 0$ in the complete solution and note that the probability $y_x$ pertaining to the direct witness T, should coincide with $p_n$ in § 38. And, assuming unity at $x = 0$ for the product of the factors included in the solution, we will have

$p_n = (1/m) + c$, $c = (mp_n - 1)/m$.

For some value of $x$ we will therefore get

$$y_x = \frac{1}{m}[1 + (mp_n - 1)X], \quad X = \frac{(mk_1 - 1)(mk_2 - 1)...(mk_x - 1)}{(m-1)^2}.$$

Note that according to some hypothesis $C_i$ differing from $C_n$ the probability $y'_x$ concerning the direct witness T should also be the probability denoted by $p_i$ in § 38 and

$$y'_x = \frac{1}{m}[1 + (mp_i - 1)X],$$

which, like $p_i$, is independent from $i$. I substitute these values in $w_n$ and obtain the probability that the number $n$, announced by the last witness $T_x$, is actually extracted

$$w_n = \frac{[1 + (mp_n - 1))X]a_n}{[1 + (mp_n - 1))X]a_n + [1 + (mp_i - 1))X](\mu - a_n)}$$

which is what was required to determine.

The product $X$ can be replaced by

$$X = h_1 h_2 \ldots h_x, \quad h_x = k_x - \frac{1-k_x}{m-1}.$$

Number $m$ is always larger than unity, and $k_x$ represents a positive fraction which can not be larger than unity, so each of the factors of $X$ can be positive or negative without ever exceeding the limits ± 1. When the number $x$ of these factors is very large, the product is negligible and even zero if that number becomes infinite. However, if those factors $h_1$, $h_2$, … form a series, continuously converging to unity, the indicated statement is excluded.

If we neglect the terms including $X$ in the formula for $w_n$, it will be reduced to $a_n/\mu$. Therefore, in general, the probability of an event transmitted to us by a traditional chain of a very large number of witnesses, does not sensibly differ from the chance proper to that event or independent from the testimonies. However, if a long number of direct witnesses attests an event, its probability will closely approach unity provided (§ 37) that it is possible to bet more than one against one on the honesty of each witness.

In the particular case in which the urn contains only one ball of each number, $a_n = 1$ and $\mu = m$. Then because of equation (38.1) the value of $w_n$ becomes

$w_n = [1 + (mp_n - 1)X]/m.$

This probability coincides with $y_x$, that is, with the probability of the announcement of number $n$ by witness $T_x$ in accord with hypothesis $C_n$ which stipulates that that number was really extracted from the urn. However, we can not admit that conclusion in advance, as Laplace (1812/1886, § 44) did when solving this problem. These probabilities, $y_x$ and $w_n$, are only identical when $(\mu - a_n)/a_n = (m - 1)$.

**40.** If desired, it is possible to express each of the magnitudes $k_1$, $k_2$, … through $m$ and the probabilities that the witness, to whom they correspond, is not mistaken and does not wish to deceive. I denote by $u_{x'}$ the probability that witness $T_{x'}$ belonging to the traditional chain is not mistaken, and by $v_{x'}$, the probability that he does not wish to deceive. If these two circumstances coexist, the witness is not mistaken. He can also be not mistaken if he is mistaken and wishes to deceive. In that second case the chance that he will announce number $n$ is $1/(m - 1)$ since $n$ is taken from $(m - 1)$ numbers in which he does not believe to have arrived. Those two cases are the only ones in which he is not mistaken when announcing that number, and the composite value of $k_{x'}$ is

$$k_{x'} = u_{x'} v_{x'} + \frac{(1 - u_{x'})(1 - v_{x'})}{m - 1}.$$

If $x' = 0$ it coincides with $p_n$ of § 38 if $u_i$ and $v_i$ (? - O.S.) are substituted in its expression. That magnitude $k_{x'}$ is the probability to be attached to the testimony of $T_{x'}$ or the value of the testimony itself, i. e., to the reason to believe that number $n$ is extracted from an urn containing $m$ kinds of different numbers when we only know that that arrival was attested by $T_{x'}$ for whom $u_{x'}$ and $v_{x'}$ are the probabilities of not being mistaken and not wishing to deceive. If $T_{x'}$ is certainly mistaken and wishes to deceive, then $u_{x'} = v_{x'} = 0$ and the probability $k_{x'}$ that number $n$ was drawn resulting from his testimony [from $T_{x'}$] will be nevertheless equal to $1/(m-1)$.

This is certain for $m = 2$. And actually in this case the witness, when announcing a number in whose arrival he does not believe and believing in the extraction of the number which was not drawn, necessarily announces the truth. If $m = 3$, even money can be bet on the arrival of the announced number. This can easily be verified by enumerating all the possible combinations. Also verified can be the value $1/(m-1)$ of probability $k_{x'}$ relative to some number $m$.

The case of a witness who is mistaken and certainly wishes to deceive should not be confused with that of an interrupted traditional chain in which witness $T_{x'-1}$ preceding $T_{x'}$ does not exist so that $T_{x'}$ certainly wishes to be mistaken since he supposes that $T_{x'-1}$ exists and therefore $v_{x'} = 0$. But the probability that $T_{x'}$ is not mistaken is not at all zero; he has no notion about the arrived event and the probability that he announces the really drawn number is $1/m$, and this is therefore the value of his testimony. And with $u_{x'} = 1/m$ and $v_{x'} = 0$ the preceding formula leads to $k_{x'} = 1/m$. This brings about the condition $h_{x'} = 0$ and reduces the probability $w_n$ of the arrival of number $n$ to $a_n/\mu$, to the chance proper to that event, as it should be evident.

**41.** By the rule concerning the probabilities of causes, we may actually complete what was stated at the end of § 7 about the tendency of our mind to believe that certain events undoubtedly have a special cause independent from chance. Suppose we observe an event which, taken by itself, has a very low probability. If it presents some symmetry or some other remarkable thing, we are naturally led to think that it is not the effect of chance or, more generally, of a unique cause attaching to it that feeble chance, but is due to a more powerful cause, such as the desire of someone who had a particular goal for producing it.

If for example[9] we find 26 printed letters of the alphabet arranged on a table in their natural order a, b, c, …, x, y, z, we will not doubt that someone had wished to dispose them in that way. But still that arrangement is not by itself less probable than any other which does not present anything remarkable and which we therefore do not hesitate to attribute to chance. If these 26 letters were successively and randomly drawn from an urn, there would have been the same chance of their arrival in the natural order or an order determined in advance, as for example b, p, w, …, q, a, t, which I have chosen arbitrarily. That chance will be slim, but not slimmer than in the other case.

Similarly, if an urn contains an equal number of white and black balls and 30 of them had to be successively extracted with replacement, the probability that all of them are white will be $(1/2)^{30} \approx 1/10^9$. However, the probability of some of them white, and some black in any order and any ratios will not be either higher or lower than in the previous case. Just the same, we could have again bet about $10^9$:1 against the appearance of that fixed arrangement. And nevertheless, when seeing the arrival of 30 white balls in succession, we can not believe that that event took place by chance whereas we easily attribute to chance the arrival of 30 balls presenting nothing regular or remarkable.

What we call *hazard* (§ 27) with the same facility produces, so to say, an event which we find remarkable and another one, unremarkable. If the equally possible events are very numerous, events of the first kind are much rarer than those of the second. It is for this reason that an arrival of events of the first kind is astonishing to us and leads us to search for their special cause. Actually, its existence is likely, but its high probability does not result from the rarity of remarkable events; it is founded on another principle which we will apply to the rules demonstrated above.

**42.** Let us denote the possible remarkable and unremarkable events by $E_1, E_2, \ldots$ and $F_1, F_2, \ldots$ When considering the 30 balls extracted from an urn containing an equal number of white and black balls, the events $E_1, E_2, \ldots$ will be the appearance of all 30 of the same colour; or alternatively white and black; or 15 of one colour followed by 15 of the other; etc. In the case of [about] 30 printed letters disposed one after another, those events will be an arrangement of those letters in the alphabetic order; or in the inverse order; or making up a phrase in the French or another language.

Denote the number of all the remarkable events by $m$, and, by $n$, the other events $F_1, F_2, \ldots$ Suppose that all are equally possible when only due to chance, then the probability of each of either kind will be $p = 1/(m + n)$. It will not be the same if those events were produced by a particular cause C independent from the probability $p$ and being, for the sake of definiteness, someone's desire and choice.

We assume that this choice is determined by diverse circumstances rendering some of the possible events remarkable. Thus, there exists a certain probability $p_1$ that the choice of that person attaches to $E_1$, attaches probability $p_2$ to $E_2$, ... If these various events are the only ones possible, it follows that $p_1 + p_2 + \ldots = 1$. And if all those probabilities are the same, their common value is $1/m$, very high compared with $p$ if the total number $(m + n)$ of the possible cases is very large by itself and when compared with $m$.

In general, these probabilities can be very unequal but we are unable to know them. For us, however, it suffices that they are very high compared with probability $p$. This should have certainly happened when that probability is extremely low, or the number $(m + n)$ excessively large, like in the examples below.

So this is the principle from which we issue when determining the probability of cause C after observing one of the events $E_1, E_2, \ldots, F_1, F_2, \ldots$ or at least when showing that it is very high if the observed event belongs to the first kind.

Suppose that that event is $E_1$. It is possible to formulate two hypotheses: it is either due to cause C, or resulted from chance. If the first hypothesis is certain, $p_1$ will be the probability of the arrival of $E_1$; if the second one is certain, that probability will be $p$. Denote by $r$ the probability of the first hypothesis existing after the observation, and consider both of them equally probable before it. Then, by the rule of § 28,

$r = p_1/(p_1+ p)$.

Suffice it that probability $p_1$ is very high compared with a very slim chance $p$ for that value of $r$ to differ very little from unity or certainty. In a preceding example, the number of possible events exceeded $10^9$, and $p$ was lower than $1/10^9$. Suppose that thousand is the number of sufficiently remarkable events for choosing one of them, and assume $1/1000$ as the value of $p_1$. Then $r$ will differ from unity less than by $1/10^6$ and much less if, as it is possible to believe, probability $p_1$ is higher than $1/1000$. If then one of those remarkable events will be observed, for example the extraction of 30 balls of the same colour from an urn containing equal numbers of balls of two colours, we should without any doubt attribute that fact to someone's desire or to some other special cause, as it is certainly done, rather than to a simple effect of chance.

Nevertheless, the probability $r$ of the cause C will be considerably lower if before the observation the latter's existence and absence were not equally possible as was supposed in the preceding formula and if its absence was initially more probable. This is what happened in the example above when many precautions were taken before the drawings for excluding the influence of any desire from the extraction of the balls. By allowing for that circumstance taking place prior to the observation, the lowering of $r$ according to the rule of § 34 becomes appreciable.

That probability will be heightened or lowered, sometimes greatly, when all the events $E_1$, $E_2$, …, $F_1$, $F_2$, …, are not at all equally possible; heightened, if the chance proper to each event is lower for those of the first kind than for the second, and lowered otherwise.

The harmony which we observe in nature is undoubtedly not occasioned by chance. By attentive examination over a very long period of time of a very large number of phenomena we came to discover their physical causes if not with absolute certainty, at least with a probability closely approaching it. Regarding phenomena presenting remarkable circumstances as things $E_1$, $E_2$, …, we will have the case in which those things possess by themselves a probability high enough for rendering interventions of the cause, which we denoted by C, very unlikely, sometimes for it being useless to consider such interventions.

It is reasonable to attribute physical phenomena whose causes are still unknown to causes similar to the known ones and to believe that they obey the same laws. With the progress of science their number decreases, however, from day to day[10]. Nowadays, we know, for example, what produces lightnings and how are the planets kept on their orbits, which was unknown to our predecessors. Those, who come after us, will discover the yet unknown causes of other phenomena.

**43.** When the number of various causes which can be attributed to an observed event E is infinite, their probabilities either before or after the arrival of E become infinitely low, and the sums included in the formulas of §§ 32 and 34 are transformed into definite integrals.

For bringing about that transformation, let us suppose that the observed event E is the extraction of a white ball from an urn containing infinitely many white or black balls. We can formulate an infinity of hypotheses about the unknown ratio $x$ of the number of white balls to the total number of balls. They can be thought to be so many mutually incompatible causes of the arrival of E.

Now, $x$ is susceptible of all values increasing by infinitely small steps, from infinitely small corresponding to the case in which the drawn ball was the only white ball in A, to $x = 1$ corresponding to the other extreme case in which that urn contained only white balls. Represent by $X$, which always is a known function of $x$, whose value is supposed certain, the probability that that ratio attaches to the arrival of E. Considering that value as a possible cause of E, determine the infinitely low probability of $x$, whether all the causes are equally possible before the observation, or different from the beginning.

In the first case, the required probability is derived from $w_n$ of § 28 by assuming an infinite $m$ and substituting the values of $x$ corresponding to $X$ instead of $p_1, p_2, \ldots$ Understanding sums in the same way as in § 32, and denoting the probability of $x$ by $w$, we will have $w = X/\sum X$. However, by the fundamental theorem of the [theory of] definite integrals

$$\sum X dx = \int_0^1 X dx.$$

Therefore, supposing that the differential $dx$ is constant, and multiplying both sides of the preceding fraction $w$ by it, we arrive at

$$w = X dx \div \int_0^1 X dx.$$

Denote by $X'$ the probability, corresponding to $x$, of a future event E' depending on the same causes as E, and by $w'$, the composite probability of the appearance of E'. Then, by the rule of § 30,

$$w' = \sum X' w = \int_0^1 XX' dx \div \int_0^1 X dx.$$

Here, $w$ was replaced by its preceding value and the sum transformed into an integral. In each example [problem] $X'$ will be a given function of $x$. If, before observing E, the diverse values of $x$ are unequally probable, denote by $Ydx$ the infinitely small probability of the chance $x$ of E and replace $q_n$ in the formulas of § 34 by $Ydx$. Then

$$w = XYdx \div \int_0^1 XYdx, \quad w' = \int_0^1 XX'Ydx \div \int_0^1 XYdx$$

will be the probabilities of the chances of the arrival of E and of the probability of the future arrival of E.

**44.** If we are certain from the very beginning that *x* can not extend from 0 to 1, but should be confined between given limits, these limits must be applied as the limits of the definite integrals included in these (? - O.S.) formulas. Alternatively, if desired to preserve their limits 0 and 1, *Y* will be a discontinuous function of *x* disappearing beyond the given limits of that variable.

Let *x* be susceptible of all values from 0 to 1, or contained in given limits of that variable. Denote by λ the probability that after the arrival of the observed event E its unknown value actually becomes contained in more narrow limits α and β. Then λ will be the sum of the values of *w* corresponding to the values of *x* thus contained:

$$\lambda = \int_\alpha^\beta XYdx \div \int_0^1 XYdx.$$

That formula can be applied in approximate calculations when the number of causes, to which the event E can be attributed, is only very considerable instead of infinite. Suppose for example that E is the extraction of *n* white balls drawn successively with replacement, without interruption, from an urn B containing a very large number of white and black balls. The probability *X* of E corresponding to ratio *x* of the number of white balls to the total number of all balls in B will be the *n*-th power of that ratio.

When required the probability that the number of white balls exceeds the number of black balls, we take α = 1/2 and β = 1 in the expression of probability λ. And if in addition, before the drawings all possible values of *x* were equally probable, *Y* will not vary with *x* and therefore disappears from that expression:

$$X = x^n, \quad \int_0^1 Xdx = \frac{1}{n+1}, \quad \int_{1/2}^1 Xdx = \frac{1}{n+1}(1 - \frac{1}{2^{n+1}})$$

and λ = (1 − 1/$2^{n+1}$) with its precision increasing with the number of black or white balls in B. Before the drawings even money could have been bet on the number of white balls exceeding the number of black balls. However, after extracting just one white ball from B, λ = 3/4 and we can bet 3 against 1 on the number of white balls being larger than the number of black balls. After some small number of white balls has arrived in succession, the probability λ that there are more white balls than black considerably approaches certainty.

**45.** As I said above (§ 30), E and E′ can be considered as events composed of the same simple event G and connected with each other by their common dependence on that event[11]. The chance of G is unknown. Before E arrived, the probability of its taking value *x* is *Ydx* and *w* after that. Since this value is certainly contained between *x* = 0 and 1, the sum of the corresponding values of *Ydx* will be unity, just as the sum of the values of *w*. So, be the given function *Y* of *x* continuous or discontinuous, it should always satisfy the condition ∫*Ydx* = 1. By the rule [the definition] of mathematical expectation (§ 23)[12] applied to the chance of G, we should assume as its value before the arrival of E the sum of all its possible values multiplied by their respective probabilities, i. e., the sum of all the products *xYdx* from *x* = 0 to *x* = 1. Denoting this chance, or, more precisely, denoting what should be assumed as its unknown chance before the arrival of E by γ, we will have

$$\gamma = \int_0^1 xYdx.$$

When considering *x* and *X* as the abscissa and ordinate of a plane curve and noting that the entire area of that curve[13] or the integral of *Ydx* is unity, γ will be the abscissa of the centre of gravity of that same area. It is this value of γ adopted as the chance of G that we should choose when betting on the first arrival of that event, but not on many successive arrivals. Indeed, depending on whether G takes place at the first trial or not, the probability of further arrivals will be heightened or lowered.

Let for example all the values of *x* be in advance equally probable, then *Y* should be independent from *x*. According to the two preceding equations, *Y* = 1 and γ = 1/2 and we have no reason to believe that G rather than the contrary event will arrive at the first trial. However, if considering E and E′ as the simple event G, we will have

$$X = x,\ X' = x \text{ and } w' = \int_0^1 XX'dx \div \int_0^1 Xdx = \frac{2}{3}$$

for the probability that, having arrived for the first time, G will arrive for the second time since the probability of that recurrence will after the first trial be heightened by 1/6 [ from 1/2 to 2/3]. That probability will be lowered by the same fraction and become 1/2 − 1/6 = 1/3 if the contrary event took place at the first trial. Indeed, assuming that event as E, and, as always, considering G as E′, that is, taking *X* = 1 − *x* and *X'* = *x*, we conclude that

$$w' = \int_0^1 x(1-x)dx \div \int_0^1 (1-x)dx = \frac{1}{3}$$

will be the probability that G, not arriving at the first trial, will appear at the second.

In advance, the probability of G occurring twice in succession is (§ 9) the product of the probability 1/2 that it takes place the first time and probability 2/3 that it arrives at the second trial, and is equal to 1/3 instead of 1/4 which is its value had the probability of G been 1/2 at the second trial. The coincidence of two events in the first two trials will have a double probability, i. e., 2/3, since it can occur by repetition of either G or of its equally probable contrary event. Comparing 2/3 = [1 + (1/3)]/2 with the probability $(1 + \delta^2)/2$ of the coincidence as determined in § 26, we have $\delta = 1/\sqrt{3}$. Therefore, when having no advance knowledge about the chance of an event G, we may have equal grounds for supposing that $x$ takes any of all possible values[14].

We will now determine the probability of coincidence in the case in which we know in advance that the values of $x$ instead of being equally possible likely differ very little from a known or unknown fraction.

**46.** We always denote by G a simple event with an unknown chance; let H be the event contrary to G with chance equal to 1 less the chance of G. We suppose that **1)** An observed event E is the arrival of G and H $m$ and $n$ times in some order. **2)** A future event E′ is the arrival of these events $m'$ and $n'$ times, again in some order.

The chances of G and H are $x$ and $1 - x$, so the probabilities $X$ and $X'$ of E and E′ are (§ 14)

$$X = Kx^m(1-x)^n, \quad X' = K'x^{m'}(1-x)^{n'},$$

where $K$ and $K'$ do not depend on $x$. Therefore

$$w' = K'\int_0^1 Yx^{m+m'}(1-x)^{n+n'}dx \div \int_0^1 Yx^m(1-x)^n dx$$

is the probability of E′ existing after E was observed. Magnitude $K$ has disappeared from this formula and

$$K' = C_{m'+n'}^{m'}. \tag{46.1}$$

If E′ is the occurrence of G and H $m'$ and $n'$ times in a fixed order, then $K' = 1$.

If before observing E we had no reason to believe that some value of $x$ is more probable than another, then $Y = 1$. And, integrating by parts, we get

$$\int_0^1 x^m(1-x)^n dx = \frac{n!}{(m+1)(m+2)...(m+n)(m+n+1)} = m!n!/(m+n+1)!.$$

Similarly,

$$\int_0^1 x^{m+m'}(1-x)^{n+n'}dx = \frac{(m+m')!(n+n')!}{(m+m'+n+n'+1)!}$$

and, see (46.1),

$$w' = \frac{(m'+n')!(m+m')!(n+n')!(m+n+1)!}{m'!n'!m!n!(m+m'+n+n'+1)!}.$$

For that formula to include the case in which one of the magnitudes *m, n, m', n'* is zero, we should assume that 0! = 1. Then, if *n* = *n'* = 0,

$$w' = \frac{m+1}{m+m'+1}.$$

This formula expresses the probability that G arrives *m'* times without interruption after it had taken place *m* times in succession without H taking place. For *m'* = 1 and *n'* = 0, the value of *w'* corresponding to *Y* = 1 is reduced to

$$w' = \frac{m+1}{m+n+2},$$

and, for *m'* = 0 and *n'* = 1,

$$w' = \frac{n+1}{m+n+2}.$$

The sum of these two fractions is unity which should have taken place since the former expresses the probability that after (*m* + *n*) trials G arrives at the next trial, whereas the latter expresses the contrary. The former is higher or lower than the latter depending on whether *m* > *n* or < *n*, i. e., on whether in the first (*m* + *n*) trials G had occurred more often or less than the contrary event H. These probabilities become equal to each other and equal to 1/2, as it was before the trials, when those two events had taken place the same number of times. However, in general this will not persist when we know in advance either by the nature of G, or by the results of previous trials that the unknown chance of G takes unequally probable values, so that *Y* ≠ 1.

The fraction γ (§ 45), which should have been assumed as the chance of G existing before the (*m'* + *n'*) new trials, will not at all be 1/2; at the next trial the probability of G, although arriving oftener than the contrary event H, can be lower than γ, or higher even in the opposite case. And this I will show in the next example.

**47.** I suppose that beforehand the chance of G likely very little deviates in either direction from a certain fraction *r*,

*x* = *r* + *z*.

Magnitude *Y* is a function of *z* and is only appreciable at very small positive or negative values of that variable. The plane curve whose current coordinates are *x* and *Y* only noticeably deviates from the *x*-axis on a very short interval extended on both sides from the ordinate corresponding to *x* = *z*. The centre of gravity of the curve's area is therefore situated on this interval, and its abscissa very little deviates from *r*. When neglecting that difference, *r* will become the value of γ from § 45.

The limits of the integrals will be $z = -r$ and $1 - r$ corresponding to $x = 0$ and 1. So, adopting $m' = 1$, $n' = 0$ and $dx = dz$ in the first expression of *w'* in § 46, we will have

$$w' = \int_{-r}^{1-r} Y x^{m+1}(1-x)^n dz \div \int_{-r}^{1-r} Y x^m (1-x)^n dz$$

as the probability that G arrives after having occurred *m* times, and H, *n* times in (*m* + *n*) trials. However, due to the nature [to the behaviour] of *Y*, it is possible, if desired, to restrict the extent of the integrals to very small values of *z*. And, when expanding the other factors in power series of *z*, in general they will converge rapidly. The only exception occurs when *r* or $1 - r$ are also very small fractions. In all other cases only the first terms of the series can be left, so that when neglecting the square of *z* (? - O.S.),

$x^m(1-x)^n = r^m(1-r)^n + [mr^{m-1}(1-r)^n - nr^m(1-r)^{n-1}]z + \frac{1}{2}[m(m-1)r^{m-2}(1-r)^n - 2mnr^{m-1}(1-r)^{n-1} + n(n-1)r^m(1-r)^{n-2}]z^2$.

Substituting (*m* + 1) instead of *m* I will derive $x^{m+1}(1-x)^n$.

Then, I insert $x^m(1-x)^n$ and $x^{m+1}(1-x)^n$ in the expression of *w'*; I also note that

$$\int_{-r}^{1-r} Y dz = 1, \quad \int_{-r}^{1-r} Y z \, dz = 0;$$

denote for the sake of brevity

$$\int_{-r}^{1-r} Y z^2 dz = h,$$

and, when neglecting $h^2$, which can only be a very small fraction, come to

$$w' = r + \left(\frac{m}{r} - \frac{n}{1-r}\right)h.$$

This proves that the probability *w'* of the arrival of G after (*m* + *n*) trials is higher if [$m/r < n/(1-r)$], or lower otherwise, than *r* or γ, which should be assumed as the chance of G existing before those trials. If *r* certainly was the chance of G and *m* and *n* were large numbers, G and H will likely appear in the ratio of their respective chances, *r* and $(1 - r)$. Equality $m/r = n/(1-r)$ renders the probability *w'* equal to the chance *r*, as it should have been.

**48.** When $m = 1$ and $n = 0$, the preceding value of *w'* becomes

$w' = r + h/r,$

the probability that G, having arrived at the first trial, also occurs at the second. The probability of the former is $r$, so $rw' = r^2 + h$ expresses the probability of the repetition at those trials. Substituting $(1 - r)$ instead of $r$, we get $(1 - r)^2 + h$ for the probability of the repetition of the contrary event, and the sum of those probabilities,

$1 - 2r + 2r^2 + 2h,$

is the probability of the coincidence of the results of those trials.

If $m = 0$ and $n = 1$, $w'$ from § 47 becomes

$w' = r - h/(1 - r),$

which is the probability that G, not arriving at the first trial, occurs at the second. When multiplied by $(1 - r)$, it will express the probability of the succession of two contrary events. The double of that product,

$2r - 2r^2 - 2h,$

is the probability of different outcomes in those two trials. By subtracting from unity the probability of coinciding outcomes it will also be derived.

The difference of the probabilities of coinciding and contrary results is therefore

$(1 - 2r)^2 + 4h,$

which shows that this excess increases when $r$ is not exactly the chance of G; we only know that this chance very little deviates from $r$. Even if knowing that $r = 1/2$, it will have still been advantageous to bet even money on coincidence. And this is what happens at playing *heads* and *tails* when first tossing the coin: equality of the chances of its two sides is physically impossible, but, in accord to its minting, the chance of each side likely very little deviates from 1/2.

**49.** Already here I announce a theorem whose demonstration follows in the next chapter. It determines the chances of an event by experience, not with certainty and rigour, but as being likely its very good approximant.

Let $g$ be the known or unknown chance of an event G, i. e., the ratio of the number of equally possible cases favourable to it, to the number of all cases that can take place and are also equally possible. Suppose that $\mu$ trials are made during which that chance proper to G and differing from its probability (§ 1) remains constant. Let $r$ be the ratio of the number of arrivals of G in that series of trials to its total number $\mu$. If that $\mu$ is not very considerable, $r$ varies with it and can much differ from $g$ in either direction. However, when $\mu$ becomes large, the difference $(r - g)$, abstracting its sign, will ever decrease with the further increasing $\mu$; if $\mu$ can become infinite, then exactly $r - g = 0$.

Denoting an arbitrarily small fraction by ε, we can always assign such a large μ that the probability of $r - g < ε$ will approach certainty as closely as desired. Below, we provide an expression of the probability of that inequality as a function of μ and ε.

And so, suppose that urn A contains *a* white and *b* black balls. When successively drawing with replacement a very large number μ of balls, α of them white and β, black, we will have

$$\frac{α}{α+β} = \frac{a}{a+b}, \quad \frac{β}{α+β} = \frac{b}{a+b}, \quad \frac{α}{β} = \frac{a}{b}$$

ever more exactly and with an ever higher probability as $μ = α + β$ becomes larger. Conversely, if the ratio of the numbers of white and black balls in A is unknown, and during a very large number of trials that ratio does not vary, we may with a very high probability take for the approximate values of that unknown ratio and the unknown chance of extracting a white ball the magnitudes α/β and α/(α + β) whatever is the number *a* of balls in that urn. However, we should remark that if the number of white balls is very small as compared with the number *b* of black balls, α will also be very small as compared with β and vice versa. And the ratio of one of the fractions α/β and *a/b* to the other can much differ from unity, at least until the series of trials does not excessively extend. When the known or unknown chance of drawing a white ball is very slim, the approximate equality α/β = *a/b* only signifies that both fractions are very small.

The announced rule equally applies to chances of diverse and mutually incompatible causes which can be attributed to an event E observed a very large number of times. If γ is the known or unknown chance of one of these causes, C, the ratio $γ/(1 − γ)$ will with a high probability closely approximate the ratio of the number of arrivals of E due to C to that number occasioned by all other causes. This is how that ratio can be determined if γ was known in advance or found out by experience.

Let the event E be the extraction of a white ball from either an urn A with *a* white and *a′* black balls, or an urn B with *b* white and *b′* black balls respectively. The values of the chance of A being the cause of E and of the contrary chance of B being that cause, are, by the rule of § 28,

$$γ = \frac{a(b+b′)}{a(b+b′)+b(a+a′)}, \quad 1 - γ = \frac{b(a+a′)}{a(b+b′)+b(a+a′)}.$$

When drawing with replacement a very large number μ of white balls from either urn, the ratio ρ of the number of them extracted from A to that of the extracted from B will likely very little deviate from $γ/(1 − γ)$:

$$ρ = \frac{γ}{1-γ} = \frac{a(b+b′)}{(b\,a+a′)}.$$

If $a + a' = b + b'$, then $\rho = a/b$. In this case all the balls can be combined in a single urn D (§ 10) without changing the ratio of white balls extracted from A and B which will be drawn from D.

For a very large number of balls of that colour the ratio of the former number to the latter will be almost equal to *a/b*. This can be verified by differently marking the balls taken from A and B and replacing each ball after its arrival.

**50.** The work of Buffon (1777, § 18) contains numerical results of an experiment on the game of *heads* and *tails* which will provide an example and a verification of the preceding rule. In that game, the chance of each side of the coin depends on its physical constitution not sufficiently known to us. And even had we known it, it would be a problem of mechanics no one can solve for finding out the chance[s] of *heads* or *tails*. So an approximate value of that chance [of those chances] should be derived by experiment for each coin in particular.

Thus, if *heads* appeared *m* times in a very large number $\mu$ of trials, $m/\mu$ should be assumed as its chance. It will also be the probability or the reason to believe that that side will arrive in a new trial made with the same coin. After that series a fair bet on the arrival of *heads* can be made by betting *m* against $(\mu - m)$. It is also by means of that probability $m/\mu$ of the simple event that the probabilities of compound events should be calculated, at least if they are not very low.

Suppose now that a very large number *m* of series of trials is made, each of them, like in the cited experiment, ending when *heads* appear. Let $a_1, a_2, a_3, \ldots$ be the number of arrivals of *heads* at the first; the second, … toss. The total number of tosses, or trials, and arrivals of *heads* will be

$$\mu = a_1 + 2a_2 + 3a_3 + \ldots, \quad m = a_1 + a_2 + \ldots$$

The chance of that side will be $p = m/\mu$ with the approximation being better and the precision higher, as $\mu$ becomes larger.

The probabilities of *heads* at the first toss, at only the second toss, at only the third one, … will be $p, p(1-p), p(1-p)^2, \ldots$ Since the numbers of those arrivals in *m* series of trials are supposed to be $a_1, a_2, \ldots$, the equalities

$$p = a_1/m, \; p(1-p) = a_2/m, \; p(1-p)^2 = a_3/m, \ldots$$

ought to take place almost exactly if that number is very large and the probabilities do not become too low. Dividing each of these equations by the preceding, we obtain different values of $(1-p)$ and therefore

$$p = a_1/m = 1 - a_2/a_1 = 1 - a_3/a_2 \ldots$$

These values, or at least a certain number of the first of them, will differ one from another and from $m/\mu$ the less, the larger are *m* and $\mu$. For those values to become certainly equal one to another, *m* and $\mu$ should be infinite. When assuming the mean of those very little differing fractions as *p* or when taking $p = m/\mu$, resulting from all the trials, then

$$a_1 = mp, a_2 = mp(1-p), a_3 = mp(1-p)^2, \ldots$$

These calculated values, or at least the first terms of that decreasing geometric progression, should very little deviate from the respective observed numbers $a_1, a_2, \ldots$ In Buffon's experiment $m = 2048$. From what he reported it follows that

$a_1 = 1061, a_2 = 494, a_3 = 232, a_4 = 137, a_5 = 56,$
$a_6 = 29, a_7 = 25, a_8 = 8, a_9 = 6.$

Numbers $a_{10}, a_{11}, \ldots$ are absent which means that the number $m$ of the series of trials was not sufficiently large for *heads* not to appear at one or many of the preceding tosses. That number is the sum of the values of $a_1, a_2, \ldots$, and $\mu = 4040$ so that $p = m/\mu = 0.50693$. Applying that value[15] and neglecting fractions, we have

$a_1 = 1038, a_2 = 512, a_3 = 252, a_4 = 124, a_5 = 61,$
$a_6 = 30, a_7 = 15, a_8 = 7, a_9 = 4, a_{10} = 1.$

The next numbers, $a_{11}, a_{12}, \ldots$ were less than unity.
   When comparing this series of calculated values with the numbers resulting from observation, we see that their first terms do not much differ one from another. Then the deviations become larger; for example, the calculated value of $a_7$ only amounts to 3/5 of the observed value. However, it corresponds to an event whose probability is lower than 1/100. When restricting the calculations to the three first terms of the series of the observed numbers, we get

$p = a_1/m = 0.51806, p = 1 - a_2/a_1 = 0.53441,$
$p = 1 - a_3/a_2 = 0.53033$

which very little differ one from another. Their mean, or the third of their sum, is $p = 0.52760$ and barely differs by 0.02 from the ratio $m/\mu$ of $p$ as resulted from all the trials.
   I have chosen this experiment owing to the name of its author and because the work where it is placed renders it authentic. Each reader can make many other experiments of the same kind either with a coin or a die having six faces. In the latter case, the number of arrivals of each in a very large number of trials will be about 1/6 of that number, at least if the die is fair and not poorly manufactured.
   **51.** The theorem, on which the preceding rule is founded, is due to Jakob Bernoulli who had been thinking about its proof for 20 years. What he had done is derived from the binomial formula by means of the following propositions.
   Let $p$ and $q$ be the given chances of contrary events E and F at each trial; suppose also that $g, h$ and $k$ are natural numbers, such that

$p = g/k, q = h/k, g + h = k, p + q = 1.$

Denote by *m, n* and μ other natural numbers connected with *g, h* and *k* by equations

$m = gk, n = hk, \mu = m + n = (g + h)k$.

Therefore, $p/q = m/n$. These numbers can be as large as desired if appropriately increasing *g, h* and *k* without changing their ratios. So,

[1] In the expansion of $(p + q)^\mu$ the largest term corresponds to $p^m q^n$. And since that term is the probability of the arrival of E and F *m* and *n* times (§ 14), it follows that that compound event, i. e., the arrival of the events in a direct proportion of their respective chances, is the most probable of all the compound events which can take place in some number μ of trials.

[2] If this number μ is very large, the ratio of the largest term of the expansion of $(p + q)^\mu$ to the sum of all the terms or to unity will be a very small fraction indefinitely decreasing as μ increases still more. Therefore, in a long series of trials that ratio for the most probable compound event will be nevertheless very small and ever smaller as the trials are continued.

[3] However, if we consider the largest term of $(p + q)^\mu$ and the *l* next terms on both of its sides, and denote by λ the sum of these $(2l + 1)$ consecutive terms, we can always, without changing either *p* or *q*[16], choose a sufficiently large μ for that fraction λ to differ from unity as little as desired. And as μ increases still more, λ will ever closer approach unity.

We conclude that in a large series of trials there always exists a high probability λ[17] that the number of arrivals of event E will be contained between $m \pm l$, and F, between $n \mp l$. Therefore, without changing the interval $2l$ of the limits of those two numbers, we can choose the number μ of the trials sufficiently large for the probability λ arbitrarily to approach certitude as well.

When assuming the ratios of these limits to the number μ of the trials; taking account of the preceding equations; and denoting

$l/\mu = \delta, p \pm \delta = p', q \mp \delta = q'$,

these ratios will be $p'$ and $q'$. And since δ indefinitely decreases as μ increases, it follows that these ratios, varying with μ, with a very high probability also indefinitely approach the chances *p* and *q* of E and F, which the excellent theorem of Jakob Bernoulli[18] has announced.

Regarding the demonstration of these properties of the terms of the expansion of $(p + q)^\mu$, we refer to the works in which they are described, pt. 4 of the *Ars Conjectandi* and the first section of Lacroix (1816). As to the theorem itself, which is placed in the next chapter, it is based on the integral calculus. At present, it should be borne in mind that that theorem essentially supposes the invariability of the chances of the simple events E and F during all the trials. However, in the applications of the calculus of probability either to various physical phenomena or moral things those chances most often vary from one trial to another and, again most often, quite irregularly. Jakob

Bernoulli's theorem is not suited for such kind of problems; however, there exist other more general propositions which are applicable for any variations of the successive chances of the events and serve as the foundation of the most important applications of the theory of probability. They are also demonstrated in the following chapters.

[Now, however,] we will derive the *law of large numbers* already considered in the Preamble as a general fact resulting from observations of all kinds.

**52.** Denote by $p_1$ the chance of event E of some nature at the first of a very large number μ of trials; by $p_2$, that chance at the second trial, …, by $p_\mu$, at the last one, and the mean of all these chances by $\bar{p}$. The mean chance of the contrary event F is $(1 - p_1) + (1 - p_2) + \ldots + (1 - p_\mu)$ divided by μ. Denote it by $\bar{q}$, $\bar{p} + \bar{q} = 1$. Here is one of the general propositions which we intend to consider. Denote the number of arrivals of E and F in that series of trials by *m* and *n*. The ratio of *m* and *n* to μ = *m* + *n* will be very likely almost equal to $\bar{p}$ and $\bar{q}$, and conversely $\bar{p}$ and $\bar{q}$ will be the approximate values of *m*/μ and *n*/μ.

**[1]** If these ratios are calculated for a long series of trials, the mean chances $\bar{p}$ and $\bar{q}$ will become known as well as, by the rule of § 49, the chances themselves, *p* and *q*, of E and F when they are constant. However, for these approximate values $\bar{p}$ and $\bar{q}$ to serve, also approximately, for evaluating the number of arrivals of E and F in a new series of a large number of trials, they should certainly, or at least very probably remain the same for that new series as well as for the first one. This is indeed the case in virtue of another general proposition which we will now formulate.

**[2]** I suppose that, by the nature of the events E and F, the arrival of one of them at each trial can be due to one of the mutually incompatible causes $C_1, C_2, \ldots, C_\nu$, which I regard at first as equally possible. I denote the chance which some cause $C_i$ attaches to the arrival of E by $c_i$, so that at some fixed trial, for example at the first one, that chance is $c_1$, and $c_2$, if the intervened cause is $C_2$, … If only one cause is possible, the chance of E will necessarily be the same at all the trials; however, by our hypothesis at each trial it can have ν equally probable values and therefore varies from one trial to another.

The mean chance $\bar{c}$ of E in a very large number of already made or future trials will likely be almost equal to γ whose magnitude is independent from the number of the trials. Therefore, the mean chance $\bar{p}$ of E can be regarded as constant for two or more series, each consisting of a very large number of trials.

When combining this second proposition with the first one, we conclude that if E arrived, or will arrive *m* times in a very large number μ of trials, and *m*′ times in another very large number μ′ of trials, then very probably and almost exactly *m*/μ = *m*′/μ′. These two ratios will be exactly equal to each other and to the unknown γ if μ and μ′ can be infinite. When their values derived from observation notably differ, there will be a reason to believe that during the time between the two series of trials some of the causes $C_1, C_2, \ldots$ ceased to be possible and other causes had appeared. This will change the chances $c_1, c_2, \ldots$ and therefore the value of γ. Nevertheless, such change will not be certain and we will provide the expression of its probability in terms of the observed difference $m/\mu - m'/\mu'$ and numbers $\mu$ and $\mu'$.

We will return to this consequence of the two preceding propositions in the [discussion of the] Jakob Bernoulli theorem itself. Note that, according to the hypothesis underlying the second proposition, $\gamma$ is the chance of E, unknown but the same in both series of trials. Actually, because each cause $C_1$, $C_2$, … has probability $1/\nu$, that event can arrive at each trial. According to the rule of § 5 the chance of its arrival due to some cause $C_i$ is $c_i/\nu$, and by the rule of § 10 its composite chance will be the sum $c_1/\nu + c_2/\nu + \ldots = \gamma$.

For simplifying my account, I assume that all the causes $C_1$, $C_2$, … are equally possible, but we can suppose that each of them enters into the total number $\nu$ of causes once or many times. They therefore become unequally probable. Denote by $\nu\gamma_i$ the number of times that some cause $C_i$ will be repeated and let $\gamma_i$ be the probability of that cause. Then

$$\gamma = \gamma_1 c_1 + \gamma_2 c_2 + \ldots + \gamma_\nu c_\nu,$$

and at the same time

$$\gamma_1 + \gamma_2 + \ldots + \gamma_\nu = 1,$$

since one of those causes should certainly take place at each trial. If the number of possible causes is infinite, the probability of each becomes infinitely low. Let then $x$ be one of the chances $c_1$, $c_2$, …, $c_\nu$ whose values extend from 0 to 1 and $Ydx$, the probability of the cause which provides that chance $x$ to event E. Then, just like in § 45,

$$\gamma = \int_0^1 Yxdx, \quad \int_0^1 Ydx = 1. \qquad (52.1a, b)$$

**53.** Suppose that instead of two possible events there actually is a given number $\lambda$ of them only one of which should arrive at each trial. Such is the case in which we consider a thing A of some nature susceptible of $\lambda$ known or unknown values $a_1$, $a_2$, … only one of which should take place or took place and is the observed or future event. Let also $c_{ij}$ be the chance that the cause $C_i$, if certain, provides to the value $a_j$ of A. The values of $c_{ij}$ for various indices from $i = 1$ to $\nu$ and from $j = 1$ to $\lambda$ are known or unknown, but for each $j$ we should have

$$c_{i1} + c_{i2} + \ldots + c_{i\lambda} = 1,$$

since, if cause $C_i$ is certain, one of the values $a_1$, $a_2$, …, $a_\lambda$ certainly arrives due to it.

Denote also by $\alpha_j$ the sum of the chances $a_j$ taking place or having taken place in a very large number $\mu$ of consecutive trials divided by that number, i. e., the mean chance of that value $a_j$ of A in that series of trials. When considering $a_j$ as an event E, and the set of the $(\lambda - 1)$ other values of A as the contrary event F, we may take, by the second general proposition of § 51,

$$\alpha_j = \gamma_1 c_{1j} + \gamma_2 c_{2j} + \ldots + \gamma_\nu c_{\nu j}.$$

Now, $\gamma_1, \gamma_2, \ldots, \gamma_\nu$ are as always the probabilities of the various causes which can bring about the events during the series of trials or in other words which can produce the observed or future values of A.

[3] The third general proposition, which I did not yet take up, consists in that the sum of these $\lambda$ values of A divided by their number, or the mean value of that thing likely differs very little from the sum of all its possible values multiplied by their mean chances[19]. And so, denoting by $s$ the sum of the actual values of A, the equality

$$s/\mu = a_1 \alpha_1 + a_2 \alpha_2 + \ldots + a_\lambda \alpha_\lambda$$

will be almost exact and have a high probability so that when denoting by $\delta$ an arbitrarily small fraction, we may always suppose that $\mu$ is sufficiently large for the difference between the sides of that equation to become less than $\delta$ with probability arbitrarily close to 1.

Note also that, due to the preceding expression of $\alpha_j$ and of the values of $\alpha_1, \alpha_2, \ldots$, the right side is independent from $\mu$. So, when that number is very large, the sum $s$ will be appreciably proportional to it. Therefore, denoting by $s'$ the sum of the values of A in another series of a very large number $\mu'$ of trials, the difference $(s/\mu - s'/\mu')$ will likely be very small. Neglecting it, we will have $s/\mu = s'/\mu'$.

In most problems the number $\lambda$ of the possible values of A is infinite. They increase by infinitely small degrees and are contained between given limits $l_1$ and $l_2$. The probability that the cause $C_i$ provides to each of these values thus becomes infinitely low. Denote by $Z_i dz$ the chance which $C_i$ provides to some value $z$ between $l_1$ and $l_2$, then

$$\int_{l_1}^{l_2} Z_i dz = 1. \qquad (53.1)$$

The total chance of this value of $z$, which is very close to its mean value during the series of trials, will be $Zdz$, where, for the sake of brevity,

$$\gamma_1 Z_1 + \gamma_2 Z_2 + \ldots + \gamma_\nu Z_\nu = Z.$$

It follows that

$$\frac{s}{\mu} = \int_{l_1}^{l_2} Zzdz.$$

Z is a known or unknown function of $z$; however, the sum of the fractions $\gamma_1$, $\gamma_2$, ... is unity as is each integral (53.1), and always

$$\int_{l_1}^{l_2} Zdz = 1,$$

whether the number ν of the possible causes is restricted or infinite.

**54.** The law of large numbers consists of two equations

$$m/\mu = m'/\mu', \ s/\mu = s'/\mu',$$

applicable to all cases of eventuality pertaining to physical and moral things. It has two different exact meanings, each of them corresponding to one of those equations. They both are continuously verified, as is seen by the various examples provided in the Preamble. These examples of every kind can not leave any doubt about its [!] generality and exactness. However, owing to the importance of that law as being the necessary basis for applying the calculus of probability in the most interesting cases, it is proper to demonstrate it à priori. Moreover, its proof, based on the propositions of §§ 52 and 53, is advantageous by indicating the very cause of its existence.

By the first equation, the number $m$ of the arrivals of event E of some nature in a very large number $\mu$ of trials can be regarded as proportional to $\mu$. For each kind of things the ratio $m/\mu$ has a special value γ which it exactly reaches if $\mu$ can become infinite. And theory proves that this value is the sum of the possible chances of E at each trial respectively multiplied by the probabilities of the corresponding causes.

What characterizes the set of these causes is the relation which exists for each of them between its probability and the chance it provides, if it [the cause] is certain, to the arrival of E. Provided that that law of probabilities does not change, we observe the permanence of the ratio $m/\mu$ in various series comprised of large numbers of trials. If, on the contrary, that law will change between two series of trials, the mean chance γ also notably changes as a similar change in the value of $m/\mu$ will reveal.

If, in the interval between two series of observations, some circumstances render more probable physical or moral causes providing greatest chances to the arrival of E, the value of γ will increase, and the ratio $m/\mu$ will become larger in the second series than it was in the first one. A contrary result will occur when circumstances increase the probabilities of causes providing least chances to the arrival of E. By the nature of that event, if all possible causes are equally probable, we will have $Y = 1$ and $\gamma = 1/2$, and the number of arrivals of E in a long series of trials will likely very little deviate from a half of its number.

Similarly, if the causes of E have probabilities proportional to the chances provided by them to its arrival, and if their number is still infinite, we will have $Y = ax$, with $a = 2$ for satisfying condition (52.1b) so that γ = 2/3. It follows that in a long series of trials there is a probability closely approaching certainty that the number of arrivals of E will almost double that of the contrary event. However, in most problems the law of probabilities of causes is unknown, the mean chance γ can not be calculated in advance, and it is experience that provides its approximate and likely value when the series of trials is sufficiently long for the ratio $m/\mu$ to become appreciably invariable and adopting it as the value of γ.

Bearing in mind all the variations of the chances during a long series of trials, the almost perfect invariability of that ratio $m/\mu$ for events of each nature is a fact worthy of remarking. We are tempted to attribute that invariability to an intervention of a mysterious power distinct from the physical or moral causes of the events and somehow aiming at order and invariability. However, theory proves that that permanence is necessary until the law of probabilities of causes relative to each kind of events does not begin to change. So in each case we should regard that permanence as a natural state of things which subsists all by itself without being aided by some alien cause. On the contrary, such a cause is needed for accomplishing a notable change. This can be compared with the state of repose of bodies which solely subsists by inertia of matter [unless and] until some alien cause does not begin to disturb it.

**55.** Before considering the second of the two preceding equations, it will be opportune to provide some examples concerning the first one, and proper for interpreting the problem.

Suppose that ν urns $C_1, C_2, \ldots, C_\nu$ contain white and black balls and denote by $c_n$ the chance of drawing a white ball from some urn $C_n$; that chance can be the same for many of these urns. One of them is chosen at random and replaced by a similar urn; the same is done with a second urn, then with a third, … so that the set of urns always remains the same. An indefinitely lengthened series of urns $B_1, B_2, \ldots$, is thus formed only consisting of the given being repeated more or less times.

Denote the chance of extracting a white ball from $B_1$ by $b_1$, from $B_2$ by $b_2, \ldots$ The unbounded sequence of $b_1, b_2, \ldots$ only contains the possibly repeated given chances $c_1, c_2, \ldots$ A ball is drawn from $B_1$, another from $B_2, \ldots$, and a ball from $B_\mu$. Let β be the mean chance of extracting a white ball in these μ successive drawings, then

$\beta = (b_1 + b_2 + \ldots + b_\mu)/\mu.$

Now, urns $C_1, C_2, \ldots$ represent the ν only possible causes of the arrival of a white ball at each trial. So, if μ is a large number, and, as above,

$\gamma = (c_1 + c_2 + \ldots + c_\nu)/\nu$

and $m$ is the number of the drawn white balls, we will get, according to the above, almost exactly and with a high probability

$m/\mu = \beta, \beta = \gamma, m = \mu\gamma.$

Therefore, *m* does not appreciably change if the drawings from the same urns $B_1, B_2, \ldots, B_\mu$ or from $\mu$ other urns are repeated. And if the balls are drawn from another very large number $\mu'$ of urns, the number of the extracted white balls will very probably be approximately $\mu'm/\mu$. If we extract at random $\mu$ balls from the set of urns $C_1, C_2, \ldots$, each time replacing it in its previous urn, the chance of drawing a white ball will be the same for all the trials, and, by the rule of § 10, equal to $\gamma$. According to the rule of § 49, when their number is very large, the number of drawn white balls, just like in the previous problem, will likely be almost equal to $\mu\gamma$. However, these two problems are essentially different and their results only coincide when $\mu$ is a very large number. Otherwise the chance of extracting a given number *m* of white balls depends, in the first problem, not only on the given system of urns $C_1, C_2, \ldots$, but on the system of urns $B_1, B_2, \ldots$ derived by chance. For example, I reduce the given system to three urns, $C_1, C_2, C_3$, and take $\mu = 2$ and $m = 1$ and will, as required, determine the chance of drawing a white ball from one of the two urns, $B_1$ and $B_2$, and a black ball from the other urn. For those two urns nine different combinations are possible:

$B_1 = B_2 = C_1$;  $B_1 = B_2 = C_2$;  $B_1 = B_2 = C_3$;  $B_1 = C_1$ and $B_2 = C_2$;
$B_1 = C_1$ and $B_2 = C_3$;  $B_1 = C_2$ and $B_2 = C_3$;  $B_1 = C_2$ and $B_2 = C_1$;
$B_1 = C_3$ and $B_2 = C_1$;  $B_1 = C_3$ and $B_2 = C_2$

For each of them the required chance has a determined value and its possible values will be

$2c_1(1 - c_1);\ 2c_2(1 - c_2);\ 2c_3(1 - c_3)$

for the first three of them,

$c_1(1 - c_2) + c_2(1 - c_1);\ c_1(1 - c_3) + c_3(1 - c_1);$
$c_2(1 - c_3) + c_3(1 - c_2)$

for the three intermediate, and the same for the three last ones.

It is easy to see that the mean value of these nine chances should be the chance of extracting one white and one black ball when drawing at random for the first time from the group of the three urns $C_1, C_2, C_3$, and then, for the second time, after replacing the extracted ball in the urn from which it was drawn. Actually, that chance will equal twice the product

$1/3(c_1 + c_2 + c_3)[1 - 1/3(c_1 + c_2 + c_3)]$,

which is also 1/9 of the sum of the preceding nine chances.

Before the system of urns $B_1, B_2, B_3, \ldots$ was formed out of the system of given urns, we had no reason to believe that urn $B_n$ will be rather one than another of $C_1, C_2, C_3 \ldots$ For us, the probability of extracting a white ball at the *n*-th drawing will be the sum of the chances $c_1, c_2, c_3, \ldots$ divided by their number, i. e., will be $\gamma$. However, although it is the same for all the drawings, and their number $\mu$ is as large as desired, we are not authorized to conclude, solely by the rule of § 49, that the number of extractions $m$ of a white ball from the urns $B_1$, $B_2$, $B_3$, … should likely very little deviate from μγ. Indeed, we should not overlook that that rule is based on the chance proper to the considered event rather than on its probability or on the reason we can have for believing in its arrival.

**56.** As a second example, I suppose that there is a very large number of 5-franc coins which I denote by $A_1$, $A_2$, … Each one is tossed and falls to the ground. Let the chance of *heads* for coin $A_i$ depending on its physical constitution be $a_i$. It is not known in advance but determined by experience after tossing the coin a very large number $m$ of times. And since this chance remains constant during that series of trials, we may, by the rule of § 49, take $a_i = n_i/m$ where $n$ is the number of arrivals of *heads*, as a likely and close approximant of that value. It will serve for calculating the probabilities of various future events concerning the same coin. We can fairly bet $m$ against $(m − n_i)$ on the appearance of *heads* at a new trial, $m^2$ against $(m^2 − n_i^2)$ on its appearance twice in succession, $2n_i(m − n_i)$ against $[m^2 − 2n_i(m − n_i)]$ on its appearance only once in two trials etc.

In a new series of a very large number $m'$ of trials the number $n'_i$ of the arrivals of *heads*, again by the rule of § 49, will likely be very close to $m'a_i$. Two ratios, $n_i/m$ and $n'_i/m'$, should very little deviate from each other. However, the value of $a_i$ given by experience is only very probable but not certain, so the probability of a small difference between these ratios will not be as high as it would have been if that value were certain and given in advance, see below.

Suppose that, instead of taking the same coin a very large number of times, we toss successively a very large number μ of such coins randomly chosen from those minted in the same way, and let $n$ be the number of arrivals of *heads*. Denote by α, the mean chance of that event, not only for those being tossed, but for all the coins of the same kind and minting. In virtue of the two general propositions of § 52 we will likely have an almost exact equality α = $n$/μ as though the unknown chances $a_1$, $a_2$, … were equal one to another.

Depending on $n$/μ > 1/2 or < 1/2, we conclude that for the 5-franc coins of that minting the chance of the arrival of *heads* is generally greater or less than that for *tails*. In particular, for coin $A_i$ the chance $a_i$ will differ from α; if $n$/μ > 1/2, it can happen that $n_i$/μ < 1/2 or vice versa.

If those coins are tossed again, or, more generally, when tossing a very large number μ′ of other coins of the same kind and minting, with the same images, the magnitude α, as defined, will not change, and if $n$′ is the number of the arrivals of *heads* in that new series of trials, we should have

$$n'/\mu' = n/\mu \qquad (56.1)$$

and, for the same coin $A_i$ in those two different series of trials, $n_i/m = n_i'/m'$. However, when the coins in those series are not of the same kind or minting, those ratios are in general unequal; just the same, the ratio $n_i/m$ varies from one coin to another.

**57.** Although the constant chances and the mean chances of the events are determined by experience in the same way and with the same probability[20], there are essential differences in their possible usages. Both the mean and the constant chance immediately lead to the probability of the arrival of the event under consideration in a new isolated trial. This, however, is not always so when dealing with the arrival of a compound event consisting of them.

As an example, I consider the coincidence of the outcomes of two consecutive tosses of a 5-franc coin. Two different cases should be examined. We may suppose that the two trials are made with two coins, either different or not, chosen at random from all the $\lambda$ coins $A_1$, $A_2$, … of the same minting or with a single coin, also selected at random. In the first case, the probability of coincidence only depends on $\alpha$ of § 56 and the same happens when the chances are constant. In the second case, the probability additionally depends on another unknown magnitude by which it differs from its value when the chances are invariable.

To explain this, I note that for two coins, $A_i$ and $A_j$, the probability of the coincidence of the outcomes will be

$$a_i a_j + (1 - a_i)(1 - a_j).$$

In the first case each of the coins $A_1$, $A_2$, … can be combined with itself or with any other coin so that the number of all the equally possible combinations is $\lambda^2$. Denoting by $s$ the composite probability of coincidence, we will have by the rule of § 10

$$s = (1/\lambda^2)[\sum a_i \sum a_j + \sum(1 - a_i)\sum(1 - a_j)],$$

where the sums extend from $i, j = 1$ to $i, j = \lambda$.

Let

$$\alpha = 1/2(1 + k),\ a_i = 1/2(1 + k + \delta_i),\ a_j = 1/2(1 + k + \delta_j),$$

where $k$, $\delta_i$, $\delta_j$, $k + \delta_i$, $k + \delta_j$ are positive or negative fractions, the first of which is derived from $n/\mu$ of § 56 and given by observations, and the others vary from one coin to another, so that

$$\sum \delta_i = 0,\ \sum \delta_j = 0.$$

At the same time

$$1 - a_i = 1/2(1 - k - \delta_i),\ 1 - a_j = 1/2(1 - k - \delta_j).$$

Because of the preceding equations, the sums which enter the expression of $s$ are

$$\sum a_i = \sum a_j = 1/2\lambda(1 + k),\ \sum(1 - a_i) = \sum(1 - a_j) = 1/2\lambda(1 - k)$$

so that

$$s = 1/2(1 + k^2).$$

This magnitude only depends on $k$ or the mean chance α of the arrival of *heads*, but not on the inequalities of the chances $δ_1$, $δ_2$, …

When repeating the tossing of the two coins chosen at random a very large number $a$ of times, $s$ will also be the mean chance of coincidences in double trials at that series. Denoting by $b$ the number of coincidences, we will, according to § 52 approximately get $b = as$, which can be verified by experience.

In the second case, each double trial should be made with the same coin; the probability of coincidence for coin $A_i$ will be

$$a_i^2 + (1 - a_i)^2.$$

Denoting by $s'$ the composite probability of coincidence, we conclude that

$$s' = (1/\lambda)[\sum a_i^2 + \sum(1 - a_i)^2] = 1/2(1 + k^2 + h^2), \; h^2 = (1/\lambda)\sum δ_i^2.$$

It is seen that this probability exceeds the probability $s$ of the first case and depends on a new unknown $h$ which in turn depends on the inequalities between $δ_1$, $δ_2$, …

When repeating a very large number $a'$ of times the double trials of the same randomly chosen coin, $s'$ will express the probability of coincidence in that series. Denote by $b'$ the number of the coincidences, then almost exactly $b' = a's'$, which will serve for determining the value of $h$, whereas $k$ is already known.

**58.** I note that, when tossing the same coin randomly chosen from $A_1$, $A_2$, …, the probability of coincidence of three outcomes is expressed by the preceding probability $s'$ and is therefore known without making new trials. For some coin $A_i$ that probability is

$$a_i^3 + (1 - a_i)^3,$$

and the composite probability

$$s'' = \frac{1}{\lambda}[\sum a_i^3 + \sum(1 - a_i^3)] = \frac{1}{4}[1 + 3(k^2 + h^2)] = \frac{1}{2}(3s' - 1).$$

Magnitude $s''$ is also the mean chance of triple coincidences in a very long series of such trials. Denote the number of triples by $a''$, and the number of coincidences by $b''$, then $b'' = a''s''$. Substituting $b'/a'$ and $b''/a''$ instead of $s'$ and $s''$, we will derive the relation

$$a'a'' = 3b'a'' - 2b''a'$$

between $a'$, $a''$, $b'$, $b''$ which will be the more precise and the more probable the greater are those numbers.

Since it is independent from the law of $a_1$, $a_2$, … it therefore subsists just as well when all these are the same, i. e., when instead of changing the coin for each double and triple trial we always toss the same one.

Therefore, if tossing the same coin a very large number 6c of times, then separating that series in double trials consisting of the first two simple trials, of the third and fourth, …; then separating that series in triple trials consisting of the first three, then of the fourth, fifth and sixth trial, … and applying the preceding equation to these two new series, we will have

$a' = 3c, a'' = 2c,$

so that $c = b' − b''$, which signifies that the number of double coincidences less that number of triple coincidences equals 1/6 of the total number of simple trials. Simple relations can be obtained for cases in which more trials are combined.

**59.** These coincidences are immediately applicable to the births of boys and girls. For this, suffice it to replace the coins $A_1, A_2, …$ by so many different marriages and assume that the chance of the birth of a boy in some marriage $A_i$ is $a_i$.

In France, the yearly number of births of both sexes amounts to almost a million and for that total number observations prove that the ratio of male to female births exceeds unity by about 1/15. During 10 years from 1817 to 1826, its mean value was 1.0656 from which its extreme values barely deviated by 0.005 in either direction. I based my memoir (1830) on the sex ratio at birth of that period. From 1817 to 1833 inclusive the mean ratio was 1.0619, not more than 0.005 different from its value during the 10 first years of those 17.

The cause of the excess of male births is unknown, and there is room to believe that it essentially differs from one marriage to another and that the chances $a_1, a_2, …$ are very unequal so that many of them are undoubtedly lower than 1/2. Nevertheless, it is seen that the yearly sex ratio very little varied during 17 years which provides a remarkable verification of the *law of large numbers*.

When taking 16/31 for that ratio for a large number of male births and the corresponding total number of them, it will also be the mean chance of a male birth, so that the magnitude $k$ from § 57 is 1/31. We do not know whether the chance of a male birth remains the same for each infant born into the same family or whether it varies, for example, as it does from one marriage to another. In the second case, the mean chance of the coincidence of the sexes of two first-born will be $(1 + k^2)/2$ which only exceeds 1/2 by about 0.0005. […] In the first case the first of these two numbers can exceed the second much more than by a half of the second because of the unknown magnitude $h$ included in the expression $(1 + k^2 + h^2)/2$ of the mean chance of coincidences.

The relation of § 58 is always is applicable to the coincidences of the sexes of the two and three first-born in a very large number of families.

**60.** By virtue of the second equation of § 54, if A is a thing, susceptible of different values at each trial, the sum of its values observed in a long series of trials will likely be almost proportional to their number. For a given thing A the ratio of that sum to that number, as it increases still more, indefinitely converges to a special value, which it reaches if that number can become infinite and which depends on the law of probabilities of the diverse possible values of A. Remarks similar to those, which we stated in § 54 when considering its first equation, can be formulated about that ratio.

The second equation, or rather

$$\frac{s}{\mu} = \int_{l_1}^{l_2} Zz\,dz,$$

just like the first one, leads to numerous and useful applications. I suppose for example that α is an angle which we desire to measure; it exists, its value is unique and fixed. However, because of variable inevitable errors of observation, a measured angle is a thing susceptible of an infinite number of different values.

I take this angle, measured many times successively, as a thing A, so that $Zdz$ expresses the chance of some value $z$ of A, resulting from the construction of the instrument and the ability of the observer. Let $k$ be the abscissa of the centre of gravity of the area of a plane curve, whose abscissa and ordinate are $z$ and $Z$ extending from $z = l_1$ to $l_2$, when denoting, like in § 53, the possible limits of A by $l_1$ and $l_2$. Let also

$$z = k + x,\ l_1 = k + h_1,\ l_2 = k + h_2$$

and represent by $X$ the value of $Z$ when $z = k + x$. Then

$$\int_{l_1}^{l_2} Z\,dz = \int_{h_1}^{h_2} X\,dx = 1,\ \int_{h_1}^{h_2} Xx\,dx = 0.$$

Therefore, by the equation above, almost exactly $s/\mu = k$ where $s$ is the sum of the values of A obtained in a large number μ of trials. It is the constant $k$ to which the mean value $s/\mu$ converges ever closer as μ increases further. However, even when that ratio becomes appreciably constant, i. e., when it becomes appreciably the same in many series of other large numbers of measurements, it can sometimes happen that that mean essentially differs from the angle α which we wish to measure. It will always be the approximate value of the constant γ, which possibly will not at all coincide with that angle[21]. Actually, let

$$z = \alpha + u,\ l_1 = \alpha + g_1,\ l_2 = \alpha + g_2.$$

Denote by $U$ the value of $Z$ when $(\alpha + u)$ is substituted instead of $z$, then

$$\int_{l_1}^{l_2} Z\,dz = \int_{g_1}^{g_2} U\,du = 1,\ k = \alpha + \int_{g_1}^{g_2} Uu\,du.$$

The difference $u$ between the angle α and a possible value $z$ of the measured angle A is a possible error of the instrument and the observer; it can be positive or negative and extends from $u = g_1$ to $g_2$ and its infinitely small probability is $Udu$. And if there is no cause connected with the construction of the instrument for attaching greater chances to positive rather than to negative errors or vice versa, and if it is the same about the method of the observer's work, the limits $g_1$ and $g_2$ will be equal in magnitude and have opposite signs and the function $U$ will take equal values when the variable $u$ has values equal in magnitude and opposite in signs. Then

$$\int_{g_1}^{g_2} Uudu = 0, \ k = \alpha.$$

In this most usual case the ratio $s/\mu$ will be the approximate value of $\alpha$. However, if the instrument, owing to its construction, or the observer, because of his manner of sighting, attach more weight to positive or negative errors, the preceding integral will not disappear, the constants $\alpha$ and $k$ will differ from each other, and the ratio $s/\mu$ will in general appreciably deviate from the veritable value of $\alpha$. We can only notice that circumstance when measuring the same angle either with another instrument or by other observers. I restrict my account by indicating this application of the theory of probability. Concerning the errors of observation and the methods of calculation proper for decreasing and evaluating their influence, I refer to Laplace (1812) and to my memoirs (1824; 1829).

**61.** As a second example of the equation cited at the beginning of § 60, let us suppose that all the causes $C_1$, $C_2$, … determine the chances of the duration of human life in a given country and at a determined period. Those causes include the diverse physical constitutions of the newborn babies, the well-being of the population, the diseases that restrict that duration, and undoubtedly some of the causes resulting from life itself which hinder its extension beyond the limits that are never crossed. Actually, there is room for believing that, had the diseases been the only causes of death, and had they been, so to say, accidental, some people from the immense number of all the living, would have escaped these dangers for more than two centuries. This, however, was never observed.

The thing A will be the lifetime of a newborn baby, $z$ will express a possible value of A, and $Zdz$, the chance of $z$ resulting from all the possible causes which can determine it, although not for a particular infant, but for mankind in the place and at the period under consideration.

Imagine that a certain physical constitution at birth provides the chance $Z'dz$ to live exactly $z$, that another constitution provides the chance $Z''dz$ to live to the same age, … Let also $\varsigma'$, $\varsigma''$, … be the probabilities of those various constitutions. Owing to those causes, the function $Z$ will be $Z'\varsigma' + Z''\varsigma'' + \ldots$ extending over all possible constitutions. And if their number is infinite, $Z$ will become a definite integral with an unknown but definite value. In a country where babies are born more robust or better constituted, that integral will undoubtedly have a greater value. In each country it can differ for the two sexes, and undoubtedly the values of $Z'$, $Z''$, … also depend on possible diseases and the well-being of the population. The function $Z$ and therefore the integral of $Zzdz$ will differ at two epochs remote from each other if, for example, during the interval between them some disease disappears or the well-being of the population betters with the progress of society.

If desired, it is possible to take the limits $l_1$ and $l_2$ of that integral as 0 and ∞ and consider $Z$ as a function disappearing beyond some value of $z$, which is unknown just as the form of $Z$. The observed values of A will be the ages at death of a very large number $\mu$ of individuals born in the same country and at the same epoch. Denoting by $s$ the sum of these ages, we will likely have almost exactly

$$\frac{s}{\mu} = \int_0^\infty Zzdz.$$

The ratio $s/\mu$ is therefore called the *mean life*, constant for each country if only some of the causes $C_1$, $C_2$, … whether known or not, do not experience any notable change.

In France, it is supposed that the mean life is about 29 years, but that evaluation was based on observations made before the usage of the [Jennerian] vaccine and is already very dated. Today, that magnitude should be much longer, and its new determination is desirable, separately for men and women, for the differing strata of society and various parts of the kingdom.

We also consider the mean life at different ages, and then $s$ is the remaining number of years of life for a very large number $\mu$ of individuals. The ratio $s/\mu$ is the mean life for that age together with which it varies, but remains constant for a given age. It is supposed to attain its maximal value at an age between 4 and 5 years when it amounts to 43 years. *Mortality tables* have another aim. For a very large number $\mu$ of individuals born in the same country and at the same epoch they make known the number of those who will continue living a year, two, three years, … until no one is left alive. Denote by $m$ the number of those living at a given age; owing to the first equation of § 54, the ratio $m/\mu$ is appreciably invariable at least until a very advanced age, and $m$ is not [yet] becoming a very small number. For the age of 100 years, that invariability consists in that $m/\mu$ is always a very small fraction.

Instead of assuming that in the integral

$$\int_0^\infty Zzdz$$

$z$ varies by infinitely small steps, let that variable increase by very small intervals, and for the sake of definiteness suppose that each of them is of unit time. Denote by $h_1$, $h_2$, … the series of the values of $z$, and by $H_1$, $H_2$, … the corresponding values of $Z$. Then the sum

$v = H_1h_1 + H_2h_2 + …$

is known to be an approximate value of that integral. The same sum will be the mean life at birth; $H_n$ is the chance of dying at age $h_n$, so that with regard to the duration of human life the mean life $v$ can be considered the mathematical expectation (§ 23) of a newborn baby with an unknown physical constitution. Nevertheless, according to mortality tables, more than half of a very large number of infants will die before attaining that age.

**62.** Suppose, for the last example that for a given place and day of year the excess of the height of the sea level caused by the simultaneous action of the Sun and the Moon was calculated. We assume as the thing A the yearly differences of that excess as calculated and observed in the same place and at the same epoch.

The values of A vary from year to year because of the winds which can blow in that place and at that epoch and which determine the chances of those diverse values. When considering all the possible directions and intensities of those winds, their respective probabilities and the chances that these causes attach to some value $z$ of A, the integral

$$\int_{l_1}^{l_2} Zz\,dz$$

will take an unknown but definite value which remains constant if the law of probabilities of each possible wind does not change.

The ratio $s/\mu$ will also be almost invariable; here, $s$ is the sum of the values of A observed during a long sequence of years. We do not know in advance whether $s/\mu$ is zero or a negligible fraction; that is, whether the influence of the winds on the general laws of the *tides* is insensible. Only experience can let us know the value of that ratio and tell us whether it varies in different periods of the year and in different places of observation, on the shorelines, in ports or at sea.

For finding out the influence of a certain wind the only applicable values of A are those observed while it was blowing. To avoid the need to have a very considerable number of years of observation, it is possible to study the values corresponding to many consecutive days during which the wind's direction changed inconsiderably. Many scholars are now occupied by that examination which demands work over a long period of time and will not fail to lead to interesting results.

**63.** The exposition of the rules of the calculus of probability and their general consequences in this and the previous chapter is actually completed and I return to the notion of *cause* and *effect* which was only indicated in § 27.

The cause proper to a thing E is, as stated in that section, another thing C possessing the *power* to produce necessarily E, whatever is its nature and the manner in which it is exercised. Thus, what is called the attraction of the Earth is a certain thing which has the power to cause unsupported bodies to fall on the surface of our planet. Just the same, in our volition resides the possibility of producing by means of our muscles and nerves a part of those movements which are for this reason called voluntary. In nature, the thing E sometimes has only one single cause C which is able to produce it, so that the observation of E always supposes the intervention of C. In other cases, that thing can be attributed to many distinct causes which act together or mutually exclude each other so that the production of E can only be due to one of them.

Such are the simplest ideas concerning the principle of causality, which I believe to be generally admitted. Nevertheless, the illustrious historian of England expressed a different opinion about that metaphysical point which deserves to be studied in more detail and on which the calculus of probability can throw much light. According to Hume (1748), we can only imagine *causality* as a *concurrence* rather than a necessary *connection* between what we call *cause* and *effect*. And for us that concurrence is only a strong presumption resulting from what we had observed a large number of times. If we had observed something an inconsiderable number of times, we are judging nature by a far too small sample and presume that the same will be reproduced in the future. Others share that opinion and try to justify it by the rules of the probability of future events determined by observing past events.

Hume, however, goes further. Even without turning to those laws [!] of probability, he thinks that the habit of seeing the effect following the cause produces in our mind some association of ideas which leads us to believe that the effect will occur when the cause takes place[22]. This is perhaps actually true for most of those who do not examine the principle of their belief or its degree of probability. For them, that association of ideas must be compared with what goes in our mind about the name of a thing and the thing itself. The name reminds us of the thing independently from our thoughts or volition.

One of the examples chosen by the author for describing his opinion is the shock of a moving body against a free body in repose and the movement of the latter after the collision. The concurrence of the shock and the movement of the hit body is actually an event which we have seen a very large number of times without ever seeing a contrary event. For us, this is sufficient, leaving aside all other considerations, to believe that there exists a very high probability that that concurrence will take place in the future as well.

The same happens about all concurrences of causes and effects which we note without exception day after day. Their probability is nourished, so to say, by that continuing experience, and reason or calculation along with habit strongly assure us in that in the future those causes will always be followed by their effects. However, by the rules stated above, in case of a phenomenon only observed an inconsiderable number of times after the cause which we assign to it, there will only be not a very high probability for the future concurrence of that cause and that effect.

Nevertheless, it often happens that we do not doubt the reproduction of that phenomenon if its cause takes place anew. That assurance assumes that our mind attributes to the cause some *power* for producing its effect and admits a necessary connection between these two things independently from the more or less large number of the observed concurrences. Thus, Oersted discovered that, when connecting the two poles of a voltaic pile by a metallic wire, a magnetic needle freely suspended nearby deviates from its natural direction. After repeating that fundamental experiment a small number of times, the illustrious physicist undoubtedly became convinced that that phenomenon will not fail to be continuously reproduced in the future. Nevertheless, if our reason to believe that that reproduction is only based on the concurrence of the voltaic pile and the deviation of the magnetic needle observed, for example, about 10 times, the probability of its occurrence at a new trial will only be 11/12 (§ 46). For a new series of 10 trials 11 can be bet against 10, which is almost even money, on that event to occur without interruption. And for a longer series of future experiments it will be reasonable to believe that that event will not be reproduced without interruption.

As another example, I will cite Biot's successful application of sorts of [the indicated] to the chemical composition of bodies by observing the *successive polarization of light*[23]. A long time ago he had established that [possibility] for homogeneous and non-crystallized media. After an inconsiderable number of thorough observations we recognize that a given substance deflected a polarized ray, let us say, to the right of the observer and that these deviations were sufficiently large for leaving no doubts about their direction. For us, this suffices to become convinced that the same substance will always deflect light in the same direction; assured as though it is a thing in which no one doubts.

But still, the concurrence of that substance and a deviation to the right observed not very many times only provides a feeble probability, even lower than 1/2, that in the same or in a somewhat larger number of new trials deviation to the right will occur. These and other easily imagined examples prove, as it seems to me, that the confidence of our mind in effects following after their causes can not be solely based on former more or less repeated observations. We actually see that, independently from any habit of our mind, the possibility alone of a certain ability of the cause to produce necessarily its effect greatly increases the reason to believe in that recurrence and, in spite of the former observations being few in number, can render its probability very close to certainty.

**64.** Before a phenomenon P was observed, if knowing that it will arrive or not in the entire series of experiments, we therefore admit that the existence of a cause C capable of necessarily producing it is not impossible. We also admit that prior to these experiments the existence of such a cause had a certain probability *p* which results from particular considerations rendering it more or less likely. Suppose also that P was observed in each of *n* experiments; after that observation the probability of the existence of C changes and becomes *w*; it is required to determine it.

However thoroughly we attempted to decrease the influence of other causes, able to produce phenomenon P at each trial in the absence of C, we may nevertheless believe that that influence did not completely disappear. Consequently, suppose that there exist certain causes $B_1$, $B_2$, …, $B_n$, known or unknown, which can also engender that phenomenon in the absence of C by combining themselves with hazard (§ 27), namely, $B_1$, at the first trial, ..., $B_n$, at the last one. Let in general $r_i$ be the probability of the existence of $B_i$, multiplied by the chance that that cause, if certain, provides to the arrival of P, and denote for the sake of brevity

$$r_1 r_2 ... r_n = \rho.$$

That product is the probability that phenomenon P arrives at all the $n$ experiments as a result of the set of the causes $B_1$, $B_2$, ... if cause C does not exist. Now, $(1 - p)$ is the probability of the non-existence of C and in that case $(1 - p)\rho$ is the probability of the observed event, of the invariable arrival of P. According to the contrary assumption, its probability is $p$; i. e., that it is just the probability of the existence of C before the observation since that cause necessarily leads to the arrival of P at all the trials. Therefore, by the rule of § 28, the probability of this second hypothesis on the existence of C after observation is

$$w = \frac{p}{p + (1-p)\rho},$$

and of its non-existence,

$$1 - w = \frac{(1-p)\rho}{p + (1-p)\rho}.$$

We can arrive at that result when considering each of the $n$ experiments one by one rather than all of them at once. Actually, according to the hypothesis, the probability of the existence of C before the first experiment was $p$; after that it became $p'$; after the second experiment, $p''$, ... so that

$$p' = \frac{p}{p + (1-p)r_1}, \quad 1 - p' = \frac{(1-p)r_1}{p + (1-p)r_1},$$

$$p'' = \frac{p'}{p' + (1-p')r_2}, \quad 1 - p'' = \frac{(1-p')r_2}{p' + (1-p')r_2}, ...$$

Eliminating at first $p'$ and $(1 - p')$ from the values of $p''$ and $(1 - p'')$, then eliminating $p''$ and $(1 - p'')$ from $p'''$ and $(1 - p''')$, ... we obtain the preceding expressions of $w$ and $(1 - w)$ of the existence and non-existence of C after the $n$-th experiment. Now let $w'$ be the probability of the arrival of P without interruption in a new series of $n'$ experiments. Whatever is this number $n'$, the probability of that event occurring owing to cause C, if certain, is the probability $w$ of the existence of cause C derived from the first $n$ experiments. In its absence, the arrival of P can also be due to the other causes $B'_1$, $B'_2$, ..., $B'_{n'}$, similar to $B_1$, $B_2$, ..., $B_n$, whose influence it was impossible to avoid completely.

Denote the new values of $r_1$, $r_2$, ..., $r_n$ by $r'_1$, $r'_2$, ..., $r'_{n'}$, so that $r'_i$ with respect to $B'_i$ is the same as was $r_i$ with respect to $B_i$, and let

$$r'_1 r'_2 \ldots r'_{n'} = \rho'.$$

If cause C does not exist, the probability of the arrival of P in these future experiments will be $(1 - w)\rho'$, so

$$w' = w + (1 - w)\rho'$$

is the composite expression of $w'$. When inserting the preceding values of $w$ and $(1 - w)$,

$$w' = \frac{p + (1-p)\rho\rho'}{p + (1-p)\rho}.$$

The expressions of $w$ and $w'$ show how the probability of the existence of C, which could have been feeble before the observation of P, can become very high after P was observed an inconsiderable number of times, and after that cause C can provide a probability very close to certainty for the constant occurrence of that phenomenon in future experiments. Suppose, for example, that, owing to some reason, and, for that matter, to a prejudice of our mind, the prior probability of C was only 1/100,000. Let us also admit that the influence of accidental causes, in spite of the precautions taken for their avoidance, can still be such, that each of the magnitudes $r_1$, $r_2$, ... equals 1/10 or less. Then, if P was only observed 10 times without interruption,

$$p = 0.00001,\ \rho < 0.00001 p,\ w > \frac{1}{1 + 0.00001(1 - p)}.$$

Therefore, the probability of the existence of C after observation will differ from unity less than by 1/100,000 and its non-existence becomes less probable than its existence before observation. Whatever is the value of $n'$, the probability $w'$ that P will invariably arrive at a new series of $n'$ experiments exceeds that of the existence of C or can not be lower.

**65.** In that application of the calculus of probability, cause C was considered in an abstract manner, independently from any theory which brings phenomenon P to conformity with more general laws, provides an exact explication according to the cause to which it is attributed and thus still more heightens the probability of the existence of that cause. We assumed that phenomenon P had taken place without interruption, and the preceding calculations aimed at proving that any belief in its future reproduction after being only observed a small number of times can not only be based on the idea that there exists just one cause capable of necessarily producing a phenomenon of that nature. And in addition the calculus of probability is unable either to inform us about the essence of that effective cause or to determine which cause is more probable among those that were able to produce necessarily the phenomenon and could have been attributed to it.

If P fails at one or many trials, but we are nevertheless certain that the cause C, if existing, is capable of producing it necessarily, and should have produced it at each trial, it becomes evident that there exists neither that cause, nor any other of the same kind. However, besides such causes there exist others which act at each trial but are only able to provide a certain chance to the arrival of a phenomenon P (§ 27), when combining themselves with hazard, or variable and irregular causes which act only sometimes (not to be, however, confused with hazard). They can influence the mean chance of the arrival of P in a long sequence of experiments, and therefore (§ 52) the number of its occurrences, past or future, divided by the total number of the trials.

And still, when thoroughly decreasing as much as possible the influence of accidental causes, so that they can be supposed insensible, and observing the phenomenon P $m$ times in a very large number $\mu$ of trials, we will have a very high probability that there exists a permanent cause favourable or contrary to the production of P depending on whether $m$ is notably larger or smaller than $\mu/2$.

Consider for example the case of a two-sided body tossed very many times in the air. The existence of a favourable or contrary cause of the arrival of a determined side can be regarded extremely probable when the numbers of the two possible outcomes notably differ as in the Buffon experiment (§ 50)[24].

What is that permanent cause? The calculus of probability only proves its necessity but is unable to indicate its nature. It is the laws of mechanics which tell us how much greater should be the weight of one of the parts of that body. But, owing to the difficulty involved in that problem, they do not determine either the effects of such a cause or the chance it provides to the arrival of one or another side. This can only be found by experience.

By suchlike means it is possible, as Laplace (1814, p. 133) suggested, to determine the existence or non-existence of certain mysterious causes which are not absolutely impossible beforehand and are unable to produce necessarily the phenomena attributed to them. To achieve this aim, long series of trials are needed with the influence of accidental causes being excluded as much as possible and the numbers of arrivals and failures of a phenomenon exactly registered. If the ratio of the former number to the latter is notably larger than unity, the probability of the existence of some cause and of the chance it provides to the production of that phenomenon becomes very high.

If gamblers A and B play a very large number $\mu$ of games with A winning $m$ times and the ratio $m/\mu$ exceeding 1/2 by not a very small fraction, the existence of a cause favouring A can be regarded as almost certain. If no gambler gives a start to the other one, that cause is the superiority of A over B with the ratio $m/\mu$ providing, so to say, its measure. When playing cards, *piquet* for example, the results of each game can only depend on the difference of the gambler's abilities and the distribution of the cards among them. If none cheats, that distribution is random; it can influence the ratio of the numbers of games won by the gamblers when these numbers are inconsiderable.

And that can be called *luck* and *misfortune*, if these notions are not attached to one or another gambler. Indeed, it is absurd to suppose that there exists some connection between them and the cards only distributed by chance. Each time those obtained by one of them could have just as well been dealt to the other. But in a sufficiently long series of games it is only the gambler's ability that can influence the chances. In a long run, only able gamblers are lucky, and vice versa. If A and B play a new long series of µ′ games, then there will be a high probability that the number of games won by A will be very close to µ ′m/µ. Otherwise, we will be justified in thinking that during the interval between those series of games the superiority of A over B either increased or decreased.

## Notes

**1.** See Note 5 to Chapter 1.

**2**. Event E and all the others introduced by Poisson were thus random.

**3.** Was a set of causes necessary? In general, this definition can be understood as randomness in a wide sense, not obeying any stochastic laws and therefore not studied by the theory of probability (Kolmogorov 1983). An example is provided by deviations from the tree of animal life (Lamarck 1815, p. 133). Tables of random numbers are however needed, at least in various applications. It occurs that, according to Poisson, randomness only acts upon random events!

**4.** Fienberg (1971) noted that it is difficult to ensure equal chances of extracting balls from an urn even independently from the circumstance mentioned by Poisson.

**5.** See Note 21.

**6.** Failing to account for that circumstance, the solution of this problem (Poisson 1830, § 17) is not exact, and I have derived from it a false corollary. Poisson

**7.** See Note 4 to the Preamble.

**8.** Poisson is considering four cases: the witness is (is not) mistaken and wishes (does not wish) to deceive. Their repetition makes his considerations dull and difficult to read. It would have been better to enumerate those cases and to refer to, say, case II.

**9.** Laplace, see Note 10 to Chapter 1, considered a similar example.

**10.** By the end of the 19th century physicists are known to have thought that soon nothing new will be left in their science. History of science had refuted such opinions.

**11.** Cf. Gauss (1823, §§ 18 and 19): independent linear functions of observations should not contain common observations.

**12.** See Note 3 to the Contents.

**13.** Poisson several times (also in Chapter 4) made use of that loose phrase.

**14.** I am only able to provide the end of this phrase in its original French and in its German translation of 1841:

*la probabilité de similitude dans deux épreuves consécutives, est la même que s'il y avait, entre les chances de G et de l'événement contraire, une différence $1/\sqrt{3}$, sans que l'on connût la chance la plus favorable.*

*so ist die Wahrscheinlichkeit der Übereinstimmung zwei auf einander folgenden Versuche dieselbe, als wenn zwischen der Wahrscheinlichkeit der Ereignisse G und der des entgegengesetzten Ereignisses eine Differenz $1/\sqrt{3}$ stattfinde, ohne daß man der günstigste Wert kennt.*

**15.** I (Note 9 to Chapter 1) have indicated that Poisson had been calculating with superfluous numbers of significant digits.

**16.** But of course any change of *p* (or *q*) results in an appropriate change of *q* (or *p*).

**17.** It is easy to understand how λ became a probability, but an explanation should not have been lacking.

**18.** Which equation did Poisson (above) mean? And the last phrase is incomprehensible. A few lines below, Poisson wrongly stated that the Jakob Bernoulli theorem is *based on the integral calculus*.

**19.** In several instances, e. g., in § 54, Poisson applied mathematical expectation, but did not mention that notion.

**20.** It is unclear how in both cases the chances are determined *in the same way and with the same probability*.

**21.** When describing the measurement of angles, Poisson made many mistakes which proves that he had been ignorant of such work. In particular, he (like other French mathematicians dealing with the treatment of observations and even like Laplace) never cited Gauss. Yes, Gauss was guilty of failing to acknowledge Legendre's formal priority, but the ensuing consequences were still hardly excusable. I also indicate that Fourier had explained the notion of *real value* of a measured constant in the best possible way (Sheynin 2007), but Poisson never mentioned that.

**22.** Here is Hume's opinion of 1739 (Hald 1998, p. 127): *One wou'd appear ridiculous who wou'd say that 'tis only probable the sun will rise to-morrow …'* The problem about the probability of that event became classical (with no references to Hume).

**23.** Here, Poisson inserted the text of Biot's letter to him (without stating its date). Biot *clearly* described the principle of the indicated application.

**25.** In that section, Poisson did not indicate this circumstance which remains doubtful.

# Chapter 3. Calculus of Probabilities Depending on Very Large Numbers

**Misprints/Mistakes Unnoticed by the Author**

**1.** In § 70, p. 179 of the original text, the formula displayed on the last line. The exponent of the first factor should be − μ/2 rather than μ/2, see next page of the original text.

**2.** In § 70, p. 181 of the original text second integral should be over [0, ∞] rather than over [∞, *x*].

**3.** In § 72, p. 188 of the original text, Poisson denoted by *a′* and *b′* the number of balls left in the urn, but then, the number of those extracted. The ratio *p′/q′* was at first equated to *a/b*, but then, to *a′/b′*.

**4.** In § 74, p. 191 of the original text, the ratio of two integrals above formula (10). The denominator of the second term of the series on the right side should be (1 + α) instead of (1 + 2). Someone corrected this mistake by hand.

**5.** In § 76, p. 195 of the original text, the last displayed line includes a superfluous and strange equality *K′* = 1/2.

**6.** In § 79, p. 202 of the original text, Poisson's remark concerning the displayed expression of *n* ended by an apparently wrong reference to formula (6).

**7.** In § 86, p. 221 of the original text, on line 2 above § 87, the sign of the inequality should be ≥ rather than ≤.

**8.** In § 88, p. 225 of the original text, formula for $v_1$. In the second term on the right side *v*μ$_1$/μ should precede the square root.

**9.** In § 90, p. 233 of the original text, the displayed formula after the words *Il en résultera*. The numerator on the right side contains factor φ(*a* – *g, b* – *g*) instead of φ(*a* – *g, b* – *h*).

**10.** In § 90, p. 234 of the original text, at the end of that section Poisson replaces the arbitrary magnitudes *a* and *b* by (*a* + *n*) and (*b* + *m*) rather than by (*a* + *n*) and (*b* + *m*).

**11.** In § 91, p. 234 of the original text, on the first line of that section numbers *a, b, a* – *m, b* – *n* should be printed rather than numbers *a, b, a* – *m, a* – *n*.

All those misprints/mistakes are corrected in the translation, but the confusion mentioned in No. 3 is left intact.

**66.** When wishing to calculate the ratio of very large powers of two given numbers, it can always be easily done by applying logarithmic tables; if necessary, tables with more significant digits than usually needed. […] But the situation differs when the ratio of two products, each consisting of a very large number of unequal factors, should be calculated. […] We have to turn to approximate methods the first example of which is due to Stirling. […]

Such ratios of products of a very large number of factors and of sums of very large numbers of rations occur in most of the most important applications of the calculus of probability. This renders the rules described in the two preceding chapters, although provided in a completed form, barely useful or simply illusory when unaided by formulas proper to calculate numerical values with sufficient precision. We are going to discuss such formulas.

**67.** First of all, let us consider the product $n!$[1]. Here and throughout this chapter $e$ will represent the base of the Naperian logarithms. Integrating by parts, we will have

$$\int_0^\infty e^{-x} x^n dx = n!. \qquad (67.1)$$

The integrand disappears at $x = 0$ and $\infty$. Within these limits, it is never infinite and only has one maximal value […] $H$ at $x = h$. We have $h = n$, $H = e^{-h} h^n$ and it is possible to take $e^{-x} x^n = H\exp(-t^2)$, where $t$ is a new variable increasing from $-\infty$ to $\infty$. The values $t = -\infty, 0, \infty$ correspond to $x = 0, h, \infty$. […] If $x = h + x'$, then […]

$$x' = h't + h''t^2 + \ldots$$

and […]

$$1 + \frac{1}{2}\frac{d^2 \ln H}{dh^2} h'^2 = 0, \quad \frac{d^2 \ln H}{dh^2} h'' + \frac{1}{6}\frac{d^3 \ln H}{dh^3} h'^2 = 0,\ldots \quad (67.2)$$

[…]

$$n! = n^n e^{-n} \sqrt{2\pi n}\, \left(1 + \frac{1}{12n} + \frac{1}{288n^2} + \ldots\right). \qquad (67.3)$$

**68.** The first terms of the series in (67.3) converge the rapidly the larger is number $n$. Nevertheless, [Poisson corrected the next phrase] the law of that series is unknown and can belong to those, which become divergent when sufficiently continued[2]. But still, if reducing it to its convergent part, we may always apply formula (67.3) […] and it is not even necessary to have a very considerable $n$ for achieving a very good approximation. [An example with $n = 10$ and a derivation of the Wallis formula follow.]

It is to Laplace that analysis is indebted for the method which we had applied for reducing integrals to series convergent in their first terms and proper for calculating approximate values when the integrands are functions raised to very large powers. And we will [also] apply that method otherwise.

**69.** Let E and F be contrary events of some nature only one of which arrives at each trial. Denote their probabilities supposed constant by $p$ and $q$, and, by $U$, the probability that in $\mu$ trials E and F occur $m$ and $n$ times (§ 14)

$$U = C_\mu^m p^m q^n, \quad m + n = \quad p + q = \qquad (69.1)$$

If $\mu$, $m$ and $n$ are very large numbers, it will be necessary to apply the Stirling formula for calculating the numerical value of $U$. And if each of these three numbers is sufficiently large for allowing to reduce that formula to its first term, then […] and approximately

$$U = \left(\frac{\mu p}{m}\right)^m \left(\frac{\mu q}{n}\right)^n \sqrt{\frac{\mu}{2\pi mn}}. \qquad (69.2)$$

It is easy to establish that the most probable compound event, or such for which this value of $U$ is maximal, corresponds to the case in which $m/n$ as closely as possible approaches $p/q$. If, on the contrary, $m$ and $n$ are given, and $p$ and $q$ are variables whose sum is unity, and can increase by infinitely small steps from zero to unity, we will find by ordinary rules that the maximal value of $U$ corresponds to $p = m/\mu$ and $q = n/\mu$. However, owing to the large number of other compound events, although less probable than the above, its probability will be inconsiderable and lowers as the number $\mu$ of trials supposed very large increases further. For example, if $p = q = 1/2$, and $\mu$ is an even number, the most probable compound event will correspond to $m = n = \mu/2$. By formula (69.2) its probability will be

$$U = \sqrt{2/\pi\mu} \qquad (69.3)$$

and lowers, as it is seen, inversely proportional to $\sqrt{\mu}$. If $\mu = 100$, $U = 0.07979$, $1 - U = 0.92021$.

It is possible to bet a bit more than 92 against 8 on the equally probable contrary events E and F not arriving the same number of times in 100 trials. When retaining the last multiplier of the Stirling formula the expression (69.3) should be multiplied by $(1 - 1/4\mu)$. […]

**70.** Not only is the compound event for which the ratio $m/n$ most closely approaches $p/q$, always the most probable, but [in addition] when the number $\mu$ of trials is very large, the probabilities of other compound events just begin to lower rapidly when $m/n$ deviates from $p/q$ in either direction beyond certain limits whose extent is inversely proportional to $\sqrt{\mu}$. Let once more $p = q = 1/2$ and $g$ be a given positive or negative magnitude less than $\sqrt{\mu}$ without taking account of the signs. Then, if in formula (69.2)

$$m = \frac{\mu}{2}\left(1 + \frac{g}{\sqrt{\mu}}\right), \quad n = \frac{\mu}{2}\left(1 - \frac{g}{\sqrt{\mu}}\right), \qquad (70.1)$$

$$U = \left(1 - \frac{g^2}{\mu}\right)^{-\mu/2} \left(1 - \frac{g}{\sqrt{\mu}}\right)^{g\sqrt{\mu}/2} \left(1 + \frac{g}{\sqrt{\mu}}\right)^{-g\sqrt{\mu}/2} \sqrt{\frac{2}{\pi(\mu - g^2)}}.$$

If $g$ is a fraction or a very small [natural] number compared with $\sqrt{\mu}$, we will have by the binomial formula (§ 8) almost exactly

$\exp(g^2/2)$, $\exp(-g^2/2)$, $\exp(-g^2/2)$

for the three first factors of $U$. Substituting $\mu$ instead of $(\mu - g^2)$, we will obtain

$$U = \sqrt{2/\pi\mu} \exp(-g^2/2)$$

for the law of the lowering of probability *U* in a small vicinity in either direction of its maximal value. For example, if μ = 200 and $g = 1/\sqrt{2}$, the ratio of the probability, that in 200 trials events E and F having equal chances will take place 105 and 95 times, to the probability of their arriving 100 times each is $e^{-1/4}:1 \approx 3:4$.

Formula (69.2) supposes that each of the three numbers, μ, *m* and *n*, is very large. If that condition is fulfilled and *m/n* much deviates from *p/q*, that formula will provide a very small *U* with respect to its maximal value. However, it is proper to note that, when applying another method of approximation, a very small value of *U* will be derived if the difference (*m/n* − *p/q*) is a very small fraction and possibly it will not coincide with the value calculated by formula (69.2). The ratio of those two approximate values can much deviate from unity.

To show it, I note that by virtue of a formula in one of my memoirs (1823, *J. Ecole Polyt.*, No. 19, p. 490) on definite integrals,

$$\frac{2}{\pi}\int_0^{\pi/2} \cos^\mu x \cos[(m-n)x]dx = \frac{1}{2^\mu}C_\mu^m.$$

Whatever are the numbers *m* and *n* and their sum μ, we will have, according to formula (69.1), if *p* = *q* = 1/2, which is sufficient to consider,

$$U = \frac{2}{\pi}\int_0^{\pi/2} \cos^\mu x \cos[(m-n)x]dx.$$

And if μ is a very large number and if, when calculating approximately, it is treated like an infinite number, the factor $\cos^\mu x$ in the integrand will disappear when *x* is finite, and cos[(*m* − *n*)*x*] will always be finite. It follows that, without changing the value of that integral, it can be considered over [0, α] where α is positive and infinitely small. Then

$$\cos x = 1 - x^2/2, \quad \cos^\mu x = \exp(-\mu x^2/2),$$
$$U = \frac{2}{\pi}\int_0^\alpha \exp(-\mu x^2/2)\cos[(m-n)x]dx$$

Actually, however, the exponential factor disappears at all finite values of *x*, and it will also be possible to extend the limits of this new integral without changing its value beyond α, and, if wished, to *x* = ∞. Then, since by a known formula

$$\int_0^\infty \exp(-\mu x^2/2)\cos[(m-n)x]dx = \sqrt{\frac{\pi}{2\mu}}\exp\left(-\frac{(m-n)^2}{2\mu}\right);$$

$$U = \sqrt{\frac{2}{\pi\mu}} \exp(-\frac{(m-n)^2}{2\mu}).$$

Now, introduce as above [as follows from formula (70.1)]

$$m - n = g\sqrt{\mu}.$$

This value of $U$ will only coincide with the result of (69.2) when $g$ is a very small number compared with $\sqrt{\mu}$. At other values of $g$ the ratio of these two values of $U$ will much differ from unity and can even become a very large number. For example, if $g = \sqrt{\mu}/2$, and $m - n = \mu/2$, the preceding formula will provide

$$U = \sqrt{2/\pi\mu}\, e^{-\mu/8}.$$

Formula (69.2) will lead to

$$U = (1 - 1/4)^{-\mu/2}(1 - 1/2)^{\mu/4}(1 + 1/2)^{-\mu/4}\sqrt{8/3\pi\mu},$$

and, since the second factor here almost equals the third,

$$U = (9/8)^{-\mu/2}\sqrt{8/3\pi\mu}.$$

The two calculated values of $U$ correspond to each other in the sense that both are very small and thus prove that in a very large number $\mu$ of trials the probability $U$ of the two events, E and F, having equal chances and arriving $3\mu/4$ and $\mu/4$ times is extremely feeble. However, when dividing the last value of $U$ by the first one, we will have

$$\frac{2}{\sqrt{3}}\left(\frac{64}{81}\sqrt{e}\right)^{\mu/4},$$

a magnitude indefinitely increasing with $\mu$ and already exceeding 800 at $\mu = 100$.

**71.** Suppose, like above, that the chances of E and F are constant but unknown. It is only known that in $\mu = m + n$ trials they arrived $m$ and $n$ times. It is required to determine the probability $U'$ that in $\mu_1 = m_1 + n_1$ trials they will occur $m_1$ and $n_1$ times. We will have (§ 46) [...]

$$U' = C_{\mu_1}^{m_1}\frac{(m+m_1)!(n+n_1)!(\mu+1)!}{m!n!(\mu+\mu_1+1)!}.$$

[After transformations]

$$U' = C_{\mu_1}^{m_1} K(1+\frac{m_1}{m})^m(1+\frac{n_1}{n})^n(1+\frac{\mu_1}{\mu})^{-\mu} \times$$

$$(\frac{m+m_1}{\mu+\mu_1})^{m_1}(\frac{n+n_1}{\mu+\mu_1})^{n_1}, \qquad (71.1)$$

$$K = \frac{\mu+1}{\mu+\mu_1+1}\sqrt{\frac{(m+m_1)(n+n_1)(\mu+1)}{mn(\mu+\mu_1+1)}}. \qquad (71.2)$$

If $m_1$ and $n_1$ are very small numbers compared with $m$ and $n$, [...] and if substituting $K$ by unity, from which it very little differs, then, almost exactly,

$$U' = C_{\mu_1}^{m_1}(m/\mu)^{m_1}(n/\mu)^{n_1}.$$

It follows from formula (69.1) that that expression coincides with the probability that events E and F arrive $m_1$ and $n_1$ times in $(m_1 + n_1)$ trials when the chances $p$ and $q$, are given beforehand and are certainly $p = m/\mu$ and $q = n/\mu$. In particular, if $m_1 = 1$ and $n_1 = 0$, we have $C_{\mu_1}^{m_1} = 1$ and $U' = m/\mu$ which is the probability that E, after taking place $m$ times in a very large number $\mu$ of trials, will arrive once more at a new trial, and this is confirmed by the rule of § 49.

However, when the numbers $m_1$ and $n_1$ are comparable to $m$ and $n$, the probability $U'$ will not be the same as when the chances of E and F are given in advance and are certainly equal to $m/\mu$ and $n/\mu$. To illustrate this by an example, I denote by $h$ a natural number or a not too little fraction and assume that $m_1 = mh$, $n_1 = nh$, and $\mu_1 = \mu h$ so that $K$ will be almost equal to $1/\sqrt{1+h}$. Since $\mu = (m + n)$, formula (71.1) will be reduced to

$$U' = \frac{1}{\sqrt{1+h}}C_{\mu_1}^{m_1}(m/\mu)^{m_1}(n/\mu)^{n_1}.$$

Denoting by $U_1$ the new form of that formula, comparing it with formula (69.1), assuming that $p = m/\mu$ and $q = n/\mu$, and substituting $m$ and $n$ by $m_1$ and $n_1$, we conclude that

$$U' = \frac{1}{\sqrt{1+h}}U_1.$$

Therefore, $U'$ is less than $U_1$ in the ratio of $1/\sqrt{1+h}$, and is very low when $h$ is a very large number.

And so, there exists an essential difference between the probabilities $p$ and $q$ of the events E and F as given by a hypothesis and their probabilities $m/\mu$ and $n/\mu$ derived from the number of the arrivals of E and F in a very large number of trials: the probability that E and F will occur given numbers of times in a new series of trials is lower in the second case than in the first. That difference is occasioned by the probabilities of E and F derived from observations, however large are the numbers $m$ and $n$ on which they are based, being only possible whereas the probabilities given beforehand are certain[3].

If, for example, it is known that an urn A contains an equal number of white and black balls, there will be a probability almost equal to 0.07979, see § 69, that in a 100 successive drawings with replacement 50 balls of each colour were extracted. However, if the ratio of the white and black balls is not given, and it is only known that 50 balls of each colour were extracted, the probability that the same happens in 100 new trials will only be $0.07979/\sqrt{2} = 0.05658$ if $h = 1$ in the preceding equation.

**72.** As an example of the chances of events E and F varying during trials, I suppose that an urn A contains $c$ balls, $a$ of them white, and $b$, black, and that $\mu$ balls were successively extracted without replacement. I denote by $V$ the probability that $m$ white and $n$ black balls will be drawn in some order, $m + n = \mu$. By the formula of § 18, denoting by $a'$ and $b'$, $a' + b' = c'$, the number of balls remaining in A in the notation of § 71 will be

$$V = C_{c'}^{a'} \frac{a!b!\mu!}{m!n!c!}.$$

If $a, b, m,$ and $n$ are very large numbers, the Stirling formula will provide an approximate value of $V$:

$$V = C_{c'}^{a'} \left(\frac{a}{c}\right)^a \left(\frac{b}{c}\right)^b \left(\frac{m}{\mu}\right)^{-m} \left(\frac{n}{\mu}\right)^{-n} \sqrt{\frac{ab\mu}{mnc'}}.$$

If $m = a, n = b, \mu = c$ so that $C_{c'}^{a'} = 1$, this probability will exactly equal unity since all $a + b$ balls will be drawn from A. And if $m/n = a/b$, then $a/c = m/\mu$ and $b/c = n/\mu$, and if $a/c = p'$ and $b/c = q'$, then

$$V = C_{c'}^{a'} p'^{a'} q'^{b'} \sqrt{c/\mu}.$$

Comparing it with formula (69.1) and denoting by $V'$ the probability that two events, having constant chances equal to the chances $a/c$ and $b/c$ of extracting a white and a black ball at the beginning of the trials, arrive $a'$ and $b'$ times in $c'$ trials, we will have

$$V = V'\sqrt{c/\mu}. \qquad (72.1)$$

This proves that, if only $\mu$ is very large, $V/V' = \sqrt{c}/\sqrt{\mu}$ for any number $c'$ of balls remaining in A after the drawings.

We may remark that $a' = p'(c - \mu), b' = q'(c - \mu)$ so that the ratio of the balls of both colours remaining in A $a'/b' = p'/q' = a/b$. If for example $p' = q' = 1/2$, so that $a' = b' = c'/2$, then (§ 69) $V' = \sqrt{2/\pi c'}$, and, since $c' = c - \mu$, $V = \sqrt{2c/[\pi\mu(c-\mu)]}$. If $\mu = c/2$, then $V = \sqrt{4/\pi\mu} = V'\sqrt{2}$.

This means that if an urn A contains very large and equal numbers of white and black balls, and a half of them is extracted without replacement, the probability of an equal number of balls of both colours having been drawn exceeds its value in case of drawings with replacement in the ratio of √2:1.

**73.** I return to the case of constant chances $p$ и $q$ of the events E and F and consider the probability $P$ that in $\mu = m + n$ trials E arrives not less than $m$ times, and F, not more than $n$ times. This probability is the sum of the first $m$ terms of the expansion $(p + q)^\mu$ in increasing powers of $q$. However, it is difficult to transform this expansion into an integral to which, if $m$ and $n$ are very large numbers, the method of § 67 can be applied. We will therefore first of all look for another expression of P better suited for that purpose.

It is also possible to say that the appropriate compound event, call it G, means that F will not occur more than $n$ times in $\mu$ trials. This can happen in $(n + 1)$ cases.

**[1]** If E arrived in each of the first $m$ trials, since only $\mu - m = n$ are left and F will not occur more than $n$ times. The probability of this case is $p^m$.

**[2]** If E appeared $m$ times, and F, once, but not at the last of those $(m + 1)$ first trials. The stated restriction is necessary for the second case not to repeat the first one. In the following $(n - 1)$ trials F obviously can not arrive more than $(n - 1)$ times and therefore can not occur more than $n$ times at all the trials. The probability of E appearing $m$ times, and F, once, in a determined trial is $p^m q$. Since F can arrive in $m$ different trials, the probability of this case favourable for G is $mp^m q$.

**[3]** If E occurred $m$ times in $(m + 2)$ first trials, and F, twice, in determined trials, but not at the second one, which is necessary and sufficient for this third case not to repeat either of the first two cases. The probability of the stated conditions is $p^m q^2$, and the composite probability of this case is $m(m + 1)p^m q^2/2$.

When continuing in such a manner, we will come to the $(n + 1)$-st case in which E occurs $m$ times and F, $n$ times, but not at the last trial so that this case will not repeat any of the preceding cases. The probability of this last case is

$$\frac{m(m+1)...(m+n-1)}{n!} p^m q^n.$$

These $(n + 1)$ cases are distinct one from another and present all the different manners in which the event G can take place. Its composite probability is the sum of their respective probabilities (§ 10):

$$P = p^m [1 + mq + \frac{m(m+1)}{2!} q^2 + ... + \frac{m(m+1)...(m+n-1)}{n!} q^n]. \quad (73.1)$$

This expression[4] should coincide with the expansion of $(p + q)^\mu$. Its advantage is that it is easily transformed into definite integrals whose numerical values can be calculated by the method of § 67 the more exactly the larger are $m$ and $n$.

**74.** [Poisson calculates

$$\int Xdx, \quad X = \frac{x^n}{(1+x)^{\mu+1}}$$

integrating it $(n + 1)$ times by parts and derives the formula]

$$P = \frac{\int_\alpha^\infty Xdx}{\int_0^\infty Xdx} = \frac{1}{(1+\alpha)^m}\left[1 + m\frac{1}{1+\alpha} + \frac{m(m+1)}{2!}\frac{\alpha^2}{(1+\alpha)^2} + \ldots\right.$$

$$\left. + \frac{m(m+1)\ldots(m+n-1)}{n!}\frac{\alpha^n}{(1+\alpha)^n}\right]. \qquad (74.1)$$

Here, α is non-negative. And if $\alpha = q/p$ with $p + q = 1$, the right side of this last expression coincides with [the right side of] formula (73.1). At that value of α the left side of that equation will be $P$, and if $n = 0$ and $m = \mu$, $P$ will be the probability that E arrives at least μ times, that is, at each trial. Therefore, $P$ should equal $p^\mu$. And actually, if $n = 0$,

$$\int_\alpha^\infty Xdx = \frac{1}{\mu(1+\alpha)^\mu} = \frac{1}{\mu}p^\mu, \quad \int_0^\infty Xdx = \frac{1}{\mu}, \quad P = p^\mu. \qquad (74.2)$$

If $n = \mu - 1$ and $m = 1$, $P$ will be the probability that E arrives at least once, or that F will not occur at each trial. Therefore, $P = 1 - q^m$, which can also be verified. Indeed, let

$x = 1/y$, $dx = -dy/y^2$, $\alpha = 1/\beta$.

For $n = \mu - 1$

$$\int_\alpha^\infty Xdx = \int_0^\beta \frac{dy}{(1+y)^{\mu+1}} = \frac{1}{\mu}\left[1 - \frac{1}{(1+\beta)^\mu}\right], \quad \int_0^\infty Xdx = \int_0^\infty \frac{dy}{(1+y)^{\mu+1}} = \frac{1}{\mu}$$

and since $\beta = 1/\alpha = p/q$ and $1/(1 + \beta) = q$ formula (74.1) coincides with the preceding value of $P$.

**75.** We begin by applying the method of § 67 to the integral (74.2). Just like in that section, we denote the value of $x$, corresponding to the maximal value $H$ of $X$, by $h$. The equation $dX/dx = 0$, which serves to determine $h$, will be

$n(1 + h) - (\mu + 1)h = 0$

so that

$$h = \frac{n}{m+1}, \quad H = \frac{n^n(m+1)^{m+1}}{(\mu+1)^{\mu+1}}. \qquad (75.1)$$

After assuming in equations (67.2)

$$H = \frac{h^n}{(1+h)^{\mu+1}},$$

differentiating $H$ with respect to $h$ and substituting the previous value of that variable, we obtain

$$h' = \sqrt{\frac{2(\mu+1)n}{(m+1)^3}}, \quad h'' = \frac{2(\mu+1+n)}{3(m+1)^2}, \ldots$$

If $m$, $n$ and $\mu$ are very large numbers of the same order, it will be easy to see that those magnitudes $h'$, $h''$, … form a very rapidly decreasing converging series whose first term, $h'$, is of the same order of smallness as the fraction $1/\sqrt{\mu}$, the second term, $h''$, is of the order of $1/\mu$, the third term, of the order of $1/(\mu\sqrt{\mu})$ etc. And so, the given integral expanded into a series is

$$\int_0^\infty X dx = H\sqrt{\pi}\left(h' + \frac{1\cdot 3}{2}h''' + \frac{1\cdot 3\cdot 5}{4}h^v + \ldots\right). \qquad (75.2)$$

**76.** The expression of the integral in the numerator of the left side of formula (74.1) depends on whether $\alpha > h$ or $< h$, where, as before, $h$ is the value of $x$, corresponding to the maximal value of $X$. Actually, the variable $t$ in the transformation of § 67 should be positive at all values of $x > h$ and negative for $x < h$. Denoting by $\theta$ and $A$ the values of $t$ и $X$ at $x = \alpha$ [here follows the derivation of $A$ at $\alpha = q/p$; taking into account the previous value of $H$, it occurs that] $\theta = \pm k$,

$$k^2 = n\ln\frac{n}{q(\mu+1)} + (m+1)\ln\frac{m+1}{p(\mu+1)}. \qquad (76.1)$$

Supposing that $k > 0$, we should choose $\theta = k$ if $q/p > h$ and $\theta = -k$ if $q/p < h$. Therefore, in accord with the transformation of § 67, we will have in the first case

$$\int_\alpha^\infty X dx = H\int_k^\infty \exp(-t^2)\frac{dx'}{dt} dt.$$

In the second case that integral is

$$H\int_{-k}^\infty \exp(-t^2)\frac{dx'}{dt} dt = H\int_{-\infty}^\infty \exp(-t^2)\frac{dx'}{dt} dt -$$
$$H\int_{-\infty}^{-k} \exp(-t^2)\frac{dx'}{dt} dt \qquad (76.2)$$

and, as in § 67,

$$dx'/dt = h' + 2h''t + 3h'''t^2 + \ldots$$

[…] Denote[5]

$$\int_k^\infty \exp(-t^2)t^{2i}dt = K_i, \quad \int_k^\infty \exp(-t^2)t^{2i+1}dt = K'_i.$$

The first term on the right side of formula (76.2) is equal to the integral (75.2), so that, if $q/p > h$, see (75.1),

$$\int_\alpha^\infty Xdx = H(h'K_0 + 3h'''K_1 + 5h^vK_2 + ...) +$$

$$H(2h''K'_0 + 4h^{IV}K'_1 + 6h^{VI}K'_2 + ...) = I + II,$$

$$\int_\alpha^\infty Xdx = \int_0^\infty Xdx - I + II \text{ if } q/p < h. \qquad (76.3a, b)$$

Each series in those formulas in general has the same degree of convergence [converges as rapidly] as the series (75.2). If $k \neq 0$, the values of the integrals $K_i$ can only be obtained by approximation, whereas the integrals $K'_i$ are always expressed in a finite form

$$K'_i = (1/2) \exp(-k^2)[k^{2i} + ik^{2i-2} + i(i-1)k^{2i-4} + ... + i!k^2 + i!].$$

If $\alpha = h$, formulas (76.3a, b) should coincide. At the same time,

$$\frac{q}{p} = \frac{n}{m+1\mu}, \quad q = \frac{n}{1+\mu}, \quad p = \frac{m+1}{1+}$$

and $k$, determined by the equation (76.1), disappears:

$$K_0 = \frac{\sqrt{\pi}}{2}, \quad K_i = 1 \cdot 3 \cdot 5 \cdot ... \cdot (2i-1)\frac{\sqrt{\pi}}{2^{i+1}}, \quad K'_i = \frac{i!}{2}$$

According to equation (75.1) the formulas (76.3a, b) are reduced to the same expression

$$\int_\alpha^\infty Xdx = \frac{H\sqrt{\pi}}{2}(h' + \frac{1 \cdot 3}{2}h''' + \frac{1 \cdot 3 \cdot 5}{4}h^v + ...) +$$

$$H(h'' + 1 \cdot 2h^{IV} + 1 \cdot 2 \cdot 3h^{VI} + ...).$$

**77.** We suppose that the numbers $m$, $n$ and $\mu$ are sufficiently large for neglecting, in various formulas, $h'''$, $h^{IV}$ etc. Taking the above values of $h'$ and $h''$, we will have

$$\frac{h''}{h'} = \frac{(\mu+1+n)\sqrt{2}}{3\sqrt{n(m+1)(\mu+1)}}.$$

By means of equation (74.1) and formulas (75.2) and (76.3) we obtain, respectively, for $q/p > h$ and $q/p < h$

$$P = \frac{1}{\sqrt{\pi}} \int_k^\infty \exp(-t^2)dt + \frac{(\mu+n)\sqrt{2}}{3\sqrt{\pi\mu mn}} \exp(-k^2), \quad (77.1a)$$

$$P = 1 - \frac{1}{\sqrt{\pi}} \int_k^\infty \exp(-t^2)dt + \frac{(\mu+n)\sqrt{2}}{3\sqrt{\pi\mu mn}} \exp(-k^2). \quad (77.1b)$$

Here, as given by equation (76.1), $k > 0$ и $k^2$ is determined by formula (76.1). For simplifying the account, $\mu$ and $m$ are substituted in the last terms of these formulas instead of $\mu + 1$ and $m + 1$. They will provide the required probability $P$ with a sufficient approximation.

If $\mu$ is even, $m = n = \mu/2$ and $q > p$, then $h = \mu/(\mu + 2)$, $q/p > h$ and we should apply formula (77.1a). It and the equation (76.1) lead to

$$P = \frac{1}{\sqrt{\pi}} \int_k^\infty \exp(-t^2)dt + \sqrt{\frac{2}{\pi\mu}} \exp(-k^2), \quad (77.2)$$

$$k^2 = \frac{\mu}{2} \ln \frac{\mu}{2q(\mu+1)} + \frac{\mu+2}{2} \ln \frac{\mu+2}{2p(\mu+1)}.$$

$P$ expresses the probability that in a very large even number of trials the more probable event F will not nevertheless arrive more often than the contrary event E. Denote by $U$ the probability that both these events occur the same number of times, then $P − U$ will be the probability that F arrives less often than E. If $p = q = 1/2$, $P − U$ will evidently be also the probability that E arrives less often than F; $2(P − U)$ added to the probability $U$ will therefore provide certainty: $2P − U = 1$, so that

$$U = \frac{2}{\sqrt{\pi}} \int_k^\infty \exp(-t^2)dt - 1 + 2\sqrt{\frac{2}{\pi\mu}} \exp(-k^2).$$

This is easy to verify.
   Owing to [the transformations performed]

$$k^2 = \frac{1}{4\mu} + \frac{1}{4(\mu+2)} + ...$$

When only preserving the terms of the same order of smallness as $1/\sqrt{\mu}$, we will have $k = 1/\sqrt{2\mu}$, $\exp(-k^2) = 1$,

$$\int_k^\infty \exp(-t^2)dt = \int_0^\infty \exp(-t^2)dt - \int_0^k \exp(-t^2)dt = \frac{\sqrt{\pi}}{2} - \frac{1}{\sqrt{2\mu}}$$

so that $U = \sqrt{2/\pi\mu}$, which actually coincides with the value deduced from formula (69.3), in turn derived from (69.2) at $m = n$ and $p = q$.

If μ is odd, let $m = (\mu - 1)/2$, and, as previously, $q > p$, then $q/p > h$. The formula (77.2) persists and from formula (76.1) it follows that

$$k^2 = \frac{\mu-1}{2} \ln \frac{\mu-1}{2q(\mu+1)} + \frac{\mu+3}{2} \ln \frac{\mu+3}{2p(\mu+1)}.$$

*P* is the probability that in a very large number μ of trials the more probable event will nevertheless occur less frequently; when μ is odd, an equal number of arrivals of E and F is impossible. If $p = q = 1/2$ this probability is 1/2, which we will verify. [Transforming the expression of $k^2$, Poisson gets]

$$k^2 = 1/(\mu-1) + 1/(\mu+3) + ...$$

If neglecting, as above, terms of the order of smallness $1/\mu$, then

$$k = \frac{\sqrt{2\pi}}{\mu}, \quad \exp(-k^2) = 1, \quad \int_k^\infty \exp(-t^2)dt = \frac{\sqrt{\pi}}{2} - \frac{\sqrt{\pi}}{\mu}$$

and the previous value of *P* becomes equal to 1/2.

**78.** Let now *n* differ from $(\mu + 1)q$ by a positive or negative magnitude ρ, very small comparing to this product. Since $p + q = 1$ and $m + n = \mu$,

$$n = (\mu + 1)q - \rho, \quad m + 1 = (\mu + 1)p + \rho.$$

The corresponding value of *h* is

$$h = \frac{(\mu+1)q - \rho}{(\mu+1)p + \rho} < \frac{q}{p} \text{ if } \rho > 0.$$

Expanding the right side of equation (76.1) in powers of ρ, we can obtain

$$k^2 = \frac{\rho^2}{2(\mu+1)pq}\left[1 + \frac{(p-q)\rho}{3(\mu+1)pq} + ...\right].$$

Introduce $r > 0$, then

$$\rho = r\sqrt{2(\mu+1)pq}, \text{ then } k = r\left[1 + \frac{(p-q)r}{3\sqrt{2(\mu+1)pq}} + ...\right],$$

and exclude very small *p* and *q*, then the obtained series in parentheses will converge rapidly since it is ordered in powers of $r/\sqrt{\mu+1}$ or $\rho/(\mu+1)$. Preserving only the first two terms, we can for the sake of brevity write down *k* in the form

$$k = r + \delta, \ \delta = \frac{(p-q)r^2}{3\sqrt{2(\mu+1)pq}}. \qquad (78.1)$$

At the same time,

$$n = (\mu+1)q - r\sqrt{2(\mu+1)pq}.$$

However, in the second term on the right side of formula (77.1a) it is sufficient to assume $k = r$ and substitute $m$ and $n$ by $p\mu$ and $q\mu$, for deriving

$$P = \frac{1}{\sqrt{\pi}} \int_{r+\delta}^{\infty} \exp(-t^2)dt + \frac{(1+q)\sqrt{2}}{3\sqrt{\pi\mu pq}} \exp(-r^2). \qquad (78.2a)$$

Let [now] $\rho < 0$, then $h > q/p$. Denote a positive magnitude by $r'$ and assume

$$\rho = - r'\sqrt{2(\mu+1)pq},$$

then

$$n = (\mu+1)q + r'\sqrt{2(\mu+1)pq}.$$

The value of $k$, derived from equation (76.1), should nevertheless invariably remain positive:

$$k = r' - \delta', \ \delta' = \frac{(p-q)r'^2}{3\sqrt{2(\mu+1)pq}}.$$

The formula (77.1b), which we should apply, will take the form

$$P = 1 - \frac{1}{\sqrt{\pi}} \int_{r'-\delta'}^{\infty} \exp(-t^2)dt + \frac{(1+q)\sqrt{2}}{3\sqrt{\pi\mu pq}} \exp(-r'^2). \qquad (78.2b)$$

If subtracting from it the previous value of $P$, the difference

$$R = 1 - \frac{1}{\sqrt{\pi}} \int_{r+\delta}^{\infty} \exp(-t^2)dt - \frac{1}{\sqrt{\pi}} \int_{r'-\delta'}^{\infty} \exp(-t^2)dt +$$

$$\frac{(1+q)\sqrt{2}}{3\sqrt{\pi\mu\ pq}}[\exp(-r'^2) - \exp(-r^2)] \qquad (78.3)$$

will be the probability that in a very large number $\mu$ of trials the number of arrivals of F will not exceed the second value of $n$ but will exceed its first value at least by unity [Poisson later added:] if that value is expressed by a natural number and not less than by unity otherwise.

**79.** Let for the sake of simplicity $N$ be the largest natural number contained in $\mu q$ and $f = \mu q - N$. Denote by $u$ such a magnitude, that $u\sqrt{2(\mu+1)pq}$ will be a natural number, very small compared with $N$ and let also

$$q + f - r\sqrt{2(\mu+1)pq} = -u\sqrt{2(\mu+1)pq} - 1,$$
$$q + f + r'\sqrt{2(\mu+1)pq} = u\sqrt{2(\mu+1)pq}.$$

The limits of the values of $n$, corresponding to the probability $R$, will become

$$n = N - u\sqrt{2(\mu+1)pq} - 1, \; n = N + u\sqrt{2(\mu+1)pq}$$

and formula (78.3) will therefore express the probability that $n$ not less than by unity exceeds the first one but will not be larger than the second. In other words, it will express the probability, that this number is contained within $N \pm u\sqrt{2\mu pq}$, equally remote from $N$ ($\mu$ was, however, substituted instead of $\mu + 1$), or that it coincides with either of them.

According to the equations above and expressions of $\delta$ and $\delta'$

$$r + \delta = u + \varepsilon + \frac{1}{\sqrt{2(\mu+1)pq}}, \; r' - \delta' = u - \varepsilon,$$

where $\varepsilon$ is a magnitude of the order of smallness $1/\sqrt{\mu}$. And if $v$ is a magnitude of the same order whose square we neglect, then

$$\int_{u+v}^{\infty} \exp(-t^2)dt = \int_{u}^{\infty} \exp(-t^2)dt - v\exp(-u^2).$$

Applying this equation to both integrals in formula (78.3) and equating $r' = r$ in the terms beyond integrals since they already included $\sqrt{\mu}$ in the denominator, we obtain

$$R = 1 - \frac{2}{\sqrt{\pi}} \int_{u}^{\infty} \exp(-t^2)dt + \frac{1}{\sqrt{2\pi\mu pq}} \exp(-u^2), \tag{79.1}$$

where, in the last term, $\mu$ is substituted instead of $\mu + 1$.

For the interval of the values of $n$ whose probability is $R$ not to include its inferior limit, the least of the two previous values of $n$ should be increased by unity. In this case the last term of the sum $(r + \delta)$, as well as the last term of the formula (79.1), will also disappear. And, for that interval not to include its superior limit, the largest value of $n$ should be decreased by unity. This will decrease the difference $(r' - \delta')$, and the last term of formula (79.1) will disappear once more. Finally, the sign of this term should be changed so that the interval of the values of *n* will not include any of its limits. Therefore, this term ought to be the probability of the exact equality

$$n = N + u\sqrt{2\mu pq},$$

where *u* is such a positive or negative magnitude, that the second term becomes very small as compared with the first one.

When neglecting the terms of the order of $1/\mu$,

$$\frac{n}{\mu} = q + u\sqrt{\frac{2pq}{\mu}}, \quad \frac{m}{\mu} = p - u\sqrt{\frac{2pq}{\mu}}$$

and [after transformations] formula (69.2) becomes

$$U = \frac{\exp(-u^2)}{\sqrt{2\pi\mu pq}},$$

which was required to verify.

The value of *P* in formula (78.2a) is the probability that *n* does not exceed the limit $\mu q - r\sqrt{2\mu pq}$, where $\mu$ is substituted instead of ($\mu$ + 1). Therefore, if, in that value of *U*, *u* = *r*, the difference $Q = P - U$ will be the probability that *n* does not reach that very limit. Just the same, if in that value of *U*, *u* = *r'* the difference $Q' = P - U$, where *P* is taken from the formula (78.2b), will be the probability of *n* being less than the limit $\mu q + r'\sqrt{2\mu pq}$. Therefore

$$Q = \frac{1}{\sqrt{\pi}} \int_{r+\delta}^{\infty} \exp(-t^2)dt + \frac{q-p}{3\sqrt{2\pi\mu pq}} \exp(-r^2), \quad (79.2a)$$

$$Q' = 1 - \frac{1}{\sqrt{\pi}} \int_{r'-\delta'}^{\infty} \exp(-t^2)dt + \frac{q-p}{3\sqrt{2\pi\mu pq}} \exp(-r'^2). \quad (79.2b)$$

Recall that *r* и *r'* are here positive magnitudes, very small as compared with $\sqrt{\mu}$, and that the limits of *n*, corresponding to *Q* and *Q'*, little deviate, in either direction respectively, from $\mu q$. At the same time, the values of $\delta$ and $\delta'$, included in those magnitudes, are very small as compared with *r* and *r'*. And, when replacing ($\mu$ + 1) by $\mu$, then, [just like in § 78]

$$\delta = \frac{(p-q)r^2}{3\sqrt{2\mu pq}}, \quad \delta' = \frac{(p-q)r'^2}{3\sqrt{2\mu pq}}.$$

**80.** When dividing the limits of *n*, corresponding to formula (79.1), by $\mu$, and noting the meaning of *U*, we will obtain the limits of the ratio $n/\mu$ corresponding to the probability *R*

$$q - \frac{f}{\mu} \mp \sqrt{\frac{2pq}{\mu}}.$$

If neglecting the fraction $f/\mu$, the magnitude $R$, determined by formula (79.1), will become the probability of the difference $(n/\mu - q)$ being contained in the limits $\pm u\sqrt{2pq/\mu}$. With changed signs and the same probability these will be the limits of the difference $(m/\mu - p)$ since the sum of those differences, $(m + n)/\mu - p - q = 0$.

It is always possible to choose such a large value of $u$ that that probability $R$ will arbitrarily little differ from certainty, and even if $u$ is not large. For example, it is sufficient to have $u = 4$ or 5, for the integral of $\exp(-t^2)$ over $[u, \infty]$ and for $(1 - R)$ as well, to become barely appreciable. If $u$ takes such a value and remains invariable, the limits of the difference $(m/\mu - p)$ will narrow as $\mu$, already being very large, will increase further. The ratio $m/\mu$ of the number of the arrivals of E to the total number of trials will ever less deviate from the chance $p$ of that event. It is always possible to increase the number $\mu$ so that the probability $R$ of the difference $(m/\mu - p)$ will become arbitrarily low. Conversely, if continuously increasing $\mu$ and assuming a constant and given magnitude $l$ for each of the indicated limits, thus ensuring that $u$ will increase in the same proportion as $\sqrt{\mu}$, the value of $R$ will ever closer approach unity. It is always possible to increase $\mu$, so that the probability $R$ of $(m/\mu - p)$ being contained within limits $\pm l$ will arbitrarily little differ from certainty. This, indeed, is the Jakob Bernoulli theorem (§ 49).

**81.** We (§ 78) have ignored the case of very small chances $p$ or $q$. Suppose now that $q$ is a very small fraction, so that the probability of the event F is very low. If the number $\mu$ of the trials is very large, the ratio $n/\mu$ of the number of the arrivals of F to $\mu$ will also be a very small fraction. Substituting $\mu - n$ instead of $m$ in formula (73.1), we will get $q\mu = w$, $q = w/\mu$. Neglecting the fraction $n/\mu$, the right side of this formula will be

$$p^m(1 + w + \frac{w^2}{2!} + \frac{w^3}{3!} + ... + \frac{w^n}{n!}).$$

At the same time,

$$p = 1 - \frac{w}{\mu}, \quad p^m = \left(1 - \frac{w}{\mu}\right)^\mu \left(1 - \frac{w}{\mu}\right)^{-n}.$$

We may replace the first factor by the exponential function $e^{-w}$, and assume that the second one is unity. In accord with equation (73.1) we will obtain almost exactly

$$P = (1 + w + \frac{w^2}{2!} + \frac{w^3}{3!} + ... + \frac{w^n}{n!})e^{-w}$$

for the probability that an event, whose chance at each trial is only a very small fraction $w/\mu$, will not appear more than $n$ times in a very large number $\mu$ of trials.

If $n = 0$, this value of $P$ becomes $e^{-w}$. This, therefore, is the probability that the event we consider will not arrive at all during $\mu$ trials. The probability that it will occur at least once is therefore $(1 - e^{-w})$, as noted in § 8. If, however, $n$ is not very small, the value of $P$ will very little deviate from unity. This is seen when $P$ is written in the form

$$P = 1 - \frac{w^{n+1}e^{-w}}{(n+1)!}[1 + \frac{w}{n+2} + \frac{w^2}{(n+2)(n+3)} + ...].$$

Let $w = 1$ and $n = 10$, then the difference $(1 - P)$ will be almost exactly equal to $1/10^8$. In other words, it is almost certain that an event, having a very slim chance $1/\mu$ of arriving at each of the $\mu$ trials, will not appear more than 10 times.

**82.** The integral in formula (79.1) can, in general, be calculated in quadratures. At the end of his book Kramp (1799) inserted a table of the values of that integral for $u = 0$ to $u = 3$, according to which

$$\int_3^\infty \exp(-t^2)dt = 0.00001957729...$$

Integrating by parts at $u > 3$ it is possible to derive [Poisson derived the formula

$$\int_0^u \exp(-t^2)dt = u - \frac{u^3}{1 \cdot 3} + \frac{u^5}{1 \cdot 2 \cdot 5} - \frac{u^7}{1 \cdot 2 \cdot 3 \cdot 7} + ...$$

and noted that the series on the right side rapidly converges at $u < 1$]. If the value of $u$ answering to $R = 1/2$ is desirable, we can equate that value of $R$ (79.1) to $1/2$ and apply that last series:

$$\frac{\sqrt{\pi}}{4} - \frac{\exp(-u^2)}{2\sqrt{2\mu pq}}.$$

Let $u = a$ be the root of the equation thus obtained. Neglecting the second term on the right side, we can put down that, to the order of smallness of $1/\mu$,

$$u = a - \frac{1}{2\sqrt{2\mu pq}}.$$

After a few attempts I found that $a \approx 0.4765$. This means that there exists the same probability that the difference $(m/\mu - p)$ will either be contained within the limits

$$\pm [0.4765\sqrt{\frac{2pq}{\mu}} - \frac{1}{2\mu}],$$

or be beyond them. For some value of $u$ there exists probability $R$, that the difference between $(m/\mu - p)$ and $(n/\mu - q)$ will not exceed $\pm 2u\sqrt{2pq}/\mu$. Therefore, if $p = q = 1/2$, the magnitude $(m - n)/\mu$ will with probability 1/2 be within limits

$$\pm [\frac{0.6739}{\sqrt{\mu}} - \frac{1}{\mu}].$$

Thus, if the chances of E and F are equal to each other, the difference of the numbers of their arrivals will with the same probability be either more or less than $0.6739\sqrt{\mu} - 1$.

Let A and B play very many fair games, for example $10^6$. Then we can bet even money that one of them will win 674 games more than the other. That difference, which can equally favour either of them, will indeed be the share of the hazard. But if, in each game, the chance $p$ of gambler A exceeds the chance $q$ of gambler B, there will be probability $R$, invariably heightening with $\mu$, that A will win $\mu(p - q) \pm 2u\sqrt{2\mu pq}$ more games than B. And since the first term, appearing because of the differing ability of the gamblers, increases as the number of games, whereas the second one only increases as the square root of that number, in a long run the more able gambler will always be the winner however small is the difference $(p - q)$.

**83.** Supposing that the chances $p$ and $q$ were known, we determined the likely good approximations of the ratios $m/\mu$ and $n/\mu$ provided that $\mu$ was very large. Conversely, if the chances are not known in advance, but those ratios are determined, our formulas established the likely values of $p$ and $q$ with a good approximation. Namely, there exists probability $R$, calculable by formula (79.1), that the chance $p$ of event E is contained within limits $m/\mu \pm u\sqrt{2pq/\mu}$. If $R$ very little differs from unity, the fraction $p$ will almost exactly be equal to $m/\mu$, and $q$, to $n/\mu$. Substituting therefore $m/\mu$ and $n/\mu$ instead of $p$ and $q$ in the double-valued terms of their limits and in the last term of the formula (79.1), which was already divided by $\sqrt{\mu}$, we will get the probability

$$R = 1 - \frac{2\mu}{\sqrt{\pi}} \int_u^\infty \exp(-t^2)dt + \sqrt{\frac{}{2\pi mn}} \exp(-u^2) \qquad (83.1)$$

of $p$ being contained within the limits

$$\frac{m}{\mu} \pm \frac{u}{\mu}\sqrt{\frac{2mn}{\mu}}.$$

If $m$, $n$ and $\mu$ are very large numbers, we can, in general, apply the approximate values $m/\mu$ and $n/\mu$ of the chances $p$ and $q$ for calculating the probability of a future event consisting of E and F; for example, of their arriving $m'$ and $n'$ times in $\mu' = m' + n'$ new trials, if $\mu'$ is very small as compared with $\mu$. We can apply formula (79.1) even if $\mu'$ is a very large number. Substitute $\mu'$, $m/\mu$ and $n/\mu$ instead of $\mu$, $p$, and $q$ in the formula itself and in the corresponding limits; then

$$R = 1 - \frac{2\mu}{\sqrt{\pi}} \int_u^\infty \exp(-t^2)dt + \frac{1}{\sqrt{2\pi\mu'mn}} \exp(-u^2) \quad (83.2)$$

will express the probability that $n'$ is contained within the limits

$$\frac{\mu'n}{\mu} \mp \frac{u}{\mu}\sqrt{2\mu'mn}. \quad (83.3)$$

Here, $\mu'n/\mu$ is introduced instead of the maximal natural number contained in that ratio. No matter how close were $m/\mu$ and $n/\mu$ to $p$ and $q$, they are only probable rather than certain, and they should not be applied, as noted in § 71, when $\mu'$ is comparable with $\mu$. This is why we will consider in another manner the problem of establishing $p$ and $q$ by observed events for applying them when studying probabilities of future events.

**84.** Suppose, like previously, that the arrival of events E and F was observed $m$ and $n$ times in a very large number $\mu = m + n$ of trials during which the chances of those events, $p$ and $q$, did not change. According to the above, there exists a very high probability that these unknown chances very little differ from the ratios $m/\mu$ and $n/\mu$, which can therefore be assumed as their approximate values. These chances can take an infinite number of values which increase by infinitely small increments and the probability of an exact value of $p$ and the corresponding value of $q$ is infinitely low. We should determine at least those values of $p$ and $q$ which little deviate from $m/\mu$ and $n/\mu$.

The magnitude $Q$, determined by formula (79.2a), is the probability of $n < \mu q - r\sqrt{2\mu pq}$ and at the same time of the unknown chance $q$ of event F which arrived $n$ times in $\mu$ events, to satisfy the condition $q > n/\mu + r\sqrt{2pq/\mu} = n/\mu + (r/\mu)\sqrt{2mn/\mu}$. In the right side $p$ and $q$ were replaced by their approximate values $m/\mu$ and $n/\mu$. If $(r - dr)$ is substituted instead of $r$ and only infinitely small magnitudes of the first order are left, then $[Q - (dQ/dr)dr]$ will also become the probability that

$$q > \frac{n}{\mu} + \frac{r}{\mu}\sqrt{\frac{2mn}{\mu}} - \sqrt{\frac{2mn}{\mu}}\frac{dr}{\mu}.$$

Therefore, $- dQ/dr$ will express the infinitely low probability that exactly

$$q = n/\mu + (r/\mu)\sqrt{2mn/\mu}$$

at all positive values of *r*, very small with respect to √μ, as it was indeed supposed in the formula for *Q*. Similarly, formula (79.2b) will express the probability *Q*′ that $q > n/\mu - (r'/\mu)\sqrt{2mn/\mu}$. Substituting (*r*′ + *dr*′) instead of *r*′, we will obtain [*Q*′ + (*dQ*′/*dr*′)*dr*′] for the probability that

$$q > \frac{n}{\mu} - \frac{r'}{\mu}\sqrt{\frac{2mn}{\mu}} - \sqrt{\frac{2mn}{\mu}}\frac{dr'}{\mu}.$$

Therefore, (*dQ*′/*dr*′)*dr*′ will be the probability that *q* exceeds the second, but not the first limit, or that exactly

$$q = n/\mu - (r'/\mu)\sqrt{2mn/\mu}.$$

Here, *r*′ is also a positive magnitude very small as compared with √μ.

[Poisson transforms the expressions of *dQ*/*dr* and *dQ*′/*dr*′ to within small magnitudes of the order of 1/√μ and notes that they turn into each other if *r* and − *r*′ are interchanged. He denotes] by *v* a positive or negative variable, very small as compared with √μ, [and introduces]

$$V = \frac{1}{\sqrt{\pi}}\exp(-v^2) - \frac{2(m-n)v^2}{3\sqrt{2\pi\mu mn}}\exp(-v^2). \qquad (84.1)$$

*Vdv* is the probability that

$$q = n/\mu + (v/\mu)\sqrt{2mn/\mu}.$$

At the same time that infinitely low probability concerns

$$p = m/\mu - (v/\mu)\sqrt{2mn/\mu}. \qquad (84.2)$$

The magnitude *V* is seen to decrease very rapidly with the increase of *v*. And even before the order of that variable *v* becomes √μ, the order of *V* will be extremely small owing to the factor exp (− *v*²). Express the values of *p* and *q* essentially deviating from *m*/μ and *n*/μ by means of *v*, and represent their probabilities by *V*′*dv*, where *V*′ is a function of *v*, differing from, and much smaller than *V*. It is quite insensible so that we are not required to determine it.

And so, let E′ be a future event comprised of E and F, and suppose that Π is the probability that it will occur at certain values of the chances of E and F and is therefore a given function of *p* and *q*. Then, let Π′ be the veritable probability of E′ when taking account of the values of *p* and *q* substituted in Π. Multiplying Π by this infinitely low probability of *p* and *q* and integrating the product from *p* = 0 and *q* = 1 to *p* = 1 and *q* = 0, we will obtain Π′. However, bearing in mind the considerations above, we can neglect the part of that integral answering to essential deviations of *p* and *q* from *m*/μ and *n*/μ. Therefore, inserting the previous values of *p* and *q* in Π, we get

$$\Pi' = \int \Pi V dv, \qquad (84.3)$$

when integrating over positive and negative values of *v*, very small as compared with √μ. This result corresponds with that, obtained more directly (Poisson 1830, § 2).

**85.** For providing the first application of the formulas (84.1) and (84.3), I will suppose that Π′ is the probability that in a very large number μ′ = *m*′ + *n*′ of new trials events E and F occurred *m*′ and *n*′ times with the ratio *m*′/*n*′ being very close to the ratio *m*/*n* of the number of the arrivals of those events in previous trials. In other words, I suppose that

*m*′ = *mh* − α√μ′, *n*′ = *nh* + α√μ′, μ′ = μ*h*,

where *h* and α are given and α is positive or negative but very small with respect to √μ′.

According to formula (69.2), for $U' = \sqrt{\mu'/2\pi m'n'}$,

$$\Pi = U'\left(\frac{\mu' p}{m'}\right)^{m'}\left(\frac{\mu' q}{n'}\right)^{n'}.$$

And now we can note that *U*′ is the probability of event E′, when assuming that *m*/μ and *n*/μ are the exact and given chances *p* and *q* of E and F and α = 0. In addition, since

$$\frac{m'}{\mu'} = \frac{m}{\mu} - \frac{\alpha}{\sqrt{\mu'}},\ \frac{n'}{\mu'} = \frac{n}{\mu} + \frac{\alpha}{\sqrt{\mu'}},$$

and denoting

$$\frac{v}{\mu}\sqrt{\frac{2mn}{\mu}} - \frac{\alpha}{\sqrt{\mu'}} = v_1,$$

we can write down the values of *p* and *q* from § 84 as

*p* = (*m*′/μ′) − *v*₁, *q* = (*n*′/μ′) + *v*₁.

Inserting these values in Π, we get

$$\Pi = U'\left(1 - \frac{\mu' v_1}{m'}\right)^{m'}\left(1 + \frac{\mu' v_1}{n'}\right)^{n'}.$$

The fractions here are of the order of 1/√μ or 1/√μ′, so that […]

$$\Pi = U'\exp\left(-\frac{\mu'^3 v_1^2}{2m'n'}\right). \qquad (85.1)$$

For the same reason we can neglect the second term in formula (84.1), and therefore [after transformations] formula (84.3) will become

$$\Pi' = \frac{1}{\sqrt{\pi}} \mu U' \int \exp(-v^2 - \frac{\mu'^3 v_1^2}{2m'n'}) dv.$$

The integral should only extend over values of $v$, very small as compared with $\sqrt{\mu}$. Owing to the exponential factor, the coefficient of $dv$ becomes quite insensible at values of $v$ comparable with $\sqrt{\mu}$ and we can, and will extend the integral over such values of $v$ and assume $-\infty$ and $\infty$ as its limits. Then, substituting $mh$ and $nh$ instead of $m'$ and $n'$ in $v_1$, we get

$$v^2 + \frac{\mu'^3 v_1^2}{2m'n'} = v^2(1+h) - \frac{2v\alpha\mu'\sqrt{h}}{\sqrt{2m'n'}} + \frac{\alpha^2\mu'^2}{2m'n'}$$

and introduce

$$v\sqrt{1+h} - \frac{\alpha\mu'\sqrt{h}}{\sqrt{2m'n'(1+h)}} = x, \quad dv = \frac{dx}{\sqrt{1+h}}.$$

With that new variable $x$, the limits of the integral will still be infinite, and the required probability becomes

$$\Pi' = \frac{1}{\sqrt{1+h}} \alpha \mu U' \exp[-\frac{\alpha^2 \mu'^2}{2m'n'(1+h)}] \qquad (85.2)$$

and is obviously simplified if $\alpha = 0$. When recalling the value of $U'$, we say that this coincides with the result obtained otherwise in § 71.

**86.** For the second example of applying formulas (84.1) and (84.3), we assume that $\Pi'$ is the probability that the difference $(n'/\mu' - n/\mu)$ does not exceed $\alpha/\sqrt{\mu'}$, which it should have reached in the preceding example of § 85.

The magnitude $\Pi$ is a function of the chances $p$ and $q$ of events E and F, and the probability that in $\mu'$ future trials F will not arrive more than $n'$ times, $n' = (\mu'/\mu)$ ± $\sqrt{\mu'}$, with E occurring $m'$, or not less than $[(\mu'/\mu)$ ± $\sqrt{\mu'}]$ times, cf. formulas (85.1). The value of that probability will be represented by one of the formulas (77.1) after substituting $\mu'$, $m'$, $n'$ instead of $\mu$, $m, n$.

If invariably calculating to within magnitudes of the order of $1/\sqrt{\mu}$ or $1/\sqrt{\mu'}$, the mentioned extreme values of $m'$ and $n'$ will satisfy the condition

$$\frac{n'}{m'+1} = \frac{n}{m}[1 + \frac{\alpha\mu^2}{mn\sqrt{\mu'}}].$$

According to the values of $p$ and $q$ from § 85, the following relation will also take place:

$$\frac{q}{p} = \frac{n}{m}[1 + v\sqrt{\frac{2\mu}{mn}}].$$

If, however, the variable $v$ is restricted, so that, without taking into account its sign

$$v < \frac{\alpha\mu^2}{\sqrt{2\mu\mu'mn}}, \text{ then } \frac{q}{p} < \frac{n'}{m'+1} \text{ or } \frac{q}{p} > \frac{n'}{m'+1}$$

for positive and negative values of the constant α respectively. It occurs that in those cases, because of formulas (77.1b) and (77.1a) respectively,

$$\Pi = 1 - \frac{1}{\sqrt{\pi}}\int_k^\infty \exp(-t^2)dt + \frac{2(\mu'+n')}{3\sqrt{2\pi\mu'm'n'}}\exp(-k^2),$$

$$\Pi = \frac{1}{\sqrt{\pi}}\int_k^\infty \exp(-t^2)dt + \frac{2(\mu'+n')}{3\sqrt{2\pi\mu'm'n'}}\exp(-k^2).$$

Magnitude $k > 0$, and its square is determined by formula (76.1). Here, the extreme values of $m'$ and $n'$, as well as of $p$ and $q$ should be applied. The result is

$$q = \frac{n'}{\mu'+1} - v', \quad p = \frac{m'+1}{\mu'+1} + v',$$

where, for the sake of brevity,

$$\frac{\alpha}{\sqrt{\mu'}} - \frac{v\sqrt{2mn}}{\mu\sqrt{\mu}} - \frac{n'}{\mu'(\mu'+1)} = v'.$$

Magnitude $v'$ is of the order of $1/\sqrt{\mu'}$. [After transformations Poisson gets]

$$k^2 = \frac{\mu'^3 v'^2}{2m'n'} - \frac{(m'-n')\mu'^4 v'^3}{3m'^2 n'^2}, \quad k = \pm k'[1 - \frac{2(m'-n')k'}{3\sqrt{2\mu'm'n'}}],$$

$$k' = \frac{\alpha\mu'}{\sqrt{2m'n'\mu}} - \frac{v\mu'\sqrt{\mu'mn}}{\sqrt{m'n'}}. \tag{86.1}$$

Because of the limit assigned for $v$, magnitude $k'$ has the same sign as α. For $k$ to be positive at $\alpha > 0$ and $\alpha < 0$ we should choose its upper and lower signs respectively.

[After transformations and, in particular, after deduction of new formulas for Π, Poisson continues:] Taking into account the formulas (84.1) and (84.3), the appropriate expressions of Π′ become:

$$\Pi' = \frac{1}{\sqrt{\pi}} \int \exp(-v^2)dv - \frac{1}{\pi}\int_{k'}^{\infty}\int \exp(-t^2 - v^2)dtdv +$$

$$\frac{2n'}{\sqrt{2\pi\mu'm'n'}}\int \exp(-k'^2 - v^2)dv - \frac{2(m-n)}{3\sqrt{2\pi\mu mn}}[\int \exp(-v^2)v^3 dv -$$

$$\frac{1}{\sqrt{\pi}}\int_{k'}^{\infty}\int \exp(-t^2 - v^2)v^3 dtdv], \qquad (86.2a)$$

$$\Pi' = \frac{1}{\pi}\int_{-k'}^{\infty}\int \exp(-t^2 - v^2)dtdv + \frac{2n'}{\sqrt{2\pi\mu'm'n'}}\int \exp(-k'^2 - v^2)dv -$$

$$\frac{2(m-n)}{3\pi\sqrt{2\mu mn}}\int_{-k'}^{\infty}\int \exp(-t^2 - v^2)v^3 dtdv]. \qquad (86.2b)$$

[Poisson simplifies these formulas; thus, just like above, he extends the limits of $v$ from $-\infty$ to $\infty$ and introduces]

$$\frac{\alpha\mu'}{\sqrt{2m'n'\mu}} = \pm\beta, \quad \frac{\mu'\sqrt{\mu'mn}}{\sqrt{m'n'}} = \gamma, \; t = \theta \pm \gamma v, \; dt = d\theta.$$

Here $\beta > 0$ and the signs correspond to those of $\alpha$. In formulas (86.2a, b), respectively, $\alpha > 0$ and $< 0$,

$$k' = \pm\beta - \gamma v, \; t = \theta \mp \gamma v$$

and in both cases the integrals are still over $\theta = \beta$ and $\infty$. Therefore

$$\Pi' = 1 - \frac{1}{\pi}\int_{\beta}^{\infty}\int_{-\infty}^{\infty} \exp[\theta^2 - 2\gamma\theta\, v (1+\gamma)^2 v^2]\, d\theta\, dv +$$

$$\frac{2n'}{\sqrt{2\pi\mu'm'n'}}\int_{-\infty}^{\infty} \exp[\beta^2 - 2\gamma\beta\, v (1+\gamma)^2 v^2]\, dv +$$

$$\frac{2(m-n)}{3\pi\sqrt{2\mu mn}}\int_{\beta}^{\infty}\int_{-\infty}^{\infty} \exp[\theta^2 - 2\gamma\theta\, v (1+\gamma)^2 v^2]\, v^3 d\theta\, dv$$

$$\Pi' = \frac{1}{\pi}\int_{\beta}^{\infty}\int_{-\infty}^{\infty} \exp[\theta^2 - 2\gamma\theta\, v (1+\gamma)^2 v^2]\, d\theta\, dv +$$

$$\frac{2n'}{\sqrt{2\pi\mu'm'n'}}\int_{-\infty}^{\infty} \exp[\beta^2 - 2\gamma\beta\, v (1+\gamma)^2 v^2]\, dv -$$

$$\frac{2(m-n)}{3\pi\sqrt{2\mu mn}}\int_{\beta}^{\infty}\int_{-\infty}^{\infty} \exp[\theta^2 - 2\gamma\theta\, v (1+\gamma)^2 v^2]\, v^3 d\theta\, dv$$

[…] The first value of Π′ is the probability that

$$n' \leq \beta\frac{n\mu'}{\mu} + \sqrt{\frac{2m'n'}{\mu'}},$$

with the sum on the right side very little exceeding its first term. The second value of Π′ is the probability that

$$n' \geq \beta\frac{n\mu'}{\mu} - \sqrt{\frac{2m'n'}{\mu'}},$$

with the difference on the right side slightly less than its first term[6].

**87.** It is possible to note that, because of the infinite limits of the variable $v$, the first two integrals for both values of Π′ coincide, and the third integrals are equal in magnitude and contrary in signs. The first value exceeds the second by [after transformations]

$$\varphi = 1 - \frac{2}{\sqrt{\pi}} \int_u^\infty \exp(-t^2)dt,$$

where $u = \beta/\sqrt{1+\gamma^2}$. When taking into account the value of $\gamma$, it is seen that $\varphi$ expresses the probability that $n'$ is contained within the limits

$$\frac{n\mu'}{\mu} \mp \frac{u\sqrt{2(\mu^3 m'n' + \mu'^3 mn)}}{\mu\sqrt{\mu\mu'}} \qquad (87.1)$$

or equal to the superior limit. If desired that the interval between the limits included the inferior limit as well, the probability that $n'$ is exactly equal to it, see formula (85.2)[7] with $\alpha\sqrt{\mu'}$ equal to the second term of (87.1) and taken with the sign *plus*, should be added to $\varphi$.

Denote by $w$ the probability that $n'$ will be within the stated limits or equal to one of it:

$$w = 1 - \frac{2}{\sqrt{\pi}} \int_u^\infty \exp(-t^2)dt +$$

$$\frac{\sqrt{\mu\mu'}}{\sqrt{2\pi m'n'(\mu+\mu')}} \exp[-\frac{u^2(\mu^3 m'n' + \mu'^3 mn)}{\mu^2 m'n'(\mu+\mu')}]. \qquad (87.2)$$

When comparing this value of $w$ with $R$ in formula (83.2), it is seen that they only differ in their last terms and are therefore almost equal to each other. However, if the number $\mu'$ of future trials is not very small as compared with $\mu$, the second terms in the expressions (87.1) and (83.3) for the limits of $n'$, corresponding to $w$ and $R$ will not coincide. The limits whose probability is $w$, can be much less narrow than in the second case.

If w little differs from certainty, the approximate values $n\mu'/\mu$ and $m\mu'/\mu$ of $n'$ and $m'$ can be substituted instead of them in the respective limits. They will become

$$\frac{n\mu'}{\mu} \mp \frac{u}{\mu}\sqrt{2\mu'mn(1+h)},$$

where $h = \mu'/\mu$. When comparing these limits with (83.3), it is seen that, for the same values of $u$ they widened being multiplied by $\sqrt{1+h}$. For equating them with each other the value of $u$ has to be decreased by the same factor and their probability will lower with respect to $R$. If $h$ is a very small fraction, formulas (83.2) and (87.2) will almost coincide as also the corresponding limits of the values of $n'$. This result is in agreement with what I found out by another method (Poisson 1830).

Formula (87.2) also expresses the probability that the difference $(n'/\mu' - n/\mu)$ is contained within the limits

$$\mp \frac{u\sqrt{2(\mu^3 m'n' + \mu'^3 mn)}}{\mu\mu'\sqrt{\mu\mu'}} \qquad (87.3)$$

or equals one of them. With opposite signs, those are also the limits of $(m'/\mu' - m/\mu)$. And, when $u = 3$ or 4, so that the probability $w$ becomes very close to certainty (§ 80)[8], but when the abovementioned observed differences notably overstep those limits, we can justifiably conclude that the unknown chances of events E and F had most likely changed either during the trials or between the two series.

Note that at the same value of $u$ and, therefore, with the same probability, the limits indicated above will be widest at $\mu = \mu'$ and narrower when one of these numbers is very large as compared with the other. If $\mu' = \mu$, $m' \approx m$ and $n' \approx n$ and the coefficient of $u$ [in the expression (87.3)] will become equal to $2\sqrt{mn}/\mu\sqrt{\mu}$. If, on the contrary, $\mu'$ is very large as compared with $\mu$, then $m' \approx m\mu'/\mu$ and $n' \approx n\mu'/\mu$ and the same coefficient becomes equal to $\sqrt{2mn}/\mu\sqrt{\mu}$, which is $\sqrt{2}$ times less.

**88.** Suppose that the contrary events E and F with unknown chances $p$ and $q$ arrived $m$ and $n$ times in a very large number $\mu$ of trials, and that other contrary events $E_1$ and $F_1$ with unknown chances $p_1$ and $q_1$ occurred $m_1$ and $n_1$ times in a very large number $\mu_1$ of trials. If the ratios $m/\mu$ and $m_1/\mu_1$, as well as $n/\mu$ and $n_1/\mu_1$, essentially differ one from another, the inequalities $p \neq p_1$ и $q \neq q_1$ should be considered certain or almost so.

However, if the differences of those ratios are small fractions, the inequalities of the chances are possibly insensible and happened because the events did not at all arrive exactly proportional to their chances. It is therefore useful to determine the probability of the inequalities $p \neq p_1$ and $q \neq q_1$, corresponding to the small differences $(m/\mu - m_1/\mu_1)$ and $(n/\mu - n_1/\mu_1)$, equal in magnitude and opposite in signs. And this is what I will do.

Just like in formula (84.2), I will denote by $p$ a magnitude little differing from $m/\mu$, so that $v$ will be a positive or negative variable, very little with respect to $\sqrt{\mu}$. Let also

$$p_1 = m_1/\mu_1 - (v_1/\mu_1)\sqrt{2m_1n_1/\mu_1}$$

only little differ from $p$, so that $v_1$ will be a positive or negative variable very small with respect to $\sqrt{\mu_1}$ and

$(m_1/\mu_1 - m/\mu) = \delta,$

where $\delta$ is a small fraction, also positive or negative. Then

$$p_1 - p = z = \delta + \frac{v}{\mu}\sqrt{\frac{2mn}{\mu}} - \frac{v_1}{\mu_1}\sqrt{\frac{2m_1n_1}{\mu_1}},$$

$$v_1 = (\delta - z)\mu_1\sqrt{\frac{\mu_1}{2m_1n_1}} + \frac{v\mu_1}{\mu}\sqrt{\frac{\mu_1 mn}{m_1n_1}}.$$

If $\varepsilon$ is a small positive fraction, then, if the probability of $p_1 \geq p + \varepsilon$ is required, the variable $z$ should take positive values not less than $\varepsilon$. The infinitely small probabilities of the previous values of $p$ and $p_1$ will be $Vdv$ and $V_1dv_1$ with $V$ determined by formula (84.1), and $V_1$, by the same formula with $\mu$, $m$, $n$ and $v$ replaced by $\mu_1$, $m_1$, $n_1$ and $v_1$. The probability of the concurrence of these two values is $VdvV_1dv_1$, and the required probability will be

$$\lambda = \iint VV_1 dv dv_1.$$

For the sake of simplification I neglect the second term of the formula (84.1) and then

$$\lambda = \frac{1}{\pi}\iint \exp(-v^2 - v_1^2)dv dv_1.$$

If desirable, $v_1$ can be replaced by the variable $z$, then $dv_1/dz$ should be calculated by issuing from the previous value of $v_1$. The variable $v_1$, as assumed here, increases, $z$ should therefore also increase, and it is necessary to change the sign of $dv_1$, so that

$$dv_1 = (\mu_1\sqrt{\mu_1}/\sqrt{2m_1n_1})dz .$$

[After troublesome work Poisson got]

$$\lambda = \frac{1}{\sqrt{\pi}}\int_u^\infty \exp(-t^2)dt, \ \lambda = 1 - \frac{1}{\sqrt{\pi}}\int_u^\infty \exp(-t^2)dt \quad (88.1a, b)$$

for ε − δ > 0 and < 0 respectively. Here,

$$\frac{(\varepsilon - \delta)\mu\mu_1\sqrt{\mu\mu_1}}{\sqrt{2(\mu^3 m_1 n_1 + \mu_1^3 mn)}} = \pm u, \qquad (88.2)$$

where $u$ is a positive magnitude whose signs correspond to ε − δ > 0 and < 0.

I ought to remark that, when neglecting the second term in the formula (84.1), the probability of the exact equality $p_1 − p = \varepsilon$ is also neglected and λ becomes the probability of $p_1 − p > \varepsilon$ rather than $\geq \varepsilon$. At ε = δ the magnitude $u = 0$, and both values of λ become equal to 1/2.

The formulas (88.1) also serve for calculating the probability that the unknown chance $p_1$ exceeds a given fraction. Indeed, I assume that in [the appropriately transformed equation] (88.2)

$$\mu = \infty, \ m/\mu = w, \ \delta = (m_1/\mu_1) - w.$$

Then […] for the sake of simplification, $w$ replaces ε + $w$, and μ, $m$, $n$ are substituted instead of $\mu_1, m_1, n_1$:

$$u = \pm (w - \frac{m}{\mu})\frac{\mu\sqrt{\mu}}{\sqrt{2mn}}. \qquad (88.3)$$

According to the difference ($w − m/\mu$) being positive or negative, the formulas (88.1a, b) express the probability that the unknown chance of an event occurring $m$ times in a very large number $\mu = m + n$ of trials exceeds a given fraction $w$.

**89.** Turning now to numerical applications of those various formulas, I choose the Buffon experiment (§ 50) as an example. The arrival of *heads* and *tails* in numerous tosses of a coin will be events E and F. According to Buffon, $m = 2048$, $n = 1992$, $\mu = 4040$ ($m$ и $n$ are the numbers of the arrivals of these events in μ tosses). Substituting these values in formula (83.1) and assuming that $u = 2$, we get

$$\frac{2}{\sqrt{\pi}}\int_u^\infty \exp(-t^2)dt = 0.00468, \ R = 0.99555.$$

At the same time, $(0.50693 \pm 0.02225)$ will be the corresponding limits of $p$. The unknown chance $p$ of the occurrence of *heads* will with probability 0.99555 be contained within 0.48468 and 0.52918. If desired to determine the probability of that chance exceeding 1/2, the preceding values of μ, $m$ and $n$ should be inserted in (88.3). If choosing $w = 1/2$ and the inferior sign, and therefore formula (88.1b), we obtain

$u = 0.62298, \lambda = 0.81043, 1 − \lambda = 0.18957.$

This proves that we can bet somewhat less than 5 against 1 on the chance of *heads* to exceed 1/2.

We can subdivide the Buffon experiment in two parts consisting of 2048 and 1992 trials. In the first one, *heads* appeared 1061 times, and *tails*, 987 times; in the second part, 987 и 1005 times respectively. When, however, considering all the trials and applying the formula (87.2), we can also calculate the probability that the number of the arrivals of *heads* and *tails* should be contained within given limits in each partial experiment. For accomplishing this aim, we assume in that formula and in the corresponding limits

$m'/\mu' = m/\mu = 0.50693$, $n'/\mu' = n/\mu = 0.49307$.

In other words, we replace $m'/\mu'$ and $n'/\mu'$, which were not supposed to be known, by their approximate values corresponding to the entire experiment. This is possible since $m'$ and $n'$ only enter in the terms of the order of smallness of $1/\sqrt{\mu}$. Instead of $\mu$ we should also assume its total value, 4040.

For the first part of the experiment, $\mu' = 2048$; at $u = 2$ (see above) $w = 0.99558$ will be the probability that $n'$, that is, the number of the arrivals of *tails*, will be contained within the limits $1001 \pm 79$. This condition was indeed fulfilled. For the second part, $\mu' = 1992$ and $w$, again with $u = 2$, was 0.99560, which is the probability that the number $n'$ of the arrivals of tails is contained within the limits $982 \pm 77$, which was also fulfilled. […]

Suppose that we do not know whether the same coin was tossed in both parts of the experiment and that, knowing their results, it is required to determine whether the probability λ of the chance of *heads* in the first part by a given fraction exceeded the same chance in the second. We should first of all insert $\mu = 1992$, $m = 987$, $n = 1005$ and $\mu_1 = 2048$, $m_1 = 1061$, $n_1 = 987$ in equation (88.2). Then, $\delta = m_1/\mu_1 - m/\mu = 0.02257$, and that equation becomes

$u = \pm 44.956 (\varepsilon - 0.02257)$.

If, for example, $\varepsilon = 0.02$, we should choose the inferior sign and apply formula (88.1b), so that

$u = 0.11553$, $\lambda = 0.56589$, $1 - \lambda = 0.43411$

and we can bet barely 4 against 3 on the chance of *heads* to be larger by 1/50 in the first part than in the second.

If $\varepsilon = 0.025$ we should choose the superior sign and apply formula (88.1a), so that

$u = 0.10925$, $\lambda = 0.43861$, $1 - \lambda = 0.56139$

and we can bet less than 1 against 1 on that excess to be larger than 1/40.

**90.** I am now solving a problem admitting an interesting application; it will be based on the preceding formulas and on a following lemma[9].

An urn contains $c = a + b$ balls, $a$ of them white and $b$, black; $l$ balls are randomly extracted without replacement, then, in the same way, $\mu = m + n$ other balls are drawn. I say that the probability of drawing $m$ white and $n$ black balls in that second series of extractions does not depend either on the number, or on the colour of those extracted previously[10] and is equal to their arrival in the case in which $l = 0$.

Suppose that the $l + \mu$ drawings are carried out consecutively; denote by $i$ the total number of combinations of the $l + \mu$ extracted balls; by $i'$ the number of combinations of the $\mu$ last balls, $m$ of them white, and $n$, black; and by $i_1$ the number of combinations of the $\mu$ first balls, $m$ of them white, and $n$, black. The chance of extracting $m$ white and $n$ black balls after the first series of drawings is $i'/i$, and the chance of the same result before all the drawings was $i_1/i$. The numbers $i'$ and $i_1$ are equal because in general combinations of $l$, then of $\mu$ definite balls coincide with the combinations of $\mu$ and then of $l$. In particular, for each combination of the $\mu$ last extracted balls, $m$ of them white and $n$, black, there always exists a combination of the $\mu$ first extracted balls with the same number of balls of both colours, and vice versa. The fractions $i'/i$ and $i_1/i$, are therefore also equal as are the probabilities which they express, which it was required to prove.

It is possible to verify this result in the following way. The urn contained $a$ white and $b$ black balls; the chance of extracting $m$ and $n$ such balls in $(m + n)$ first drawings is a function $f(a, b, m, n)$. Similarly, the chance of $g$ and $h$ white and black balls appearing in $(g + h)$ first drawings will be $f(a, b, g, h)$. The number of balls left in the urn will become $(a - g)$ and $(b - h)$[11], and the chance of $m$ white and $n$ black balls arriving in $\mu = m + n$ new drawings will be $f(a - g, b - h, m, n)$. The product of the two last values of that function will be equal to the chance of drawing $m$ white and $n$ black balls after the appearance of those $g$ white and $h$ black balls.

Therefore, adding up the $(l + 1)$ values[12] of that product corresponding to all the $l$ natural and zero values of $g$ and $h$, we will obtain the complete expression of the chance of extracting $m$ white and $n$ black balls after the appearance of $l$ balls of some colour. It is required to establish that that chance does not depend on $l$ and is equal to $f(a, b, m, n)$, i. e., to establish that

$$f(a, b, m, n) = \sum f(a, b, g, h) f(a - g, b - h, m, n),$$

where the sum extends from $g = 0$ and $h = l$ to $g = l$ and $h = 0$.

I note that, according to § 18,

$$f(a,b,m,n) = \frac{\varphi(m,n)\, \varphi(a-m,b-n)}{\varphi(a,b)}, \quad \varphi(a,b) = C_{a+b}^{a}$$

so that

$$f(a,b,g,h) f(a-g,b-h,m,n) =$$

$$\frac{\varphi(g,h)\varphi(a-g,b-h)}{\varphi(a,b)} \frac{\varphi(m,n)\varphi(a-g-m,b-h-n)}{\varphi(a-g,b-h)}$$

or, which is the same,

$$f(a,b,g,h)f(a-g,b-h,m,n) =$$
$$\frac{\varphi(m,n)}{\varphi(a,b)}\varphi(g,h)\varphi(a-g-m,b-h-n).$$

And, taking into account the value of $f(a, b, m, n)$, the equation to be verified becomes

$$\varphi(a-m, b-n) = \sum \varphi(g,h) \cdot \varphi(a-g-m, b-h-n),$$

where $\varphi(m, n)/\varphi(a, b)$ was cancelled from both sides. And, since $a$ and $b$ are arbitrary, they can be replaced by $(a + m)$ and $(b + n)$. We will obtain

$$\varphi(a, b) = \sum \varphi(g,h) \cdot \varphi(a-g, b-h).$$

The left side is the coefficient of $x^a y^b$ in the expansion of $(x + y)^c$, and the right side, is again the coefficient of $x^a y^b$ in the product of the expansions of $(x + y)^l$ and $(x + y)^{c-l}$, that is, of $(x + y)^c$. These sides are identical, which it was required to verify.

**91.** Suppose that the numbers $a, b, a - m, b - n$ are very large. The approximate values of $\varphi(m, n)$, $\varphi(a - m, b - n)$, $\varphi(a, b)$, as well as of $f(a, b, m, n)$ are calculated by means of the Stirling formula. If only taking into account its first term and calculating $f(a, b, m, n)$, this probability can be represented [after very difficult transformations] in the form

$$f(a,b,m,n) = H\exp(-t^2)[1 - \frac{4t^3(a-b)(c-2\mu)}{3\sqrt{2(\mu)\mu}} \quad abc], \quad (91.1)$$

$$H = \sqrt{\frac{ab\mu(c-\mu)}{2\pi cmn(a-m)(b-n)}},$$

as the chance of extracting $m$ white and $n$ black balls

$$m = \frac{\mu a}{c} - \frac{t\sqrt{2(\mu)\mu} \quad abc}{c^2}, \quad (91.2a)$$

$$n = \frac{\mu b}{c} + \frac{t\sqrt{2(\mu)\mu} \quad abc}{c^2}. \quad (91.2b)$$

The difference $(n - m)$ will be even if $\mu$ is even, and odd otherwise. Let $i$ be a positive natural number and $n - m = 2i$ or $2i - 1$, then $t$ will be

$$t = 2\delta + \gamma, \quad \delta = \frac{c^2}{2\sqrt{2(c-\mu)\mu\ abc}},$$

$$\gamma = \frac{(a-b)\mu c}{2\sqrt{2(c-\mu)\mu\ abc}}, \quad \gamma_2 = \frac{(a-b)\mu c}{2(\sqrt{(c-\mu)\mu\ abc}} - \delta$$

for even and odd values of μ. When inserting that *t* in formula (91.1), it will represent the probability that after μ consecutive drawings the number of the extracted black balls exceeds the number of the arrived white balls by 2*i* or 2*i* − 1. Taking *i* = 1, 2, 3, … until exp(− $t^2$) becomes insensible, or, if desired, until *i* = ∞, the sum of the results will be the probability *s* that after those μ drawings the number of the black balls will exceed the number of the white balls by some even or odd numbers of unities:

$$s = \sum H \exp(-t^2)[1 - \frac{4t^3(a-b)(c-2\mu)}{3\sqrt{2(c-\mu)\mu\ abc}}],$$

where the sum extends over all the values of *t* from γ + 2δ to ∞, increasing by increments of 2δ. According to the hypotheses, 2δ is a very small fraction, and the sum can express a rapidly converging series arranged by the powers of such an increment.
 [After long transformations]

$$s = 1 - \frac{1}{\sqrt{\pi}} \int_v^\infty \exp(-t^2)dt - \Gamma \exp(\gamma)^2, \qquad (91.3a)$$

$$s = \frac{1}{\sqrt{\pi}} \int_v^\infty \exp(-t^2)dt - \Gamma \exp(\gamma)^2 \qquad (91.3b)$$

for γ < 0 and > 0 and positive *v* equal in magnitude to γ, where

$$\Gamma = \frac{(a-b)(c-2\mu)(7+4\gamma^2)+3c^2}{6\sqrt{2\pi(c-\mu)\mu abc}}.$$

The probability that after μ drawings *m* = *n* = μ/2, which is only possible for even values of μ, is

$$\sigma = \frac{c^2 \exp(\gamma)^2}{\sqrt{2\pi(c-\mu)\mu abc}}. \qquad (91.4)$$

**92.** Suppose that after μ drawings μ′, and then μ″, … other balls were extracted, until all the *c* balls in the urn are drawn. Then

*c* = μ + μ′ + μ″ + …

Suppose also that μ′, μ″, … as well as μ are very large numbers. Denote by *s*′, *s*″, … the new values of *s* obtained after inserting μ′, μ″, … instead of μ when applying formulas (91.3a, b) depending on whether before the drawings the number *b* of black balls in the urn was more or less than the number *a* of the white balls and, consequently, whether γ became negative or positive. By the lemma of § 90 the chances of extracting more black than white balls in μ, μ′, μ″, … consecutive drawings are *s, s′, s″*, … These chances only vary owing to the inequalities of μ, μ′, μ″, …; they would have been identical had those numbers been equal.

Denote the mean value of *s, s′, s″*, … by *r*. Suppose that the total number α of the drawings is very large[13], and that in *j* of them the number of black balls exceeded the number of white balls. Then, by the first proposition of § 52, the probability that *j* is within given limits will be the same as when all the chances *s, s′, s″*, … are equal one to another and equal to their mean, *r*. Therefore, substituting α, *r*, (1 − *r*) instead of μ, *q, p* in formula (79.1), we will get the probability

$$R = 1 - \frac{2}{\sqrt{\pi}} \int_u^\infty \exp(-t^2)dt + \frac{1}{\sqrt{2\pi\alpha r(1-r)}} \exp(-u^2)$$

that *j* is within the limits $\alpha r \mp u\sqrt{2\alpha r(1-r)}$ or equal to one of them. Here, *u* is small as compared with √α.

So this is the solution of the problem I proposed. It can be applied to the election of deputies in a large country, in France for example. Denote the number of the electors in France by *c*; among them, *a* hold one opinion, and *b* = *c* − *a* are of the opposite view. They are randomly distributed over α electoral colleges, and a deputy is elected in each by a majority vote. It is required to determine the probability *R* that the number *j* of deputies holding the second opinion will be contained within given limits if there are μ, μ′, μ″, … voters in those colleges.

Let the limits of *j* be those just indicated, then the required probability *R* will be represented by the previous formula. Each electoral college consists of voters living in the same locality rather then selected randomly from their general list, as we supposed above. And still it will be useful to find out what happens by our hypothesis and to show it by examples.

**93.** In France, the number of electoral colleges, as also the number of deputies, is 459, and the total number of voters can be estimated as approximately 200,000[14]. I suppose that the numbers μ, μ′, μ″, … coincide, and that μ is odd:

α = 459, μ = 435, *c* = αμ = 199,665.

Let also *a* = 94,835 and *b* = 104,830, so that their difference amounts to almost 1/20 of the number of voters. Then γ < 0 and *v* = − γ. Assuming the second value of γ from § 91, we get

$$v = 0.77396, \quad \frac{1}{\sqrt{\pi}} \int_v^\infty \exp(-t^2)dt = 0.13684,$$

and, owing to formula (91.3a), $s = 0.85426$, $1 − s = 0.14574$.

The chance of electing a deputy by the voters of the more numerous party will therefore exceed 21/25, whereas the minority of the electors, although not much differing from their majority, can never hope to elect more than 4/25 of the deputies. Substituting the values of $s$ and $(1 − s)$ instead of $r$ and $(1 − r)$ in the expression of $R$ from § 92, we will find that, at $\alpha = 459$ and $u = 2$, $R = 0.99682$ for the probability that the number of the deputies elected by the more numerous party will be contained within $392 \mp 21$ and within $67 \pm 21$ by the other party. As compared with $\alpha$, these limits are wide since that $\alpha$ is not extremely large.

I invariably suppose that the difference $(b − a)$ is almost $c/20$. But let now $\mu$ be even:

$\alpha = 459$, $\mu = 436$, $c = \alpha\mu = 200{,}124$.

Suppose also that $a = 95{,}064$ and $b = 105{,}060$. As previously, $v = -\gamma$, but the first value of $\gamma$ should now be selected from § 91, so that

$$v = 0.74006, \quad \frac{1}{\sqrt{\pi}} \int_v^\infty \exp(-t^2)dt = 0.14764,$$

and $s = 0.84279$, $1 − s = 0.15721$.

Since $\mu$ is even, the case of $m = n$ is possible. By formula (91.4) its chance is $\sigma = 0.02218$. Adding $\sigma/2$ to $s$, we have $s = 0.85388$, very little less than in the case of an odd $\mu$.

To show the influence of the inequality of the number of voters in the colleges, I will suppose that a half of them are equally distributed over 1/3 of the colleges, and the other half, over the other colleges. For the first third

$\alpha/3 = 153$, $\mu = 654$, $\alpha\mu/3 = 100{,}062$,

and for the rest,

$2\alpha/3 = 306$, $\mu = 327$, $2\alpha\mu/3 = 100{,}062$.

I will also suppose that $a = 95{,}062$, $b = 105{,}062$, $c = 200{,}124$, so that both parties are as unequally strong as previously, with the inequality being almost 1/20 of the voters. In the first case $\mu$ is even, and odd in the second case. Respectively,

$s = 0.89429$, $\sigma = 0.01376$, $s + \sigma/2 = 0.90117$; $s = 0.81981$,

and the stronger party has a mean chance of electing a deputy equal to

$r = (0.90117 + 0.81981)/2 = 0.86049$.

It little exceeds the corresponding chance in the case in which the number of the voters in each college was the same. However, if the difference $(b − a)$ begins to increase, the chance of electing a deputy by the minority will decrease very rapidly and soon almost disappears. For proving this statement, I will suppose that the voters are equally distributed over the colleges and assume the values of α, μ и c from the first example. In addition, let $a = 89,835$ and $b = 109,830$, so that the difference $b - a \approx c/10$, twice larger than in that example. Then

$s = 0.98176, 1 - s = 0.01824$

and the chance of electing a deputy by the minority party becomes not greater than about $1/60^{15}$. Since $s$ is small, the probability $P$ that in such elections the number of deputies from the minority party will not exceed a given number $n$, should be calculated by the formula of § 81. It occurs that

$w = α(1 - s) = 8.3713, n = 15, P = 0.98713, 1 - P = 0.01287$.

And if the difference between the strengths of the parties is increased to 30,000, or to 3/20 of the total number of voters, the chance $(1 - s)$ becomes less than 1/1000. An election of one single deputy by the minority party will then be very unlikely.

But then, a representative government becomes nothing but a deception, since a minority of 90,000 voters out of their total number 200,000 will only be represented by a very small number of deputies, while a minority of 85,000 will only have a very slim chance of electing one single representative of their interests to the chamber of deputies. And if, between two sessions, only 3/20 of the voters change their opinion, the entire chamber will follow suit.

Voters in each electoral college are not selected randomly from their general list for the entire France, as I have supposed. In each district the prevailing opinion is formed and maintained by particular causes, such as local interests, the influence of the Government and some citizens. But still it is useful to indicate the extreme variability that randomness can effect in the composition of the chamber of deputies by very small changes in the ratio of voters keeping to opposite opinions.

## Notes

**1.** Notation $n!$ was not yet known but Poisson introduced two other appropriate symbols. In this section, only the Stirling formula (67.3) is mainly needed, but Poisson referred below to some intermediate results.

**2.** The Stirling series is divergent (Fichtenholz 1947/1950, vol. 2, chapter 8, § 5/501, p. 820).

**3.** I would say: $p$ and $q$ are given in the direct, but not in the inverse problem. This circumstance explains why the variances of the respective random variables are different. In the inverse problem, the variance is larger, and more trials are needed for achieving the same precision as in the direct problem. Bayes, who did not yet possess the notion of variance (introduced by Gauss in 1823), understood this fact (Sheynin 2010). Laplace (1814/1994, p. 120) approvingly mentioned Bayes, but did not even refer to his appropriate memoir. Poisson (paper of 1837, p. 73) apparently only once actually applied the variance when estimating the quality of guns.

**4.** That formula is due to Montmort (1708/1713, p. 245), see also Todhunter (1865/1965, p. 97).

**5.** Note that $\int_{-\infty}^{t} (x - a)^r \varphi(x) dx$, where φ(x) is the appropriate density, is called the *r*-th incomplete moment. Below, Poisson calculated the integral of $t^{2i}\exp(-t^2)$ over [0, ∞] at various values of *i*, whereas Gauss (1816) had established the usual (complete) moments of the normal distribution. The relation between Poisson's § 82 and that memoir of 1816 deserves to be studied, but in any case Poisson's persistent refusal to apply the results of that great scholar turned out against him.

**6.** As printed, both values of the probability were thus lower than approximately one and the same expression. Corrected in the translation.

**7.** Formula (85.2) does not contain $\alpha\sqrt{\mu'}$.

**8.** At the stated place Poisson mentioned *u* equal to 4 or 5.

**9.** After the publication of my note [see note 18 to Chapter 1], I was informed that the proposition included there was already contained in this lemma which I (1825, p. 70) had applied for the solution of the problem of *thirty-and-forty*. Poisson. The page number is wrong. O. S.

**10.** Poisson had thus considered subjective probabilities.

**11.** The explanation, both here and below, is not clear enough.

**12.** It is seen below that *l* = *g* + *h*.

**13.** It can be concluded that α = *c*. The exposition is unclear.

**14.** Voters constituted only a small part of the population.

**15.** More precisely, 1/55.

# Chapter 4. Calculus of Probabilities Depending on Very Large Numbers, Continued

## Misprints/Mistakes Unnoticed by the Author

**1.** In § 94, p. 247 of the original text, the first of the three integrals in this section lacks *dx*.

**2.** In § 95, p. 251 of the original text. Both integrals in the first displayed formula lack limits (and left without limits in the translation). The first of them also lacks the differential *dz*.

**3.** In § 96, p. 253 of the original text. The integrand in the third displayed formula should be exp $(-t^2)$; the $t^2$ is missing.

**4.** In § 103, p. 274 of the original, the first integral representation of $\gamma_i$. The result of integration should include $\gamma_i$ rather than γ.

**5.** In § 104, p. 277 of the original text. On line 10 the chance of A is mentioned instead of chance E.

**6.** In § 106, p. 285 of the original text, after the first displayed formula; μ in the denominator should have been mentioned rather than 1/μ.

**7.** In § 106, p. 286 of the original. Poisson introduced a magnitude λ and included it in three equalities. On the left side of the third equality appeared the nonsensical difference $s/\mu - s^2/\mu^2$. There also, on p. 287, line 3 after first displayed formula. The integral should be over [− *u*, *u*] rather than over [*u*, − *u*].

**8.** In § 107, p. 288 of the original text. The second equation in the second displayed formula. Replace the denominator μ on the left side by μ′.

**9**. In § 109, p. 295 of the original text, line 4 from bottom. Magnitude δ should be very small as compared with $1/\sqrt{\mu}$ rather than with μ. The correct statement is on p. 296, again on line 4 from bottom.

**10.** In § 110, p. 301 of the original text, line 5 from bottom. Replace *maximal natural number included in s* by … *included in* α.

**11.** In § 112/2, p. 308 of the original text, line 2 after formula (b). Replace *chances F and F* by *chances E and F*.

**12.** In § 112/8, p. 312 of the original text. The last term in the formula of Γ lacks the exponent. There also, the first bracket lacks the term − 4*g*.

**13.** In § 112/12, p. 315 of the original text, line 2 after first displayed formula. Replace *probability that magnitude γ* by … *magnitude A*.

**14.** At the beginning of § 113, p. 316 of the original text. Laplace's initial equations are only partly, and, for that matter, mistakenly mentioned. See Note 18.

If not stated otherwise, all these misprints/mistakes are corrected in the translation.

**94.** We will now consider formulas which are related to variable chances. This will lead us to prove the three main propositions, indicated in §§ 52 and 53 and therefore to the *law of large numbers*.

Suppose that a series of $\mu = m + n$ successive trials is made with the somehow variable chances of contrary events E and F being $p_1$ and $q_1$ at the first trial, $p_2$ and $q_2$ in the second, ..., and $p_\mu$ and $q_\mu$, in the last one. Then

$$p_1 + q_1 = p_2 + q_2 = \ldots = p_\mu + q_\mu = 1.$$

Denote by $U$ the probability that those events will arrive $m$ and $n$ times in some order. By the rule of § 20 $U$ will be the coefficient of $u^m v^n$ in the expansion of the product[1]

$$X = (up_1 + vq_1)(up_2 + vq_2) \ldots (up_\mu + vq_\mu). \qquad (94.1)$$

If $u = \exp(x\sqrt{-1})$, $v = \exp(-\sqrt{-1})$, the term $Uu^m v^n$ of that product will be $U\exp[(m-n)x\sqrt{-1}]$. All other terms will have exponents differing from the indicated so that, when multiplying the product by $\exp[(m-n)x\sqrt{-1}]$ and integrating over $[-\pi, \pi]$, all the other factors will disappear:

$$\int_{-\pi}^{\pi} X \exp[-(m-n)x\sqrt{-1}]dx = 2\pi U.$$

[...] The factors in (94.1) are

$$up_i + vq_i = \cos x + (p_i - q_i)\sin x \sqrt{-1}.$$

Let

$$\cos^2 x + (p_i - q_i)^2 \sin^2 x = \rho_i^2$$

and introduce a real angle $r_i$

$$(1/\rho_i)\cos x = \cos r_i, \quad (1/\rho_i)(p_i - q_i)\sin x = \sin r_i$$

so that

$$up_i + vq_i = \rho_i \exp(r_i \sqrt{-1}).$$

The magnitude $\rho_i$ is two-valued; we will suppose it positive. For the sake of brevity introduce

$$\rho_1 \rho_2 \ldots \rho_\mu = Y, \quad r_1 + r_2 + \ldots + r_\mu = y. \qquad (94.2\text{a, b})$$

Then [...]

$$U = \frac{1}{2\pi} \int_{-\pi}^{\pi} Y \cos[y - (m-n)x]dx + \frac{\sqrt{-1}}{2\pi} \int_{-\pi}^{\pi} Y \sin[y - (m-n)x]dx.$$

[The right side is reduced to the first integral which is denoted as formula (94.3).]

Usual integration in a finite form is possible. However, if μ is not a large number, that formula will be useless, but if μ is very large, an arbitrarily approximate value of $U$, as shown below, can be derived from it.

**95.** At $x = 0$ each factor in formula (94.2a) equals unity, and is still less at any other values of $x$ within the limits of integration. Therefore, when μ is very large, the product in formula (94.2a), is, in general, very small except for very small values of $x$, and if μ becomes infinite, $Y$ will disappear at any finite values of $x$.

There exists an exceptional case in which all the factors of $Y$ indefinitely tend to unity; the product of an infinite number of such factors is known to be possibly finite. Since

$$\rho_i^2 = 1 - 4p_i q_i \sin^2 x,$$

one of the chances of events E and F or their product $p_i q_i$ indefinitely decreases during the trials, but excluding that particular case, and having a very large μ we can consider the variable $x$ as a very small magnitude and neglect the part of the preceding integral corresponding to other values of $x$.

Therefore, the following series arranged in powers of $x^2$ rapidly converge:

$$\rho_i = 1 - 2p_i q_i x^2 + [(2/3)p_i q_i - 2p_i^2 q_i^2]x^4 - \ldots,$$
$$\ln\rho_i = -2p_i q_i x^2 + [(2/3)p_i q_i - 4p_i^2 q_i^2]x^4 - \ldots,$$
$$\ln Y = -\mu k^2 x^2 + \mu[(1/3)k^2 - k'^2]x^4 - \ldots$$

Here

$$\mu k^2 = 2\sum p_i q_i, \quad \mu k'^2 = 4\sum p_i^2 q_i^2$$

with the sums extending from $i = 1$ to $i = \mu$.

Supposing also that $x = z/\sqrt{\mu}$, considering the new variable $z$ as a very small magnitude as compared with $\sqrt{\mu}$, and neglecting magnitudes of the order of smallness of $1/\mu$, we get

$$Y = \exp(-k^2 z^2).$$

In addition, by the values of $\rho_i$ and $\sin r_i$ it occurs that

$$r_i = (p_i - q_i)x + (4/3)(p_i - q_i)p_i q_i x^3 + \ldots$$

Denote

$$h = (4/3\mu)\sum(p_i - q_i)p_i q_i$$

and let the mean chances of E and F be $p$ and $q$, $p + q = 1$. When only preserving magnitudes of the order of smallness of $1/\sqrt{\mu}$,

$$y = z(p - q)\sqrt{\mu} + z^3 h/\sqrt{\mu},$$
$$\cos[y - (m - n)x] = \cos(zg\sqrt{\mu}) - [z^3 h/\sqrt{\mu}]\sin(zg\sqrt{\mu}),$$

where

$$g = (p - m/\mu) - (q - n/\mu).$$

I substitute these values of $Y$ and $\cos[y - (m - n)x]$ in formula (94.3), and note that $dx = dz/\sqrt{\mu}$ [see above]:

$$U = \frac{2}{\pi\sqrt{\mu}}[\int \exp(-k^2 z^2)\cos(\mu g\sqrt{\ }) \, dz -$$

$$\frac{h}{\sqrt{\mu}}\int \exp(-k^2 z^2) z^3 \sin(\mu g\sqrt{\ }] \cdot dz$$

[After transformations] that probability becomes

$$U = \frac{1\theta\exp(\theta)}{k\sqrt{\pi\mu}}\exp(\theta) - \frac{h}{2k^4\mu\sqrt{\pi}}(3 - 2\theta)^2, \qquad (95.1)$$

where

$$p - m/\mu = \theta k\sqrt{\mu}, \; q - n/\mu = -\theta k\sqrt{\mu}, \text{ so that } g = 2\theta k\sqrt{\mu},$$

and that probability thus corresponds to

$$m = p\mu - \theta k\sqrt{\mu}, \; n = q\mu + \theta k\sqrt{\mu},$$

that is, to numbers, almost proportional to the mean values $p$ and $q$ and the number $\mu$ of trials[2].

**96.** Those $m$ and $n$ are natural numbers, so $\theta$ should be a multiple of $\delta = 1/k\sqrt{\mu}$ or 0. At $\theta = 0$ the formula (95.1) will indicate the probability $1/\pi\mu\sqrt{\ }$ that exactly $m/n = p/q$. Denote by $t$ a positive multiple of $\delta$; substitute $\theta = -t$ and then $t$ in that formula. Their sum

$$\frac{2}{k\sqrt{\pi\mu}}\exp(-t^2) \qquad (96.1)$$

is the probability that $m$ and $n$ will be respectively equal to one of the values

$$p\mu \mp kt\sqrt{\mu}, \; q\mu \pm kt\sqrt{\mu}.$$

Denote by $u$ a given multiple of $\delta$ and assume in (96.1) $t = \delta, 2\delta, \ldots, u$. The sum of the thus obtained results increased by the value of $U$ at $\theta = 0$, will be

$$R = \frac{1}{k\sqrt{\pi\mu}} + \frac{2}{k\sqrt{\pi\mu}} \sum \exp(-t^2).$$

This $R$ is the probability that $m$ and $n$ are contained within the limits

$$p\mu \mp uk\sqrt{\mu}, \quad q\mu \pm uk\sqrt{\mu}$$

or equal to one of them. The sum extends from $t = \delta$ to $u$ and increases by increments equal to $\delta$. However, it can be replaced by the difference of the sums of the exponential function extending from $t = \delta$ to $\infty$ and from $t = u + \delta$ to $\infty$. By the Euler formula, which we had applied in § 91, this latter sum multiplied by $\delta$, when approximating as we should, i. e., when neglecting the square of $\delta$, will be

$$\int_u^\infty \exp(-t^2)dt - \frac{\delta}{2}\exp(-u^2).$$

At $u = 0$ the first sum extending from $t = \delta$ to $\infty$ and multiplied by $\delta$, will be $(\sqrt{\pi} - \delta)/2$. Therefore, subtracting this magnitude from the previous and dividing the difference by $\delta$, we get the sum in the expression for $R$

$$\sum \exp(-t^2) = \frac{1}{2\delta}\sqrt{\pi} \exp(\frac{1}{\delta}\int_u^\infty) \quad -t^2 \; dt \; \exp(\frac{1}{2} \frac{1}{2}). \quad -u^2$$

When taking into account the value of $\delta$, this expression becomes

$$R = 1 - \frac{2}{\sqrt{\pi}} \int_u^\infty \exp(-t^2)dt + \frac{\exp(-u^2)}{k\sqrt{\pi\mu}}. \qquad (96.2)$$

If the values of $p_i$ and $q_i$ are constant and therefore equal to their mean values $p$ and $q$, then $k = \sqrt{2pq}$. Formula (96.2) and the previous limits of $m$ and $n$ will coincide respectively with the formula (79.1) and the limits corresponding to it. If necessary, this coincidence of results obtained by such differing methods can confirm our calculations.

At an insignificant value of $u$, such as 3 or 4, the value of $R$ will be very close to unity. It is therefore almost certain that, if the number $\mu$ of trials is very large, the ratios $m/\mu$ and $n/\mu$ will very little differ from the mean chances $p$ and $q$, and will approach them all the more as $\mu$ increases further, finally coinciding with them if $\mu$ can become infinite. This is the first of the two general propositions of § 52.

**97.** Suppose now that A is a thing able to take many positive and negative values, multiples of a given magnitude $w$. These values are contained within limits $\alpha w$ and $\beta w$ inclusive, so that $\beta - \alpha + 1$ is their number. Here, $\alpha$ and $\beta$ are natural numbers or zeros, and, without considering their signs, the second of them exceeds the first. If, however, A can only take one value, then $\beta = \alpha$. At each trial, made for establishing A, its possible values are not equally probable. In addition, for the sake of greater generality, I suppose that the chances of one and the same value vary from one trial to another.

Let $n$ be some number contained between $\alpha$ and $\beta$ or equal to one of these magnitudes. Then, denote by $N_1$ the chance of the value $nw$ of A at the first trial, by $N_2$, that chance at the second trial, etc., and let $s$ be the sum of those values at $\mu$ successive trials. It is required to determine the probability that that sum is contained within given limits. Denote first of all the probability of the exact equality $s = mw$ by $\Pi$. Here, $m$ is a given number *situated between $\alpha$ u $\beta$* or equal to one of these magnitudes.[Poisson deleted the words which I italicized but after that the phrase became incomprehensible.]

Compile the product

$$\sum N_1 t^{nw} \sum N_2 t^{nw} \ldots \sum N_\mu t^{nw},$$

where $t$ is an indefinite magnitude and the sums extend over all the values of $n$ from $\alpha$ to $\beta$. Expand this product in powers of $t^w$, and it will be seen that $\Pi$ becomes the coefficient of $t^{mw}$. This is evident for $\mu = 1$. Let now $\mu = 2$ and denote by $n'w$ and $n''w$ the exponents of $t$, which that magnitude will take in the corresponding sums. The thing A can evidently take the value $mw$ in as many different ways as there are different solutions of the equation $n' + n'' = m$ with $n'$ and $n''$ being contained between $\alpha$ and $\beta$. The probability of each way is equal to the product of the values of $N_1$ and $N_2$ corresponding to each pair of numbers $n'$ and $n''$. Therefore, the composite probability of the equality $s = mw$ is expressed by the coefficient of $t^{mw}$ in the product of those first two sums. This consideration can easily be extended to $\mu = 3, 4, \ldots$ If all the magnitudes $N_1, N_2, \ldots$ coincide, their product will be the power of $\mu$ of one of the polynomials corresponding to the sums $\sum$, and this case was studied in § 17.

By a similar reasoning, when assuming that $t^w = \exp(\theta \sqrt{-1})$ and denoting by $X$ the product of $\mu$ sums $\sum$, we get[3]

$$\Pi = \frac{1}{2\pi} \int_{-\pi}^{\pi} X \exp(\theta m \sqrt{-\theta}).d$$

Let $i$ and $i_1$ be two given numbers and $P$, the probability that the sum $s$ is contained within $iw$ and $i_1 w$ or equal to one of those limits. Then the value of $P$ can be established by $\Pi$ when assuming that $m = i, i + 1, \ldots, i_1$. And the sum of the values, corresponding to $\exp(\theta m \sqrt{-})$ will be expressed in the following way […].

For simplifying, I suppose that $w$ is infinitely small, $i$ and $i_1$ are infinite numbers, and

$iw = c - \varepsilon$, $i_1 w = c + \varepsilon$, $\theta = wx$, $d\theta = wdx$,

where $c$ и $\varepsilon$ are given constants and $\varepsilon$ is positive so that the inequality $i_1 > i$, implied in the expression of $P$, persists. The limits of the integral with the new variable $x$ are infinite. Then, $\sin(\theta/2) = wx/2$, and, when neglecting the numbers $\pm 1/2$ as compared with $i$ and $i_1$, $P$ becomes

$$P = \frac{1}{\pi} \int_{-\infty}^{\infty} X \exp(-cx\sqrt{-1}) \sin\varepsilon x \, \frac{dx}{x}. \qquad (97.1)$$

The possible values of A increase by infinitely small increments and their number ought to be supposed infinite with the probability of each being infinitely low. Denote given constants by $a$ and $b$, and the continuous variable by $z$, and let

$$\alpha w = a, \; \beta w = b, \; nw = z, \; t^{nw} = \exp(xz\sqrt{-1}).$$

At the same time let

$$N_1 = wf_1 z, \; N_2 = wf_2 z, \; \ldots$$

Each sum included in $X$, becomes a definite integral over $[a, b]$. Assuming $w = dz$, we conclude that

$$X = \int_a^b e^{xz\sqrt{-1}} f_1 z \, dz \int_a^b e^{xz\sqrt{-1}} f_2 z \, dz \ldots \int_a^b e^{xz\sqrt{-1}} f_\mu z \, dz \qquad (97.2)$$

becomes the product of $\mu$ factors and should be substituted in formula (97.1) instead of $X$.

**98.** That formula expresses the probability that after $\mu$ trials the sum of the values of A will be contained within the given magnitudes $c - \varepsilon$ and $c + \varepsilon$. At the $n$-th trial the infinitely small chance of value $z$ will be $f_n z \, dz$, with all possible values of A being contained, by the hypothesis, within the limits $a$ and $b$, and one of these values certainly taking place at each trial. Therefore

$$\int_a^b f_n z \, dz = 1.$$

The function $f_n z$ can be continuous or discontinuous, but positive within the limits $a$ and $b$.

If the chance of each value $z$ remains constant during the trials, this function is independent from $n$. Denoting it by $fz$, we have

$$X = \left( \int_a^b e^{xz\sqrt{-1}} fz \, dz \right)^\mu, \; \int_a^b fz \, dz = 1.$$

And if all the values of A are equally probable, $fz$ will be a constant, and, since it should satisfy the last equation, equal to $1/(a - b)$.

Let $a = h - g, \; b = h + g$, then

$$fz = \frac{1}{2g}, \quad \int_a^b e^{xz\sqrt{-1}} fz\, dz = \frac{\sin gx}{gx} e^{hx\sqrt{-1}},$$

and formula (97.1) becomes [...]

$$P = \frac{2}{\pi} \int_0^\infty \left(\frac{\sin gx}{gx}\right)^\mu \frac{\sin\varepsilon x}{x} \cos[(\mu h - c)x]\, dx.$$

[After long transformations Poisson obtained in § 99]

$$2(2g)^\mu P = \frac{\Gamma - \Gamma_1}{\mu!}, \tag{99.1}$$

$$\Gamma = \pm(\gamma + \mu g + \varepsilon)^\mu \mp \mu(\gamma + \mu g - 2g + \varepsilon)^\mu \pm$$
$$C_\mu^2(\gamma + \mu g - 4g + \varepsilon)^\mu \mp C_\mu^3(\gamma + \mu g - 6g + \varepsilon)^\mu + \dots,$$

with $\Gamma_1$ equal to $\Gamma$ with a changed sign of $\varepsilon$. The equation (99.1), as was required, represents the value of $P$ in a finite form.

**100.** If, when there is only one observation, $\mu = 1$, $P$ is the probability that the value of A, which should by the hypothesis be contained within the given limits $a$ and $b$ or $h - g$ and $h + g$, will after the observation be within limits $c - \varepsilon$ and $c + \varepsilon$. If these latter include the former, the equality $P = 1$ should be satisfied. If, however, the former include the latter, $P$ ought to be the ratio of the intervals of the latter, $2\varepsilon$, to the former, $2g$. Then, if both the latter limits are beyond the interval of the former, then necessarily $P = 0$; if $c - \varepsilon$ is within the interval $h - g$ and $h + g$, and $c + \varepsilon$ is beyond it, then $P$ should be equal to the ratio of the difference $[(h + g) - (c - \varepsilon)]$ to the interval $2g$. Finally, if $c + \varepsilon$ is within the interval $h - g$, $h + g$, and $c - \varepsilon$ is beyond it, $P$ should be equal to the ratio of the difference $[(c + \varepsilon) - (h - g)]$ to the same interval. Here are these values of $P$:

$$P = 1, \quad P = \frac{\varepsilon}{g}, \quad P = 0, \quad P = \frac{h + g - c + \varepsilon}{2g}, \quad P = \frac{c + \varepsilon - h + g}{2g}.$$

They are derived from equation (99.1), which at $\mu = 1$ becomes

$$P = \frac{1}{4g}(\Gamma - \Gamma_1).$$

In addition, $\gamma = h - c$, so that

$$\Gamma = \pm(h) + g - c + \quad \mp (h) - g - c + \tag{100.1a}$$
$$\Gamma_1 = \pm(h) + g - c - \quad \mp (h) - g - c - \tag{100.1b}$$

In the *first* of the five cases $c + \varepsilon > h + g$ and $c - \varepsilon < h - g$. In formulas (100.1a, b) the magnitudes in brackets are positive and negative respectively. Therefore, in the first formula we ought to choose the superior signs, and inferior signs in the second.

In the *second* case $h + g > c + \varepsilon$, $h - g < c - \varepsilon$. Superior signs ought to precede the first terms of both formulas, and inferior signs, to precede the second terms.

In the *third* case $h - g > c + \varepsilon$ and we ought to choose superior signs in both formulas. Here, however, the condition $h + g < c - \varepsilon$ is also possible, and then we ought to choose the inferior signs, and $P = 0$.

In the *fourth* case $c - \varepsilon > h - g$, $c - \varepsilon < h + g$, $c + \varepsilon > h + g$. We ought to choose the inferior signs in the second formula and, in the first formula, the superior sign before the first term and the inferior sign before the second.

In the *fifth* case $c - \varepsilon < h - g$, $c + \varepsilon > h - g$, $c + \varepsilon < h + g$. We ought to choose the superior signs in the first formula and the superior sign before the first term of the second formula and the inferior sign before its second term.

Here are the results of all those cases.

1. $\Gamma = 2g$, $\Gamma_1 = -2g$, $P = 1$.
2. $\Gamma = 2h - 2c + 2\varepsilon$, $\Gamma_1 = 2h - 2c - 2\varepsilon$, $P = \varepsilon/g$.
3. $\Gamma = 2g$, $\Gamma_1 = 2g$, $P = 0$.
4. $\Gamma = 2h - 2c + 2\varepsilon$, $\Gamma_1 = -2g$, $P = \dfrac{h + g - c + \varepsilon}{2g}$.
5. $\Gamma = 2g$, $\Gamma_1 = 2h - 2c - 2\varepsilon$, $P = \dfrac{c + \varepsilon - h + g}{2g}$.

The values of $P$ when only one observation is available can also be verified by formula (97.1) for the general case. When considering $f_1 z$ as a discontinuous function disappearing at all values of $z$ beyond the given limits $a$ and $b$, the probability $P$ that the value of A should be contained within the limits $c - \varepsilon$ and $c + \varepsilon$ will evidently be

$$P = \int_{c-\varepsilon}^{c+\varepsilon} f_1 z \, dz.$$

At $\mu = 1$, by formulas (97.2) and (97.1),

$$X = \int_a^b e^{xz\sqrt{-1}} f_1 z \, dz,$$

and [...]

$$P = \frac{1}{\pi} \int_a^b \left( \int_0^\infty \frac{\sin[(\varepsilon + \ )\mp\ x}{x} dx - \int_0^\infty \frac{\sin[(\ \varepsilon - \ )\mp\ x}{x} dx \right) f_1 z \, dz.$$

However, $\int_0^\infty \frac{\sin\gamma x}{x}dx = \pm\frac{\pi}{2}$

for a positive and negative constant γ. Therefore, when the signs of $c + \varepsilon - z$ and $c - \varepsilon - z$ coincide or are contrary, the indicated difference will be equal to 0 or π, and the integral with variable $z$ will disappear at all values of $z$ either larger than $c + \varepsilon$, or smaller than $c - \varepsilon$. It ought only to extend between $a$ and $b$, and between $c - \varepsilon$ and $c + \varepsilon$. And since we suppose that $f_1 z$ is zero beyond the first limits, the value of $P$ is reduced to the integral of $f_1 z$ in the limits from $z = c - \varepsilon$ to $c + \varepsilon$, which was indeed required to verify.

**101.** If μ is a very large number, we can apply transformations similar to those made in § 95 for replacing formula (97.1) by another one which will establish an approximate value of $P$. We note first of all that it is possible to write formula (97.2) in the form

$$X = \int_a^b e^{xz_1\sqrt{-1}} f_1 z_1 dz_1 \int_a^b e^{xz_2\sqrt{-1}} f_2 z_2 dz_2 \cdots \int_a^b e^{xz_\mu\sqrt{-1}} f_\mu z \, dz.$$

[After transformations]

$$X = Ye^{y\sqrt{-1}},$$
$$\rho_1 \rho_2 \ldots \rho_\mu = Y, \quad r_1 + r_2 + \ldots + r_\mu = y.$$

Substituting $X$ in formula (97.1), we get […]

$$P = \frac{2}{\pi} \int_{-\infty}^\infty Y \cos(y - cx) \sin\varepsilon x \frac{dx}{x}. \tag{101.1}$$

[…] When μ is a very large number and very small values of $x$ are excluded, the product $Y$, equal to unity at $x = 0$, in general becomes a very small fraction and disappears if μ can be infinite. Just like in § 95, without considering the special case[4], in which $Y$ approaches a non-zero magnitude, we assign only very small values to $x$, see the integral in formula (101.1). Just beyond these values $Y$ becomes insignificant, so that if $Y = \exp(-\theta^2)$, the variable θ can be supposed infinite there.

Therefore, when replacing the variable $x$ by θ, the limits of the integral ought to be $\theta = 0$ and ∞. For representing $x$ and $dx$ through θ and $d\theta$, I expand the previous expressions[5] of $\rho_n \cos r_n$ and $\rho_n \sin r_n$ in powers of $x$ and replace $z_n$ by $z$ in the integral. Assuming

$$\int_a^b z f_n z \, dz = k_n, \quad \int_a^b z^2 f_n z \, dz = k'_n, \quad \int_a^b z^3 f_n z \, dz = k''_n, \ldots,$$

we get converging series

$$\rho_n \cos r_n = 1 - \frac{x^2}{2!} k'_n + \frac{x^4}{4!} k'''_n - \ldots, \quad \rho_n \sin r_n = x k_n - \frac{x^3}{3!} k''_n + \ldots$$

Then, I introduce

$(k'_n - k^2_n)/2 = h_n$, $(k''_n - 3k_n k'_n + 2k_n^3)/6 = g_n$, ... [...]
$\sum k_n = \mu k$, $\sum h_n = \mu h$, $\sum g_n = \mu g$, $\sum(l_n - h_n^2/2) = \mu l$, ...

Here and below the sums extend from $n = 1$ to $\mu$.

The following formula takes place

$\ln Y = -\theta^2 = -x^2 \mu h + x^4 \mu l - ...,$

so that

$$x = \frac{\theta}{\sqrt{\mu h}} + \frac{l\theta^3}{2\mu h^2 \sqrt{\mu h}} + ..., \quad \frac{dx}{x} = \frac{d\theta}{\theta} + \frac{l\theta d\theta}{\mu h^2} + ...$$

[After transformations it occurs that]

$$P = \frac{2\theta}{\pi} \int_0^\infty \exp(\theta^2)\cos[(\mu k - c)\sin\epsilon \ x\frac{d}{\theta} +$$

$$\frac{2g}{\pi h\sqrt{\mu h}} \int_0^\infty \exp(\theta^2)\sin[(\mu k - c)\sin\epsilon \ x\theta^2 d\theta, \quad (101.2)$$

and

$$P = 1 - \frac{2}{\sqrt{\pi}} \int_u^\infty \exp(-t^2) dt \quad (101.3)$$

expresses the probability that, having a very large number $\mu$ of trials, the sum $s$ of the values of A is contained within the limits

$\mu k \mp 2u\sqrt{\mu h}, \ k = c/\mu.$ \quad (101.4)

Dividing that sum by $\mu$, we can determine the corresponding probability for the mean $s/\mu$.

**102.** Even if $u$ is insignificant, the probability (101.3) will very little differ from unity. We conclude therefore that the ratio $s/\mu$ probably very little differs from $k$. That magnitude is the sum of the possible values of A multiplied by their chances at each corresponding trial and divided by the number $\mu$ of these trials, i. e., the sum of the indicated values, multiplied by their corresponding mean chances; our conclusion coincides with the proposition of § 53, which is thus proved in all generality.

And so, given a very large number $\mu$ of trials, there invariably exists a probability, very close to certainty, that the mean value of A very little differs from $k$. The difference $(s/\mu - k)$ indefinitely decreases with the increase of $\mu$ and will become exactly zero if that number is infinite.

A plane curve with $z$ and $f_n z$ as its current coordinates represents the law of probabilities of the values of A at the $n$-th trial, so that the element $f_n z dz$ of the area of that curve will be the infinitely low probability of the value of A expressed by the abscissa $z$. A curve with current coordinates $z$ and $(1/\mu)\sum f_n z$ will express the law of probabilities of the mean value of A in a series of $\mu$ trials. Like the complete area of that curve from $a$ до $b$ equal to unity, the integral

$$\int_a^b f_n z dz = 1. \qquad (102.1)$$

Denote by $\varsigma$ the abscissa of its centre of gravity, then

$$\varsigma = \frac{1}{\mu}\sum \int_a^b z f_n z dz =$$

It is indeed equal to $k$, to which the mean value of A invariably converges. It disappears if at each trial the values of A equal in magnitude and contrary in sign are equally probable, that is, if $f_n(-z) = f_n z$ for all values of $n$ and $z$.

The constant $h$ should be positive.
[Poisson next proves the equality]

$$4h_n = \int_a^b \int_a^b (z - z')^2 f_n z f_n z' dz dz'.$$

The magnitude $4h_n$ is obviously positive and can not disappear since all the elements of the double integral are positive. The same holds for $\sum h_n$ and $h$.

In the simplest case all the possible values of A during the trials remain equally probable. Then the equality $f_n z = 1/(b - a)$ will take place at any $n$ and therefore[6]

$$k_n = k = (a + b)/2, \quad h_n = h = (a^2 + ab + b^2)/6 - (a + b)^2/8$$

and the limits of $s/\mu$ corresponding to probability $P$ will be

$$\frac{1}{2}(a+b) \mp \frac{u(b-a)}{\sqrt{6\mu}}, \text{ or } \mp \frac{2ub}{\sqrt{6\mu}} \text{ if } a = -b. \qquad (102.2a, b)$$

Let for example (§ 82) $u = 0.4765$, then the mean $s/\mu$ will with the same probability be within or beyond the limits $0.389b/\sqrt{\mu}$. For $\mu = 600$ we can bet even money on $s/\mu$ not to deviate from zero more than by $0.4765b/3 \cdot 10 \approx 0.016b$.

Such is the case in which at each trial point M with equal probability arrives anywhere on the segment of length $2b$. For a very large number $\mu$ of trials the mean distance of M from the segment's midpoint has probability $P$ of not exceeding the fraction $2u/\sqrt{6\mu}$ of $b$.

And if at each trial the point M arrives on a circle of radius $b$ with its equal distances from the circle's centre being equally probable, then, evidently, the probability $f_n z dz$ of distance $z$ will be proportional to that $z$. Suppose that that function is constant at all trials and note that all the possible distances are contained within 0 and $b$. Then, for satisfying condition (102.1), it will be necessary to assume that $f_n z = 2z/b^2$. Consequently[7],

$$k_n = k = 2b/3, \quad 2h_n = 2h = b^2/2 - 4b^2/9$$

and $P$ will be the probability that for $\mu$ trials the mean distance of M from the centre of the circle will be contained within the limits

$$\frac{2b}{3} \mp \frac{ub}{3\sqrt{\mu}}.$$

**103.** We (§ 97) supposed that the thing A can take all, even if not equally probable values between $a$ and $b$, but the formulas derived above are just as applicable to the case in which the number of those possible values is restricted. For proving this, it is sufficient to assume that the functions $f_1 z, f_2 z, \ldots$, which express the laws of the probabilities of the values of A in $\mu$ successive trials, are discontinuous[8].

Suppose that $c_1, c_2, \ldots, c_v$ are $v$ values of $z$ contained within $a$ and $b$, and that the function $f_n z$ disappears at any value of $z$, not infinitely little deviating from one of the magnitudes $c_1, c_2, \ldots, c_v$. Denote an infinitely small magnitude by $\delta$ and suppose also that

$$\int_{c_1-\delta}^{c_1+\delta} f_n z\, dz = \gamma_1, \quad \int_{c_2-\delta}^{c_2+\delta} f_n z\, dz = \gamma_2, \ldots, \quad \int_{c_v-\delta}^{c_v+\delta} f_n z\, dz = \gamma_v.$$

The thing A will therefore only take $v$ given values $c_1, c_2, \ldots, c_v$, whose probabilities at the $n$-th trial are $\gamma_1, \gamma_2, \ldots, \gamma_v$ and can vary from one trial to another, that is, vary with $n$. However, one of those values certainly takes place at the $n$-th trial, and the equality

$$\gamma_1 + \gamma_2 + \ldots + \gamma_v = 1$$

should hold for all values of $n$ from 1 to $\mu$.

In addition, that sum is equal to the integral of $f_n z$ over $[a, b]$ and the derived equation replaces the condition (102.1). [After transformations it occurs that]

$$\int_a^b z f_n z\, dz = \gamma_1 c_1 + \gamma_2 c_2 + \ldots + \gamma_v c_v,$$

$$\int_a^b z^2 f_n z\, dz = \gamma_1 c_1^2 + \gamma_2 c_2^2 + \ldots + \gamma_v c_v^2,$$

so that the magnitudes *k* and *h* from § 101 become

$$k = \frac{1}{\mu}\sum(\gamma_1 c_1 + \gamma_2 c_{2v} + \ldots + \gamma_v c_v),$$

$$h = \frac{1}{2\mu}\sum[(\gamma_1 c_1^2 + \gamma_2 c_{2v}^2 + \ldots + \gamma_v c_v^2) - (\gamma_1 c_1 + \gamma_2 c_{2v} + \ldots + \gamma_v c_v)^2].$$

The sums are extended over the µ trials and the formula (101.3) will express the probability that the sum *s* of the values of A in that series of trials will be contained within the limits (101.4), where *k* и *h* should be their determined values which are easy to calculate if the *v* possible values of A and their probabilities are given for each trial.

And if these probabilities are constant, and moreover, equal one to another, their common value will be $1/v$ and then simply

$$k = (c_1 + c_2 + \ldots + c_v)/v,$$
$$h = [v(c_1^2 + c_2^2 + \ldots + c_v^2) - (c_1 + c_2 + \ldots + c_v)^2]/2v^2.$$

Suppose that the possible values of A are the 6 numbers indicated on the faces of a usual die thrown a very large number µ of times successively. Neglecting the possible small differences between the chances of those faces, we have

$v = 6, c_1 = 1, c_2 = 2, c_3 = 3, c_4 = 4, c_5 = 5, c_6 = 6,$
$k = 7/2, h = 35/24,$

and the formula (101.3) expresses the probability that the sum *s* of the numbers arrived after those µ trials will be contained within the limits

$$\frac{1}{2}\left(7\mu \mp u\sqrt{\frac{70\mu}{3}}\right).$$

For $u = 0.4765$ and $\mu = 100$, the sum *s* will with equal probabilities be contained within or beyond the limits $350 \pm 11.5$.

**104.** Consider now, like in § 52, an event E of some nature which can only arrive due to *v* different and incompatible causes $C_1, C_2, \ldots, C_v$. Suppose that cause $C_i$ provides chance $c_i$ to the occurrence of E if it acted, and that $\gamma_i$ is the probability of that action. The chance of E can therefore vary from one trial to another and take *v* different values $c_1, c_2, \ldots, c_v$ if only their probabilities $\gamma_1, \gamma_2, \ldots, \gamma_v$ and the causes $C_1, C_2, \ldots, C_v$ do not change. Assuming such a chance of E, there will exist probability *P* (101.3) that for a very large number µ of trials its mean value will be contained within the limits determined by (101.4). There, *k* and *h* take their first values of § 103 provided that $C_1, C_2, \ldots$ and $\gamma_1, \gamma_2, \ldots$ remain constant during the trials. Therefore,

$$k = \gamma_1 c_1 + \gamma_2 c_2 + \ldots + \gamma_v c_v, \qquad (104.1)$$
$$h = [(\gamma_1 c_1^2 + \gamma_2 c_2^2 + \ldots + \gamma_v c_v^2) - (\gamma_1 c_1 + \gamma_2 c_2 + \ldots + \gamma_v c_v)^2]/2$$

and it is seen that they are independent from μ, whichever is the number of the magnitudes included in these formulas or how unequal they are.

With a small value of *u* the probability *P*, that the mean chance of the event E probably very little differs from the sum of the products in formula (104.1), indefinitely approaches that sum as μ increases further, and will be very close to certainty. In this conclusion consisted the second of the two general propositions of § 52, which we still had to prove.

In two series of very large numbers of μ and μ′ trials in which the event E arrived *m* and *m*′ times, the ratios *m*/μ and *m*′/μ′ probably very little differ from the corresponding mean chances of E (§ 96). It is therefore likely that they very little differ from the previous value of *k* (104.1) and therefore from each other since the values of *k* in those two series coincide if only the causes $C_1$, $C_2$, … did not change in the interval between those series. But how probable is a small given difference between *m*/μ and *m*′/μ′? We will consider this important problem below.

**105.** In most problems to which the formula (101.3) is applicable, the law of probabilities of the values of A is unknown, and the magnitudes *k* и *h*, included in the limits of the mean value of A can not be determined in advance. However, a long series of trials can serve for eliminating those unknowns, which are included in the limits of the mean value of A, from another such series, also consisting of a large number of trials and subjected to the action of the same causes attaching the same chance to each value of A and having the same probability themselves. A complete solution of this problem is the subject of the calculations below.

I suppose that in formula (101.2) *c* = ε […] and determine under that condition the probability that the sum *s* of the values of A in μ trials is contained within the limits of 0 and 2ε. Then, the derivative of *P* with respect to ε […] expresses the infinitely low probability that exactly *s* = 2ε. Let also

$$2\varepsilon = \mu k + 2v\sqrt{\mu h}, \quad d\varepsilon = \sqrt{\mu h}\, dv.$$

Denote by *wdv* the corresponding value of (*dP*/*d*ε)*d*ε, where magnitudes of the order of smallness of 1/μ are neglected. Then *x* can be reduced to the first term $\theta/\sqrt{\mu h}$ of its expression in § 101. And [after transformations]

$$\frac{s}{\mu} = k + \frac{2v\sqrt{h}}{\sqrt{\mu}},$$

where[9] *v* is positive or negative, but very small as compared with √μ.

Denote by $C_1$, $C_2$, …, $C_v$ all the mutually incompatible causes, whether known or not, which are able to provide the thing A one of its possible values, and by $\gamma_1$, $\gamma_2$, …, $\gamma_v$ their respective probabilities whose sum is unity. If the number of the causes is infinite, each of those probabilities is infinitely low. All the values of A are contained between *a* and *b*; and if their number is infinite, the chance of each provided by each of the indicated causes will be infinitely low. Denote by $Z_i dz$ the chance attached by cause $C_i$, if certain, to the value $z$ of the thing A. The integral

$$\int_a^b z f_n z \, dz, \qquad (105.1)$$

concerning the *n*-th trial can take $v$ values corresponding to the integrals of $zZ_1, zZ_2, \ldots, zZ_v$ over $[a, b]$ with the probabilities of these values being the probabilities of their causes. For some trial the probability $\gamma_i$ expresses the chance of the integral of $zZ_i$. The infinitely low probability of the mean value

$$\frac{1}{\mu} \sum \int_a^b z f_n z \, dz$$

is therefore determined by the preceding rule about the mean value $s/\mu$ of a thing in a very large number $\mu$ of trials: $s$ is the sum of $\mu$ unknown values of integral (105.1) in that series of trials, and $k$ and $h$ should be determined by its $v$ possible values.

Adopt the $v$ values of the integrals of $zZ_1, zZ_2, \ldots$ as $c_1, c_2, \ldots, c_v$ of § 103 and introduce for the sake of brevity

$$\gamma = \sum \gamma_i \int_a^b zZ_i \, dz,$$

$$\beta = \frac{1}{2} \sum \gamma_i \left( \int_a^b zZ_i \, dz \right)^2 - \frac{1}{2} \left( \sum \gamma_i \int_a^b zZ_i \, dz \right)^2, \qquad (105.2)$$

where[10] the sums extend from $i = 1$ to $v$. By the formulas of that section, these $\gamma$ and $\beta$ independent from $\mu$ should indeed be assumed as $k$ and $h$. Denote now by $v_1$ a positive or negative magnitude, very small as compared with $\sqrt{\mu}$, and by $V_1$, a polynomial only containing odd powers of $v_1$, and introduce an infinitely small magnitude $w_1 dv_1$

$$w_1 dv_1 = \frac{1}{\sqrt{\pi}} [1 - \frac{V_1}{\sqrt{\mu}}] \exp(-v_1^2) dv_1,$$

which is the probability of the equation

$$\frac{1}{\mu} \sum \int_a^b z f_n z \, dz = \gamma + \frac{2\beta_1 \sqrt{\phantom{x}}}{\sqrt{\mu}}. \qquad (105.3)$$

Consider also the magnitude

$$\frac{1}{2}\int_a^b z^2 f_n z dz - \frac{1}{2}\left(\int_a^b z f_n z dz\right)^2$$

as a thing that can take $v$ values in accordance with causes $C_1$, $C_2$, …, $C_v$ with probabilities at each trial coinciding with the probabilities of these causes. In addition, denote by $\tilde{v}$ such a positive or negative magnitude that $\tilde{v}/\sqrt{\mu}$ will be a very small fraction, and by $\tilde{V}$, a polynomial only containing odd powers of $\tilde{v}$, introduce magnitude

$$\tilde{w}d\tilde{v} = \frac{1}{\sqrt{\pi}}[1 - \frac{\tilde{V}}{\sqrt{\mu}}]\exp(-\tilde{v}^2)d\tilde{v}$$

and, for the sake of brevity, let

$$\alpha = \frac{1}{2}\sum \gamma_i \int_a^b z^2 Z_i dz - \frac{1}{2}\sum \gamma_i \left(\int_a^b z Z_i dz\right)^2.$$

The expression $\tilde{w}d\tilde{v}$ will be the probability that the mean of $\mu$ values of the magnitude under consideration,

$$\frac{1}{2\mu}\sum \int_a^b z^2 f_n z dz - \frac{1}{2\mu}\sum \left(\int_a^b z f_n z dz\right)^2,$$

will only differ from $\alpha$ by a definite magnitude of the order of smallness of $1/\sqrt{\mu}$, whose determination is meaningless. This mean, however, simply coincides with $h$ from § 101. And, when neglecting magnitudes of the order of $1/\mu$, it will be sufficient to assume $\alpha$ instead of $h$ in the second term of the previous value of $s/\mu$, already of order $1/\sqrt{\mu}$. Thus,

$$\frac{s}{\mu} = k + \frac{2\alpha\sqrt{\phantom{x}}}{\sqrt{\mu}},$$

and, if only the assumed value of $h$ is certain, the probability of this equality is again $wdv$. However, this value only has probability $\tilde{w}d\tilde{v}$, depending on the variable $\tilde{v}$, which does not enter $s/\mu$. Therefore, the composite probability of that last value is the product of $wdv$ and the sum of the values of $\tilde{w}d\tilde{v}$, corresponding to all possible values of $\tilde{v}$. These values should be very small as compared with $\sqrt{\mu}$, but, because of the exponential factor of $\tilde{w}d\tilde{v}$, it is possible to extend its integral from $\tilde{v} = -\infty$ до $\infty$ without appreciably changing it. Its part depending on $\tilde{V}$ will disappear since it consists of pairs of elements equal in magnitude and contrary in sign, so that simply

$$\int_{-\infty}^{\infty} \tilde{w}d\tilde{v} = 1.$$

Therefore, the probability of the preceding equation will invariably be *wdv*, as though the assumed approximate value of *h* is certain. Note also that the mean making up the left side of equation (105.3) coincides with magnitude *k* from § 101, and $w_1 dv_1$ is the probability that

$$k = \gamma + \frac{2\beta_1 \sqrt{\ }}{\sqrt{\mu}}.$$

Substituting this value in *s*/μ, we get

$$\frac{s}{\mu} = \gamma + \frac{2\beta_1 \sqrt{\ }}{\sqrt{\mu}} + \frac{2\alpha \sqrt{\ }}{\sqrt{\mu}},$$

that is, the probability that for each pair of the values of *v* и $v_1$ this last equation is the product $wdv w_1 dv_1$. Neglecting the term with μ in the denominator, I denote it by σ:

$$\sigma = \frac{1}{\pi}[1 - \frac{1}{\sqrt{\mu}}(V + V_1)]\exp(-v^2 - v_1^2) dv dv_1.$$

In addition, denote by θ a positive or negative magnitude, very small, like *v* and $v_1$, as compared with √μ. Then we will be able to suppose that

$$v_1 \sqrt{\beta} + v\sqrt{\alpha} = \theta\sqrt{\alpha + \beta}.$$

If desired, $v_1$ и $dv_1$ can be replaced in the previous differential formula by this new variable:

$$v_1 = \frac{\theta\sqrt{\alpha + \beta}}{\sqrt{\beta}} - \frac{v\sqrt{\alpha}}{\sqrt{\beta}}, \quad dv_1 = \frac{\sqrt{\alpha + \beta}}{\sqrt{\beta}} d$$

[After transformations it occurred that the equation]

$$\frac{s}{\mu} = \gamma + \frac{2\theta\sqrt{\alpha + \beta}}{\sqrt{\mu}} \qquad (105.4)$$

only contains the variable θ so that its composite probability is the sum of the values of σ at all positive and negative values which it is possible to provide for the other variable, *v*. Furthermore, owing to the exponential included in σ, it will be possible to extend this integral (? - O.S.) from $v = -\infty$ to ∞ without appreciably changing its value.

And so […]

$$\eta d\theta = \frac{1}{\sqrt{\pi}}\exp(-\theta^2)d\theta - \frac{1}{\sqrt{\pi\mu}}\chi\exp(-\theta^2)d\theta$$

is the probability of equation (105.4); χ is a polynomial only containing odd powers of θ. It is required to eliminate the unknown (α + β) from that equation. This is possible since the expression of (α + β) is reduced to

$$\alpha + \beta = \frac{1}{2}\sum\gamma_i\int_a^b z^2 Z_i dz - \frac{1}{2}\left(\sum\gamma_i\int_a^b z Z_i dz\right)^2 \qquad (105.5)$$

and is independent from the sum $\sum\gamma_i[\int_a^b zZ_i dz]^2$ which is contained in both α and β. We can therefore calculate (α + β) independently from its second term.

**106.** [Poisson introduces]

$$\varphi = \frac{1}{\mu}\sum\int_a^b z^2 f_n z dz$$

[and proves that] there exists probability $\tilde{w}d\tilde{v}$ that

$$\frac{1}{2}\sum\gamma_i\int_a^b z^2 Z dz$$

only differs from φ/2 by a determined magnitude of the order of smallness of $1/\sqrt{\mu}$.

Moreover, if invariably neglecting terms containing μ in the denominator, then, like in § 105, it is possible to replace the first term of the expression α + β by φ/2, without changing at all the probability η*d*θ of this expression. The other term of the value of α + β is exactly $\gamma^2/2$, so that

$$\alpha + \beta = \varphi/2 - \gamma^2/2,$$

and the equation (105.4) becomes

$$\frac{s}{\mu} = \gamma + \frac{\theta\sqrt{2\varphi - 2\gamma^2}}{\sqrt{\mu}}.$$

Let now *Z* be a given function of *z*. The analysis in §§ 97 and 101, as also the expression *wdv* from § 105, are easily extended on the sum of the values of *Z* in μ trials under our consideration. Suffice it to replace A by another thing $A_1$ taking the same values as the function *Z*. An infinitely low probability of some value of $A_1$ is equal to that of the corresponding value of *z,* and at the *n*-th trial it will be $f_n z dz$. Denote by $k_1, h_1, g_1, \ldots$ the magnitudes now pertaining to $A_1$, which in

§ 101 concerned A and were denoted by *k, h, g,* …, then

$$\frac{s_1}{\mu} = k_1 + \frac{2\nu\sqrt{h_1}}{\sqrt{\mu}}$$

and if $Z = z^2$

$$k_1 = \frac{1}{\mu}\sum\int_a^b z^2 f_n z \, dz = \varphi.$$

In the preceding expression of $s/\mu$ we may assume to within our degree of approximation that $\varphi = s_1/\mu$. And we are assured, just like in § 105, that the probability of that expression does not change, and that $\eta d\theta$ remains the infinitely low probability of the equation

$$\frac{s}{\mu} = \gamma + \frac{\theta\sqrt{2s_1/\mu - 2\gamma^2}}{\sqrt{\mu}} \quad \text{or} \quad \frac{s}{\mu} = \gamma + \frac{\theta\sqrt{2s_1/\mu - 2s^2/\mu^2}}{\sqrt{\mu}}.$$

This can be derived from the previous when once more neglecting magnitudes of the order of smallness of $1/\mu$.

I denote by $\lambda_n$ the value of A at the *n*-th trial and introduce for the sake of brevity

$$\sum \lambda_n/\mu = \lambda, \quad \sum(\lambda_n - \lambda)^2/\mu = l^2/2.$$

The following equalities are identities:

$$\frac{s_1}{\mu} = \frac{\sum \lambda_n^2}{\mu}, \quad \frac{s}{\mu} = \frac{\sum \lambda_n}{\mu}, \quad \frac{s_1}{\mu} - \frac{s}{\mu} = \frac{\sum(\lambda_n - \lambda)^2}{\mu}$$

and the preceding equation becomes

$$\frac{s}{\mu} = \gamma + \frac{\theta l}{\sqrt{\mu}}.$$

Therefore, when denoting by *u* a given positive magnitude, the integral of the probability $\eta d\theta$ of this equation over $[-u, u]$ will express the probability that $s/\mu$ is contained within the limits $\gamma \mp ul/\sqrt{\mu}$. Denoting this probability by $\Gamma$ and taking account of the expression $\eta d\theta$, we get

$$\Gamma = \frac{1}{\sqrt{\pi}} \int_{-u}^{u} \exp(\theta)^2 \, d\theta - \frac{1}{\sqrt{\pi\mu}} \int_{-u}^{u} \chi \exp(\theta)^2 \, d\theta.$$

However, $\chi$ is a polynomial only containing odd powers of $\theta$, so the second integral disappears, and $\Gamma$ coincides with the probability *P* as derived in formula (101.3).

This formula thus expresses the probability that the limits $\mp u l/\sqrt{\mu}$, which after the trials do not anymore depend on any unknowns, contain the difference between the mean value $s/\mu$ of A and the special magnitude $\gamma$. This mean indefinitely approaches $\gamma$ and reaches it if $\mu$ becomes infinite with the causes $C_1, C_2, \ldots, C_v$ of the possible values of A remaining invariable.

**107.** Suppose that two series of large numbers of trials $\mu$ and $\mu'$ are made. Denote the sums of the values of A in these series by $s$ and $s'$, with $\lambda_n$ and $\lambda'_n$ being the values of A at the $n$-th trials. Let

$$\sum \lambda_n/\mu = \lambda, \ \sum(\lambda_n - \lambda)^2/\mu = l^2/2, \ \sum \lambda'_n/\mu' = \lambda', \ \sum(\lambda'_n - \lambda')^2/\mu' = l'^2/2.$$

The sums extend over all the trials of each series, i. e., from $n = 1$ to $\mu$ in the first series, and from $n = 1$ to $\mu'$ in the second. If the causes $C_1, C_2, \ldots, C_v$ did not change between the series, $\gamma$ (105.2) will also remain invariable. Denote by $\theta$ and $\theta'$ positive or negative variables, very small as compared with $\sqrt{\mu}$ and $\sqrt{\mu'}$. Equations concerning the mean values of A in these series will be[11]

$$\frac{s}{\mu} = \gamma + \frac{\theta l}{\sqrt{\mu}}, \ \frac{s'}{\mu'} = \gamma + \frac{\theta' l'}{\sqrt{\mu'}}, \qquad (107.1)$$

and their probabilities $\eta d\theta$ and $\eta' d\theta'$

$$\eta d\theta = \frac{1}{\sqrt{\pi}} \chi [1 - \frac{1}{\sqrt{\mu}}] \exp(-\theta^2) d\theta,$$

$$\eta' d\theta' = \frac{1}{\sqrt{\pi}} \chi [1 - \frac{1}{\sqrt{\mu'}}] \exp(-\theta'^2) d\theta'.$$

Here, $\chi$ and $\chi'$ are polynomials only containing odd powers of $\theta$ and $\theta'$. And if the series are composed of different trials, $s/\mu$ and $s'/\mu'$ can be considered as independent events. By the rule of § 5 the probability of their simultaneous occurrence is the product of $\eta d\theta$ and $\eta' d\theta'$. The same takes place for any combination of equations (107.1) and in particular for their difference. […] And so, neglecting terms with $\sqrt{\mu\mu'}$ in the denominator, we obtain the probability of that difference for each pair of the values of $\theta$ and $\theta'$

$$\psi \equiv \eta \eta' d\theta d\theta' = \frac{1}{\pi} \chi [1 - \frac{\chi}{\sqrt{\mu}} - \frac{\chi'}{\sqrt{\mu'}}] \exp(-\theta^2 - \theta'^2) d\theta d\theta'. \qquad (107.2)$$

Just like in § 105, I assume that

$$\frac{\theta' l'}{\sqrt{\mu'}} - \frac{\theta l}{\sqrt{\mu}} = \frac{t \sqrt{l'^2 \mu + l^2 \mu'}}{\sqrt{\mu\mu'}},$$

so that

$$\frac{s'}{\mu'} - \frac{s}{\mu} = \frac{t\sqrt{l'^2\mu + l^2\mu'}}{\sqrt{\mu\mu'}}.$$

[…] Then, I replace θ in equation (107.2) by a new variable $t$, so let

$$\theta' = \frac{t\sqrt{l'^2\mu + l^2\mu'}}{l'\sqrt{\mu}} + \frac{\theta l\sqrt{\mu'}}{l'\sqrt{\mu}}, \quad d\theta' = \frac{\sqrt{l'^2\mu + l^2\mu'}}{l'\sqrt{\mu}}dt.$$

Therefore

$$\psi = \frac{dt d\theta\sqrt{l'^2\mu + l^2\mu'}}{\pi l'\sqrt{\mu}}(1-\Pi)\exp(-\theta'^2 - t^2),$$

where $\Pi$ is a polynomial whose each term contains an odd power of $t$ or θ. The difference $s'/\mu' - s/\mu$ only contains the variable $t$, and its probability is equal to the integral of ψ over all the values of the other variable, θ. Because of the exponential included in ψ, that integral can be extended from $\theta = -\infty$ to $\infty$ without appreciably changing its value. Assuming therefore that $\theta' = t'$, $d\theta' = dt'$ and denoting by $\Pi'$ the appropriately changed $\Pi$, we get

$$\psi = \frac{1}{\pi}(1-\Pi')\exp(-t'^2 - t^2)dt'dt.$$

[After transformations it occurs that]

$$\mp\frac{u\sqrt{l'^2\mu + l^2\mu'}}{\sqrt{\mu\mu'}}$$

are the limits within which the difference $s'/\mu' - s/\mu$ is contained with probability

$$\Delta = \frac{2}{\sqrt{\pi}}\int_0^u \exp(-t^2)dt.$$

This probability coincides with the value of $P$ as derived in formula (101.3), so $P$ is the probability that the difference of the mean values of A in two long series of trials will be contained within limits including nothing unknown.

Assigning a value for $u$ sufficient for $P$ to differ very little from unity, and discovering that the difference mentioned is beyond the established limits, we can justifiably conclude that between the series of trials the causes $C_1, C_2, \ldots, C_v$ of the possible values of A did not remain invariable, so that the probabilities $\gamma_1, \gamma_2, \ldots, \gamma_v$ of these causes or chances which they provide for the different values of A underwent some changes.

In accord with what was said in § 106, each of the magnitudes *l* and *l′* should likely very little differ from the unknown magnitude $2\sqrt{\alpha+\beta}$, constant in both series of trials. It is therefore very probable that they very little differ from each other. Without appreciably changing either the limits indicated above or their probability, we can therefore assume that *l′* = *l*. In a future series of trials there will be probability *P*, determined by the formula (101.3), that the limits of the mean value *s′*/μ′ of A are

$$\frac{s}{\mu} \mp \frac{ul\sqrt{\mu+\mu'}}{\sqrt{\mu\mu'}}.$$

At each given value of *u* these limits only depend on the results of the accomplished first series of trials.

At one and the same value of *u*, that is, with the same level of probability, the interval between these limits is larger, in the ratio of $\sqrt{\mu+\mu'}$ to $\sqrt{\mu'}$ than for the difference γ − *s*/μ. If μ′ is a very large number as compared with the very large number μ, the intervals for both series almost coincide.

**108.** If two series of μ and μ′ trials are made for measuring the same thing by different instruments having equally probable errors of equal magnitudes and contrary signs, the mean values *s*/μ and *s′*/μ′ derived from those series will indefinitely tend to one and the same magnitude which will be the veritable value of A (§ 60)[12]. In this case, the unknown γ will therefore be the same in both series and the means *s*/μ and *s′*/μ′ will likely differ only very little.

However, the values of the unknown α + β in these series can very much differ one from another and then *l* and *l′* will be quite unequal. The values of these magnitudes are known and it is possible to require how best to combine the means *s*/μ and *s′*/μ′ for establishing the limits γ or the veritable value of A.

So I denote by *g* and *g′* indefinite magnitudes whose sum is unity and add up the equations (107.1), at first multiplying them by *g* and *g′* respectively:

$$\gamma = \frac{gs}{\mu} + \frac{g's'}{\mu'} - \frac{gl\theta}{\sqrt{\mu}} - \frac{g'l'\theta'}{\sqrt{\mu'}}.$$

In accord with the above, the probability of the derived equation at any pair of values of θ and θ′ is ψ, and from calculations similar to those just made it follows that *P*, as given by (101.3), expresses the probability that the unknown γ is contained within the limits

$$\frac{gs}{\mu} + \frac{g's'}{\mu'} \mp \frac{u\sqrt{g'^2l'^2\mu + g^2l^2\mu'}}{\sqrt{\mu\mu'}}.$$

If desired that for the same probability *P*, i. e., for each given value of *u*, the interval between these limits becomes as short as possible, *g* and *g'* should be determined by equating to zero the derivative of the coefficient of *u* with respect to these magnitudes. Since $g + g' = 1$ and $dg' = - dg$, it turns out that

$$g = \frac{l'^2 \mu}{l'^2 \mu + l^2 \mu'}, \quad g' = \frac{l^2 \mu'}{l'^2 \mu + l^2 \mu'},$$

and the tightest limits of γ will be

$$\frac{sl'^2 + s'l^2}{l'^2 \mu + l^2 \mu'} \mp \frac{ull'}{\sqrt{l'^2 \mu + l^2 \mu'}},$$

with formula (101.3) as previously indicating their probability.

It is easy to generalize this result on any number of series of numerous observations made by different instruments for measuring one and the same thing A. Suppose that the magnitudes μ, *s*, *l* concerning the first series are denoted by μ′, *s*′, *l*′, μ″, *s*″, *l*″, … in the second, the third, … series. Suppose also that

$$\frac{\mu}{l^2} + \frac{\mu'}{l'^2} + \frac{\mu''}{l''^2} + \ldots = D^2, \quad \frac{\mu}{D^2 l^2} = q, \quad \frac{\mu'}{D^2 l'^2} = q', \quad \frac{\mu''}{D^2 l''^2} = q'', \ldots$$

Formula (101.3) will then express the probability that the unknown value of A is contained within the limits

$$\frac{sq}{\mu} + \frac{s'q'}{\mu'} + \frac{s''q''}{\mu''} + \ldots \mp \frac{u}{D},$$

derived from the most favourable combination of the observations. With a small value of *u*, the magnitude *P*, see formula (101.3), can be rendered very close to unity, so that the value of A will likely very little differ from the sum of the means *s*/μ, *s'*/μ′, *s″*/μ″, …, multiplied by *q*, *q'*, *q″*, … respectively.

The result of each series of observations influences that approximate value of A and the interval between its limits ∓*u*/*D* the more, the larger is the corresponding ratio μ/*l*², or μ′/*l*′², or μ″/*l*″², … If all the series are made with the same instrument, they can be considered as a single series of μ + μ′ + μ″ + … observations. And, as stated above, *l*, *l*′, *l*″, … will likely be almost equal one to another. Suppose that the sums extend over the entire single series from *n* = 1 to μ + μ′ + μ″ + … and introduce

$$\frac{\sum \lambda_n}{\mu + \mu' + \mu'' + \ldots} = \lambda, \quad \frac{\sum (\lambda_n - \lambda)^2}{\mu + \mu' + \mu'' + \ldots} = \frac{l_1^2}{2}.$$

Then we can consider $l_1$ as the common value of $l, l', l'', \ldots$ The previous limits of the unknown $\gamma$, whose probability indicates formula (101.3), become

$$\frac{s+s'+s''+\ldots}{\mu+\mu'+\mu''+\ldots} \mp \frac{ul_1}{\sqrt{\mu+\mu'+\mu''+\ldots}},$$

and coincide with the result of § 106 concerning one single series of trials.

**109**. The problem formulated at the end of § 104 is solved by the method similar to the just applied. Suppose that in a very large number $\mu$ of trials an event E of some nature arrived $m$ times. Its chance is variable and equals $p_n$ at the $n$-th trial. Let $\sum p_n/\mu = p$ and $\sum p_n^2 = q$ and denote by $v$ a positive or negative magnitude, very small as compared with $\sqrt{\mu}$, and suppose that $U$ is the probability of the equation

$$\frac{m}{\mu} = p - \frac{v\sqrt{2p-2q}}{\sqrt{\mu}}.$$

For simplifying, neglect the second term of formula (95.1), take into account the included there magnitude $k$ and replace $\theta$ by $v$, then

$$U = \frac{\exp(-v^2)}{\sqrt{2\pi\mu(p-q)}}.$$

Just like in § 104, denote all the possible causes of the event E, whose number can be either finite or infinite, by $C_1, C_2, \ldots, C_\nu$ with $\gamma_1, \gamma_2, \ldots, \gamma_\nu$ being their probabilities and $c_1, c_2, \ldots, c_\nu$, the chances they provide to the occurrence of E. Supposing that $p_n$ can take these $\nu$ values whose probabilities are $\gamma_1, \gamma_2, \ldots, \gamma_\nu$, we introduce

$$r = \gamma_1 c_1 + \gamma_2 c_2 + \ldots + \gamma_\nu c_\nu, \; \rho = \gamma_1 c_1^2 + \gamma_2 c_2^2 + \ldots + \gamma_\nu c_\nu^2.$$

Denote by $v_1$ a positive or negative variable, very small as compared with $\sqrt{\mu}$. Then the infinitely low probability of the exact equality

$$p = r + \frac{v_1\sqrt{2\rho - 2r^2}}{\sqrt{\mu}} \qquad (109.1)$$

will be $w_1 dv_1$ (§ 105) or simply $(1/\sqrt{\pi})\exp(-v_1^2)dv_1$ when neglecting the second term of its expression. Denote also by $\tilde{v}$ a variable, very small as compared with $\sqrt{\mu}$. Then the probability that $(p - q)$ only differs from $(r - \rho)$ by a definite magnitude proportional to $\tilde{v}$ and of the order of smallness of $1/\sqrt{\mu}$, is $\tilde{w}d\tilde{v}$, or simply $(1/\sqrt{\pi})\exp(-\tilde{v}^2)d\tilde{v}$, see § 105.

Then, it is seen that, when neglecting magnitudes of the order of $1/\mu$, it is possible to replace $(p - q)$ by $(r - \rho)$ without changing the probability $U$ of the previous value of the ratio $m/\mu$. That ratio then becomes

$$\frac{m}{\mu} = p - \frac{v\sqrt{2r-2\rho}}{\sqrt{\mu}}. \qquad (109.2)$$

Introduce now

$$\delta = \frac{1}{\sqrt{2\mu(r-\rho)}}.$$

Since $m$ is a natural number, $v$ can only be a positive or negative multiple of $\delta$, very small as compared with $1/\sqrt{\mu}$.

It follows from formulas (109.1) and (109.2) that

$$\frac{m}{\mu} = r + \frac{v_1\sqrt{2\rho-2r^2}}{\sqrt{\mu}} - \frac{v\sqrt{2r-2\rho}}{\sqrt{\mu}}.$$

For any pair of values of $v$ and $v_1$ the probability of that equation is the product of $U$ and $(1/\sqrt{\pi})\exp(-v_1^2)dv_1$, and I denote it by $\varepsilon$. When $(p-q)$ in the expression of $U$ is replaced by $(r-\rho)$, this product becomes

$$\varepsilon = \frac{1}{\pi\sqrt{2\mu(r-\rho)}} \exp(-v^2 - v_1^2)dv_1.$$

Let

$$v_1 = \theta\sqrt{\frac{r-r^2}{\rho-r^2}} + v\sqrt{\frac{r-\rho}{\rho-r^2}}, \quad dv_1 = \sqrt{\frac{r-r^2}{\rho-r^2}}d\theta,$$

then, neglecting terms of the order of smallness of $1/\mu$,

$$\frac{m}{\mu} = r + \frac{\theta\sqrt{2r-2r^2}}{\sqrt{\mu}}, \quad r = \frac{m}{\mu} - \frac{\theta\sqrt{2m(\mu-m)}}{\mu\sqrt{\mu}}.$$

[Poisson next writes down the formula for $\varepsilon$, after including there $\delta$ and $\theta$.]

However, the expression for $r$, does not include $v$ and its probability is also independent from $v$ and is equal to the sum of the values of $\varepsilon$, corresponding to all those which $v$ can take and which ought to increase by increments $\delta$, whose multiple is $v$. Since $\delta$ is small, an approximate value of this sum can be derived by substituting $dv$ instead of $\delta$ in $\varepsilon$ and replacing the sum by an integral. The calculated value will be exact to within magnitudes of the order of $\delta$ or $1/\sqrt{\mu}$. The variable $v$ should be very small as compared with $\sqrt{\mu}$, but owing to the exponential included in $\varepsilon$, it is possible to extend that integral from

$v = -\infty$ to $\infty$ without appreciably changing its value. And when assuming that

$$v\sqrt{\frac{r-r^2}{\rho-r^2}} + \theta\sqrt{\frac{r-\rho}{\rho-r^2}} = \theta_1, \; d\theta_1 = \sqrt{\frac{r-r^2}{\rho-r^2}}dv,$$

the integral with variable $\theta_1$ will also have infinite limits. Denote now by $\chi d\theta$ the infinitely low probability of the expression for $r$, then

$$\chi d\theta = \frac{d\theta}{\pi}\exp(-\theta^2)\int_{-\infty}^{\infty}\exp(-\theta_1^2)d\theta_1 = \frac{1}{\sqrt{\pi}}\exp(-\theta^2)d\theta.$$

And so, $u$ is a positive and given magnitude, and the probability that the unknown value of $r$ will be within the limits

$$\frac{m}{\mu} \mp \frac{u\sqrt{2m(\mu-m)}}{\mu\sqrt{\mu}} \tag{109.3}$$

is

$$\int_{-u}^{u}\chi d\theta = \frac{2}{\sqrt{\pi}}\int_{0}^{u}\exp(-\theta^2)d\theta$$

and coincides with $P$ from formula (101.3). Thus, $P$ is the probability that the special magnitude $r$, to which the ratio $m/\mu$ indefinitely tends as the large number $\mu$ increases still further, and does not differ from that ratio more than by the magnitude contained in the limits which are indicated by formula (109.3) and do not include any unknowns.

Suppose that in the second series of a very large number $\mu'$ of trials the event E arrived $m'$ times. Denote by $\theta'$ a positive or negative variable, very small as compared with $\sqrt{\mu'}$. The infinitely low probability of the equation

$$r = \frac{m'}{\mu'} - \frac{\theta'\sqrt{2m'(\mu'-m')}}{\mu'\sqrt{\mu'}}$$

is $(1/\sqrt{\pi})\exp(-\theta'^2)d\theta'$, and the probability of the equation

$$\frac{m'}{\mu'} - \frac{m}{\mu} = \frac{\theta'\sqrt{2m'(\mu'-m')}}{\mu'\sqrt{\mu'}} - \frac{\theta\sqrt{2m(\mu-m)}}{\mu\sqrt{\mu}},$$

derived for each pair of the values of $\theta$ and $\theta'$ [of $\theta'$ and $\theta$] by subtracting that value of $r$ from its previous value, equals the product of $(1/\sqrt{\pi})\exp(-\theta'^2)d\theta'$ and $(1/\sqrt{\pi})\exp(-\theta^2)d\theta$.

[Poisson replaces] the variable $\theta'$ by $t$ without changing $\theta$, then replaces $\theta$ by $t'$ without changing $t$. The probability of the previous equation becomes $(1/\pi)\exp(-t^2 - t'^2)dtdt'$, and the equation itself is now

$$\frac{m'}{\mu'} - \frac{m}{\mu} = \frac{t\sqrt{2\mu^3 m'(\mu' - m') + 2\mu'^3 m(\mu - m)}}{\mu\mu'\sqrt{\mu\mu'}}.$$

It only contains the variable $t$, and its composite probability is equal to the integral with variable $t'$ of the differential expression above. Without appreciably changing its value this integral can be extended from $t' = -\infty$ to $\infty$, and that expression will become $(1/\sqrt{\pi})\exp(-t^2)dt$, so that $P$, see formula (101.3), will express the probability that the difference $(m'/\mu' - m/\mu)$ is contained within the limits

$$\mp \frac{u\sqrt{2\mu^3 m'(\mu' - m') + 2\mu'^3 m(\mu - m)}}{\mu\mu'\sqrt{\mu\mu'}}.$$

Here $u$ is positive and given, and this formula only includes known numbers. The derived limits coincide with those obtained much simpler in § 87, although only for a constant chance of event E being the same in both series of trials. Formula (87.2) contains a term of the order of $1/\sqrt{\mu}$ or $1/\sqrt{\mu'}$, lacking in formula (101.3), since terms of that order of smallness were neglected in the probabilities under consideration.

**110.** In this work, I do not suppose to study the numerous problems to which the preceding formulas are applicable; I indicated such main problems in § 60 and subsequent sections[13]. As an example, I choose the known problem about the planetary and cometary orbits.
[After long transformations it occurred that Poisson's derived formula]

$$P = \frac{1}{\mu!}[\beta^\mu - \mu(\beta-1)^\mu + C_\mu^2(\beta-2)^\mu - C_\mu^3(\beta-3)^\mu + ...$$
$$- \alpha^\mu + \mu(\alpha-1)^\mu - C_\mu^2(\alpha-2)^\mu + C_\mu^3(\alpha-3)^\mu - ...], \quad (110.1)$$

coincides with Laplace's formula (Laplace 1812/1886, p. 261) obtained with the same aim but by quite another way.
[Poisson determined magnitudes α and β from equalities introduced above:]

$$h = g, \gamma = \mu g - c, c - \varepsilon = 2g\alpha, c + \varepsilon = 2g\beta.$$

If all the values of A within the limits from 0 to $2g$ are equally possible and impossible beyond them, formula (110.1) will express the probability that after some number μ of trials the sum of the values of a thing A is contained between $2g\alpha$ and $2g\beta$. The series containing β and α in formula (110.1) are extended until the magnitudes raised to the power of μ are still positive. If $n$ represents the maximal natural number included in β, the appropriate series stops at the $(n + 1)$-th term or earlier depending on the inequalities $\mu > n$ or $< n$. The same takes place for the second series for the maximal natural number included in α.

Whatever was the cause which determined the formation of planets, it is supposed that initially all possible inclinations from 0 to 90° of the planes of their orbits to the ecliptic were equally probable. And it is required to determine the probability that, under that hypothesis, the sum of the inclinations of the 10 known planets excluding the Earth[14] should be contained in given limits, for example between 0 и 90°. When assuming the inclination of a planet as a thing A to which formula (110.1) is applicable, 90° ought to be chosen as the interval $2g$ of the possible values of A, and, in the formula itself, $α = 0$, $β = 1$ and $μ = 10$, so that $P = 1/10!$.

This fraction is approximately a quarter of a millionth and therefore a sum of inclinations smaller than a right angle is absolutely unlikely. We should conclude that it certainly exceeds a right angle. However, it only amounts to approximately 82°, and, since it only experiences very small periodic variations, the hypothesis of equal probability of all inclinations at the formation of the planets is inadmissible. The unknown cause of their formation must have rendered the least inclinations much more probable.

The inclinations of the planets are supposed here to be independent from the direction of their movement; that is, from either the direction of the Earth's movement round the Sun or the contrary movement. If both these directions at the formation of the planets were equally probable, the probability of the movement of the 10 other planets in the same direction as the Earth would be 1/2 to the power of 10, i. e. lower than 1/1000. And an equal chance of both contrary movements is very unlikely and proves that the unknown cause of the formation of the planets must have rendered a high probability to their movement in one and the same direction.

Let now the thing A be the eccentricity of a planetary orbit and suppose that initially all its values from 0 to 1 were equally probable. Then the probability that, at $α = 0$, $β = 1.25$ and $μ = 11$, the sum of the known eccentricities should be contained, for example, between 0 and 5/4, will by formula (110.1) be $P = (1/11!)(1.25^{11} − 0.25^{11})$. This probability is lower than three millionth; on the contrary, it is extremely probable that the sum of the 11 eccentricities ought to exceed 1.25. And still, that sum, which only experiences small periodic variations, is a little smaller than 1.15. The hypothesis of equal probability of all possible values of A is thus inadmissible; the unknown cause of the formation of the planets undoubtedly rendered the least eccentricities, just as the least inclinations, much more probable.

**111.** For 138 comets of those observed since AD 240 [since BC 240], astronomers calculated parabolic elements as thoroughly as possible. 71 of them move directly and 67 have retrograde motions. The small difference between those numbers already proves that the unknown cause of the formation of the comets did not lead to their more probable motion in some common direction. The sum of the inclinations of those 138 comets to the ecliptic amounts to almost 6752°, which[15] exceeds 75 right angles almost by 2°. For establishing whether it should so little differ from that magnitude, had all the possible inclinations from 0 to 90° been equally probable, we ought to assume α and β in formula (110.1) little differing in either direction from 75. This, however, will prevent numerical calculations. For determining the probability *P* that, under the same hypothesis, the sum of the inclinations of the orbits of all the observed comets should be contained within given limits, we ought to turn to formula (101.3).

I assume the inclination of a cometary orbit to the plane of the ecliptic as the thing A. The limits of its possible values, usually denoted by *a* and *b*, will be *a* = 0 and *b* = 90°, with all these values supposed to be equally probable. Formula (101.3) will express the probability *P* that the mean of a large number μ of the observed inclinations is, in degrees (cf. § 102), between 45 $\mp 90u/\sqrt{6\mu}$. Assuming *u* = 1.92 and having μ = 138, we will get *P* = 0.99338 for the probability that, under the hypothesis of equal chances of all possible inclinations, the mean inclination of the 138 observed comets will not be beyond the limits of 45° ± 6°. We can bet almost 150 against 1 on that mean to be between 39° и 51°. Actually, it is 48°55′. We can not therefore believe that the unknown cause of the formation of the comets rendered their differing inclinations unequally probable.

Without introducing any hypotheses about the law of probabilities of these inclinations, formula (101.3) also expresses the probability that the mean inclination of a large number μ of comets observed in the future will only differ from the indicated value 48°55′ by a number of degrees contained within the limits (cf. § 107)

$$\mp \frac{ul\sqrt{138\mu + \mu'}}{\sqrt{138\mu'}}.$$

By the calculated inclinations of 138 comets, the nephew of Bouvard[16] established that the value of *l* included in these limits was 34°49′. Suppose for example that μ′ = μ and assume, as above, that *u* = 1.92, and it will be possible to bet 150 against 1 on the difference of the mean inclinations of the 138 observed comets and the same number of new comets to be included within ± 8°21′. The number μ′ of the existing comets is undoubtedly extremely large as compared with the comets whose orbits astronomers were able to calculate. After assuming μ′ as the number of the unknown comets, the preceding limits will become almost $\mp ul/\sqrt{138}$, narrower in the ratio 1: √2, than if μ′ = μ. Assuming as previously that *u* = 1.92, we will obtain probability almost equal to 150/151 that the difference mentioned above will be contained within the limits ± 5°42′.

Divide the observed comets in two groups, consisting of 69 earlier and another 69 more modern of them. The mean inclination in the first group will be 49°12′, and 48°38′ in the second, differing barely more than by half a degree. This example is really proper for proving that the mean values of one and the same thing correspond to each other even if the available numbers of observations are not extremely large and the observed values essentially differ one from another. Indeed, the maximal and least inclinations are here 1°41′ and 89°48′. The mean inclinations of the 71 and 67 comets with direct and retrograde motion are more different; the former is 47°3′, and the latter, 50°54′.

Erect in the northern hemisphere a perpendicular from the centre of the Sun to the plane of the ecliptic and it will cut the celestial sphere at the northern pole of the ecliptic. And if a perpendicular is erected from that centre in the same hemisphere to the plane of a cometary orbit, it will cut the celestial sphere at the northern pole of that orbit. The angular distance between those poles will be the inclination of that orbit to the ecliptic. However we should not confuse, as did the respected translator of Herschel (1834), the assumption that all points of the celestial sphere can with the same probability be the poles of the cometary orbits with the hypothesis of equal probability of all possible inclinations of the comets.

Let *a* and *b* be circular zones situated in the northern hemisphere of one and the same infinitely small width with a common centre at the northern pole of the ecliptic and with their angular distances from that pole being α and β. Denote by *p* and *q* the probabilities that a random point of that hemisphere belongs to zones *a* and *b*. The fractions *p* and *q* are evidently in the ratio *a*:*b*, and therefore in the ratio of the sines of angles α and β. By the hypothesis of equal possibility for all the points of the celestial sphere to be the poles of cometary orbits, *p* and *q* express the chances of the distances α and β of these poles from the ecliptic; in other words, the chances of the inclinations of the two cometary orbits to equal these distances.

And so, according to the adopted hypothesis, the chances of the various inclinations are proportional to the sines of the inclinations themselves but they are not equal to each other. The chance of an inclination of 90° becomes twice larger than of an inclination of 30°, and they both are infinitely large as compared with the chance of an infinitely small inclination[17].

**112.** For concluding this chapter, we adduce a summary of the formulas of probability proved here and in the previous chapter. The number μ of trials is supposed very large. It consists of two parts, *m* and *n*, also supposed very large. The formulas are the more exact the larger is the number μ, and they will become absolutely exact at infinite μ.

**112/1.** Let *p* and *q, p + q = 1*, be the chances of contrary events E and F, invariable during the entire series of trials. Denote by *U* the probability that in μ = *m* + *n* trials E arrives *m* times, a F, *n* times. According to § 69

$$U = \left(\frac{\mu p}{m}\right)^m \left(\frac{\mu q}{n}\right)^n \sqrt{\frac{\mu}{2\pi mn}}. \qquad (a)$$

In § 79 this formula was reduced to

$$U = \frac{\exp(-v^2)}{\sqrt{2\pi\mu pq}},$$

with

$$m = \mu p - v\sqrt{2\mu pq}, \ n = \mu q + v\sqrt{2\mu pq}.$$

Here $v$ is a positive or negative magnitude, very small as compared with $\sqrt{\mu}$. In this form, the formula above is equally applicable if the chances of E and F vary from one trial to another, it will only be necessary, in accord with formula (95.1), to assume as $p$ and $q$ their the mean values over all the series of $\mu$ successive trials.

**112/2.** Suppose that events E and F with unknown chances $p$ and $q$ arrived $m$ and $n$ times in $\mu$ trials. Denote by $U'$ the probability that they will occur $m'$ and $n'$ times in $\mu' = m' + n'$ future trials. The numbers $m'$ and $n'$ are proportional to $m$ and $n$, so that

$$m' = \mu'm/\mu, \ n' = \mu'n/\mu.$$

Formula (§ 71)

$$U' = \sqrt{\frac{\mu}{\mu + \mu'}}\ U_1 \qquad\qquad (b)$$

will take place whatever is $\mu'$. Here, $U'$ is the probability of a future event provided that $m/\mu$ and $n/\mu$ certainly are the chances of E and F, i. e., provided that

$$U' = C_{\mu'}^{m'}(m/\mu)^{m'}(n/\mu)^{n'}.$$

**112/3.** Invariable chances $p$ and $q$ of events E and F are given, and $P$ is the probability that in $\mu = m + n$ trials E arrives not less than $m$ times, and F, not more than $n$ times. And (§ 77), for $q/p > n/(m + 1)$ and $q/p < n/(m + 1)$,

$$P = \frac{1}{\sqrt{\pi}}\int_k^\infty \exp(-t^2)dt + \frac{(\mu+n)\sqrt{2}}{3\sqrt{\pi\mu mn}}\exp(-k^2)\ ,\qquad (c_1)$$

$$P = 1 - \frac{1}{\sqrt{\pi}}\int_k^\infty \exp(-t^2)dt + \frac{(\mu+n)\sqrt{2}}{3\sqrt{\pi\mu mn}}\exp(-k^2).\qquad (c_2)$$

Here, $k$ is a positive magnitude whose square is

$$k^2 = n\ln\frac{n}{q(\mu+1)} + (m+1)\ln\frac{m+1}{p(\mu+1)}.$$

**112/4.** Denote by $R$ the probability that the number of occurrences of the events E and F in $\mu$ trials will not be beyond the limits

$$\mu p \mp u\sqrt{2\mu pq}, \ \mu q \pm u\sqrt{2\mu pq}.$$

Here, $u$ is a positive magnitude, very small as compared with $\sqrt{\mu}$. Then (§ 79)[18]

$$R = 1 - \frac{2}{\sqrt{\pi}} \int_u^\infty \exp(-t^2)dt + \frac{1}{\sqrt{2\pi\mu pq}} \exp(-u^2). \qquad (d)$$

Conversely, if the chances $p$ and $q$ are unknown, and E and F arrived $m$ and $n$ times in $\mu = m + n$ trials, then (§ 83)

$$R = 1 - \frac{2\mu}{\sqrt{\pi}} \int_u^\infty \exp(-t^2)dt + \sqrt{\frac{\mu}{2\pi mn}} \exp(-u^2) \qquad (e)$$

is the probability that the values of $p$ and $q$ will not be beyond the limits

$$\frac{m}{\mu} \pm \frac{u}{\mu}\sqrt{\frac{2mn}{\mu}}, \quad \frac{n}{\mu} \mp \frac{u}{\mu}\sqrt{\frac{2mn}{\mu}}.$$

**112/5.** In two series of very large numbers $\mu$ and $\mu'$ of trials event E occurred or will occur $m$ and $m'$ times, and event F arrived or will arrive $n$ and $n'$ times. Denote by $u$ a positive magnitude, very small as compared with $\sqrt{\mu}$ and $\sqrt{\mu'}$. Then the probability $w$ that the difference $(m/\mu - m'/\mu')$ will not be beyond the limits

$$\mp \frac{u\sqrt{2(\mu^3 m'n' + \mu'^3 mn)}}{\mu\mu'\sqrt{\mu\mu'}},$$

and $(n/\mu - n'/\mu')$ will not be beyond those limits taken with contrary signs, − then that probability will by § 87 be

$$w = 1 - \frac{2}{\sqrt{\pi}} \int_u^\infty \exp(-t^2)dt +$$
$$\frac{\sqrt{\mu\mu'}}{\sqrt{2\pi m'n'(\mu + \mu')}} \exp[-\frac{u^2(\mu^3 m'n' + \mu'^3 mn)}{\mu^2 m'n'(\mu+\mu')}]. \qquad (f)$$

Since $m/\mu \approx m'/\mu'$, and $n/\mu \approx n'/\mu'$, it is possible, without appreciably changing the value of $w$, to interchange the magnitudes $\mu'$, $m'$, $n'$ and $\mu$, $m$, $n$ in its last term which will always be a small fraction. This formula, if at least neglecting its last term (§ 109), is applicable to the general case in which the chances of E and F vary from one trial to another provided that the known or unknown causes of these events did not undergo any changes in these two series, − if the existence of these causes retains the same probabilities, and each cause invariably provides the same previous chances to the arrival of E and F.

**112/6.** Events E and F occur $m$ and $n$ times in $\mu$ trials. Suppose that two other contrary events $E_1$ and $F_1$ arrive $m_1$ and $n_1$ times in $\mu_1$ trials. Suppose that $\mu$ and $\mu_1$ are large numbers, and that

$$m_1/\mu_1 - m/\mu = \delta,$$

where δ is a very small positive or negative fraction and denote by $p$ and $p_1$ the unknown and constant chances of the arrival of E and $E_1$, and by λ, the probability that $p_1$ exceeds $p$ at least by a small and given positive fraction ε. Then, let $u$ be a positive magnitude and suppose that

$$\pm \frac{(\varepsilon - \delta)\mu\mu_1\sqrt{\mu\mu_1}}{\sqrt{2(\mu^3 m_1 n_1 + \mu_1^3 mn)}} = u.$$

The sign coincides with the sign of the factor (ε − δ). By § 88

$$\lambda = \frac{1}{\sqrt{\pi}}\int_u^\infty \exp(-t^2)dt, \; \lambda = 1 - \frac{1}{\sqrt{\pi}}\int_u^\infty \exp(-t^2)dt. \quad (g)$$

These formulas correspond to (ε − δ) > 0 and < 0. They also express the probability that the unknown chance $p$ of the arrival of the event E exceeds the ratio $m/\mu$ given by observations by a given fraction $w$, for which it is sufficient that

$$u = \pm (w - \frac{m}{\mu})\frac{\mu\sqrt{\mu}}{\sqrt{2mn}} \text{ for } (w - m/\mu) > 0 \text{ or } < 0.$$

**112/7.** Suppose that the chances of the contrary events E and F vary from one trial to another and are $p_i$ and $q_i$ at trial $i$, so that $p_i + q_i = 1$. Denote for the sake of brevity that

$$\sum p_i/\mu = p, \; \sum q_i/\mu = q, \; 2\sum p_i q_i/\mu = k^2,$$

where the sums extend from $i = 1$ to μ.

Suppose also that the events E and F arrived $m$ and $n$ times in μ trials and denote by $u$ a positive magnitude, very small as compared with $\sqrt{\mu}$. Then (§ 96)

$$R = 1 - \frac{2}{\sqrt{\pi}}\int_u^\infty \exp(-t^2)dt + \frac{\exp(-u^2)}{k\sqrt{\pi\mu}} \quad (h)$$

expresses the probability that the ratios $m/\mu$ and $n/\mu$ will be contained within the limits $p \mp uk/\sqrt{\mu}, \; q \pm uk/\sqrt{\mu}$. This conclusion coincides with the formula (d) in the particular case of constant chances.

**112/8.** Suppose that a thing A can take all values within the limits $(h \pm g)$ and that they are equally possible and the only possible. Denote by $P$ the probability that for some number $i$ of trials the sum of the values of A is contained within also given limits $(c \pm \varepsilon)$. Then (§ 99)

$$2(2g)^i P = \frac{\Gamma - \Gamma_1}{i!}, \qquad (i)$$

where

$$\Gamma = \pm(ih + ig - 4g) - c + ({}^i \mp i \; i/2 + ig - g) - c + {}^i \pm$$
$$C_i^2(ih + ig - 4g) - c + ({}^i \mp C_i^3 \; i/6 + ig - g) - c \pm {}^i \pm$$

whereas $\Gamma_1$ is equal to $\Gamma$ with a changed sign of ε. In each term of this formula the superior sign should be chosen if the magnitude in brackets is positive and the inferior sign otherwise. Magnitudes $g$ and ε are positive, and $h$ and $c$ can be both positive and negative.

**112/9.** Whatever is the law of probabilities of the possible values of a thing A at each trial, and however it varied from one trial to another, if $s$ is the sum of those values, for a very large number μ of trials we have (§ 101)

$$P = 1 - \frac{2}{\sqrt{\pi}} \int_u^\infty \exp(-t^2) dt \qquad (k)$$

for the probability of the mean value $s/\mu$ of A to be within the limits $k \mp 2u\sqrt{h/\mu}$. Here, $u$ is a positive magnitude, very small as compared with $\sqrt{\mu}$, $h$ is positive and $k$ and $h$ depend on the probabilities of the values of A in the entire series of trials. If these probabilities are constant, and the same for all possible values of A within given limits $a$ and $b$ and disappear beyond these limits, then

$$k = (a + b)/2, \; h = (b - a)/2\sqrt{6}.$$

If A only takes a finite number $v$ of constant and equally probable values $c_1, c_2, \ldots, c_v$, then

$$k = \frac{1}{v} \sum (c_1 + c_{2v} + \ldots + c\;),$$
$$h = \frac{1}{2v^2} \sum v[(c_1^2 + c_{2v}^2 + \ldots + c^2) - (c_1 + c_{2v} + \ldots + c)^2].$$

**112/10.** Suppose that $\lambda_n$ is the value of A at the $n$-th trial and

$$\sum \lambda_n/\mu = \lambda, \; \sum (\lambda_n - \lambda)^2/\mu = l^2/2,$$

with the sums extending from $n = 1$ to μ. Suppose also that not a single cause of all the possible values of A had experienced any changes either in its probability or chances that it provides to each of those values. Then, there exists a special magnitude γ, which the mean value $s/\mu$ of A indefinitely approaches as μ ever increases and which it reaches if μ becomes infinite. And formula (k) expresses the probability that γ is contained within the limits (§ 106)

$$\frac{s}{\mu} \mp \frac{ul}{\sqrt{\mu}},$$

which do not include any unknowns.

**112/11.** Suppose that in the second series of a very large number $\mu'$ of trials $s'$ is the sum of the values of A, and $l'$, the new value of $l$ from the first series. The formula (k) will also express the probability that the difference $(s'/\mu' - s/\mu)$ of the two means is contained within limits (§ 107)

$$\mp \frac{u\sqrt{l'^2\mu + l^2\mu'}}{\sqrt{\mu\mu'}}.$$

If $l' \approx l$ the same formula will also express the probability that the mean $s'/\mu'$ of the second series is contained within the limits

$$\frac{s}{\mu} \mp \frac{ul\sqrt{\mu + \mu'}}{\sqrt{\mu\mu'}},$$

only depending on the results of the first series and the given magnitude $u$ and becoming the narrower the more $\mu'$ is exceeding $\mu$.

**112/12.** Suppose that many series of very large numbers $\mu$, $\mu'$, $\mu''$, … of trials are made for establishing the values of one and the same thing A. I denote the sums of its values obtained in those series by $s$, $s'$, $s''$, … As before, the preceding magnitude $l$ applies to the first series, and the corresponding magnitudes for the other series are $l'$, $l''$, … Suppose that the causes of the errors of measurement vary from one series to another, but that, as the numbers $\mu$, $\mu'$, $\mu''$, …, increase further, all the means $s/\mu$, $s'/\mu'$, $s''/\mu''$, … nevertheless indefinitely tend to one and the same unknown magnitude $\gamma$, the veritable value of A. This indeed takes place if the indicated causes do not lead in any series to unequal probabilities of errors equal in magnitude and contrary in sign[19]. And the formula (k) will also express the probability that A is contained within the limits (§ 108)

$$\frac{sq}{\mu} + \frac{s'q'}{\mu'} + \frac{s''q''}{\mu''} + \ldots \mp \frac{u}{D}, \qquad (112.1)$$

where

$$\frac{\mu}{l^2} + \frac{\mu'}{l'^2} + \frac{\mu''}{l''^2} + \ldots = D^2, \quad \frac{\mu}{D^2 l^2} = q, \quad \frac{\mu'}{D^2 l'^2} = q', \quad \frac{\mu''}{D^2 l''^2} = q'', \ldots$$

Furthermore, the sum in (112.1), or of the products of the means $s/\mu$, $s'/\mu'$, $s''/\mu''$, …, respectively multiplied by the magnitudes $q$, $q'$, $q''$, …, is an approximate value of $\gamma$, the most favourable of all those, which can be derived from combining all the series of observations and is that value of $\gamma$ whose limits of errors $\mp u/D$ will be the narrowest possible for a given value of *u* or for the same level of probability.

**112/13.** Finally, suppose that the causes of the arrival of the event E, as formula (f) indeed expressed it, do not change during the trials. Then the ratio $m/\mu$ of the number of occurrences of E to the number $\mu$ of trials indefinitely tends to a special magnitude *r* and exactly reaches it if $\mu$ becomes infinite. And this formula (f), if neglecting its last term, or formula (k) will be the probability that the unknown value of *r* is contained within the limits (§ 109)

$$\frac{m}{\mu} \mp \frac{u\sqrt{2m(\mu-m)}}{\mu\sqrt{\mu}}.$$

**113.** For completing the survey of the formulas, it would have been necessary to add those which apply to the probabilities of the values of one or many magnitudes derived from a very large number of linear equations corresponding to the results of the same number of observations. However, concerning these additional formulas I refer the readers to Laplace (1812/1886, pp. 516 – 519). Issuing from the system of 126 equations compiled by Bouvard and concerning the motion of Saturn in longitude[20], and applying to them the *method of least squares*, he concluded that we can bet a million against one that the mass of Jupiter, when assuming the mass of the Sun as unity, does not differ in any direction by more than a hundredth part from the fraction 1/1070.

Nevertheless, later observations of another nature led to the magnitude almost equal to 1/1050, which exceeds fraction 1/1070 by approximately 1/50 of its value, and, as it seemed, established a defect in the calculus of probability. The value 1/1050 of the mass can not be doubted, as Encke concluded by the perturbations of a comet having a period [of returning] of 1204 days, and as Gauss and Nicolai [21] decided by the perturbations of Vesta and Junon, and Airy by the recently measured elongations of the satellites of Jupiter.

And still, even if Laplace's calculations determined the mass of that planet with a probability very close to certainty but 1/50 less than its veritable value, we can not conclude that the intensity of Jupiter's force of attraction could have acted on Saturn less than on its own satellites, on comets or minor planets. And that result did not occur because of some inaccuracy in the formulas of probability applied by Laplace. We can believe that the mass of Jupiter as obtained by him was a bit too small owing to some faulty terms in the expression of such complicated perturbations of Jupiter. They were somewhat corrected, but other terms can still demand new corrections. A complete alteration of the tables of motion of Saturn and Jupiter, even now so precise, is an important point of the *Mécanique céleste,* and the present work of Bouvard[22] will certainly clear it up.

## Notes

**1.** The sum $\sum U(m, n)u^m v^n$ is called the double generating function of $u$ and $v$. See also § 20.

**2.** The formula (95.1) expresses a limit theorem for the pattern of trials with variable probability. Then comes the integral limit theorem, formulas (97.1) and (101.1). Poisson also studied the central limit theorem, only proved (if disregarding its particular case, the theorem of De Moivre − Laplace) by Markov and Liapunov. See §§ 102 − 103, and Poisson (1824, §§ 1 – 7; 1837). Issuing from those sections, Poisson proved the *second main proposition* of the law of large numbers (§ 52), but he did not sufficiently describe the conditions of that limit theorem, and his account is methodically bad.

Hald (1998, pp. 317 – 327) described in detail Poisson's proof of the central limit theorem. Without substantiation he stated that the lack of explicitly formulated conditions of the theorem conformed to the *then prevailing custom*. Then, Hald did not indicate whether Poisson's proof was rigorous or not; it is known that only Markov and Liapunov were able to provide a rigorous proof.

**3.** Poisson derives the inversion formula for sums of differently distributed lattice random variables. Just above he applied a characteristic function.

**4.** Concerning the study of that particular case, I refer readers to my memoir (1824) earlier mentioned in § 60. Poisson

**5.** See these previous expressions in § 94.

**6.** The second formula below could have been simplified in an obvious way.

**7.** Poisson did not explain the formulas below and the second one, just as stated in Note 6 about a previous formula, could have been simplified.

**8.** Here, Poisson derives the Dirac function which he (1811/1833, p. 637) had applied even earlier.

**9.** This phrase should have been inserted somewhat above.

**10.** I (Note 3 to Chapter 3) indicated that Poisson barely applied the notion of variance.

**11.** In this section, Poisson studied the admissible difference between the mean values of a random variable in two series of trials with a variable probability; in § 109 he investigated the difference of the mean frequencies of the arrivals of a random variable. He (§ 88 of Chapter 3) only considered the difference of probabilities for series of Bernoulli trials.

**12.** See Note 17 to Chapter 2.

**13.** I can also mention the probability of target firing which I considered in a memoir written before this book. It will appear […] (1837). Poisson

**14.** Laplace also issued from data on 10 planets. It was Cournot (1843/1984, p. 175) who named them (including 4 minor planets).

**15.** Poisson did not explain why should the sum of the inclinations approximately equal 75 right angles. Above, he did not indicate that the eccentricities of the planetary orbits are conditioned by the velocity of the planets' rotation round the Sun. It would perhaps be more proper to discuss these velocities.

**16.** A. Bouvard (1767 − 1843), a French astronomer. Poisson should have named his nephew.

**17.** In his appended Note, Poisson qualitatively explains the origin of the aérolithes (meteorites) and shooting stars by the existence of an *inexhaustible* multitude of very small bodies rotating round the Sun, the planets *or even satellites*.

**18.** A comparison of the results (d) and (e) shows that Poisson believed that the direct and the inverse problem are equally precise. So also had thought Jakob Bernoulli and De Moivre, but Bayes proved that the inverse problem was less precise (Sheynin 2010). However, Poisson was not altogether consistent, see the end of § 71 in Chapter 3 and the description of §§ 69 – 72 in his Contents. Anyway, he did not study quantitatively the difference between those two problems.

**19.** See Note 17 to Chapter 2. In § 112, just below, Poisson adduced a summary of his formulas some of which are provided in a changed form.

**20.** Laplace would have been unable to determine the mass of Jupiter by the observations of Saturn. Actually, he had at his disposal 126 equations concerning Jupiter and 129, concerning Saturn, see Supplement 1 of 1816 to his *Théor. Anal. Prob*. (Laplace (1812/1886, p. 516 – 519). He concluded that, if observations will be treated the same way, his estimate of the mass of Jupiter will not change during the next century even by 1/100. Like Poisson after him, he did not account for the unavoidable existence of systematic errors. Note that Laplace almost never, and Poisson never referred to the fundamental achievements of Gauss in the theory of errors.

**21.** F. B. G. Nicolai (1793 – 1846), a German astronomer.

**22.** Poisson mentioned celestial mechanics as a scientific subject rather than as Laplace's classic, but italicized both these words.

## Bibliogrpahy

# Chapter 5. Application of the General Rules of Probability to the Decisions of Jury Panels and Judgements of Tribunals[1]

## Mistakes/Misprints Unnoticed by the Author

**1**. In § 116, p. 327 of the original text. The last bracket in the denominator of the displayed formula should be $(1 - u')$ rather than $(1 - u'')$.

**2**. In § 127, p. 350 of the original text. The left side of the second formula (14) should be $Z_i$ rather than $Z$.

**3.** In § 130, p. 358 of the original text. The last integral in the second displayed formula should be over $[-\infty, 0]$ rather than over $[\infty, 0]$.

**4.** In § 131, p. 361 of the original text. The integrand on the left side of the formula on line 2 from bottom should contain $U_i$ rather than $U'_i$.

**5.** In § 134, p. 370 of the original text, line 7, the difference of two ratios. In one of them the numerator $a_i$ is missing.

**6.** In § 135, p. 374 of the original text, line 2 after the last displayed formula. The numerator of the ratio should be $a_4$ rather than $a^4$. Several similar mistakes are in §§ 140 and 143.

**7.** In § 139, p. 383 of the original text, the last but one term on the right side of the last displayed formula should contain $t^4$ rather than $t^5$.

**8.** In § 149, p. 409 of the original text. The term $35v^4(1 - v)$ in the first line of the last displayed formula should be $35v^4(1 - v)^3$.

All those misprints/mistakes are corrected in the translation.

**114.** When studying such a delicate matter, it is convenient to begin with the simplest cases, then treating the problem in all generality. And I suppose first of all that there is only one juryman. Let $k$ be the probability of the defendant's guilt prior to the court hearing (prior probability). This probability is established by preliminary information and subsequent charges. I also denote by $u$ the probability that the juryman's decision is faultless, and by $\gamma$, the probability that the accused will be convicted. This will happen if he was guilty, and the juryman was not mistaken or innocent and the juryman mistaken,

By the rule of § 5, the probability of the first case is $ku$, and $(1 - k)(1 - u)$, of the second. Therefore (§ 10), the composite probability of the defendant's conviction is

$$\gamma = ku + (1 - k)(1 - u), \qquad (114.1)$$

and of his acquittal, $(1 - \gamma)$. This latter happens if the accused is guilty and the juryman was mistaken or innocent and the juryman was not mistaken. Therefore,

$$1 - \gamma = k(1 - u) + u(1 - k),$$

which could have been derived from the precedent. When subtracting the latter equality from the first one, we get

$$2\gamma - 1 = (2k - 1)(2u - 1),$$

so that $2\gamma - 1 = 0$ when either $(2k - 1)$, or $(2u - 1)$ is zero. It will be positive and negative when the signs of $(2k - 1)$ and $(2u - 1)$ are the same or contrary. Also,

$\gamma = 1/2 + (2k - 1)(2u - 1)/2,$

so that $\gamma$ exceeds 1/2 by a half of the positive or negative product $(2k - 1)(2u - 1)$.

After the juryman's decision we can formulate the two only possible hypotheses: was the accused guilty or not? Their probabilities, like in all hypothetic cases, are determined by the rule of § 34 and their sum is unity, so that it is sufficient to establish one of them. Suppose that the accused is convicted, and the probability of the first hypothesis or his guilt is $p$. According to the mentioned rule,

$$p = \frac{ku}{ku + (1-k)(1-u)} \qquad (114.2)$$

since the conviction is the observed event whose probability in accordance with both hypotheses, as indicated just above, is $ku$ or $(1 - k)(1 - u)$.

Suppose that the accused is acquitted, and denote by $q$ the probability of the second hypothesis, i. e., of his innocence. The observed event is the acquittal of the accused, whose probability in accordance with that hypothesis is $(1 - k)u$ or $k(1 - u)$ under the contrary assumption, so therefore

$$q = \frac{(1-k)u}{(1-k)u + k(1-u)}. \qquad (114.3)$$

The denominators of $p$ and $q$ are the expressions of $\gamma$ and $1 - \gamma$, and consequently

$p = ku/\gamma, q = (1 - k)u/(1 - \gamma).$

The probability that the juryman was not mistaken is

$u = p\gamma + q(1 - \gamma),$

which is easily verified. This happens in two different cases: the accused was guilty and will be convicted, or he was innocent and will be acquitted.

By the rule of § 9 concerning the probability of an event composed of two simple events whose chances influence each other, the probability of the first case is $\gamma p$, and of the second, $(1 - \gamma)q$, and (§ 10) the composite value of $u$ is the sum of these two products. After the juryman's decision is announced, the probability of its faultlessness is $p$ or $q$, if the accused is convicted or acquitted. And if $k \neq 1/2$, it, can only equal $u$ if $u = 0$ or 1.

The formulas above provide a complete solution of the problem in the case of only one juryman, of the same problem about the probability of an event testified by one witness (§ 36). The event, either veritable of false, is here the defendant's guilt. Before the juryman's decision was announced, there existed some reason for believing that that event is veritable; its probability was $k$ with $(1 − k)$ being the probability of the defendant's innocence. New information appeared after that decision, and $k$ became another probability, $p$, if the juryman decided, or attested that the accused is guilty; otherwise $(1 − k)$ becomes probability $q$.

In each case the preliminary probabilities $k$ and $(1 − k)$ should have evidently heightened if the chance of the juryman's faultlessness exceeded 1/2, and lowered otherwise. In other words, heightened or lowered depending on whether $u > 1/2$ or $< 1/2$. This follows from the expressions for $p$ and $q$ leading to

$$p = k + \frac{k(1-k)(2u-1)}{\gamma}, \; q = 1 - k - \frac{k(1-k)(2u-1)}{1-\gamma},$$

so that $p > k$, $p < k$, $q < (1− k)$, $q > (1− k)$ depending on whether $u > 1/2$ and $u < 1/2$. At $u = 1/2$ the preliminary probabilities $k$ and $(1 − k)$ do not change at all.

The just provided expressions of $p$ and $q$ lead to

$p\gamma + q(1 − \gamma) = u = k\gamma + (1 − k)(1 − \gamma)$.

If in addition to the probability of guilt $k$ the chance $\gamma$ of conviction somehow becomes known, these formulas can serve for calculating the probability that the juryman is not mistaken. This can also be verified when observing that he is not mistaken if the accused is guilty and convicted or innocent and acquitted. Before the juryman's decision, the probabilities of these two cases are $k\gamma$ and $(1 − k)(1 − \gamma)$, and their sum is the composite value of $u$.

If $k = 1/2$, the preliminary values of $p$ and $q$ at once become $p = q = u$. Actually, since there is no reason to believe in advance that the accused is rather guilty than innocent, after the juryman's decision our grounds for believing either possibility can not differ from the probability that the juryman is faultless. If $k = 1$, so that the probability of guilt is thought in advance to be beyond doubt, then $p = 1$ and $q = 0$. Whichever was this decision and the chance of its faultlessness, the guilt of the accused will become still more certain after that decision. The same takes place concerning the innocence of the accused when $k = 0$, i. e., when it is certain in advance. However, it is unknown whether the accused will be convicted or acquitted, and in the first case the chance of conviction is $\gamma = u$ and $1 − u$ in the second. These chances will be equal, as it should have been, to the probability that, at $k = 1$ and 0, the juryman is not mistaken or mistaken.

**115.** Suppose that after the juryman's decision the accused is judged by another juryman whose probability of faultlessness is $u'$. It is required to determine the probabilities $c, b, a$ that the accused will be convicted by both; acquitted by one of them and convicted by the other; acquitted by both. Suppose that $\gamma'$ is the probability of conviction of the accused by the second juryman after he was convicted by the first one. The chance of the first conviction is $\gamma$, so that the probability of conviction in both cases will be $c = \gamma\gamma'$. However, when the accused is brought before the second juryman, probability $p$ that he is guilty has already appeared after the decision of the first one. Therefore, the value of $\gamma'$ is derived from formula (114.1) when $p$ and $u'$ are substituted instead of $k$ and $u$:

$$\gamma' = pu' + (1 - p)(1 - u')$$

and formulas (114.1) and (114.2) lead to

$$c = kuu' + (1 - k)(1 - u)(1 - u').$$

By a similar reasoning it follows that

$$a = k(1 - u)(1 - u') + (1 - k)uu'.$$

Combining these two formulas, we obtain the probability that the decisions of both jurymen, whether to convict or acquit, coincide:

$$a + c = uu' + (1 - u)(1 - u').$$

Note that this composite probability does not depend on the prior probability of the defendant's guilt.

Suppose that the first juryman acquitted the accused and let the probability of his conviction by the second be $\gamma_1$, then the product $(1 - \gamma)\gamma_1$ will express the probability that these contrary decisions occurred in the indicated order. In addition, $(1 - q)$ will be the probability of the defendant's guilt when he, being acquitted by the first juryman, appears before the second. The value of $\gamma_1$ is determined by formula (114.1) if $k$ and $u$ are replaced by $(1 - q)$ and $u'$, so that

$$\gamma_1 = (1 - q)u' + q(1 - u').$$

Taking into account the values of $(1 - \gamma)$ and $q$, given by formulas (114.1) and (114.3), we have

$$(1 - \gamma)\gamma_1 = k(1 - u)u' + (1 - k)(1 - u')u.$$

Interchanging $u$ and $u'$ in this expression, we will obviously get the probability that the decisions of the jurymen are contrary but were made in the inverse order to the precedent. Adding this probability to the previous, we obtain the composite probability of contrary decisions made by the jurymen in any order

$$b = (1 - u)\,u' + (1 - u')u.$$

Just like in the case of coinciding decisions, it does not depend on $k$. At $u = u' = 1/2$ both these probabilities also equal 1/2, and their sum, $a + b + c$, invariably equals unity, as it should have been.

The probability $p'$ of guilt of the accused convicted by both jurymen is expressed by formula (114.2) when $k$ and $u$ are replaced by $p$ and $u'$, and the probability $q'$ of his innocence after his being acquitted by both is derived from formula (114.3) when those magnitudes are replaced by $(1 - q)$ and $u'$:

$$p' = \frac{pu'}{pu' + (1-p)(1-u')}, \quad q' = \frac{qu'}{qu' + (1-q)(1-u')}.$$

Taking account of the values of $p$ and $q$ given by the same formulas (114.2) and (114.3), these $p'$ и $q'$ become

$$p' = \frac{kuu'}{kuu' + (1-k)(1-u)(1-u')}, \quad q' = \frac{(1-k)uu'}{(1-k)uu' + k(1-u)(1-u')}.$$

Suppose also that $p_1$ is the probability of guilt of the accused after his acquittal by the first juryman and conviction by the second, and $q_1$, the probability of his innocence after conviction by the first juryman and acquittal by the second. The value of probability $p_1$ that the accused, acquitted by the first juryman is not innocent, is derived from formula (114.2) when $u$ and $k$ are replaced by $u'$ and $(1 - q)$. The value of $q_1$ is derived from formula (114.3) when those magnitudes are replaced by $u'$ и $p$. Thus,

$$p_1 = \frac{(1-q)u'}{(1-q)u' + q(1-u')}, \quad q_1 = \frac{(1-p)u'}{(1-p)u' + p(1-u')}$$

and, taking into account the same formulas (114.2) and (114.3),

$$p_1 = \frac{k(1-u)u'}{k(1-u)u' + (1-k)(1-u')u}, \quad q_1 = \frac{(1-k)(1-u)u'}{(1-k)(1-u)u' + k(1-u')u}.$$

The probability that the accused convicted by the first juryman and acquitted by the second is guilty is $(1 - q_1)$. It can evidently be derived from $p_1$ by interchanging $u$ and $u'$, which indeed takes place: then the probability of the innocence of the accused acquitted by the first juryman and convicted by the second, or $(1 - p_1)$, is derived from $q_1$ by transposing $u$ and $u'$.

At $u' = u$ we have $p_1 = k$ and $q_1 = (1 - k)$, which should have occurred since contrary decisions by jurymen having the same chance of faultlessness can not change anything in our previous reason in believing that the accused is guilty or innocent.

**116.** It is not difficult to generalize these considerations on successive decisions of some number of jurymen having a known chance of faultlessness. However, it is simpler to arrive at the proper results by the following way.

For the sake of definiteness I suppose that there are three jurymen. Denote their probabilities of faultlessness by $u$, $u'$ и $u''$ and, as before, by $k$, the prior probability that the accused is guilty. Unanimous conviction requires that he is indeed guilty and that not a single juryman is mistaken; or, that he is innocent and all the jurymen are mistaken. The composite probability of such a conviction will therefore be

$$kuu'u'' + k(1-u)(1-u')(1-u'').$$

It is also seen that the probability of a unanimous acquittal is

$$k(1-u)(1-u')(1-u'') + (1-k)uu'u''.$$

The sum of these two expressions will be the probability of a unanimous decision whether to convict or to acquit

$$uu'u'' + (1-u)(1-u')(1-u'').$$

It does not depend on $k$, which is just as true for any number of jurymen.

The accused can be convicted by two jurymen and acquitted by the third in three different ways depending on whether the chance of faultlessness of that third juryman is $u$, $u'$ or $u''$. He can also be acquitted by two jurymen and convicted by the third, again in those same three ways. It is not difficult to see that the probabilities of all the six combinations are

$$ku'u''(1-u) + (1-k)(1-u')(1-u'')u$$
$$kuu''(1-u') + (1-k)(1-u)(1-u'')u'$$
$$kuu'(1-u'') + (1-k)(1-u)(1-u')u''$$
$$k(1-u')(1-u'')u + (1-k)(1-u)u'u''$$
$$k(1-u)(1-u'')u' + (1-k)(1-u')uu''$$
$$k(1-u)(1-u')u'' + (1-k)(1-u'')uu'$$

The sum of these six magnitudes is the composite probability that the decisions were not unanimous:

$$u'u''(1-u) + uu''(1-u') + uu'(1-u'') + (1-u')(1-u'')u +$$
$$(1-u)(1-u'')u' + (1-u)(1-u')u''.$$

It does not depend on $k$. The sum of the composite probabilities of the unanimous and majority decisions should be unity which is indeed fulfilled.

After the decisions are announced, it is easy to derive the probability of the defendant's guilt, which, generally speaking, will differ from its previous value. If, for example, the accused is convicted by jurymen whose chances of faultlessness were $u$ and $u'$ and acquitted by the third, the probability of this result will be $kuu'(1-u'')$ or $(1-k)(1-u)(1-u')u''$ in accordance with the hypotheses of his guilt and innocence. Therefore, by the rule of § 34, the probability of his guilt will be

$$\frac{kuu'(1-u'')}{kuu'(1-u'')+(1-k)(1-u)(1-u')u''}.$$

If $u'= u''$ it becomes independent from the common value of $u'$ and $u''$. The same will take place if the accused is only convicted by one juryman whose chance of faultlessness is $u$. And actually after the decision is announced, various versions of the decision by the two other jurymen will not influence the grounds of my belief in the defendant's guilt or innocence. Indeed, there will be no reason to heighten rather than to lower the probability of his guilt since the chances of faultlessness of the two last jurymen are supposed to be the same.

The formulas above are equally applicable to the case in which the jurymen, instead of judging successively and without communicating with each other, are united and judge after deliberation. Discussion can clear up the situation for them and, in general, heighten their probability of faultlessness[2]. The values of $u$, $u'$ and $u''$ in those two cases can differ and less deviate from unity in the second case.

**117.** Consider in particular the case in which the chances of faultlessness are the same for all the jurymen. We will then reduce to it the general case of determining the probability of the number of convictions in a very large number of decisions.

Denote by $u$ this given probability of the faultlessness, by $n$, the number of those jurymen, and, finally, by $k$, the prior probability of the defendant's guilt. Let also $i = 0, 1, 2, …, n$ and $\gamma_i$ be the probability of his conviction by a majority verdict of $(n - i)$ votes against $i$.

For that compound event to take place it was necessary that the accused was guilty, and that $(n - i)$ jurymen were not mistaken whereas $i$ of them had mistook; or, that the accused was innocent and that $(n - i)$ jurymen were mistaken and $i$ of them made no mistake. The probability of the first case is the product of $ku^{n-i}(1-u)^i$ by the number of ways for choosing $i$ jurymen from their total number, $n$. The second case has probability equal to the product $(1-k)u^i(1-u)^{n-i}$ by the number of ways for choosing $(n - i)$ jurymen from their total number. […] As a result, the composite value of $\gamma_i$ is

$$\gamma_i = N_i[ku^{n-i}(1-u)^i + (1-k)u^i(1-u)^{n-i}] \qquad (117.1)$$

[$N_i$ was Poisson's notation for $C_n^i$].

If $n - i > i$ and $n - 2i = m$, the accused is convicted by a majority of $m$ votes. If, however, $i$ jurymen convict him, and $(n - i)$ absolve, he is acquitted by the same majority. The probability of such an acquittal which I denote by $\delta_i$, is determined by the value of $\gamma_i$ after transposing $(n - i)$ and $i$ without changing the value of $N_i$. Thus,

$$\delta_i = N_i[ku^i(1-u)^{n-i} + (1-k)u^{n-i}(1-u)^i]. \qquad (117.2)$$

The sum of the two last equations is

$$\gamma_i + \delta_i = N_i[u^{n-i}(1-u)^i + u^i(1-u)^{n-i}],$$

and does not depend on $k$. The probability of either conviction or acquittal by a given majority does not therefore depend on the prior supposed guilt of the accused. In the particular case of $u = 1/2$ the probabilities $\gamma_i$ и $\delta_i$ considered separately are also independent of $k$, and their common value is

$$\gamma_i = \delta_i = N_i/2^n.$$

At $k = 1/2$ they are again the same for any value of $u$.

**118.** Suppose now that $c_i$ is the probability of the defendant's conviction by a majority verdict of not less than $(n - i)$ votes against not more than $i$ of the jurymen, i. e., of a conviction by a majority of at least $m$ votes. Denote also by $d_i$ the probability that the accused is acquitted by a majority verdict of not less than $(n - i)$ votes against not more than $i$. Then, by the rule of § 10,

$$c_i = \gamma_0 + \gamma_1 + \gamma_2 + \ldots + \gamma_i, \; d_i = \delta_0 + \delta_1 + \delta_2 + \ldots + \delta_i,$$

and, because of the previous formulas,

$$c_i = kU_i + (1 - k)U_i, \; d_i = kV_i + (1 - k)U_i. \qquad (118.1)$$

Here, for the sake of brevity

$$N_0 u^n + N_1 u^{n-1}(1 - u) + N_2 u^{n-2}(1 - u)^2 + \ldots + N_i u^{n-i}(1 - u)^i = U_i,$$
$$N_0(1 - u)^n + N_1(1 - u)^{n-1}u + N_2(1 - u)^{n-2}u^2 + \ldots + N_i(1 - u)^{n-i}u^i = V_i,$$

where $U_i$ is a given function of $u$, which becomes $V_i$ when $u$ is replaced by $(1 - u)$. At the same time

$$c_i + d_i = U_i + V_i$$

is the probability not depending on $k$ that the accused will be either convicted or acquitted by a majority of not less than $m$ votes.

When replacing $i$ by $(n - i - 1)$ in the expression of $d_i$,

$$U_i + V_{n-i-1} = 1, \; c_i + d_{n-i-1} = 1.$$

Suppose that at least $(n - i)$ votes are needed for conviction, then the accused will be acquitted if not more than $(n - i - 1)$ jurymen decide to convict. This means that one of the two events with probabilities $c_i$ and $d_{n-i-1}$ will certainly occur.

For an odd $n$, supposing that $n = 2i + 1$ and therefore $m = 1$,

$$U_i + V_i = [u + (1 - u)]^n = 1, \; c_i + d_i = 1.$$

Consequently, which is evident, the accused will certainly be either convicted or acquitted by a majority of at least one vote. For an even $n$, the least majority is $m = 2$, so that $n = 2i + 2$. Then

$$U_i + V_i = [u + (1-u)]^n - N_{i+1}u^{i+1}(1-u)^{i+1},$$

$$c_i + d_i = 1 - \frac{(2i+2)(2i+1)2i\ldots(i+2)}{(i+1)!}[u(1-u)]^{i+1}.$$

Conviction or acquittal by majority of at least two votes is not however certain since the votes for acquittal and conviction can be equally strong. The probability of this unique case, independent from $k$, is determined by subtracting the previous value of $c_i + d_i$ from unity. Denoting it by $H_i$, we can write it in the form

$$H_i = \frac{(2i+2)![u(1-u)^{i+1}]}{[(i+1)!]^2}.$$

The function $u(1-u)$ reaches its maximal value of 1/4 at $u = 1/2$ and $H_i$ decreases as $u$ further deviates from 1/2. It also continuously decreases as $i$ increases. This conclusion is indeed derived from the expression

$$H_{i+1} = \frac{(2i+3)(2i+4)u(1-u)}{(i+1)^2}H_i.$$

After passing the maximal value indicated above, the ratio $H_{i+1}:H_i$ will invariably be less than unity; the largest value of $H_i$ corresponds to $u = 1/2$ и $i = 0$ and is equal to 1/2.

If $i + 1$ is a very large number, then (§ 67) […] and

$$H_i = \frac{[4u(1-u)]^{i+1}}{\sqrt{\pi(i+1)}}[1 - \frac{1}{8(i+1)} + \ldots]$$

will be an approximate value of that magnitude and, as it is seen, if $u$ appreciably deviates from 1/2, or $4u(1-u)$, from unity, a very small fraction. If $u = 1/2$ and, for example, $i + 1 = 6$ or $n = 12$, this formula, if restricting its second factor to its first two terms, provides 230.94 …/1024. Although $i + 1$ is not a very large number, this very little differs from the exact value, 231/1024.

The sum $U_i + V_i = G_i$ takes values not exceeding unity, and the difference $1 - G_i$ will be positive or zero. And since the expression for $c_i$ can be written as

$$c_i = k - k(1 - G_i) - (2k - 1)V_i,$$

$c_i < k$ if $2k - 1 > 0$ or $k > 1/2$. Therefore, in usual cases in which there are grounds for believing rather in the prior defendant's guilt than his innocence, the chance of his conviction by a majority of at least one vote, i. e. by any majority, will always be lower than this prior probability. Suppose for example that it was possible to bet in advance 4 against 1 on his guilt, but the ratio of stakes on his conviction will be less.

This proposition is obviously independent from either the chances of the jurymen's mistakes or the values of $u$ differing from unity. At $u = 1$ and any $i$, $U_i = 1$, $V_i = 0$, $c_i = k$, $d_i = 1 - k$, and at $u = 0$, $U_i = 0$, $V_i = 1$, $c_i = (1 - k)$, $d_i = k$. At these two extreme values of $u$ both conviction and acquittal can obviously only occur unanimously. This follows from formulas (117.1) and (117.2), which, excluding the case $i = 0$, lead to $\gamma_i = \delta_i = 0$.

**119.** With the previous notation being preserved, we will additionally express by $p_i$ the probability that the accused is guilty if convicted by a majority verdict of $(n - i)$ against $i$, that is, by a majority of $m$ votes, and by $q_i$, the probability of his innocence when acquitted by the same majority verdict. In other words, $p_i$ and $q_i$ are the probabilities that both decisions by that majority verdict are proper. In the first case the probability of the observed event, i. e., of conviction, is $N_i k u^{n-i}(1-u)^i$ or $N_i(1-k)(1-u)^{n-i}u^i$ depending on whether the accused is guilty or not. According to the rule of § 34, after cancelling the factor $N_i$ in the numerator and denominator,

$$p_i = \frac{ku^{n-i}(1-u)^i}{ku^{n-i}(1-u)^i + (1-k)(1-u)^{n-i}u^i}. \qquad (119.1)$$

In case of acquittal

$$q_i = \frac{(1-k)u^{n-i}(1-u)^i}{(1-k)u^{n-i}(1-u)^i + k(1-u)^{n-i}u^i}. \qquad (119.2)$$

If $k = 1/2$, $p_i = q_i$. Indeed, if there are no grounds for believing in advance that the accused is rather guilty than innocent, a proper decision of the case by the same majority verdict has the same probability both when convicting and acquitting. At $u = 1/2$, $u^{n-i}(1-u)^i = u^i(1-u)^{n-i}$ and, as it should be, $p_i = k$ and $q_i = 1 - k$ whichever are $n$ and $i$.

Suppose that in formulas (119.1) and (119.2) $u = t/(1 + t)$, $1 - u = 1/(1 + t)$. Noting that $n = m + 2i$, we have

$$p_i = \frac{kt^m}{kt^m + 1 - k}, \quad q_i = \frac{(1-k)t^m}{(1-k)t^m + k},$$

which proves that the probability of a proper decision, other things being equal, only depends on the appropriate majority vote but not on the total number $n$ of jurymen. Actually, a contrary vote by the same majority with an equal chance of a mistake for all the jurymen can not either heighten or lower the grounds for believing that the decision is either proper or mistaken. This conclusion, however, essentially supposes that the chance $u$ of the jurymen's faultlessness is known before their decision. As shown below, this will not take place if that chance is derived after the decision by the number of votes for and against.

With a given $u$ a decision returned by an odd number of jurymen by a majority of, say, only 1 vote, does not deserve either more or less trust than the decision made by a single juryman. However, the probability of such a decision, whether to convict or acquit, lowers as the total number of jurymen increases. This probability is equal to the sum $w_i$ of the right sides of formulas (117.1) and (117.2) at $n = 2i + 1$. Taking into account the value of $N_i$, we get

$$w_i = \frac{(2i+1)![u(1-u)]^i}{(i!)^2(i+1)}, \quad w_{i+1} = \frac{2i+3}{2i+4} 4u(1-u)w_i.$$

Since $4u(1-u)$ can not exceed unity, invariably $w_{i+1} < w_i$. When comparing that value of $w_i$ with $H_i$ from § 118, we see that the former exceeds the latter in the ratio $1:2u(1-u)$, whatever was $i$.

**120.** If it is only known that the accused is convicted unanimously or by at least a majority of $m$ votes, i. e., of $m$, $m + 2$, $m + 4$, …, $m + 2i$ votes, it occurs that the probability of his guilt $P_i$ is higher than $p_i$. Supposing that the accused is guilty, the probability of his conviction, or of the observed event, will be, in accordance with the above, $kU_i$. It will be $(1 - k)V_i$, when supposing that he is innocent. Therefore, when adding in this latter case the probability of an acquittal by [the same] majority of at least $m$ votes we will have

$$P_i = \frac{kU_i}{kU_i + (1-k)V_i}, \quad Q_i = \frac{(1-k)U_i}{(1-k)U_i + kV_i}. \qquad (120.1, 2)$$

The probabilities of proper decisions by the indicated majority verdict are also expressed by the magnitudes $P_i$ and $Q_i$. Unlike $p_i$ and $q_i$, they are not independent from the total number $n$ of the jurymen and dependent only on $m$ or $n - 2i$.[3] For quantitatively comparing them with each other I assume that $k = 1/2$, which equates $P_i$ and $Q_i$ as well as $p_i$ and $q_i$, and suppose that the prior probabilities of the defendant's innocence and guilt are equal to each other. Finally, I assume that $u = 3/4$, so that you can bet 3 against 1 on the faultlessness of each juryman. For the usual number $n = 12$ of jurymen and $i = 5$ I find first of all that $p_i = 9/10$ and $1 - p_i = 1/10$ and then, that $U_i = 7254 \cdot 3^7/4^{12}$ and $V_i = 239{,}122/4^{12}$. It occurs that almost exactly $P_i = 403/409$, $1 - P_i = 6/409$.

This proves that in the provided example the probability $(1 - P_i)$ of a mistaken conviction by a majority of at least 2 votes barely amounts to $1/7$ of the probability $(1 - p_i)$ of the error to be feared when deciding by a majority verdict of 2 votes or by 7 votes against 5.

The formulas (117.1, 117.2; 118.1; 119.1, 119.2; 120.1, 120.2) are easily applicable in the case in which the accused brought before the same number of jurymen was already convicted or acquitted by another court. The value of $k$, included in the indicated formulas, will here be the probability of the defendant's guilt resulting from the first decision and determined by one of the formulas mentioned above.

**121.** If $(n - i)$ and $i$ are very large numbers, the values of $U_i$ and $V_i$ should be calculated approximately. For this aim I assume that $(1 - u) = v$, so that $U_i$ will be the sum of the first $(i + 1)$ terms of the expansion of $(u + v)^n$, arranged in increasing powers of $v$. In accordance with formulas (77.1a, b), when replacing *p, q, n, μ* by *u, v, i, n* we obtain

$$U_i = \frac{1}{\sqrt{\pi}} \int_\theta^\infty \exp(-t^2)dt + \frac{(n+i)\sqrt{2}}{3\sqrt{\pi n i(n-i)}} \exp(\theta^2), \qquad (121.1a)$$

$$U_i = 1 - \frac{1}{\sqrt{\pi}} \int_\theta^\infty \exp(-t^2)dt + \frac{(n+i)\sqrt{2}}{3\sqrt{\pi n i(n-i)}} \exp(\theta^2), \qquad (121.1b)$$

with the square of the positive magnitude θ being

$$\theta^2 = i \ln \frac{i}{v(n+1)} + (n+1-i) \ln \frac{n+1-i}{u(n+1)}.$$

These expressions for $U_i$ correspond to the cases in which *v/u* exceeds the ratio *i/(n + 1 − i)* or is less than it.

If the accused was convicted by any majority verdict from 1 or 2 votes or unanimously, the number (*n* + 1) and the indicated ratio will be almost 2*i* and 1, and the formulas (121.1a) and (121.1b) should be applied at *v > u* и *v < u,* and if *u* and *v* essentially differ from 1/2, or 4*uv* from unity, then almost exactly θ = *i*ln(1/4*uv*).

And so, since *i* is a very large number, the value of $\theta^2$ will be sufficiently large for the integrals and the exponentials in formulas (121.1) to become inappreciable. The magnitude $U_i$ will be 1 or 0 at *u > v* and *u < v*. And since in our case the sum ($U_i + V_i$) is almost or exactly 1 at odd and even values of *i*, $V_i$ will be 0 or 1 at *u > v* and *u < v*.

It follows that if the prior probability *k* of the defendant's guilt is not a very small fraction and the chance *v* of a mistake of each of a very large number *n* jurymen is appreciably less than the chance *u* of the contrary, the probability $P_i$ of guilt after his conviction will be very close to certainty. Indeed, it is very unlikely that an indicated large number *n* of jurymen will return a decision by a feeble majority verdict. On the contrary, if *u* is appreciably less than *v*, and *k* is not very close to unity, the probability $P_i$ of a proper decision becomes a very small fraction and the defendant's innocence, extremely probable. The probabilities of conviction and acquittal, $c_i$ and $d_i$, indicated by formulas (118.1), will very little differ from *k* and (1 − *k*) if *u > v* and, on the contrary, from (1 − *k*) and *k* if *v > u*.

When *u* = *v* = 1/2 and *n* = 2*i* + 1 or 2*i* + 2, the ratio *i/(n + 1 − i)* becomes a little less than *v/u* or 1 so that formula (121.1a) should be applied. And since θ will be a very small fraction, then, if neglecting the square of θ and the terms having *i* as the divisor, and noting that in formulas (121.1)

$$\int_\theta^\infty \exp(-x^2)dx = \frac{\sqrt{\pi}}{2} -$$

the equality

$$U_i = \frac{1}{2} + \frac{1\theta}{\sqrt{\pi i}} - \frac{1}{\sqrt{\pi}}$$

will become almost exact.

For an odd $n = 2i + 1$, substituting 1/2 instead of $u$ and $v$, we get

$$\theta^2 = -i\ln[1 + \frac{1}{i}] - (i+2)\ln[1 - \frac{1}{i+2}].$$

Expanding the logarithms in series, we can obtain to within our adopted precision $\theta = 1/\sqrt{i}$ and $U_i = 1/2$. The sum $(V_i + U_i) = 1$, so that $V_i$ also equals 1/2 and, as it should have been, $P_i = k$ whichever was the number of the jurymen provided that their chances of being mistaken or not are the same. At $n = 2i + 2$ and $u = v = 1/2$

$$\theta^2 = -i\ln[1 + \frac{3}{2i}] - (i+3)\ln[1 - \frac{3}{2i+6}], \quad \theta = \frac{3}{2\sqrt{i}}, \quad U_i = \frac{1}{2} - \frac{1}{2\sqrt{\pi i}}.$$

However, taking into account the value of $H_i$ from § 118, we have in this case

$$U_i + V_i = 1 - \frac{1}{\sqrt{\pi i}},$$

but

$$V_i = \frac{1}{2} - \frac{1}{2\sqrt{\pi i}}$$

and, just like in the previous case, $P_i = k$ since $U_i = V_i$. The probability $c_i$ of conviction which we consider will be independent from $k$ and equal to $U_i$ or become a bit lower than 1/2.

**122.** Invariably assuming that there are some $n$ jurymen, suppose now that for each of them the chance of faultlessness can take $v$ different and unequally probable values, $x_1, x_2, …, x_v$ for the first of them, $x'_1, x'_2, …, x'_v$ for the second, $x''_1, x''_2, …, x''_v$ for the third etc. Denote by $X_i, X'_i, X''_i, …$ the probabilities that the chances $x_i, x'_i, x''_i …$ will indeed occur. They will also be the probabilities of the corresponding chances $1 - x_i, 1 - x'_i, 1 - x''_i, …$ One of the chances $x_1, x_2, …, x_v$ will certainly take place, the same as one of the chances $x'_1, x'_2, …, x'_v$ and $x''_1, x''_2, …, x''_v$ etc, and therefore

$$X_1 + X_2 + … + X_v = 1,$$
$$X'_1 + X'_2 + … + X'_v = 1, \quad (121.1)$$
$$X''_1 + X''_2 + … + X''_v = 1, …$$

Suppose that

$$X_1 x_1 + X_2 x_2 + … + X_v x_v = u,$$

$$X'_1 x'_1 + X'_2 x'_2 + \ldots + X'_\nu x'_\nu = u', \qquad (121.2)$$
$$X''_1 x''_1 + X''_2 x''_2 + \ldots + X''_\nu x''_\nu = u'',$$

then at the same time

$$X_1(1 - x_1) + X_2(1 - x_2) + \ldots + X_\nu (1 - x_\nu) = 1 - u,$$
$$X'_1(1 - x'_1) + X'_2(1 - x'_2) + \ldots + X'_\nu (1 - x'_\nu) = 1 - u',$$
$$X''_1(1 - x''_1) + X''_2(1 - x''_2) + \ldots + X''_\nu (1 - x''_\nu) = 1 - u''$$

with $u$, $u'$, $u''$ and $(1 - u)$, $(1 - u')$, $(1 - u'')$ being the mean chances of faultlessness and mistake of the respective jurymen.

Then the probability $\Pi$ that no juryman having chances $x_i$, $x'_{i'}$, $x''_{i''}$ … of faultlessness is mistaken will be the product of these chances by their respective probabilities $X_i$, $X'_{i'}$, $X''_{i''}$, …:

$$\Pi = X_i X'_{i'} X''_{i''} \ldots x_i x'_{i'} x''_{i''} \ldots$$

Denote by $P$ the probability that none of the $n$ jurymen having the indicated chances of faultlessness is mistaken whichever are the probabilities of those possible chances. By the rule of § 10 $P$ will be the sum of the $n\nu$ values of $\Pi$, corresponding to a successive substitution of numbers 1, 2, …, $\nu$ instead of each of the $n$ numbers $i$, $i'$, $i''$, … It is easy to see that for any $n$ and $\nu$ this sum is equal to the product of the $n$ means $u$, $u'$, $u''$, …, so that $P = uu'u'' \ldots$

Relative to the chances $x$, $x'_{i'}$, $x''_{i''}$ … of faultlessness, the probability that only one single juryman is mistaken can be calculated when replacing $x_{i'}$ by $(1 - x_{i'})$ in the product $\Pi$, if that was the first juryman, when replacing $x'_i$ by $(1 - x'_i)$, if that was the second one etc. Denote by $\Pi'$ the composite probability that only one single juryman was mistaken:

$$\Pi' = X_i X'_{i'} X''_{i''} \ldots [(1 - x_i) x'_{i'} x''_{i''} \ldots + $$
$$x_i (1 - x'_{i'}) x''_{i''} \ldots + x_i x'_{i'} (1 - x''_{i''}) \ldots + \ldots].$$

Taking into account all the possible chances of faultlessness for each juryman, denote now by $P'$ the probability that only one of them is mistaken. It will be the sum of the $n\nu$ values of $\Pi'$ obtained by successively substituting all the numbers 1, 2, …, $\nu$ for each of the indices $i$, $i'$, $i''$, … It is easy to see that this sum only depends on the means $u$, $u'$, $u''$, … and equals

$$P' = (1 - u) u' u'' \ldots + u(1 - u') u'' \ldots + uu'(1 - u'') \ldots + \ldots$$

When continuing such transformations, we can arrive at the following general proposition: The probability that $(n - i)$ jurymen out of $n$ are not mistaken, and $i$ are mistaken, does not change if for each of them the chance of faultlessness can only take one single value, $u$ for the first one, $u'$ for the second etc. Therefore, if the mean chances $u$, $u'$, $u''$, … are not equal one to another, the differing probabilities of conviction by a given majority verdict and guilt of the accused are determined by the rules of § 116 generalized on some number $n$ of jurymen. If, however, they are equal to each other, their probabilities are expressed by formulas (117.1, 117.2; 118.1; 119.1, 119.2; 120.1, 120.2) when the mean chance, the same for all the jurymen, is substituted there instead of $u$.

We can precisely represent the possibility of each juryman having a large number of unequally probable chances of faultlessness by imagining that the list from which they should be selected is separated into v classes, so that all those contained in one of them have the same chance of faultlessness. Suppose that the first juryman should be selected from a list in which $x_i$ is that chance for one of its classes, and $X_i$ is the ratio of the number of people of that class to their total number in the list.

Then the chance of faultlessness for the first juryman, provided that he belongs to the indicated class, will be $x_i$ with $X_i$ being the probability of his belonging there, i. e., of his chance $x_i$. The same takes place for the other jurymen selected from other lists.

When the jurymen for a session of some assize court[4] are randomly selected from the list of all people who could have been chosen, the mean chances $u, u', u'', \ldots$ are in advance equal one to another. However, for other assize courts their common value can differ. […] In any case, those mean chances existing before the selection of jurymen should not be confused with the chances of faultlessness for the randomly selected jurymen. We return to this essential difference just below.

**123.** If the number v of the possible chances $x_1, x_2, \ldots$ is infinite, the probability of each will be infinitely small. Suppose that $Xdx$ is the probability that the chance of faultlessness of a juryman randomly selected from a given list is $x$. Let also $u$ be the mean of all possible chances when taking into account their respective probabilities. The sum[5], which should be unity, and that, which forms $u$, in accordance with the preceding become definite integrals

$$\int_0^1 Xdx = 1, \quad \int_0^1 xXdx = u. \qquad (123.1, 2)$$

The positive magnitude $X$ can be a continuous or discontinuous function of $x$, absolutely arbitrary but obeying equation (123.1). To each given $X$ there will correspond a completely definite numerical value of $u$ which however corresponds to an infinitely differing expressions of $X$ or different laws of probabilities.

If all the values of $x$ from 0 to 1 are equally possible $X$ will be independent from $x$, and equal to unity for satisfying equation (123.1). In this case, because of equation (123.2) $u = 1/2$. If $X$ increases in the indicated interval in such a manner that the chance of a juryman's faultlessness by itself becomes the more probable the closer it is to certainty and if, in addition, $X$ increases uniformly[6], then

$X = \alpha x + \beta, \alpha > 0, \beta > 0$

and

$$\int_0^1 Xdx = \alpha/2 + \beta = 1, \ \beta = 1 - \alpha/2, \ X = 1 - \alpha/2 + \alpha x,$$

which excludes the inequality α > 2. As a result, $u = 1/2 + \alpha/12$, and the mean chance can not either exceed 2/3, or be less than 1/2; these extreme cases correspond to the values of α = 2 and 0.

Suppose also that, when *x* increases by equal increments, X changes by a geometric progression and obeys equation (123.1) and that

$$X = \frac{\alpha}{e^\alpha - 1} e^{\alpha x}$$

with an arbitrary value of α. As usual, *e* is the base of the Naperian logarithms. We will get

$$u = \frac{1}{1 - e^{-\alpha}} - \frac{1}{\alpha}.$$

Therefore, if α increases from − ∞ to ∞, the mean chance *u* will take all the possible values from 0 to 1: at α = − ∞, 0, and ∞, *u* = 0, 1/2 and 1 respectively.

Suppose that the different chances of faultlessness vary not from 0 to 1, but are contained in narrower limits; for example, the chance *x* can not be lower than 1/2; and, in addition, that if it is higher, all its values are equally possible. Then *X* should be a discontinuous function to be determined in the following way. I denote by ε a positive and finite but quite insensible magnitude and let *fx* be a function very rapidly changing in the interval $x = (1/2 − \varepsilon)$ and 1/2, disappearing within limits $x = 0$ and $(1/2 − \varepsilon)$ and taking a given constant value *g* within the limits $x = 1/2$ and 1.

Under these conditions I assume that $X = fx$; by the nature of this function *fx*

$$\int_0^1 Xdx = g/2 + \int_{1/2-\varepsilon}^{1/2} fxdx.$$

However, owing to the condition (123.1)

$$\int_{1/2-\varepsilon}^{1/2} fxdx = 1 - g/2,$$

so that $g \leq 2$, since *fx* can only be positive.

We can regard the multiplier *x* in the integral

$$\int_{1/2-\varepsilon}^{1/2} xfxdx$$

as a constant equal to 1/2. The integral will therefore be equal to

$1/2 - g/4$. And, noting that

$$\int_0^1 xfxdx = \int_0^{1/2\varepsilon} xfxdx + \int_{1/2\varepsilon}^{1/2} xfxdx + \int_{1/2}^1 xfxdx,$$

we conclude that the left side is equal to $1/2 - g/4 + g/2 - g/8$, so that $u = 1/2 + g/8$.

Therefore, the mean chance can not in this case either exceed $u = 3/4$, which corresponds to the value of $g = 2$, or be less than $1/2$, which corresponds to the value of $g = 0$.

We can thus assume an infinite number of various hypotheses about the type of the function $X$. If one of them is certain, the corresponding value of $u$ will also be doubtless. If, on the contrary, all of them are possible, their respectable probabilities will be infinitely low as also the different values of the mean chance which appears because of those hypotheses. The last case takes place when the various possible values of the chance of a juryman's faultlessness remain unknown and we do not even know the law of their probabilities. All possible hypotheses can be formulated about that law so that the mean chance will take unequally probable values.

Denote by $\varphi u du$ the infinitely low probability that that chance exactly equals $u$. The function $\varphi u$ can be continuous or discontinuous, its integral over $[0, 1]$ should be equal to unity, and all remarks made about $X$ are applicable to it as well.

**124.** The preceding formulas completely answer all the questions pertaining to the subject of this chapter if only the prior probability $k$ of the defendant's guilt is known as well as the probability of faultlessness of each juryman and kind of cases. If this chance takes many possible values, all of them and their probabilities should be given. Again, if there are infinitely many of those values, and the probability of each is infinitely low, we ought to know the function which expresses the law of their probabilities.

However, none of these necessary elements is known in advance. Before the accused appears in court, he is committed for trial, and the pertinent procedures doubtless lead to his guilt being more probable than his innocence. Consequently, there are grounds to believe that $k > 1/2$, but by how much? We are unable to know it in advance. All depends on the ability and severity of the magistrates charged with conducting the preliminary investigation and can change with the change of the kind of cases. Neither, either before or after their random selection, can we find out the chances of the jurymen's faultlessness. For each of them it depends on his enlightenment, on the advisability which he attaches to repression of one or another type of criminality, on pity inspired in him by the age or sex of the accused etc. None of these circumstances is known to us, neither do we know their quantitative influence on the votes. For applying the previous formulas, we ought to eliminate the unknown elements from them, and with that problem we will occupy ourselves now.

**125.** Consider the case in which the chance of faultlessness is the same for all the jurymen. Suppose that before the judgement it is unknown and can take all the possible values from 0 to 1 and that the infinitely low probability of its value *u*, is denoted by φ*udu*. If this value is certain, that is, if the chance of faultlessness of each juryman is indeed equal to *u*, the formula (117.1) expresses the probability that the accused, whether guilty or not, will be convicted by a majority verdict of ($n - i$) votes against *i*. In this formula, *n* is the total number of jurymen, and *k*, the prior probability of the defendant's guilt. The probability of such a division of votes will indeed by determined by the right side of that formula multiplied by φ*udu*. And if such a division took place, the probability that the common chance of faultlessness is equal to *u* will be the mentioned product divided by the sum of its values for all of the *u* from 0 to 1 (§ 43). Thus, when denoting by $w_i du$ this infinitely low probability, we get

$$w_i = \frac{[ku^{n-i}(1-u)^i + (1-k)u^i(1-u)^{n-i}]\varphi u}{\int_0^1 [ku^{n-i}(1-u)^i + (1-k)u^i(1-u)^{n-i}]\varphi u du}.$$

The factor $N_i$ in formula (117.1) does not depend on *u* and is cancelled here from the numerator and denominator. If $\lambda_i$ is the probability that the chance *u* of faultlessness is contained within given limits *l* and *l'*, then

$$\lambda_i = \frac{k\int_l^{l'} u^{n-i}(1-u)^i \varphi\, udu + (1-k)\int_l^{l'} u^i(1-u)^{n-i}\varphi\, udu}{k\int_0^1 u^{n-i}(1-u)^i \varphi\, udu + (1-k)\int_0^1 u^i(1-u)^{n-i}\varphi\, udu}. \quad (125.1)$$

When *n* is even and the votes are equally divided, $n = 2i$ and

$$\lambda_i = \frac{\int_l^{l'} u^i(1-u)^i \varphi\, udu}{\int_0^1 u^i(1-u)^i \varphi\, udu}.$$

Thus, $\lambda_i$ does not depend on *k*, but depends on it when the votes are not divided equally. When any two values of *u*, equally remote from the extreme values 0 and 1 or from the mean $u = 1/2$, are equally probable, so that $\varphi(1-u) = \varphi u$, then

$$\int_0^1 u^{n-i}(1-u)^i \varphi\, udu = \int_0^1 u^i(1-u)^{n-i}\varphi\, udu$$

If in addition $l < 1/2$ and $l' = 1 - l$, then also

$$\int_l^{l'} u^{n-i}(1-u)^i \varphi\, udu = \int_l^{l'} u^i(1-u)^{n-i}\varphi\, udu,$$

and formula (125.1) becomes

$$\lambda_i = \frac{\int_l^{1-l} u^{n-i}(1-u)^i \varphi\, udu}{\int_0^1 u^{n-i}(1-u)^i \varphi\, udu}$$

so that $\lambda_i$ is again independent from $k$, whichever was the division of the votes, $(n-i)$ and $i$.

Assume that in formula (125.1) $l = 1/2$ and $l' = 1$ and denote the result by $\lambda'_i$, then

$$\lambda'_i = \frac{k\int_{1/2}^1 u^{n-i}(1-u)^i \varphi\, udu + (1-k)\int_{1/2}^1 u^i(1-u)^{n-i}\varphi\, udu}{k\int_0^1 u^{n-i}(1-u)^i \varphi\, udu + (1-k)\int_0^1 u^i(1-u)^{n-i}\varphi\, udu}$$

will be the probability that the chance $u$ is contained within the limits 1/2 and 1, and thus exceeds 1/2. Suppose now that $l = 0$ and $l' = 1/2$, and denote the appearing form of (125.1) by $\lambda''_i$, then the probability that $u < 1/2$ is [Poisson repeated formula (125.1) with the integral in the numerator being over [1/2, 1].] And, since the probability of an exact equality $u = 1/2$ is infinitely low, $\lambda'_i + \lambda''_i = 1$. This can be verified at once by noting that the denominators of both magnitudes coincide and the sum of the integrals multiplied by $k$ in the numerator is equal to the integral multiplied by $k$ in the denominator, and that the same is true for the integrals multiplied by $(1-k)$.

**126.** Suppose that the chance of faultlessness of each juryman is $u$ and that the probability that the accused is convicted by a majority verdict of $(n-i)$ votes against $i$ is expressed by the magnitude $w_i du$. Then the probability that that accused is guilty, provided that the indicated chances are certainly $u$, is expressed by the magnitude $p_i$ from § 119.

By the rules of §§ 5 and 10 the composite probability of the defendant's guilt is the integral of $p_i w_i du$ over [0, 1]. Denoting it by $\chi_i$ we have

$$\chi_i = \frac{k\int_0^1 u^{n-i}(1-u)^i \varphi\, udu}{k\int_0^1 u^{n-i}(1-u)^i \varphi\, udu + (1-k)\int_0^1 u^i(1-u)^{n-i}\varphi\, udu} \quad (126.1)$$

where account was taken of the expressions of $w_i$ и $p_i$. This probability equals 0 and 1 together with $k$. When representing it as

$$\chi_i = k + \frac{k(1-k)\left(\int_0^1 u^{n-i}(1-u)^i \varphi\, udu - \int_0^1 u^i(1-u)^{n-i}\varphi\, udu\right)}{k\int_0^1 u^{n-i}(1-u)^i \varphi\, udu + (1-k)\int_0^1 u^i(1-u)^{n-i}\varphi\, udu},$$

it is seen that at all other values of $k$ the probability of the defendant's guilt after the verdict was returned becomes higher or lower than before when the first integrals in the preceding formula are larger or smaller than the second. If they are equal to each other, which invariably occurs when $n = 2i$ and $\varphi(1 - u) = \varphi u$, then $\chi_i = k$. Indeed, the probability of the defendant's guilt can not change at all when the votes were equally divided or when the values of $u$ and $(1 - u)$ or the expressions $1/2 \pm (1 - 2u)/2$ of the chance of faultlessness are equally probable.

In all other cases $\chi_i$ depends not only on the majority of $m$ or $(n - 2i)$ votes and on the value of $k$, as $p_i$ does, but also on the total number of the jurymen $n$ and the law of probabilities of the chances of faultlessness as expressed by the function $\varphi u$.

Suppose that conviction by 201 jurymen was returned by a majority verdict of only one single vote, or, in another case, by one single juryman, and that the chance of faultlessness was the same for all of them. The probabilities of the propriety of the decisions will be exactly equal to each other. However, if in the first case that chance essentially differed from 1/2, the observed event will be an exceptional and very rare fact; and if that chance is 1/2, the probability of the first event, according to the expression of $w_i$ from § 119, will be somewhat higher than 1/9. And, finally, if the chance of faultlessness was unknown in advance, and is determined by the decision itself, the defendant's guilt is much less probable in the first case than in the second. Indeed, the result of voting, 101:100, considered by itself is not less correct than the decision of one single juryman, but almost an equal division of votes means that the chance of faultlessness of the jurymen likely not much differed from 1/2, certainly because the case was difficult.

**127.** For expressing a precise idea about the meaning of formulas (125.1) and (126.1), we ought to suppose that before a case was decided, someone had grounds to believe with probability $k$ that the accused is guilty without knowing either the case, or anyone of the $n$ jurymen except that they were selected from the general list. For him, the probability of their faultlessness is the same for all of them (§ 122), but unknown. Before the verdict is returned, it was possible to suppose that that unknown $u$ took any possible value from 0 to 1. By some considerations, which we do not at all examine, the infinitely low probability that that person attaches to the variable $u$, is $\varphi udu$ with a given function $\varphi u$. Its integral over [0, 1] should be unity, since the value of $u$ is doubtless contained within those very limits.

Then it becomes known that the accused was acquitted by $i$ votes and convicted by $(n - i)$ of them. This is additional information, and for the mentioned person the probability $\lambda_i$ that the chance $u$ of faultlessness of each juryman will now be contained within the limits $l$ and $l'$. The grounds for believing that the accused is guilty is also strengthened or weakened because its prior probability $k$ became equal to $\chi_i$. For another person having other information those grounds will differ, but they should not be confused with the chance itself of guilt. It depends on $k$ and on the chance of faultlessness proper to each voting juryman according to his capacities and the essence of the case under consideration. Had the numerical values $u, u', u'', …$ of that chance been known for each juryman as well as $k$, the veritable chance of the defendant's guilt after the decision could have been calculated by the rules of § 116 generalized on the case of $n$ jurymen. However, a determination of these values in advance is impossible, and the indicated rules are inapplicable.

Suppose that it is only known that the accused was convicted by a majority of at least $m$ or $(n - 2i)$ votes, that is, by $m, m + 2, m + 4, …$ votes or convicted unanimously, and denote by $Y_i$ the probability that the chance of faultlessness, common for all jurymen, is contained within the limits $l'$ and $l$, and by $Z_i$, the probability that the convicted accused is guilty. Then it will be possible to derive these two magnitudes by the same reasoning as made use of for deriving $\lambda_i$ и $\chi_i$, but applying $c_i$ and $P_i$ (§§ 118 and 120) rather than $\gamma_i$ and $p_i$, as it was done when establishing formulas (125.1) and (126.1). We will thus get

$$Y_i = \frac{k\int_0^1 U_i \varphi u du + (1-k)\int_l^{l'} V_i \varphi u du}{k\int_0^1 U_i \varphi u du + (1-k)\int_0^1 V_i \varphi u du}, \qquad (127.1a)$$

$$Z_i = \frac{k\int_0^1 U_i \varphi u du}{k\int_0^1 U_i \varphi u du + (1-k)\int_0^1 V_i \varphi u du}. \qquad (127.1b)$$

These formulas, like formulas (125.1) и (126.1), can be generalized on the case in which it is known that $n'$ jurymen out of $n$ were randomly selected from the first list, $n''$, from another etc. if supposing also that the mean chance $u'$ of faultlessness in the first list had probability $\varphi' u' du'$, that in the second list, $u''$ and $\varphi'' u'' du''$ respectively etc. However, such a generalization is both easy and useless, and we will not write out the pertinent complicated formulas[7].

**128.** When $i$ and $(n – i)$ are very large numbers, we ought to apply the method of § 67 for calculating approximate values of the integrals contained in formulas (125.1, 126.1, 127.1). I begin by considering the first two. In the interval between $u = 0$ and 1 the product $u^{n-1}(1 – u)^i$ has only one maximal value. I denote it by $\beta$, and let $\alpha$ be the corresponding value of $u$. Then

$$\alpha = \frac{n-i}{n}, \quad \beta = \frac{i^i(n-i)^{n-i}}{n^n}, \quad u^{n-i}(1-u)^i = \beta\exp(-x^2).$$

Passing over to logarithms, I get

$$x^2 = \lg\beta - (n-i)\lg u - i\lg(1-u).$$

The variable $x$ continuously increases from $-\infty$ to $\infty$ and the values $x = -\infty, 0, \infty$ correspond to $u = 0, \alpha, 1$. The limits of the integral are $x = \pm\infty$ at $u = 0$ and $1$. In general, denoting by $\lambda$ и $\lambda'$ the limits of $x$ corresponding to $u = l$ и $l$, and accounting for the preceding values of $\beta$ and $x^2$ we will obtain

$$\lambda = \pm\sqrt{(n-i)\ln[(n-i)/ln] + i\ln[i/n(1-l)]},$$
$$\lambda' = \pm\sqrt{(n-i)\ln[(n-i)/l'n] + i\ln[i/n(1-l')]}.$$

If $l$ и $l'$ exceed $\alpha$, then $\lambda$ и $\lambda'$ will become positive, and the superior signs should be applied, and otherwise, the inferior. And if $l < \alpha$ and $l' > \alpha$, the superior sign should precede the second radical, and the inferior sign, the first of them, so that $\lambda$ will be negative, and $\lambda'$, positive.

For expanding $u$ in power series of $x$, we assume that the coefficients $\gamma, \gamma', \gamma'', \ldots$ are constant and that

$$u = \alpha + \gamma x + \gamma' x^2 + \gamma'' x^3 + \ldots$$

Taking into account the values of $\alpha$, $\beta$ and $x^2$, we will get

$$x^2 = \frac{n^3}{2i(n-i)}(\gamma x + \gamma' x^2 + \gamma'' x^3 + \ldots)^2 + \frac{n^4(n-2i)}{3i^2(n-i)^2}(\gamma x + \gamma' x^2 + \gamma'' x^3 + \ldots)^3 + \ldots$$

Equating the coefficients of the same powers of $x$ on both sides of this equation (? - O.S.), we can derive the values of $\gamma, \gamma', \gamma'', \ldots$, so that

$$u = \alpha + x\sqrt{\frac{2i(n-i)}{n^3}} - \frac{2x^2(n-2i)}{3n^2} + \ldots,$$
$$du = \sqrt{\frac{2i(n-i)}{n^3}}dx - \frac{4x(n-2i)}{3n^2} + \ldots$$

If the function $\varphi u$ does not decrease very rapidly when $u$ deviates in some direction from its particular value $u = \alpha$, it will be possible, after substituting $u$ as expressed by the series, to expand $\varphi u$ as well in powers of $(u - \alpha)$, then in powers of $x$. We will thus arrive at

$$\varphi u = \varphi\alpha + [x\sqrt{\frac{2i(n-i)\varphi\alpha}{n^3}}\ldots]\frac{d}{d\alpha} + \frac{2}{2}(x\sqrt{\frac{)i\ n-i\ \varphi\alpha}{n^3}}\ldots]\frac{d^2}{d\alpha^2} + \ldots$$

By applying these various values, the series expansion of the integral of $u^{n-1}(1-u)^i\varphi u du$ over $[0, 1]$ from formula (126.1) will contain integrals over $[-\infty, \infty]$ of $\exp(-x^2)dx$, multiplied by even and odd powers of *x*. The value of the integrals containing the even powers of *x*, are known, whereas the others will disappear. The numbers *i* and (*n* – *i*) are of the same order as *n* is, and the series under consideration will be arranged by magnitudes of the order of smallness of $1/\sqrt{n}$, $1/(n\sqrt{n})$, $1/(n^2\sqrt{n})$, … Restricting our calculations to its first term and noting that

$$\int_{-\infty}^{\infty} \exp(-x^2)dx = \sqrt{\phantom{n}}$$

we establish that

$$\int_0^1 u^{n-i}(1-u)^i \varphi\, udu = \frac{i^i(n-i)^{n-i}\sqrt{2\pi i(n-i)}}{n^{n+1}\sqrt{n}}\varphi(\frac{n-i}{n}).$$

In addition, interchanging *i* and (*n* − *i*), we get

$$\int_0^1 u^i(1-u)^{n-i}\varphi\, udu = \frac{i^i(n-i)^{n-i}\sqrt{2\pi i(n-i)}}{n^{n+1}\sqrt{n}}\varphi(\frac{i}{n}).$$

Denote by δ a positive magnitude, very small as compared with $\sqrt{n}$, let

$$l = \frac{n-i}{n} - \delta\sqrt{\frac{2i(n-i)}{n^3}},\ l' = \frac{n-i}{n} + \delta\sqrt{\frac{2i(n-i)}{n^3}}$$

and, neglecting terms of the order of smallness of $1/\sqrt{n}$, expand the logarithms in the expressions for λ and λ′. It will then occur that λ = − δ and λ′ = δ. After that we have

$$\int_l^{l'} u^{n-i}(1-u)^i \varphi\, udu = \frac{i^i(n-i)^{n-i}\sqrt{2i(n-i)}}{n^{n+1}\sqrt{n}}\varphi(\frac{n-i}{n})\int_{-\delta}^{\delta}\exp(-x^2)dx$$

to within terms of the order of 1/*n*. As δ increases, the integral on the right side will tend to $\sqrt{\pi}$. For its being very little different from that value, it is sufficient for δ to equal such numbers as 2 or 3. And if the limits *l* и *l'* become appreciably larger or smaller than (*n* – *i*)/*n*, the integral on the left side will almost disappear.

Denote by ε a positive magnitude, very small as compared with $\sqrt{n}$, and let

$$l = \frac{i}{n} - \varepsilon\sqrt{\frac{2i(n-i)}{n^3}},\ l' = \frac{i}{n} + \varepsilon\sqrt{\frac{2i(n-i)}{n^3}}.$$

Then

$$\int_l^{l'} u^i(1-u)^{n-i}\varphi\, udu =$$

$$\frac{i^i(n-i)^{n-i}\sqrt{2i(n-i)}}{n^{n+1}\sqrt{n}}\varphi(\frac{i}{n})\int\limits_{-\varepsilon}^{\varepsilon}\exp(-x^2)dx\,. \qquad (128.1)$$

For the limits $l$ и $l'$ becoming appreciably larger or smaller than $i/n$, the integral on the left side will almost disappear.

If the fractions $(n – i)/n$ and $i/n$ noticeably differ from each other, the first of the preceding values of $l$ and $l'$ will also differ from the value $u = i/n$, corresponding to the maximal value of $u^i(1 – u)^{n–i}$, and the integral on the left side of equation (128.1) will little differ from zero. The subsequent values of $l$ and $l'$ will also appreciably differ from the value $u = (n – i)/n$, which corresponds to the maximal value of $u^{n–i}(1 – u)^i$, and the integral of $u^{n–i}(1 – u)^i \varphi u$ over $[l, l']$ will also almost disappear.

**129.** When substituting in formula (126.1) approximate values of the integrals included there and cancelling common factors in the numerator and denominator, we will obtain the probability of the guilt of the accused convicted by a majority of $m = n - 2i$ votes if the number $n$ of jurymen is very large:

$$\chi_i = \frac{k\varphi[(n-i)/n]}{k\varphi[(n-i)/n]+(1-k)\varphi(i/n)}.$$

As is seen, it depends on the ratio $i/n$, or, we can also say, on the ratio $(n – i)/n$, but not on the difference of these numbers as the probability $p_i$ does if (§ 119) the chance $u$ of faultlessness was known in advance. Thus, if the accused is convicted by 1000 votes against 500 when the total number of jurymen is 1500 or by 100 votes against 50 when that number is 150, the probability $\chi_i$ will be the same, but the probabilities $p_i$ will very much differ. On the contrary, suppose that the second case persists but that in the first case 550 jurymen convicted the defendant and 500 acquitted him. Then the probabilities $p_i$ will be equal to each other, but the probabilities $\chi_i$ can much differ one from another.

If the accused is convicted; then $(n – i)/n > 1/2$ and $i/n < 1/2$. Supposing that at $u < 1/2$ the function $\varphi u$ is almost zero; in other words, if the mean chance of faultlessness, when lower than 1/2 or lower than the chance of error, is considered unlikely; and if, in addition, the fraction $k$ is not too close to zero, so that it will be possible to neglect the second term of the denominator of $\chi_i$ as compared with the first term, − then $\chi_i = 1$ or at least is very close to certainty.

When applying approximate values of the integrals in formula (125.1) and assuming that the fractions $(n – i)/n$ and $i/n$ do not differ too little from each other, we will get the probability

$$\lambda_i = \frac{\frac{1}{\sqrt{\pi}}k\varphi(\frac{n-i}{n})\int\limits_{-\delta}^{\delta}\exp(-x^2)dx}{k\varphi(\frac{n-i}{n})+(1-k)\varphi(\frac{i}{n})}$$

that, when deciding to convict the accused, the common chance $u$ of faultlessness is contained within the limits

$$\frac{n-i}{n} \mp \delta\sqrt{\frac{2i(n-i)}{n^3}}.$$

Retaining the hypothesis under which one of the two integrals in the numerator of formula (125.1) disappears, we have

$$\lambda_i = \frac{\frac{1}{\sqrt{\pi}} k\varphi(\frac{i}{n}) \int_{-\varepsilon}^{\varepsilon} \exp(-x^2) dx}{k\varphi(\frac{n-i}{n}) + (1-k)\varphi(\frac{i}{n})}$$

for the probability that that chance is contained within the limits

$$\frac{i}{n} \mp \varepsilon\sqrt{\frac{2i(n-i)}{n^3}}.$$

Sufficient, but not excessively large values can be assigned to the magnitudes $\delta$ and $\varepsilon$ so that the integrals with variable $x$ very little differ from $\sqrt{\pi}$. Then the sum of both values of $\lambda_i$ will also very little differ from unity. The mean chance will almost certainly be contained either within the first limits little differing from the fraction $(n-i)/n$ which exceeds 1/2, or within the second limits little differing from $i/n$ which is less than 1/2.

When supposing that $\varphi(i/n)$ is insensible or negligible as compared with $\varphi[(n-i)/n]$, the second case is eliminated, and we can assume that almost certainly the value of $u$ very little differs from $(n-i)/n$. In other words, the chances $u$ and $(1-u)$ of faultlessness and error are in the ratio of the convicting and acquitting votes, $(n-i)$ and $i$. And now it seems that the probability $\chi_i$ is not sensibly reduced to unity, but very little differs from $p_i$ at $u = (n-i)/n$. However, it should be remarked that the probability $p_i$ corresponds to the case in which the chance $u$ definitely takes only one possible value. For including this case in that to which $\chi_i$ is corresponding, we ought to suppose that $\varphi u$ only disappears within infinitely short intervals on both sides of the possible value of $u$, and that this function very rapidly decreases in the vicinity of that value. However, the analysis in § 128 essentially supposed, as it was seen, that $\varphi u$ does not change at all in any direction from the value of $u = (n-i)/n$. The probability $\chi_i$ derived there is not at all applicable to the case to which corresponds $p_i$ from § 119. It can also be remarked that formula (126.1) describes the latter.

When in general denoting by $v$ the only possible value of $u$, and by $\eta$, a positive infinitely small magnitude, and supposing that $\varphi u$ is a function disappearing at all values of $u$ beyond the limits $v \pm \eta$, the limits of the integrals in formula (126.1) will coincide with these extreme points. Within them the factors $u^{n-i}(1-u)^i$ and $u^i(1-u)^{n-i}$ are constant; after withdrawing them beyond the sign of the integrals and cancelling the integral of φ*udu* over [*v* − η, *v* + η], which will be a common multiplier in the numerator and denominator of formula (126.1), this formula will coincide with formula (119.1) at *u* = *v*.

Supposing that the fractions (*n* – *i*)/*n* and *i*/*n* do not essentially differ from each other, and assuming that ε = δ, the preceding values of λ$_i$ will correspond to the same limits of the chance *u*. However, their common value will differ from the preceding and become independent from *k* and equal to

$$\frac{1}{\sqrt{\pi}} \int_{-\delta}^{\delta} \exp(-x^2) dx$$

since in this particular case both integrals in the numerator of formula (125.1), as well as the integrals in its denominator, are almost equal to each other.

**130.** For determining approximate values of the integrals in formulas (127.1) we should express the values of $U_i$ and $V_i$ by means of formulas (121.1). The first of them concerns the case in which (1 − *u*)/*u* > *i*/(*n* + 1 − *i*), and the second, the contrary case. The first formula therefore subsists within the interval *u* = 0 and α, and the second, within *u* = α and 1, when (*n* + 1) is replaced by *n* and (*n* – *i*)/*n* = α. In accordance with the equation which determines θ in formulas (121.1),

$$u^{n-i}(1-u)^i = \frac{i^i(n-i)^{n-i}}{n^n} \exp(\theta^2).$$

This coincides with the preceding equation connecting *u* and *x* [§ 128] from which it follows that

$$u = \alpha + \theta\sqrt{\frac{2i(n-i)}{n^3}} - \frac{2\theta^2(n-2i)}{3n^2} + ...$$

However, θ should invariably be positive (§ 121), and, at *u* = 0, α and 1, it equals ∞, 0 и ∞. The variable *u* increases from 0 to α, and θ decreases from ∞ to 0. Then *u* increases once more from α to 1, and θ increases from 0 to ∞. Therefore, in accordance with formulas (121.1),

$$\int_0^\alpha U_i \varphi u\, du = \frac{1}{\sqrt{\pi}} \int_0^\infty \left( \int_\theta^\infty \exp(-x^2) dx \right) \varphi u \frac{du}{d\theta} d\theta +$$

$$\frac{(n+i)\sqrt{2}}{3\sqrt{n i(n-i)}} \int_{-\infty}^0 \exp(\theta^2) \varphi u \frac{du}{d\theta} d\theta,$$

$$\int_{0\alpha}^1 U_i \varphi u\, du = \int_0^1 \varphi u\, du - \frac{1}{\sqrt{\pi}} \int_0^\infty \left( \int_\theta^\infty \exp(-x^2) dx \right) \varphi u \frac{du}{d\theta} d\theta +$$

$$\frac{(n+i)\sqrt{2}}{3\pi\sqrt{n(i\ n-i)}}\int_0^\infty \exp(\theta^2)\varphi\ u\frac{du}{d\theta}d\theta.$$

The values of these integrals as well as the double integrals can be obtained in the form of converging series of $\theta$ when substituting the preceding series instead of $u$ and its differential coefficient instead of $du/d\theta$, and expanding $\varphi u$ into a series. This supposes that that function does not change very rapidly in any direction from the particular value of the argument $u = \alpha$. When neglecting terms of the order of smallness of $1/n$, then, simply,

$$u = \alpha + \theta\sqrt{\frac{2i(n-i)}{n^3}},\quad \frac{du}{d\theta} = \sqrt{\frac{2i(n-i)}{n^3}},\quad \varphi u = \varphi\alpha.$$

The square root can be preceded here by two signs. Plus should be chosen in the integrals [above] if the variable $\theta$ increases, and minus otherwise. Then the signs of the integrals are changed, and their limits interchanged so that

$$\int_0^\alpha U_i\varphi u\,du = \varphi\alpha\sqrt{\frac{2i(n-i)}{\pi n^3}}\int_0^\infty\left(\int_\theta^\infty \exp(-x^2)dx\right)d\theta,$$

$$\int_\alpha^1 U_i\varphi u\,du = \int_\alpha^1 \varphi u\,du - \varphi\alpha\sqrt{\frac{2i(n-i)}{\pi n^3}}\int_0^\infty\left(\int_\theta^\infty \exp(-x^2)dx\right)d\theta.$$

Adding up these formulas, we get

$$\int_0^1 U_i\varphi u\,du = \int_\alpha^1 \varphi u\,du.$$

In the general case, denoting by $a$ and $a_1$ two such values of $u$ that $a < \alpha$ and $a_1 > \alpha$, and by $b$ and $b_1$, positive values of $\theta$, corresponding to values $u = a$ and $a_1$, it will occur to within our order of approximation[8] that

$$\int_a^\alpha U_i\varphi u\,du = \varphi\alpha\sqrt{\frac{2i(n-i)}{\pi n^3}}\int_0^b\left(\int_\theta^\infty \exp(-x^2)dx\right)d\theta,$$

$$\int_\alpha^{a_1} U_i\varphi u\,du = \int_\alpha^{a_1}\varphi u\,du - \varphi\alpha\sqrt{\frac{2i(n-i)}{\pi n^3}}\int_0^{b_1}\left(\int_\theta^\infty \exp(-x^2)dx\right)d\theta.$$

When integrating by parts, we can obtain […] and therefore

$$\int_a^\alpha U_i\varphi u\,du = \varphi\alpha\sqrt{\frac{2i(n-i)}{\pi n^3}}\left(b_1\int_b^\infty \exp(-x^2)dx + \frac{1}{2} - \frac{1}{2}\exp(-b^2)\right),$$

$$\int_{\alpha}^{a_1} U_i \varphi u\, du = \int_{\alpha}^{a_1} \varphi u\, du -$$

$$\varphi\alpha\sqrt{\frac{2i(n-i)}{\pi n^3}}\left(b_1\int_{b_1}^{\infty} \exp(-x^2)dx + \frac{1}{2} - \frac{1}{2}\exp(-b_1^2)\right).$$

In formulas (121.1) I replaced $U_i$ by $V_i$ and therefore replaced $u$ by $(1-u)$ (§ 118). The first substitution takes place if $u/(1-u)$ exceeds $i/(n+1-i)$, i. e. in the interval from $u = 1 - \alpha$ to $1$; I had also replaced $(n+1)$ by $n$ and again assumed that $\alpha = (n-i)/n$.

The second substitution concerned the case in which $u$ changed from $0$ to $1 - \alpha$. Denoting by $\theta'$ the value of $\theta$ when $u$ is replaced by $(1-u)$ and neglecting, as previously, terms of the order of smallness of $1/n$, I first arrive at [...], so that

$$\int_0^1 V_i \varphi u\, du = \int_0^{1-\alpha} \varphi u\, du.$$

Now, if $a'$ and $a'_1$ are such values of $u$, that $a' < 1 - \alpha$ and $a'_1 > 1 - \alpha$, and denoting by $b'$ and $b'_1$ the positive values of $\theta'$, which are derived from the equation

$$(1-u)^{n-i}u^i = \frac{i^i(n-i)^{n-i}}{n^3}\exp(\theta'^2),$$

and correspond to values of $u = a'$ and $a'_1$, then also

$$\int_{1-\alpha}^{a'_1} V_i \varphi u\, du =$$

$$\varphi(1-\alpha)\sqrt{\frac{2i(n-i)}{\pi n^3}}\left(b'_1\int_{b'_1}^{\infty} \exp(-\theta'^2)d\theta' + \frac{1}{2} - \frac{1}{2}\exp(-b'^2_1)\right),$$

$$\int_{a'}^{1-\alpha} V_i \varphi u\, du =$$

$$\int_{a'_1}^{1-\alpha} \varphi u\, du - \varphi(1-\alpha)\sqrt{\frac{2i(n-i)}{\pi n^3}}\left(b'\int_{b'}^{\infty} \exp(-\theta'^2)d\theta' + \frac{1}{2} - \frac{1}{2}\exp(-b'^2)\right).$$

**131.** The approximate values of the integrals in formulas (127.1) are thus calculated and we obtain

$$Z_i = \frac{k\int_{\alpha}^{1} \varphi u\, du}{k\int_{\alpha}^{1} \varphi u\, du + (1-k)\int_0^{1-\alpha} \varphi u\, du} \qquad (131.1)$$

for the probability of the guilt of the accused convicted by at least
$(n - i)$ jurymen out of their very large number $n$. The ratio $\alpha = (n - i)/n$
exceeded 1/2. If at $u < 1/2$ the function $\varphi u$ is insensible or disappears,
the same will happen with the second integral in the denominator of
(131.1). And if $k$ is not a very small fraction, $Z_i$ will almost equal
unity. And if $\varphi(1 - u) = \varphi u$ for all values of $u$, then

$$\int_0^{1-\alpha} \varphi u \, du = -\int_\alpha^{0} \varphi(1-u) \, du = \int_\alpha^{1} \varphi u \, du$$

and $Z_i$ becomes equal to $k$, which should have indeed taken place.

If $a = 1 - \alpha$ and $a'_i = \alpha$, the corresponding values $b$ and $b'_1$ of $\theta$ and $\theta'$ will be equal to each other. Denote them by $c$ and take into account the meaning of $\alpha$, then $c$ will be a positive magnitude determined by the equation

$$(n - i)^i i^{n-1} = i^i (n - i)^{n-i} \exp(-c^2).$$

We will also get

$$\int_{1-\alpha}^{\alpha} U_i \varphi u \, du = \varphi \alpha \sqrt{\frac{2i(n-i)}{\pi n^3}} \left( c \int_c^\infty \exp(-x^2) \, dx + \frac{1}{2} - \exp(-c^2) \right),$$

and therefore

$$\int_{1-\alpha}^{\alpha} V_i \varphi u \, du =$$

$$\varphi(1-\alpha) \sqrt{\frac{2i(n-i)}{\pi n^3}} \left( c \int_c^\infty \exp(-x^2) \, dx + \frac{1}{2} - \exp(-c^2) \right).$$

If the accused had been convicted, the probability that the chance $u$ of faultlessness, common for all the jurymen, is contained within the limits $(1 - \alpha)$ and $\alpha$, or within $i/n$ and $(n - i)/n$, is equal to

$$Y_i = \frac{k \varphi \alpha + (1-k) \varphi(1-\alpha)}{k \int_\alpha^1 \varphi u \, du + (1-k) \int_0^{1-\alpha} \varphi u \, du} \times$$

$$\left( c \int_c^\infty \exp(-x^2) \, dx + \frac{1}{2\pi} - \exp(-c^2) \right) \sqrt{\frac{2i(n-i)}{n^3}}.$$

It is feeble owing to the very small radical. It is therefore likely that the chance $u$ was either greater than $\alpha$ or smaller than $(1 - \alpha)$. For verifying this, I assume that $a_1 = a'_1 = 1$. The corresponding values of $\theta$ and $\theta'$ will be $b_1 = b'_1 = \infty$, and

$$\int_{\alpha}^{1} U_i \varphi u\, du = \int_{\alpha}^{1} \varphi u\, du - \frac{1}{2\pi}\varphi\alpha\sqrt{\frac{2i(n-i)}{n^3}},$$

$$\int_{1\alpha}^{1} V_i \varphi u\, du = \frac{1}{2\pi}\varphi(1-\alpha)\sqrt{\frac{2i(n-i)}{n^3}}.$$

Subtract from this latter integral its previous value over [(1 − α), α], and then

$$\int_{\alpha}^{1} V_i \varphi u\, du = \varphi(1-\alpha)\sqrt{\frac{2i(n-i)}{\pi n^3}}\left(\exp(-c^2) - c\int_{c}^{\infty}\exp(-x^2)dx\right).$$

When taking into account the integrals of $U_i \varphi u\, du$ и $V_i \varphi u\, du$ over [α, 1]

$$Y_i = \frac{k\int_{\alpha}^{1}\varphi u\, du - \left(\frac{1}{2}k\varphi\alpha - (1-k)\varphi(1-\alpha)\Omega\right)}{k\int_{\alpha}^{1}\varphi u\, du + (1-k)\int_{0}^{1\alpha}\varphi u\, du},$$

$$\Omega = [\exp(-c^2) - c\int_{c}^{\infty}\exp(-x^2)dx]\sqrt{\frac{2i(n-i)}{\pi n^3}}$$

will be the probability that the chance $u$ is contained within the limits $u = \alpha$ and 1, and thus exceeds α. I also assume that $a = a' = 0$, then $b = b' = \infty$, so that […]

$$Y_i = \frac{(1-k)\int_{0}^{1\alpha} \varphi u\, du - \sqrt{\frac{2i(n-i)}{\pi n^3}}\Gamma}{k\int_{\alpha}^{1}\varphi u\, du + (1-k)\int_{\alpha}^{1\alpha}\varphi u\, du},$$

$$\Gamma = \left(\frac{1}{2}(1-k)\varphi(1-\alpha) - k\varphi\alpha[\exp(-c^2) - c\int_{c}^{\infty}\exp(-x^2)dx]\right)$$

will be the probability that the chance $u$ is contained within the limits $u = 0$ and $(1 - \alpha)$, and is thus less than $(1 - \alpha)$. The sum of the two last values of $Y_i$ is almost unity, which should be verified. If the values of $\varphi u$ at $u < 1/2$ are zero or almost zero the last value of $Y_i$ will be very low, and the preceding value very little differing from certainty. In any case, the sum of the three just calculated values of $Y_i$ is unity, as it should have indeed been.

**132.** Owing to the preceding, even if the number $n$ of the jurymen is very large, for deriving the probability of the guilt of the accused convicted by $(n - i)$ votes against $i$ we ought to formulate a hypothesis about the type of the function $\varphi u$ or the law of probabilities of the chances of faultlessness. In usual cases, when the number $n$ is not very large, this is all more necessary. Laplace's pertinent hypothesis consisted in assuming that the function $\varphi u$ is zero at all values of $u < 1/2$ and constant at all $u > 1/2$. It means that each chance of faultlessness lower than the contrary chance is assumed impossible, and that all the other chances of faultlessness are considered equally probable. His hypothesis is admissible since it satisfies the condition of the integral of φ$u$d$u$ over [0, 1] being equal to unity; in addition (§ 123), the mean of the possible values of $u$, or the same integral of $u$φ$u$d$u$, will be contained within the limits of 1/2 and 3/4, and, at $u > 1/2$, will depend on the value of φ$u$.

By that hypothesis φ$u$ = 0 at $u < 1/2$ and is constant at $u > 1/2$, and therefore the interval between the limits of the integral in formula (126.1) shortens and is contained between $u = 0$ and 1/2. The function φ$u$ can be taken out beyond the sign of the integral, and since

$$\int_{1/2}^{1} u^i (1-u)^{n-i} du = \int_{0}^{1/2} u^{n-i} (1-u)^i du,$$

that formula becomes

$$\chi_i = \frac{k \int_{1/2}^{1} u^{n-i}(1-u)^i du}{k \int_{1/2}^{1} u^{n-i}(1-u)^i du + (1-k) \int_{0}^{1/2} u^{n-i}(1-u)^i du}.$$

Here, the constant common multiplier φ$u$ in the numerator and denominator is cancelled.

Laplace did not at all account for the prior probability of the defendant's guilt. Accordingly, for the derived formula to coincide with his result we ought to suppose that this guilt is not either more, or less probable than his innocence, i. e., to assume that $k = 1/2$. Then

$$\chi_i = \frac{\int_{1/2}^{1} u^{n-i}(1-u)^i du}{\int_{0}^{1} u^{n-i}(1-u)^i du}, \quad 1 - \chi_i = \frac{\int_{0}^{1/2} u^{n-i}(1-u)^i du}{\int_{0}^{1} u^{n-i}(1-u)^i du}.$$

After the integration it occurs that the probability the innocence of the accused convicted by ($n - 2i$) jurymen out of their total number $n$ is

$$1 - \chi_i = \frac{1}{2^{n+1}}[1 + \frac{n+1}{1} + \frac{(n+1)n}{2!} + \ldots + \frac{(n+1)n\ldots(n-i+2)}{i!}]. \quad (132.1)$$

Laplace (1816/1886, c. 527) actually derived this formula. The sum in brackets consists of ($i + 1$) terms and is unity at $i = 0$ and therefore the probability of a mistaken unanimous conviction is $1/2^{n+1}$. Without admitting the value $k = 1/2$ and assuming that $i = 0$, we will get

$$\lambda_i = \frac{1-k}{k2^{n+1}-(2k-1)} = \frac{1}{2^{n+1}}\left[1 - \frac{(2k-1)(2^{n+1}-1)}{k2^{n+1}-(2k-1)}\right]$$

which is larger or smaller than $1/2^{n+1}$ at $k$ larger or smaller than 1/2.

In ordinary cases $n = 12$; from formula (132.1), supposing successively that $i = 0, 1, 2, …, 5$, we will obtain fractions

1/8192, 14/8192, 92/8192, 378/8192, 1093/8192, 2380/8192, −

the probabilities of mistaken conviction by 11 jurymen against 1; by 10 against 2; …; by 7 against 5. When the majority is least, the probability of a mistake is almost 2/7, which means that, among a very large number of accused convicted by such a majority, 2/7 were likely convicted mistakenly; for a majority of 8 against 4 the indicated probability is almost 1/8.

Let us apply the Laplace hypothesis to formula (125.1) and denote by δ a positive magnitude not exceeding 1/2. Supposing that $k = 1/2$, $l = 1/2$ and $l' = 1/2 + δ$, we get

$$\lambda_i = \frac{\int_{1/2\delta}^{1/2\delta} u^{n-i}(1-u)^i du}{\int_0^1 u^{n-i}(1-u)^i du}$$

for the probability that the chance $u$ of faultlessness which, by that hypothesis, can not be lower than 1/2, is contained within the limits 1/2 and 1/2 + δ when conviction was decided by $(n – i)$ votes against $i$. Integration is not difficult. At $i = 0$ or unanimity it occurs that

$\lambda_i = [(1/2) + \delta]^{n+1} - [(1/2) - \delta]^{n+1}$.

If, for example, $n = 12$ and $δ = 0.448$, then almost exactly $\lambda_i = 1/2$, so that even money can be bet on the chance $u$ to be contained between the limits 0.5 and 0.948. At $δ = 1/4$, without supposing that $i = 0$, we obtain

$$\lambda_i = \frac{1}{4^{n+1}}[3^{n+1} - 1 + \frac{n+1}{1}(3^n - 3) + C_{n+1}^2(3^{n-1} - 3^2) + ...$$
$$+ C_{n+1}^i(3^{n-i+1} - 3^i)]$$

for the probability that the chance $u$ is contained within the limits of 1/2 and 3/4, or that it is closer to 1/2, than to unity. At $n = 12$ and $i = 5$ it equals 0.915 … In cases of the least majority we can bet a little more than 10 against 1 on that chance being lower than 3/4.

**133.** Formula (132.1) was derived from another one in which the chance $u$ of faultlessness had been assumed the same for all the jurymen. Therefore, although Laplace did not mention it, that chance can not be admitted for each of them as their proper chance. That should be the mean chance for the entire list from which the jurymen are randomly selected (§ 122). People whose chance of faultlessness, at least in difficult cases, is lower than 1/2, i. e., lower than the contrary chance will doubtless be included in such a list. Laplace's hypothesis demands, however, that there will be an inconsiderable number of such people and that, consequently, they will not prevent the mean chance of faultlessness to be invariably higher than 1/2. The illustrious geometer also supposed that all values of that chance from 1/2 to 1 were equally probable.

He only justified these two assumptions by stating that *The opinion of the judge more strongly tends to the truth than to error*. However, issuing from this principle, he only concluded that the values of the function φ$u$, expressing the law of probabilities of the values of the mean chance, should be larger at $u > 1/2$ than otherwise. This condition can be fulfilled in an infinite number of ways without requiring that either φ$u$ = 0 at $u < 1/2$, or its invariability at $u > 1/2$. His hypothesis is not sufficiently justified in advance, whereas its corollaries, as it will be shown, render it quite inadmissible.

Actually, the formula (132.1), which is one of its necessary corollaries, does not include anything depending on the abilities of the people from the general list of jurymen. Someone who only knows, for example, that two convictions were returned by the same majority verdict and the same total number of jurymen selected from two different lists, could have had the same grounds for believing that both convictions were mistaken in spite of knowing that people included in one list were much more able than those selected from the other one. And it is impossible to agree anymore with this conclusion.

If the accused is convicted by a majority verdict of ($n - i$) votes against $i$, so that $i/n < 1/2$, then, by the hypothesis under consideration, the value of φ($i/n$) is zero or almost zero. When the number $n$ of the jurymen is very large, the probability $\chi_i$ of the defendant's guilt is therefore very close to unity whatever is the difference $[(n - i) - i]$ (§ 129). If, for example, 520 jurymen against 480 convicted him, his guilt should be considered almost certain, although 480 jurymen denied it with their chance of faultlessness, as can be supposed, not differing from that of the 520 others. This conclusion is sufficient for rejecting the hypothesis on which it was based. Indeed, no one will assign considerable confidence to such a decision, and the less so, the same confidence as to an almost unanimous decision of a 1000 jurymen.

According to that hypothesis, when the ability of people included in the general list change or when they differ in different regions or different kinds of cases, the probability of those chances of faultlessness which are closer to unity or differ from 1/2 less than the others, will heighten in the same ratio. This does not really occur; when for some reason this ability strengthens, the probability of the chances of faultlessness closest to certainty heightens, whereas the contrary happens to chances remotest from unity. By adopting a function φ$u$, which can satisfy these conditions, and in addition should not disappear or become almost zero at $u < 1/2$, it will be possible to overcome the indicated difficulties. However, that is not sufficient: an infinite number of continuous and discontinuous kinds of φ$u$ satisfy those conditions; for the same number $n$ of jurymen and the same difference $[(n - i) - i]$ they lead to very unequal values of the probability $\chi_i$ expressed by formula (126.1).

And so, when knowing those numbers for one single convicting decision and supposing that the prior probability $k$ is 1/2 or any other fraction, we will be unable, as said above, to determine the real probability of the virtue of that decision. It depends on the chance of faultlessness proper to each juryman which we can not know. In addition, for someone, who only knows that the $n$ jurymen were randomly selected from the general list and has grounds to believe that the virtue of a decision only depends on the mean chance of faultlessness common to all the jurymen from that list (§ 122), it is impossible to calculate the value of that probability. Indeed, such a calculation requires a formulation of a particular hypothesis about the law of probabilities of the values of the mean chance when it changes from 0 to 1. It should not be the Laplace hypothesis, or any other if insufficiently justified.

If the jurymen selected from the general list make only one decision, the previous formulas are useless. The same will happen if the number of decisions is inconsiderable; on the contrary, a very large number of convictions and acquittals in known ratios are pronounced by jurymen successively and randomly selected from the same general list. And we will show at once that the formulas (117.1, 117.2; 118.1; 119.1, 119.2; 120.1, 120.2) were based on these grounds. They contain only two unknown magnitudes, $k$ and $u$, and therefore only require results of two observations. And we will first of all determine these results.

**134.** The general list of citizens who can become jurymen contains a number of people. Each jury panel consists of $n$ jurymen randomly selected from that list for one year or many years, and they judge a very large number $\mu$ of accused. Denote by $a_i$ the number of those accused convicted by a majority verdict of at least $(n - i)$ votes against $i$. This means that $i$ is zero or a positive number less than $n/2$. In advance, the prior chance of such a conviction should change from one decision to another, but anyway the mean of its unknown values calculated for $\mu$ decisions will likely be almost equal to the ratio of $a_i/\mu$ (§ 95).

Furthermore, the values of this mean chance and that ratio very little change with $\mu$, supposed to be very large. And if that number increases further and further, they will indefinitely tend to a special constant $R_i$ and reach it if $\mu$ can become infinite without any changes in the various causes of convictions by the indicated majority verdict. This special constant is the sum of the chances provided by all possible causes for conviction, or the event under consideration, multiplied by the respective probabilities of those very causes (§ 104).

It is impossible to enumerate them and calculate their influence in advance, but we do not need to know them. It is sufficient to suppose that neither their probabilities, nor the chances they attach to convictions are changing, and the observation itself will let us know whether this assumption conforms to reality. If so, then, denoting by $a'_i$ the number of convictions decided by a majority verdict of at least $(n - i)$ votes against $i$ for another very large number $\mu'$ of accused, the difference $(a'_i/\mu' - a_i/\mu)$ will likely be a very small fraction (§ 109). If, on the contrary, this difference is not very small, we can justifiably believe that in the interval between the two series of decisions occurred some essential change in the causes of conviction. Calculation can only tell us about the existence of such changes, but will not say anything about its essence.

What we say about convictions by at least $(n - i)$ votes against $i$ equally concerns convictions by $(n - i)$ votes against $i$. Denote by $b_i$ the number of those when there are $\mu$ accused. There exists a special constant $r_i$, to which the ratio $b_i/\mu$ indefinitely tends as $\mu$ increases and reaches it if $\mu$ can become infinite without any changes in the causes of conviction. And if $b'_i$ is the number of convictions for $\mu'$ accused, the difference $(b'_i/\mu - b_i/\mu)$ will likely be a very small fraction. The following equalities will evidently take place:

$a_i = b_i + b_{i-1} + b_{i-2} + \ldots + b_0,$
$a'_i = b'_i + b'_{i-1} + b'_{i-2} + \ldots + b'_0,$
$R_i = r_i + r_{i-1} + r_{i-2} + \ldots + r_0.$

Suppose that $\alpha$ is a positive magnitude very small as compared with $\sqrt{\mu}$ and $\sqrt{\mu'}$ and

$$P = 1 - \frac{2}{\sqrt{\pi}} \int_{\alpha}^{\infty} \exp(-x^2)dx.$$

In accordance with the formulas of § 112 the magnitude $P$ will also be the probability that the limits of the differences $(a'_i/\mu' - a_i/\mu)$ and $(b'_i/\mu' - b_i/\mu)$ for the unknowns $R_i$ and $r_i$ are

$$\frac{a_i}{\mu} \mp \alpha \sqrt{\frac{2a_i(\mu - a_i)}{\mu^3}}, \quad \frac{b_i}{\mu} \mp \alpha \sqrt{\frac{2b_i(\mu - b_i)}{\mu^3}} \qquad (134.1a, b)$$

and again for them

$$\mp \alpha \sqrt{\frac{2a_i(\mu - a_i)}{\mu^3} + \frac{2a'_i(\mu' - a'_i)}{\mu'^3}}, \qquad (134.1c)$$

$$\mp \alpha \sqrt{\frac{2b_i(\mu - b_i)}{\mu^3} + \frac{2b'_i(\mu' - b'_i)}{\mu'^3}}. \qquad (134.1d)$$

Other things being equal, the extent of these limits shortens as $\mu$ and $\mu'$ increase almost proportional to the square roots of these large numbers. Indeed, $a_i$ and $b_i$ increase almost like $\mu$, and $a'_i$ and $b'_i$, almost like $\mu'$. And this extent narrows as much as $\alpha$ decreases, but their probability $P$ then lowers.

**135.** As I indicated in § 7 of the Preamble, all the numerical data below were picked out from the *Comptes généraux de l'Administration de la justice criminelle*[9] published by the government. From 1825 to 1830 inclusive the number of cases yearly passed over for consideration by jury panels in the entirety of France was

   5121, 5301, 5287, 5721, 5506, 5068,

and the number of the accused in these criminal cases was

   6652, 6988, 6929, 7396, 7373, 6962,

or almost 7 defendants in 5 cases yearly. During those years the number of convicted by majority verdicts of at least 7 votes against 5 was

   4037, 4348, 4236, 4551, 4475, 4130,

so that the conviction rates were

   0.6068, 0.6222, 0.6113, 0.6153, 0.6069, 0.5932.

It is already seen that during those six years with an invariable criminal legislation these yearly rates changed very little.

I accept the total number of the accused during those six years as $\mu$ and as $a_5$, the corresponding number of convictions: $\mu = 42{,}300$, $a_5 = 25{,}777$. The corresponding limits (134.1a) are $0.6094 \pm 0.00335\alpha$; assuming, for example, that $\alpha = 2$, we get $P = 0.9953$, the probability, very close to certainty, that the unknown $R_5$ only differs from 0.6094 by 0.0067.

When separating the period under consideration in two equal parts, then, for them,

   $\mu = 20{,}569$, $\mu' = 21{,}731$, $a_5 = 12{,}621$, $a'_5 = 13{,}156$,
   $a_5/\mu = 0.6136$, $a'_5/\mu' = 0.6054$, $a_5/\mu - a'_5/\mu' = 0.0082$.

The limits (134.1c) of the derived difference will be $\pm 0.00671\alpha$; supposing that $\alpha = 1.2$, we obtain the limits $\pm 0.00805$ and $P = 0.9103$, $(1 - P) = 0.0897$. We can therefore bet almost 10 against 1 on that difference to be contained within limits $\pm 0.00805$. Although the actual difference, $\pm 0.0082$, somewhat exceeds the derived limits, both the deviations and the probability that this should not have taken place are not sufficiently essential for justifiably believing that some appreciable change occurred in the pertinent causes.

In 1831 the number of the accused increased to 7606, and the number of those convicted reached 4098. The law demanded that the majority verdict for conviction should not be less than 8 votes against 4. And it accordingly occurred that

   $\mu = 7606$, $a_4 = 4098$, $a_4/\mu = 0.5388$.

If the required majority is excluded and the other causes influencing the decisions of the jurymen remained during that year as they were before, the ratio $b_5/\mu$ will be obtained by subtracting $a_4/\mu$ from $a_5/\mu$, or by subtracting 0.5388 from the derived above fraction 0.6094, and then $b_5/\mu = 0.0706$.

For verifying this result I note that in 1825 − 1830 the law stipulated an intervention of judges composing the assize[10] courts each time when the jurymen's decision was adopted by the least majority verdict of 7 votes against 5. And the *Comptes généraux* state that during the five last years of those six the number of such interventions

398, 373, 373, 395, 372

barely changed. It occurred 1911 times in all, but the corresponding number of the accused was not indicated. Accordingly, for those years, we ought to compare the number of cases rather than of the accused. Their total number amounted to 26,883, so that $\mu = 26{,}883$, $b_5 = 1911$ and $b_5/\mu = 0.0711$, very little differing from the previous result[11].

This accord of two values of $b_5/\mu$ proves that in 1831 the probabilities $u$ and $k$, on which depends the indicated ratio, remained almost the same as in the previous years. We should nevertheless point out that the calculation of the last value was based on the hypothesis that the number of convictions by 7 votes against 5 is to the number of the accused as the number of cases in which this majority took place is to the total number of cases. This can not be verified in advance, since the *Comptes généraux* lacks necessary data.

In 1832 and 1833 the number of the accused without those in political cases amounted to 7555 and 6964. An essential difference between these numbers took place owing to the new legislative measure according to which in 1833 many kinds of cases were given over from the assizes to the police courts. The number of convictions decided, as in 1831, by majority verdicts of at least 8 votes against 4 increased to 4105 and 4448, so that $a_4/\mu = 0.5887$ and $0.5895$.

It is seen that these ratios little differ from each other, but their mean, 05888, exceeds $a_4/\mu = 0.5388$ for 1831 by 0.05 or by about 1/10 of this value. If no changes occurred in the causes which could have influenced the jurymen's voting, then, after taking into account the limits (134.1c) and their probability $P$, this would have been absolutely unlikely. And indeed, from 1832 criminal legislation underwent such a change: the question of *mitigating circumstances* was brought before the jury panels. If present, they led to diminution of penalty which rendered convictions easier and increased their number.

**136.** Various ratios which we provided just above for the entirety of France are not the same in all parts of the kingdom. However, it turned out that, excluding the department of Seine and some others, the number of criminal cases judged during a few years was not enough for deriving with a sufficiently high probability and for the jurisdiction of each assize court that constant to which the rate of conviction should tend. Here are the results for the Paris court.

In 1825 − 1830 the yearly numbers of the accused were

802, 824, 675, 868, 908, 804,                                  (136.1)

and the numbers of those convicted

   567, 527, 436, 559, 604, 484,                         (136.2)

so that the rates of conviction were

   0.7070, 0.6396, 0.6459, 0.6440, 0.6652, 0.6020.

Let $\mu$ and $a_5$ be the sums of the former, (136.1), and the latter, (136.2), six numbers, then $\mu = 4881$, $a_5 = 3177$, $a_5/\mu = 0.6509$. For France in its entirety the first two of these derived numbers during the same years were 42,300 and 25,779, and we established that that rate should very little differ from 0.6094 [= 25,779: 42,300]. This is less than the previous value [0.6509 for the department of Seine] by 0.0416 or about 1/15 of its value. However, taking into account the limits (134.1c) and their probability $P$, this would have been absolutely unlikely if only some particular cause did not lead there to easier convictions than in other parts of France. What was this cause? This is what calculus can not tell us. Anyway, for the same interval of time in that department, whose population barely amounts to 1/36 part of the population of France, the number of the accused exceeded 1/9 of them for France as a whole, which means 4 times as many. So the repression of criminality is more necessary there, and perhaps for this reason local jurymen are more severe[12].

When issuing from those values of $\mu$ and $a_5$ the limits (134.1a) become $0.6509 \pm 0.00965\alpha$, and, at $\alpha = 2$, $P = 0.99532$, $(1 − P) = 0.00468$. You can bet more than 200 against 1 on the unknown $R_5$ only differing in any direction from 0.6509 by 0.0193.

The last of the four ratios quoted above, 0.6020, is appreciably less than the mean of the five others. It makes sense to investigate whether this difference indicates in a sufficient measure some particular cause that could have in 1830 led the jurymen to be less severe than in the previous years.

Denote the numbers of the accused and the convicted in the department of Seine for 1825 – 1829 by $\mu$ and $a_5$, and the same numbers for 1830 by $\mu'$ and $a'_5$. Then

   $\mu = 4077$, $a_5 = 2693$, $\mu' = 804$, $a'_5 = 484$,
   $a_5/\mu = 0.6605$, $a'_5/\mu' = 0.6019$, $a_5/\mu − a'_5/\mu' = 0.0585$.

The limits (134.1c) become $\pm 0.02657\alpha$; if $\alpha = 2$, you can bet more than 200 against 1 on the difference $(a_5/\mu − a'_5/\mu')$ not to exceed 0.05314. Actually, however, it exceeds that fraction almost by 1/10 of its value, and we can believe that at that time [1830] some anomaly happened in the jurymen's votes. The cause that rendered them a bit less severe could have been the Revolution of 1830. Whatever it was, it apparently influenced the jurymen over all France, since the rate of conviction in the kingdom lowered almost to 0.59 whereas its mean value for the 5 previous years was 0.61.

During 1826 – 1830 inclusive the number of criminal cases in the department of Seine amounted to 2963. The majority verdict of 7 votes against 5 was returned 194 times, and the [assize] court had to intervene. With $b_5/\mu = 194/2963 = 0.0655$ we obtain a value somewhat smaller than for France in its entirety.

**137.** When, like in the *Comptes généraux*, separately considering all the kinds of crimes brought before the assize courts, the numbers of the accused and the convicted will not be large enough for the ratios to become invariable and serve as the basis of our calculations. However, all criminal cases in those *Comptes* are also separated in two categories, *crimes against the person*, and *against property*. These large categories yearly lead to very different ratios, almost invariable by themselves. We will cite them.

During the six years 1825 – 1830 the number of the accused in the entirety of France in those categories amounted to

1897, 1907, 1911, 1844, 1791, 1666;
4755, 5081, 5018, 5552, 5582, 5296.

The corresponding numbers of convictions for the same criminal legislation amounted to

882, 967, 948, 871, 834, 766;
3155, 3381, 3288, 3680, 3641, 3364,

and the rates of conviction were

0.4649, 0.5071, 0.4961, 0.4723, 0.4657, 0.4598;
0.6635, 0.6654, 0.6552, 0.6628, 0.6523, 0.6352.

It is seen that in each category these rates did not change much from year to year, but that the latter appreciably exceeded the former. Denote by $\mu$ and $\mu'$ and $a_5$ and $a'_5$ the total number of the accused and the convicted, then

$\mu = 11{,}016$, $a_5 = 5268$, $\mu' = 31{,}284$, $a'_5 = 20{,}509$,
$a_5/\mu = 0.4782$, $a'_5/\mu' = 0.6556$.

The second rate exceeds the first one by a little more than 1/3 of its value. Issuing from these numbers, we can conclude that the limits (134.1a) of the unknown $R_5$ are respectively

$0.4782 \pm 0.00675\alpha$ and $0.6556 \pm 0.00380\alpha$.

It follows that for $\alpha = 2$ with probability very close to certainty $R_5$ does not differ from 0.4782 by more than 0.0135, and from 0.6556 by more than 0.0076.

In 1831 a majority verdict of at least 8 votes against 4 was required for conviction. And, accordingly,

$\mu = 2046$, $a_4 = 743$, $\mu' = 5560$, $a'_4 = 3355$,
$a_4/\mu = 0.3631$, $a'_4/\mu' = 0.6034$.

When subtracting these ratios from the previous it occurs that the rates of conviction with a majority verdict of at least 7 votes against 5 were $b_5/\mu = 0.1151$ and $b'_5/\mu' = 0.0522$. It is remarkable that the former concerning crimes against the person is almost twice larger than the latter, whereas, on the contrary, $a'_5/\mu'$ exceeds $a_5/\mu$ almost by 1/3. Thus, convictions in the second category of crimes are not only proportionally more numerous, but are also returned by a larger majority.

For both sexes, the rates under consideration are not the same at all. In the assize courts women almost invariably composed 0.18 of all the yearly accused. During the five years 1826 – 1830 inclusive, the respective numbers for women amounted to

$\mu = 1305$, $\mu' = 5465$, $a_5 = 586$, $a'_5 = 3312$,
$a_5/\mu = 0.4490$, $a'_5/\mu' = 0.6061$.

When comparing the calculated rates with their preceding values, it occurs that they are smaller, although only by about 1/16 and 1/12 of their values.

In 1832 and 1833, when convictions were returned not less than by 8 votes against 4 with *mitigating circumstances* being allowed for, the numbers of the accused and the convicted men and women were

$\mu = 4108$, $\mu' = 10,421$, $a_4 = 1889$, $a'_4 = 6664$,
$a_4/\mu = 0.4598$, $a'_4/\mu' = 0.6395$.

Like above, the letters with strokes concern crimes against property. At $\alpha = 2$ the limits (134.1a) indicate that with probability very close to certainty the unknown $R_5$ deviates from 0.4598 not more than by 0.022, and from 0.6395 not more than by 0.0133 respectively. We can remark that the ratio $a_4/\mu : a'_4/\mu'$ was almost equal to the derived above ratio $a_5/\mu : a'_5/\mu'$. When comparing $a_4/\mu$ and $a'_4/\mu'$ with those magnitudes for 1831, we can also indicate that the influence of the *mitigating circumstances* only increased the rate $a'_4/\mu'$ for crimes against property by 1/15, but $a_4/\mu$, for crimes against the person, was increased almost by 1/3 of its value.

**138.** In § 122 we established that the chance of an accused to be convicted by jurymen randomly selected from the general list of a department or an assize court would have been the same had all of them a common chance of faultlessness. Therefore, when convictions are decided by majority verdicts of at least $(n - i)$ votes against $i$, that chance is expressed by formula (118.1a), and by formula (117.1) for majority verdicts of $(n - i)$ votes against $i$. For each department and category of crime the magnitudes $c_i$ and $\gamma_i$, expressed by these formulas are those to which the rates $a_i/\mu$ and $b_i/\mu$ indefinitely tend as $\mu$, proposed very large, increases further. In other words, when considering cases of the same category and the same department separately for the accused men and women, $c_i$ and $\gamma_i$ will coincide with the unknowns $R_i$ and $r_i$ (§ 134).

As done above, we will separate all kinds of criminal cases into two distinct categories, crimes against the person and against property, and letters with a stroke will concern the latter. However, to prevent the calculations from being too complicated, we will not at all account for the sex of the accused. Its influence on the rate of conviction can be neglected since more than 5/6 of the total number of the accused consisted of men. For each department we will get

$$a_i/\mu = c_i,\ b_i/\mu = \gamma_i,\ a'_i/\mu' = c'_i,\ b'_i/\mu' = \gamma'_i, \qquad (138.1\text{a, b, c, d})$$

with the approximation being the better and the probability the higher, the larger were the numbers $\mu$ and $\mu'$.

If for different departments the ratios on the left sides of equations (138.1) are known, those four equations will suffice for determining the unknowns $k$ and $u$, included in the magnitudes $c_i$ and $\gamma_i$ and their analogues $k'$ and $u'$ in $c'_i$ and $\gamma'_i$. At present however, the need to have very considerable $\mu$ and $\mu'$ makes it impossible to apply separately the equations (138.1) to each department. We have to suppose that as a rule the unknowns $u, u', k, k'$ little change from one of them to another so that the numbers concerning the whole of France can be assumed as the left sides of those equations.

The magnitudes $u$ and $u'$ thus determined will exactly coincide with the chances of faultlessness had the lists of the jurymen of all the departments been combined together for randomly selecting each juryman from it. Since the magnitudes $k$ and $k'$ also depend on the ability of magistrates directing preliminary investigations, they can differ in different departments. However, the equations (138.1) are linear with respect to these unknowns and their derived values are the means of those actually taking place in all the departments. Finally, I should note that the need to be content with these general values of $u$, $u', k$, and $k'$ is only due to the lack of complete observational data rather than some imperfection in the described theory.

The magnitudes $c_i$ and $\gamma_i$ do not change if $k$ and $u$ are replaced by $(1 - k)$ and $(1 - u)$, see §§ 117 и 118. Therefore, if $k$ and $u$ larger than 1/2 correspond to the given $a_i/\mu$ and $b_i/\mu$ and satisfy the equations (138.1a, b), then there exist other $k$ and $u$, also satisfying them but smaller than 1/2. However, we ought to suppose that the prior probability of the guilt of the accused exceeds the probability of his innocence, and that the mean chance of a juryman's faultlessness exceeds 1/2. We should therefore apply the values of $k$ and $u$ which are larger than 1/2 and reject the other values as alien for the problem.

The same remark concerns the equations (138.1c, d) and the values of $k'$ and $\mu'$ derived from them. Nevertheless, when applying these equations to the judgements returned in the numerous political trials during the ill-fated years of the Revolution, it will be possible, as described in § 12 of the Preamble, to apply their roots smaller than 1/2, because in those times the prior legal innocence of the accused could have been more probable than their guilt, and the probability of voluntary mistakes made by the jurymen could have exceeded their chance of faultlessness.

**139.** I suppose that in the formulas (117.1) and (118.1) $n = 12$ and $i = 5$. The coefficients included there will be

$N_0 = 1$, $N_1 = 12$, $N_2 = 66$, $N_3 = 220$, $N_4 = 495$, $N_5 = 792$.

I also suppose that

$$\frac{a_5}{\mu} = c, \quad \frac{b_5}{\mu} = 792\gamma, \quad u = \frac{t}{1+t}, \quad 1-u = \frac{1}{1+t}$$

so that the equation (138.1b) becomes

$$\gamma = \frac{(t^2+1)t^5}{2(1+t)^{12}} + \frac{(2k-1)(t^2-1)t^5}{2(1+t)^{12}}. \tag{139.1}$$

Since

$$U_5 = 1 - 924u^6(1-u)^6 - V^5,$$

the equation (138.1) can be written as

$$c = k[1 - \frac{924t^6}{(1+t)^{12}}] -$$
$$\frac{(2k-1)}{(1+t)^{12}}[1 + 12t + 66t^2 + 220t^3 + 495t^4 + 792t^5]. \tag{139.2}$$

Equations (139.1) and (139.2) apply to crimes against the person. The corresponding equations for the second category of crimes are established when replacing magnitudes $c$, $\gamma$, $k$, $t$ by $c'$, $\gamma'$, $k'$, $t'$. The unknown $t$ can take all values from $t = 0$ when $u = 0$ to $\infty$ when $u = 1$. However, we should only admit the values $t > 1$ corresponding to $u > 1/2$. Then, the unknown $k$ should be contained within limits $1/2$ and $1$, and therefore, see equation (139.1), $t$ should be such that its limits could be determined by the inequalities

$$\gamma > \frac{(t^2+1)t^5}{2(1+t)^{12}}, \quad \gamma < \frac{t^7}{(1+t)^{12}}. \tag{139.3a, b}$$

The right side of the inequality (139.3a) continuously decreases from $t = 0$ to $\infty$. In inequality (139.3b) it increases from $t = 1$ to $7/5$, then decreases to $t = \infty$.

By eliminating $k$ from equations (139.1) and (139.2), a reciprocal equation of the 24[th] degree in $t$ can be derived and reduced to an equation of the 12[th] degree. However, it is much simpler directly to calculate simultaneously $k$ and $u$ satisfying equations (139.1) and (139.2) by successive approximations.

**140.** For the six years 1825 – 1830 we have

$c = 0.4782$, $\gamma = 0.1151/792 = 0.0001453$.

At $t = 2$ the fraction in inequality (139.3a) exceeds this value of $\gamma$, and at $t = 3$ $\gamma$ exceeds the fraction in (139.3b). The value of $t$ should therefore be larger than 2 and smaller than 3. It is easy to become assured that within those indicated limits that unknown has only one possible value. After a few attempts, I adopted the value 2.112, and then equation (139.1) led to $k = 0.5354$. Inserting these values in the right side of equation (139.2), I find that it is 0.4783, only differing from the left side by 0.0001. Therefore, with a very good approximation

$k = 0.5354, t = 2.112$.

For the same years

$c' = 0.6556, \gamma' = 0.0523/792 = 0.00006604$.

I insert these values in equations (139.1) and (139.2) instead of $c$ and $\gamma$ and replace $t$ and $k$ by $t'$ and $k'$. Solving them as previously, I find with the same degree of approximation that $k' = 0.6744$, $t' = 3.4865$, so that for those years

$$u = \frac{t}{1+t} = 0.6786, \; u' = \frac{t'}{1+t'} = 0.7771$$

will be some juryman's chances of faultlessness in both categories of crimes.

Before the decision of a case someone who does not know either the jurymen or even the place where it is judged, can bet a little more than 2 against 1 and almost 7 against 2 on every juryman to vote properly in those categories. Here, it is usually said, *parier tant contre tant* (bet something against something) for rendering more sensible the signification to be attached to $u$ and $u'$, and also because the proposed bet is illusory since it will never be known who won. The person who knew nothing about the case could have also bet by issuing from the previous values of $k$ and $k'$ and staked somewhat less than 7 against 6 and a little more than 2 against 1 on the guilt of the accused in cases of those categories. Below, we will see the probability of the defendant's guilt after his case was decided.

When considering the number of the accused and convicted without distinguishing between those categories, it will necessary to adopt, again for those years and for France in its entirety,

$c = 0.6094, \gamma = 0.0706/792 = 0.00008914$.

When solving equations (139.1) and (139.2), we will then establish that

$k = 0.6391, t = 2.99, u = 0.7494$.

If, however, separately studying the department of the Seine, we will have to adopt the values of $c$ и $\gamma$ (§ 136)

$c = 0.6509, \gamma = 0.0655/792 = 0.00008267$,

which lead to

$k = 0.678, t = 3.168, u = 0.7778.$

Again in those years, without distinguishing between those categories of crime, the probabilities $k$ and $u$ in the Paris assize court were a little higher than in the rest of the kingdom and somewhat exceeded 2/3 and 3/4. However, the differences of the two values of each of these magnitudes are inconsiderable, which is one of the causes to believe that the value of each of them in one part of France is the same as in another one. This justifies as much as possible the hypothesis about their equality over all the kingdom which we indeed adopted for being able to calculate their approximate values by issuing from a sufficiently large number of observations.

Thus, as I said before, in 1831 the values of $k$ and $u$ or $k'$ and $u'$ remained invariable, but they had to change later along with the ratios derived from them. We know the ratios $a_4/\mu$ и $a'_4/\mu'$ only for 1832 and 1833 which is not sufficient for determining the unknowns $u$ and $k$ or $u'$ and $k'$. These magnitudes possibly changed for the second time after the latest law which, maintaining the problem of *mitigating circumstances*, stipulated secret voting which could have influenced their [!] chance of faultlessness.

After 1831, we therefore are unable to determine the values of $u$ and $k$ or $u'$ and $k'$. However, that law, while establishing that majority verdicts of 7 votes against 5 are sufficient for convictions, required the jury panels to report whether they had decided a case by the minimal majority. If, in future, the *Comptes généraux* will indicate the number of convicted rather than only the number of cases decided by that least majority; if it will become known how many men and women were there among the accused and how many of them were in each category of crime, − if these conditions will be met, in a few years it will become possible to determine quite precisely $u$ and $k$ for the different parts of the kingdom, for men and women and for those categories of crime.

**141.** When knowing the values of $u$ and $k$, formulas (117.1, 117.2, 118.1) allow to determine the probabilities of either conviction or acquittal given a required majority verdict or the least majority verdict. For $n = 12$ and $i = 0$

$\gamma_0 = ku^{12} + (1 - k)(1 - u)^{12}, \delta_0 = (1 - k)u^{12} + k(1 - u)^{12}$

will be those probabilities for unanimous conviction or acquittal. Therefore

$\gamma_0 + \delta_0 = u^{12} + (1 - u)^{12}$

will indeed be the indicated probability. In addition, the magnitude

$\gamma_0 - \delta_0 = (2k - 1)[u^{12} - (1 - u)^{12}]$

is positive since $k > 1/2$ and $u > 1 − u$, so that an unanimous decision of a case is less probable for acquittals than for convictions. If the chance $u$ of faultlessness essentially differs from 0 and 1, these various probabilities will be very low. Adopting, for example, the values of $u$ and $k$ for France in its entirety without distinguishing between the categories of crime, which means that $k = 0.6391$ and $u = 0.7494$, we get

$\gamma_0 = 0.0201$, $\delta_0 = 0.0113$, $\gamma_0 + \delta_0 = 0.0314$.

This suffices for proving how rare should be unanimous decisions by 12 jurymen. If, however, unanimity is required for both conviction and acquittal, then, in accord with this value of ($\gamma_0 + \delta_0$), almost 32 can be bet against 1 on the lack of any decision. And if the jurymen do not communicate with each other, and do not agree to decide by majority, this will indeed happen in about 32 cases out of 33.

Denote by $M$ the probability that one out of $\mu$ cases was decided either unanimously or not. Then

$M = (1 − \gamma_0 − \delta_0)^\mu$,

and if desired that $M = 1/2$, then, maintaining the previous value of the sum ($\gamma_0 + \delta_0$),

$$\mu = \frac{-\lg 2}{\lg(1 - \gamma_0 - \delta_0)} = 21.73.$$

This means that it is possible to bet a little more than even money on at least 1 case out of 22 to be decided unanimously. And that bet will be disadvantageous for 21 cases.

**142.** Before going ahead, it is necessary to recall what was said in § 2 of the Preamble about the sense attached in the decision of a case to the word *guilty* and to derive a few important corollaries.

When announcing that the accused is guilty, the juryman declares that, in his opinion, there is sufficient proof for convicting him. Therefore, the decision that the accused is innocent means that the probability of guilt is insufficiently high for conviction. This latter consideration does not mean that the juryman believes that the accused is innocent; he doubtless oftener thinks that the accused is rather guilty. The probability of guilt sometimes exceeds 1/2 but remains lower than the value at which both the juryman's conscience and public security would have demanded conviction. The real sense of one or another decision made by a juryman consists in that the accused *is, or is not convictable*. The probabilities $P_i$ and $Q_i$ that conviction or acquittal are virtuous (§ 120) also express our grounds for believing that the convicted accused was subject to conviction and the acquitted accused was not.

Therefore, $P_i$ is doubtless lower than the real guilt of the convicted whereas, on the contrary, $Q_i$ is higher than the probability of the innocence of the acquitted accused. Unlike $P_i$ and $Q_i$, which are thus determined and considered for a very large number of decisions of cases of the same kind, those other probabilities can not at all be calculated. Nor should we believe that $P_i$ and $Q_i$ express the general opinion; they express the probabilities of conviction or acquittal by a court consisting of all citizens included in the general list from which the 12 jurymen are randomly selected.

The chance $c_i$ of being convicted by a certain number of jurymen is lower than the fraction denoted by $k$ (§ 118) which as a rule is much lower than $P_i$. Similarly, the chance $d_i$ of acquittal is always lower than the fraction $(1 - k)$, which in turn is much lower than $Q_i$. For the jurymen of each assize court and each of the two categories of crime we should therefore imagine that there exists a definite probability $z$, considered sufficient and necessary for conviction. And the chance $u$ of faultlessness for a juryman randomly selected from a list of his department is equal to the probability with which he decides whether the guilt of the accused is equal to, or exceeds $z$ or not. The probability $u$ mainly depends on the degree of instruction of the class of people included in the general list of jurymen and on the probability $z$; that is, on the opinion formed by them on the necessity of a more or less strong repression of different kinds of crimes. These two differing probabilities can therefore change in time and from on department to another. It is understandable how $u$ can be derived from observation, but we are quite unable to establish $z$. We can only conclude that, other things being equal, it heightens and lowers when noting that the rate of conviction remarkably lowered or heightened. Thus, when the problem of *mitigating circumstances* was posed before the jurymen it occurred that that rate increased from 0.54 to 0.59 (§ 135), and it was necessarily concluded that aside from positive decisions they assumed a lower than previously probability $z$, since the penalties became less severe.

The prior probability of guilt of the accused doubtless much exceeds that which we denoted by $k$. Its maximal value derived by us was almost 3/4, but nevertheless no one will hesitate to bet much more than 3 against 1 on the real guilt of someone brought before an assize court. However, what was stated about $P_i$ equally concerns $k$; and it should also be understood that $k$ only expresses the prior probability that the accused is *convictable*. This probability therefore can depend on the probability which the jurymen require for conviction but which by its nature does not depend on the probability $u$ that some juryman is faultless. It follows that $k$ can change with $z$ even when the forms of preliminary investigation and the ability of those who directed them remain invariable, whatever was the probability $u$. Here is an example of such a change.

In 1814 – 1830, criminal cases in Belgium[13] were decided by tribunals of five judges, and majority verdicts of 3 votes against 2 were sufficient for conviction. In 1830, the composition of the tribunals changed and in 1831 jury panels existing under the French domination were re-established with majority verdicts of 7 votes against 5 becoming necessary for conviction. The forms of preliminary investigation remained as they were.

It follows from the *Comptes de l'administration de la justice criminelle* recently published by the [Belgian] government, that in 1832 − 1834 the rates of conviction in that kingdom were 0.59, 0.60 and 0.61. It is seen that they changed very little and that their mean was almost equal to the mean for France before 1830. However, these *Comptes* say nothing about the number of convictions returned either by the least majority of 7 votes against 5, or by some other definite majority. The quoted rate (? - O.S.) is therefore insufficient for deriving the values of $u$ and $k$ in Belgium. However, the total rate which we denoted by $a_5/\mu$ so little differs in Belgium from the French value that we can believe that the partial rate $b_5/\mu$ in both countries is almost the same. Consequently, $u$ and $k$ in these countries are also almost the same. Therefore, we can admit that in Belgium the value of $k$ does not much differ from the fraction 0.64, which was previously obtained for France in its entirety without distinguishing between the categories of crime.

In the same *Comptes* we find that in 1826 − 1829 those rates increased to the almost equal values of 0.84, 0.85, 0.83 and 0.81 with their mean being a little larger than 0.83. However, in accordance with § 118 the probability whose approximate value is this mean, should invariably be lower than $k$. Therefore, in those years $k$ should have been much higher than in 1832 – 1833. This can only be attributed to the change of the unknown $z$ in those two periods, i. e., to the jurymen demanding for conviction a higher probability that the accused is guilty than the judges thought sufficient. This conclusion is in addition independent from the chance $u$ of faultlessness which possibly changed, and was possibly higher for the judges than for the jurymen or vice versa, a question that remains indecisive owing to the lack of necessary observational data.

Magnitude $k$ depends on the probability $z$, and the inequality of its values in the two categories of crime could have resulted from two different causes. It is more difficult to establish an assumption of prior guilt in cases of crimes against the person and in those cases the jurymen require a higher probability $z$ for conviction. It can be thought that, acting jointly, these different causes had indeed led to the indicated inequality.

It follows from the dependence between $z$ and $k$ that, if in 1832 and 1833 the allowance for *mitigating circumstances* led to a noticeable lowering of the probability which the jurymen believed sufficient for conviction, the probability $k$, on the contrary, should have heightened. The change in $u$ and $k$ in contrary directions should have also led to the heightening of $u$. Indeed, we can suppose that the chance of the jurymen's faultlessness lowered since, on the one hand, they required a lower probability for conviction, and, on the other hand, there existed a higher prior probability of the defendant being convictable.

**143.** We only have to calculate the probabilities of the guilt of the convicted and of the innocence of the acquitted by means of formulas (120.1) and (120.2) and the derived values of $u$ and $k$ or $u'$ and $k'$; more precisely, to calculate the probabilities that in the former but not in the latter case the accused is convictable. At first, however, we will transform those formulas into equations more convenient for calculation and add other formulas whose numerical values [parameters] are also very important to know.

Owing to equation (118.1a) formula (120.1) can be replaced by the equation

$$P_i c_i = kU_i,$$

where the ratio $a_i/\mu$ known from observations is adopted as the approximate value of $c_i$. The magnitude $(1 - P_i)$ is the probability that the accused convicted by a majority verdict of not less than $(n - i)$ votes against $i$ is innocent; $c_i$ is the probability that the accused, whether guilty or not, is convicted by this majority. The product $c_i(1 - P_i)$ therefore expresses the chance of an innocent accused to be nevertheless convicted.

Denoting it by $D_i$ and taking into account the preceding equation and equation (118.1a), we get

$$D_i = (1 - k)V_i.$$

This result was also possible to obtain by considerations applied in § 120 for deriving $P_i$. If the number of votes required for conviction is not less than $(n - i)$, the probability $\Pi_i$ that an acquitted accused is innocent is determined by $Q_i$ or by formula (120.2) when $i$ is replaced by $(n - i - 1)$. Allowing for the equation (118.1b), we conclude that

$$\Pi_i d_{n-i-1} = (1 - k)U_{n-i-1},$$

or, which is the same, in accordance with § 118

$$\Pi_i(1 - c_i) = (1 - k)(1 - V_i).$$

The probability that an acquitted accused is guilty is $(1 - \Pi_i)$, and the probability for the accused not to be convicted is $(1 - c_i)$. Therefore, the product $(1 - \Pi_i)(1 - c_i)$ expresses the chance $\Delta_i$ of a guilty accused to be nevertheless acquitted:

$$\Delta_i = 1 - c_i - (1 - k)(1 - V_i)$$

or, owing to the equation (118.1a),

$$\Delta_i = k(1 - U_i).$$

The chances $D_i$ and $\Delta_i$ are, so to say, measures of the danger to convict an accused not subject to be convicted, and the danger to the society to see a convictable accused acquitted. With regard to the veritable guilt or innocence of the accused we should not forget that $D_i$, unlike $P_i$, is only the superior limit and that $\Delta_i$, unlike $Q_i$, is only the inferior limit. After calculating $P_i$ and $\Pi_i$ the magnitudes $D_i$ and $\Delta_i$ are determined at once since, owing to the previous equations

$$D_i = 1 - k - \Pi_i(1 - c_i), \quad \Delta_i = k - P_i c_i,$$

and it is seen that the chances $\Delta_i$ and $D_i$ are always lower than, respectively, the probabilities $k$ and $(1 − k)$ of the prior guilt and innocence. When the number $\mu$ of the accused is very large, the numbers of convictions and acquittals will be $a_i$ and $(\mu − a_i)$. The numbers of the innocent convicted and acquitted guilty accused will likely almost equal the products $D_i a_i$ and $\Delta_i(\mu − a_i)$.

For $n = 12$ and $i = 5$ and 4, adopting $a_5/\mu$ and $a_4/\mu$ as the approximate values of $c_5$ and $c_4$ and substituting, as done before, $t/(1 + t)$ and $1/(1 + t)$ instead of $u$ and $(1 − u)$, we will derive from the previous equations that

$$\frac{a_5}{\mu} P_5 = k(\Omega + \frac{924t^6}{(1+t)^{12}}), \quad \frac{a_4}{\mu} P_4 = k(\Omega + \frac{924t^6 + 792t^7}{(1+t)^{12}}),$$

$$[1 - \frac{a_5}{\mu}]\Pi_5 = (1-k)\Omega, \quad [1 - \frac{a_4}{\mu}]\Pi_4 = (1-k)(\Omega - \frac{792t^5}{(1+t)^{12}}),$$

$$\Omega = [1 - \frac{1 + 12t + 66t^2 + 220t^3 + 495t^4 + 792t^5}{(1+t)^{12}}].$$

At the same time,

$D_5 = 1 − k − (1 − a_5/\mu)\Pi_5$, $\Delta_5 = k − (a_5/\mu)P_5$,
$D_4 = 1 − k − (1 − a_4/\mu)\Pi_4$, $\Delta_4 = k − (a_4/\mu)P_4$.

So these are the various formulas which [whose parameters] should have been reduced to numbers. The magnitudes included there concern crimes against the person; similar magnitudes concerning the second category of crimes are denoted by the same letters with strokes.

**144.** During 1831 the majority verdict required for conviction was 8 votes against 4 whereas the problem of *mitigating circumstances* did not exist. We had

$a_4/\mu = 0.3632$, $t = 2.112$, $k = 0.5354$,

so that

$P_4 = 0.9811$, $\Pi_4 = 0.7186$, $D_4 = 0.00689$, $\Delta_4 = 0.1791$.

Out of 743 then convicted accused almost 5, as that value of $D_4$ shows, should not have been convicted, and $\Delta_4$ establishes that approximately 233 accused out of 1303 should not have been acquitted. The chance of being convicted although not being convictable a bit exceeded 1/150, and of being acquitted although not being subject to acquittal was contained within the interval of 1/6 and 1/5. Finally, the probability of guilt of the accused did not differ by 1/50 from certainty, and the probability of innocence of the acquitted, that is, of insufficiently proved guilt, a little exceeded 2/3.

These results concerned crimes against the person. For crimes against property the results for the same year were

$a'_4/\mu' = 0.6034$, $t' = 3.4865$, $k' = 0.6744$

so that

$P'_4 = 0.9981$, $\Pi'_4 = 0.8199$, $D'_4 = 0.0004$, $\Delta'_4 = 0.0721$.

Mistakenly convicted were thus only 4 out of 10,000, or less than 2 out of the 3355 announced convictions. The proportion of the convictable acquitted exceeded 7/100, and they should have numbered about 159 out of the 2205 acquittals. The probability of guilt of the convicted differed from certainty less than by 2/1000, and the probability of innocence of the acquitted was a bit higher than 4/5. It is seen that these results are more satisfactory than those concerning the first category of crimes; indeed, the convictions, although proportionally more numerous, were likely returned by a stronger majority (§ 141).

The eight just calculated probabilities $P_4$, $P'_4$, … are based on the ratios $a_4/\mu$, $a'_4/\mu'$, $b_4/\mu$, $b'_4/\mu'$ derived from observation and applied previously for calculating $t$, $t'$, $k$, $k'$. All these eight magnitudes are fractions less than unity which all the more remarkably confirms the theory since, when issuing from arbitrary $t$ and $k$ and $t'$ and $k'$, although not much differing from those derived from observation, this result is not anymore repeated with the same generality.

In the years preceding 1831, a majority verdict of 7 votes against 5 had been sufficient for conviction; however, the assize courts intervened when the majority was [thus] minimal. Conviction was only finalized if a majority of the 5 judges who then composed those courts agreed with the majority of the jurymen. For this reason it is necessary to consider separately the convictions returned by this least majority and by majorities of at least 8 votes against 4. In this second case the values of the probabilities $P_4$ and $P'_4$, $\Pi_4$ and $\Pi'_4$, $D_4$ and $D'_4$, $\Delta_4$ and $\Delta'_4$ are those calculated just above since before and during 1831 the values of $t$ and $k$ and $t'$ and $k'$ had been the same (§ 137).

And so, in 1825 – 1830 about 5000 and almost 20,000 accused were convicted by this majority of at least 8 votes against 4 for crimes in the two categories. According to the previous values of $D_4$ and $D'_4$, about 35 and 8 of them were likely not convictable, doubtless too many if wishing to say that they were really innocent.

Concerning the other convictions announced by the least majority verdicts of 7 votes against 5, the probability of guilt of the accused, see formula (119.1), for $n = 12$, $i = 5$, $u = t/(1 + t)$, $1 - u = 1/(1 + t)$ is

$$p_5 = \frac{kt^2}{kt^2 + 1 - k}.$$

For crimes against the person we have, just like above, $t = 2.112$, $k = 0.5354$ and $p_5 = 0.8372$. For crimes against property, after substituting $p'_5$, $k'$, $t'$ instead of $p_5$, $k$, $t$ and adopting as above $t' = 3.4865$, $k' = 0.6744$, we get $p'_5 = 0.9618$.

Finally, without distinguishing between these categories of crimes, and again for France in its entirety, we should assume that $k = 0.6391$ and $t = 2.99$ (§ 140). Denoting by $w_5$ the corresponding value of $p_5$ of the probability of guilt of the convicted, we have $w_5 = 0.9406$.

Subtracting $p_5$, $p'_5$, $w_5$ from unity, we get almost exactly 0.16, 0.04 and 0.06 for the probability of a mistaken decision in the three cases under consideration. According to the Laplace formula (§ 132), it should have been 0.29, the same in all those three cases, that is, almost twice higher than $(1 - p_5)$ and five times higher than $(1 - w_5)$. In the next section we will see what reduction experiences this probability $(1 - w_5)$ of the innocence of the accused, if his conviction was confirmed by an assize court by a majority of not less than 3 votes against 2.

After combining the data sufficient for determining, in the manner stated in § 140, the values of $k$ and $u$ or $k'$ and $u'$ which existed at the considered time, calculating the corresponding probabilities $P_5$, $\Pi_5$, $D_5$, $\Delta_5$ or $P'_5$, $\Pi'_5$, $D'_5$, $\Delta'_5$ as previously, and comparing them with the probabilities $P_4$, $\Pi_4$, $D_4$, $\Delta_4$ or $P'_4$, $\Pi'_4$, $D'_4$, $\Delta'_4$, which we determined for 1831, we can ascertain without any illusions, for example, the relative advantages of the criminal legislations in those two epochs for both public security and guaranties due to the accused.

As remarked in § 138, two different pairs of values of $k$ and $u$ or $k'$ and $u'$ either larger than 1/2 or smaller than 1/2 correspond to the same observational data and complement the first values to unity. Thus, we found out that in 1831 for crimes against property $k' = 0.6744$, $u' = 0.7771$, but, when applying the same observational data, we can also determine that

$$k' = 1 - 0.6744 = 0.3256, \quad u' = 1 - 0.7771 = 0.2229.$$

The value of $u'$ became $(1 - u')$ and at the same time $t'$ changed into $1/t'$, so that $t' = 1/3.4865 = 0.2868$. Invariably adopting $a'_4/\mu' = 0.6034$, we find that $P'_4 = 0.000675$, since conviction instead of heightening the guilt of the accused, lowered it almost to zero. However, as indicated in § 138, we should in general reject the values of the unknowns $k$ and $u$ or $k'$ and $u'$ smaller than the contrary probabilities. They are nevertheless indicated by the calculations for including the case in which among a very large number of extraordinary judicial decisions the legal guilt of the convicted is less probable than his innocence.

**145.** If adopting in formula (118.1a) $n = 5$, $i = 2$, $u = t/(1 + t)$, $1 - u = 1/(1 + t)$, the probability that the accused will be convicted by a tribunal of five judges by a majority verdict not less than of 3 votes against two is

$$c_2 = k - \frac{(2k-1)(1 + 5t + 10t^2)}{(1+t)^5}.$$

Here, as always, $k$ denotes the prior probability of guilt of that accused, and $u$ is the chance that no judge is mistaken. Because of formula (120.1) we also have

$$c_2 P_2 = k[1 - \frac{1 + 5t + 10t^2}{(1+t)^5}]$$

or simpler, when applying the previous equation,

$$(2k - 1)c_2 P_2 = k(k - 1 + c_2)$$

for determining the probability $P_2$ of the guilt of the accused after his conviction.

Let us apply these equations to an accused convicted by the least majority verdict of 7 votes against 5 and then brought before an assize court, as it was done before 1831. Here, magnitude $k$ is the probability of his guilt resulting from the jurymen's decision; an approximate and likely value of $c_2$ is established by observation and is equal to the rate of conviction by an assize court for a very large number of convictions.

We see in the *Comptes généraux* that in 1826 – 1830 the assize courts over the kingdom received 1911 cases after the accused were convicted by a majority verdict of 7 votes against 5. Those convictions were confirmed in 1597 cases, but the *Comptes* do not indicate how these numbers are distributed among the categories of crime, and we have to determine $P_2$ and the unknown $t$ without distinguishing between them. Thus, $c_2 = 1597/1911 = 0.8357$. Adopting the value of $w_5$ from § 144 as $k$, that is, assuming that $k = 0.9406$ which exceeded, as it should be (§ 118), the proportion $c_2$ of convictions, we get $P_2 = 0.9916$. The probability of innocence of the accused convicted both by the jurymen by the least majority verdict of 7 votes against 5 and by the judges at least by 3 votes against 2 therefore very little differed from 1/100. From the 1597 convicted about 15 likely were not convictable.

The same values of $k$ and $c_2$ lead to

$$\frac{k - c_2}{2k - 1} = 0.1188,$$

and the equation from which we should determine the unknown $t$, becomes

$$1 + 5t + 10t^2 = 0.1188(1 + t)^5,$$

so that $t = 2.789$, $u = 0.7361$. This proves that the chance $u$ of the judges' faultlessness very little differs from 0.7494 as determined in § 140 for jurymen without distinguishing between the categories of crime[14].

**146.** The formulas which we applied to various problems of decisions in criminal cases equally concern all other very numerous decisions as for example made by police and military courts. However, for considering them the observations should provide the necessary data for determining the elements included in those formulas.

The *Comptes généraux de l'administration de la justice criminelle* also contain the results pertaining to the police courts. For the nine years 1825 – 1833 1,710,174 people all over France were brought before them and 1,464,500, or 0.8563 were convicted. This rate little varied from year to year and always remained within the limits of 0.84 − 0.87. The number of judges in police courts was not invariable, but not less than 3 were required, and most often there were 3 indeed. Thus, 2 votes against 1 were sufficient for conviction. Supposing that in equation (118.1a) $n = 3$ and $i = 1$, and replacing $u$ by $t/(1 + t)$, we get

$$c_1 = \frac{k(t^3 + 3t^2) + (1-k)(3t+1)}{(1+t)^3}.$$

We can adopt the approximate and likely value of $c_1 = 0.8563$ obtained from observations, but this is not sufficient for determining the two unknowns, $k$ and $t$. We ought to know in addition how many convictions out of those 1,464,500 were unanimous and how many were returned by majority verdict of 2 votes against 1. This, however, is lacking in the *Comptes généraux*. Suppose that in those police courts the chance of the judges' faultlessness, just as of judges in general, is 3/4, and adopt, in the preceding equation, $c_1 = 0.8563$ and $t = 3$, then $k$ will exceed unity which renders that hypothesis inadmissible. We can suggest that for judges that chance is higher than for the jurymen, but without the necessary observations it is impossible to say by how much.

Military courts consist of 7 judges, and for conviction the law requires a majority verdict at least of 5 votes against 2. Probability $c_2$ of conviction of an accused is derived from equation (118.1a) when $n = 7$, $i = 2$ and $u$ is replaced by $t/(1 + t)$. It occurs that

$$c_2 = \frac{k(t^7 + 7t^6 + 21t^5) + (1-k)(1 + 7t + 21t^2)}{(1+t)^7}.$$

According to the *Comptes généraux de l'administration de la justice militaire* published by the War Minister, the number of the convicted is estimated as 2/3 of the accused. This rate was derived from a large number of decisions and we can therefore adopt it an approximate and likely value of $c_2$. But this is not sufficient for determining the two unknowns included in the preceding equation. When supposing that the chance of faultlessness of military judges and judges of the assize courts very little differ from each other, assuming therefore that it is 3/4 for the former, and supposing that $t = 3$ and $c_2 = 2/3$, it follows from that equation that $k = 0.8793$, $(1 − k) = 0.1207$. We can bet somewhat more than 7 against 1 on the guilt of a serviceman brought before a military court.

Formula (120.1) and equation (118.1a) lead to the expression

$$(1 + t)^7 c_2 P_2 = k(t^7 + 7t^6 + 21t^5)$$

for determining the probability $P_2$ of the guilt of the accused after his conviction. For the preceding values of $c_2$, $t$, $k$ it occurred that $P_2 = 0.9976$ which shows how little the indicated probability differs from certainty. The obtained result is based on a hypothetical value of *t* or *u*, whose degree of precision we do not know, but it would be interesting to be able to compare positively the justice in military courts and assize courts from the standpoint of the probability of their decisions [being proper]. For achieving this goal, it is necessary to know, in addition to the rate of conviction equal to 2/3, the same rates for unanimous convictions and convictions by majority verdicts of 6 votes against 1 and of 5 votes against 2. These data are regrettably lacking in the observations and we are unable to provide any somewhat probable assumption.

**147.** For completing our work, we still have to consider the probability of decisions by tribunals dealing with civil cases. In a civil case it is required to decide who of the litigants has the right on his side. This would have been decided in a certain manner by judges who never err, and their decision will always be unanimous whatever their number.

This, however, never happens. Two equally enlightened judges who most attentively investigate the same case are often led to contrary decisions. We should admit that each judge has a chance of being mistaken at voting, or not to judge as an ideal judge for whom any cause of erring is impossible. This chance depends on the degree of enlightenment and the integrity of the judge, and it is not known in advance. If possible, its value should be derived from observations by methods which we will indicate. And if this, or the contrary chance is established for each judge of some tribunal, it will be possible to derive the probability of the virtue of their decision; or, in other words, of its conformity to the decision that would have been pronounced by ideal judges. It is also possible to establish the probability that other judges, again with a known chance of faultlessness, will confirm the decision of the former.

This second problem is similar to that which we presented in criminal cases. The magnitude formerly denoted by *k* is replaced by the probability that the right is on the side of the litigant determined by the first decision favourable to him. However, if the case is considered by a tribunal for the first time, there is no preliminary probability that the decision will favour one or another side. It is then senseless to consider a probability similar to *k*, and the only unknowns which should be determined by observations are the probability of the judges' faultlessness.

**148.** Consider first of all a tribunal of first instance consisting of three judges, A, A′ и A″. Denote by $u, u'$ и $u''$ their probability of faultlessness and by *c*, the probability of their unanimous decision. Such decisions will occur if either no judge is mistaken, or they are all mistaken. The pertinent probabilities are $uu'u''$ and $(1-u)(1-u')(1-u'')$, so that the composite value of *c* is

$$c = uu'u'' + (1-u)(1-u')(1-u'').$$

After a unanimous decision is returned, two hypotheses can be formulated, either that it was proper or not. The former means that no judge was mistaken, the latter, that all three of them were mistaken. The probability of the observed event, of a unanimous decision, will be $uu'u''$, if the first hypothesis is correct, or $(1-u)(1-u')(1-u'')$, if it is mistaken. When applying to these hypotheses the rule about the probability of causes (§ 28), and denoting by $p$ the probability of the first cause, of the correctness of the decision, we obtain

$$p = \frac{uu'u''}{uu'u'' + (1-u)(1-u')(1-u'')},$$

or, otherwise, $cp = uu'u''$.

If the decision was not unanimous, one judge voted in favour of one side, and the two others, in favour of the other side. Denote by $a$, $a'$ and $a''$ the probabilities of such a decision when the first judge was A, A′ or A″. Then

$a = u'u''(1-u) + u(1-u')(1-u'')$,
$a' = uu''(1-u') + u'(1-u)(1-u'')$,
$a'' = uu'(1-u'') + u''(1-u)(1-u')$.

For example, the first equation conforms to the case in which A′ and A″ were not mistaken, and A was mistaken, or vice versa and the same for the other equations. Now denote by $b$ the probability of some non-unanimous decision. Then

$b = a + a' + a''$,

and since a decision can only be either unanimous or not, $b + c = 1$, which is easy to verify. As a result, it occurs that simply

$b = 1 - uu'u'' - (1-u)(1-u')(1-u'')$.

For the decision to be proper, it is required that the two judges forming the majority vote the same way without erring; for the decision to be mistaken, it is necessary for them to err. Therefore, when denoting by $q$ the probability of a proper non-unanimous decision, and conforming with the rule of the probability of causes or hypotheses,

$bq = (1-u)u'u'' + (1-u')uu'' + (1-u'')uu'$.

Then, having a very large number μ of decisions announced by the same judges A, A′, A″, let γ and β be the numbers of decisions returned unanimously and otherwise. Among the latter let α, α′, α″ be the numbers of decisions in which judges A, A′, or A″ voted contrary to the two others. Then with a very close approximation it will likely occur that γ/μ = $c$, β/μ = $b$, α/μ = $a$, α′/μ′ = $a'$, α″/μ″ = $a''$.

The number β is the sum of α, α′, and α″, and the number $b$, the sum of $a$, $a'$ and $a''$ [see above]. Thus, the second of those equations is the sum of the 3 last ones, and the 5 equations are reduced to 4. If the numbers α, α′, and α″ are known from observations, and the preceding expressions for *a, a′, a″* are substituted in the last three equations, it will be possible to derive the values of *u, u′, u″*, and γ is determined by substituting the expression of *c* in the first equation. And if γ is also known from observations, the comparison of its two values will serve for confirming the theory. Since *u, u′, u″* are also determined, the probabilities *p* and *q* that the unanimous and majority decisions are proper is easily established by the preceding equations.

Observations do not indicate the numbers γ, α, α′, or α″ for any tribunal. However, for showing how to apply these formulas I have arbitrarily chosen the probabilities $u = 4/5$, $u' = 3/5$ and $u'' = 3/5$ and supposed that the chance of faultlessness of each of the three judges exceeds the contrary chance, that judges A′ and A″ are equally educated and their chances of faultlessness are the same whereas A is better educated and his chance of error is lower. Then $c = 8/25$, $b = 17/25$ and $p = 9/10$ and $q = 57/85$. We can bet 17 against 8 or somewhat more than 2 against 1 on the decision of the judges to be not unanimous, and 9 against 1 on a proper unanimous decision and only 57 against 28 or almost 2 against 1 on a proper non-unanimous decision.

For those three judges the mean chance of faultlessness is $(u + u' + u'')/3 = 2/3$. When considering them equally educated and adopting that fraction, 2/3, as the common value of *u, u′, u″*, then

$c = 1/3$, $b = 2/3$, $p = 8/9$, $q = 2/3$.

These values of *p* and *q* are a little smaller than the previous so that in our example an equal education of the three judges lowered the probability of a proper decision both unanimous or not. On the other hand, the latter value of *c* occurred to be higher than the former whereas the former value of *b* exceeded the latter value so that the equal education of the judges heightened the probability of unanimity and therefore lowered the probability of a majority decision.

If, however, we do not know whether the decision was unanimous or not, the grounds for believing it proper will differ from *p* and *q*. Denoting it in this case by *r*, we have

$r = uu'u'' + (1 − u) u'u'' + (1 − u') uu'' + (1 − u'') uu'$,

since under the hypothesis of the decision being proper it, or the observed event can occur in 4 different ways whose probabilities are the 4 terms of this formula. Under the contrary hypothesis the probability of this event will be

$(1 − u)(1 − u')(1 − u'') + u(1 − u')(1 − u'') +$
$u'(1 − u)(1 − u'') + u''(1 − u)(1 − u')$.

The sum of the probabilities of the considered event constitutes certainty or unity, whereas the divisor of *r*, which appears by the rule of § 28, is also unity. Note that, as it was possible to show directly,

$r = cp + bq$.

Taking the previous values of *u, u', u''*, we find that $r = 93/125$. The mean of this value and 2/3 is $r = 20/27$ [19/27]. This second value of *r* is a little less than the previous and therefore the virtue of the decision, as in the previous case, became less probable when all the three judges were equally educated.

**149.** It is not difficult to generalize these formulas on the case of a tribunal consisting of an arbitrary number of judges, but such a result will be impossible to apply owing to the insufficient observational data necessary for determining the chance of faultlessness of the different judges. When supposing that those chances are equal one to another and that there are three judges, then, in previous notation,

$c = u^3 + (1 - u)^3$, $b = 1 - u^3 - (1 - u)^3$, $cp = u^3$,
$bq = 3(1 - u)u^2$, $r = u^3 + 3(1 - u)u^2$.

Adopting now approximate and likely values of γ/μ or β/μ, it will be possible to determine *u* from either of the two first equations. For this determination it will suffice to know the numbers γ and β of the unanimous and non-unanimous decisions out of their total and very large number μ, but we do not have these data. If, for example, we suppose that γ = β, then

$u^3 + (1 - u)^3 = 1 - 3u + 3u^2 = 1/2$

and therefore

$u = [1 \pm \sqrt{3}/3]/2$.

One of the two values of *u* is larger, and the second is smaller than 1/2. We should believe that the judge's chance of faultlessness exceeds the contrary chance and, choosing the former value, $u = 0.7888$, we get

$p = 0.9815$, $q = 0.7885$, $r = 0.8850$.

Suppose that the decision of the 3 judges, whether unanimous or not, is revised by a tribunal of appeal consisting of, say, 7 other judges whose common chance of faultlessness is *v*. Denote by *C* the probability that the decision of the first tribunal will be confirmed by a majority verdict with not less than 4 votes against 3. The value of *C* will be provided by formula (118.1a) with *k* and *u* replaced by *r* and *v* and, for $n = 7$ and $i = 3$,

$C = r[v^7 + 7v^6(1 - v) + 21v^5(1 - v)^2 + 35v^4(1 - v)^3] +$
$(1 - r)[(1 - v)^7 + 7v(1 - v)^6 + 21v^2(1 - v)^5 + 35v^3(1 - v)^4]$.

Suppose that the first decision was proper since confirmed by the second tribunal, then out of the 7 judges of the tribunal of appeal none was mistaken, or 1, 2, or 3 were mistaken. The pertinent probabilities will conform to the 4 terms of the first square bracket. Their sum multiplied by *r* is the probability that the decision was proper and will be confirmed. It is also seen that the part of this expression of *C* which is multiplied by (1 − *r*) expresses the probability that the decision was mistaken but will nevertheless be confirmed and the sum of both parts is the composite expression of *C*. It is seen just as well that, when denoting by *C′* the probability that the second tribunal repealed the decision of the first, then

$$C' = (1 - r)[v^7 + 7v^6(1-v) + 21v^5(1-v)^2 + 35v^4(1-v)^3] + r[(1-v)^7 + 7v(1-v)^6 + 21v^2(1-v)^5 + 35v^3(1-v)^4].$$

The decision of the first tribunal should have been either confirmed or repealed, and therefore *C* + *C′* = 1 which can be verified by noting that the sum of the square brackets in the last formula equals [*v* + (1 − *v*)]⁷ = 1. Then, for *r* = 1/2 and arbitrary value of *v*, or for *v* = 1/2 and arbitrary value of *r* it occurs that C = C′ = 1/2. These results are obvious by themselves.

When separately considering the two parts of *C* and *C′* we can also say that the first part of *C* is the probability that both tribunals decided properly; its second part, the probability that both were mistaken. The same parts of *C′* are the probabilities that the first tribunal was mistaken and the second decided properly whereas its second part expresses the contrary. Denote by ρ the probability that the tribunal of appeal decided properly whether the first tribunal was, or was not mistaken, then that magnitude will be the sum of the first parts of *C* and *C′*, and (1 − ρ), the sum of their second parts which was possible to establish directly.

Denote also by Γ the probability that the decision of that [!] court was confirmed by the second Royal court composed of 7 judges as well, by Γ′, the probability of the contrary event, and by *w*, the chance of faultlessness of each of these 7 judges. Then Γ and Γ′ can be derived from *C* and *C′* when replacing *r* and *v* by ρ and *w*, so that if *w* = *v*

$$\Gamma = \rho^2 + (1 - \rho)^2, \Gamma' = 2\rho(1 - \rho).$$

These values satisfy the condition Γ + Γ′ = 1. The expressions of ρ and ρ′ can also be written as

$$\rho = \frac{r - C'}{2r - 1}, \quad 1 - \rho = \frac{r - C}{2r - 1}.$$

Denote also by *P* and *P′* the probabilities that the decision of the first court of appeal was proper if it confirmed or repealed the decision of the court of first instance. In the first case, when supposing that the decision of the court of first instance was proper or not, the probability of the observed event, i. e., of the coincidence of both decisions will be, respectively, equal to the first and the second part of *C*, and the probability *P* of the first hypothesis will be equal to that first part divided by the sum of both parts. Therefore,

$$CP = r[v^7 + 7v^6(1 - v) + 21v^5(1 - v)^2 + 35v^4(1 - v)^3].$$

In addition, $C'P'$ is equal to all the first part of $C'$, which can also be derived, as it should have been, from formulas (120.1) and (120.2) if $k = r$, $n = 7$, $i = 3$. The obtained results can be replaced by

$$CP = r\rho, \quad C'P' = (1 - r)\rho.$$

**150.** The decision in the court of first instance should be made at least by 3 judges, and at least by 7 in the court of appeal. Usually, these least numbers are not exceeded which is why I adopted them. If $r$ is replaced in my formulas by its value as a function of $u$, they will include the chances $u$ and $v$, which can only be determined by observations. Only one magnitude regrettably is known, the rate of the number of decisions in the court of first instance confirmed by the Royal courts. For applying those formulas it is therefore necessary to reduce both unknowns, $u$ и $v$, to one single magnitude by introducing a particular hypothesis. Most natural, as it seems to me, is to suppose that $v = u$, that is, to suppose that the chance of faultlessness is the same for the judges of both tribunals.

Suppose then that $m$ decisions of the court of first instance out of a very large number $\mu$ of them were confirmed and that, therefore, $(\mu - m)$ were not. The ratio $m/\mu$ can be adopted as an approximate and likely value of the probability which we denoted by $C$, so that

$$C = \frac{m}{\mu}, \quad v = u, \quad u = \frac{t}{1+t}, \quad 1 - u = \frac{1}{1+t}.$$

Therefore,

$$\frac{m}{\mu} = r - \frac{(2r-1)(1 + 7t + 21t^2 + 35t^3)}{(1+t)^7},$$

$$r = 1 - \frac{1+3t}{(1+t)^3}, \quad 2r - 1 = 1 - \frac{2(1+3t)}{(1+t)^3}.$$

By substituting these values in the ratio $m/\mu$, we obtain an equation of the 10[th] degree for determining $t$ and then $u$. With $v = u$ the expression of $C$ will not change if $u$ and $r$ are replaced by $(1 - u)$ and $(1 - r)$, which conforms to replacing $t$ by $1/t$. This means that if $t < 1$ satisfies the given value of $m/\mu$, then $t > 1$ will satisfy that value as well. The equation in the unknown $t$ is *reciprocal* and does not change when $t$ is replaced by $1/t$. We should adopt the value $t > 1$ since it conforms to the value of $u > 1/2$, to the chance of faultlessness exceeding the contrary chance which should take place if the magistrates are honest and educated.

**151.** During the three last months of 1831 and in 1832 and 1833 the *Compte généraux de l'administration de la justice civile* published by the government indicates the numbers $m$ and $(\mu - m)$ of the confirmed and not confirmed decisions for each Royal court. However, for determining $t$ the number $\mu$ is only sufficiently large in the jurisdiction of the Paris Royal court. At present, we are therefore obliged to suppose, just as we did in the case of jurymen, that the chance $u$ of faultlessness is almost the same for all the kingdom's judges. This will enable us to determine $t$ by the values of $m$ and ($\mu - m$) for all the Royal courts taken together.

For the period stated above and for France in its entirety

$m$ = 976, 5301, 5470; $\mu - m$ = 388, 2405, 2617;
$m/\mu$ = 0.7155, 0.6879, 0.6764

The last ratios for the two complete years differ one from another not more than by 1/70 of their mean, which is a very remarkable example of the action of the law of large numbers[15]. When adopting the sums of the numbers for all three periods as $m$ and $\mu$, we get

$m$ = 11,747, $\mu$ = 17,157, $m/\mu$ = 0.6847.

Separately for the Paris Royal court

$m$ = 2510, $\mu$ = 3297, $m/\mu$ = 0.7613.

The obtained ratio $m/\mu$ almost by 1/9 exceeds its mean value for France as a whole. And when adopting the value 0.6847 for France, we arrive at

$t$ = 2.157, $u$ = 0.6832, $r$ = 0.7626.

According to the obtained $r$, without knowing either the tribunal, or the nature of the case, we can bet a little more than 3 against 1 on the decision of the court of first instance to be proper. It is also seen that the chance $u$ of faultlessness for judges of civil cases very little differs from the fraction 0.6788 expressing that chance for the jurymen existing before 1832, i. e., before the law stipulated that *mitigating circumstances* ought to be considered.

When adopting this value of $r$ and the ratios $m/\mu$ и ($\mu - m$)/$\mu$ as the values of $C$ and $C'$, it follows from the formulas of § 150 that

$P$ = 0.9479, $P'$ = 0.6409, $\Gamma$ = 0.7466,

which proves that we can bet almost 19 against 1 on the proper decision of the court of appeal confirming the decision of the court of first instance, and less than 2 against 1 if those decisions were contrary.

Note also that without knowing whether the decision of the court of first instance was confirmed or not, the probability $\Gamma$ that it will be confirmed by a second Royal court when issuing from the same data is a little less than 3/4. The values of the four parts constituting the expressions of $C$ and $C'$ are

$rp$ = 0.6495, $(1 - r)\rho$ = 0.2022,
$r(1 - \rho)$ = 0.1131, $(1 - r)(1 - \rho)$=0.0352.

These fractions, whose sum is unity, express the probabilities that the decisions of the court of first instance and then the court of appeal were proper; that the former was mistaken and the latter proper; that the former was proper and the latter mistaken; that both were mistaken.

**Notes**

**1.** This problem was treated in a memoir, read at the St. Petersburg Academy in June 1834 by Mr. Ostrogradsky, a member of that Academy. Judging by the published extract which the author sent me, he considered the problem in a manner quite different from which I am following in this chapter and which is indicated in the Preamble. Poisson

**2.** Speculative considerations are unfounded. The chance of a juryman's mistake undoubtedly changes in time and depends on his physical and moral state. My witness is an attorney from *The Posthumous Papers of the Pickwick Club*, Chapter 34, who stated that *Discontented or hungry jurymen always find for the plaintiff.*

**3.** This is difficult to understand.

**4.** Concerning the assize courts see Note 7 to the Preamble.

**5.** Apparently, sums (121.1) and (121.2).

**6.** Uniformity alone would have been sufficient for the existence of a linear relation between *X* and *x*.

**7.** An easy generalization led, however, to complicated formulas.

**8.** Poisson did not always calculate with the same degree of approximation.

**9.** At the end of the chapter Poisson also issued from other sources.

**10.** In the sequence, Poisson studied data from several types of courts, and it is sometimes difficult to understand which type he referred to.

**11.** In the same section above, Poisson stated that there were 7 accused in 5 cases.

**12.** Poisson's statement is doubtful. The relative number of crimes in a certain region could have been connected with its economic situation.

**13.** Quetelet repeatedly discussed legal proceedings in Belgium (Sheynin 1986, especially pp. 302 – 303). His mathematics was pedestrian, but, for instance, in 1846 he clearly stated that after the introduction of jurymen the proportion of acquittals doubled. True, in 1833 and later he stated the contrary, but perhaps jurymen were initially selected from a narrow layer. In 1832 and later Quetelet published a table of the rates of conviction of the accused depending on their personality, sex included.

**14.** Assuming *t* = 2.789, I concluded that the left side of that equation was equal to 91.73, but that its right side was 92.78.

**15.** This law of large numbers was confirmed anew by the value of the ratio *m*/μ in 1824. Poisson

[i] (? - O.S.) symbols mean here and later, without additional notice, translator's concern that the original text allows some ambiguity or unclearness.